\documentclass[a4paper,11pt,fleqn]{scrbook} 

\usepackage{amsfonts}
\usepackage{amsmath}
\usepackage{amssymb}
\usepackage{ mathrsfs }

\usepackage[plainpages=false, pdfpagelabels]{hyperref}
\usepackage{amsthm}
\usepackage{cleveref}
\usepackage{fixltx2e}

\newcommand{\A}{{\mathbb A}}
\newcommand{\Z}{{\mathbb Z}}
\newcommand{\F}{{\mathbb F}}
\newcommand{\Q}{{\mathbb Q}}
\newcommand{\R}{{\mathbb R}}
\newcommand{\N}{{\mathbb N}}
\newcommand{\C}{{\mathbb C}}

\newcommand{\G}{\mathbb{G}}
\renewcommand{\H}{\mathbb{H}}
\newcommand{\I}{\mathbb{I}}
\newcommand{\Ac}{\mathcal A}

\newcommand{\Cc}{\mathcal C}
\newcommand{\Dc}{\mathcal D}
\newcommand{\Ec}{\mathcal E}
\newcommand{\Fc}{\mathcal F}
\newcommand{\Gc}{\mathcal G}
\newcommand{\Hc}{\mathcal H}
\newcommand{\Ic}{\mathcal I}

\newcommand{\Oc}{\mathcal O}

\newcommand{\Lc}{\mathcal L}
\newcommand{\Mc}{\mathcal M}

\newcommand{\Rc}{\mathcal R}
\newcommand{\Sc}{\mathcal S}
\newcommand{\Tc}{\mathcal T}

\newcommand{\Vc}{\mathcal V}

\newcommand{\id}{\operatorname{id}}

\newcommand{\Real}{\operatorname{Re}}
\newcommand{\image}{\operatorname{im}}

\newcommand{\vol}{\operatorname{vol}}
\newcommand{\res}{\operatorname{res}}
\newcommand{\pol}{\operatorname{pol}}
\newcommand{\Eis}{\operatorname{Eis}}
\newcommand{\aug}{\operatorname{aug}}
\newcommand{\sgn}{\operatorname{sgn}}
\newcommand{\Norm}{\operatorname{N}}
\newcommand{\Tr}{\operatorname{Tr}}
\newcommand{\Ind}{\operatorname{Ind}}
\newcommand{\Sym}{\operatorname{Sym}}

\newcommand{\df}{{\mathfrak{d}}}

\renewcommand{\lim}[1]{\underset{#1}{lim} \;}

\theoremstyle{plain}
\newtheorem{proposition}{Proposition}[section] 
\newtheorem{lemma}[proposition]{Lemma} 
\newtheorem{theorem}[proposition]{Theorem}
\newtheorem*{thm}{Theorem}  
\newtheorem{corollary}[proposition]{Corollary} 

\theoremstyle{definition}
\newtheorem{definition}{Definition} 
\theoremstyle{remark}
\newtheorem{example}[proposition]{Example} 
\newtheorem{remark}[proposition]{Remark} 

\hypersetup{pageanchor=false}

\input xy
\xyoption {all}
\usepackage[all, knot]{xy}

  \setkomafont{footnote}{\linespread{1.1}\selectfont} 

  \setkomafont{caption}{\small\linespread{1.1}\selectfont} 

\begin{document}

\begin{center}

\begin{Huge}
\textbf{Polylogarithms for \textsl{GL}\textsubscript{2}\\ over totally
real fields}
\end{Huge}

\vspace{1cm}

Philipp Graf

\vspace{1cm}

\textbf{Abstract}
\end{center}

We give a new, purely topological construction of Eisenstein cohomology classes for Hilbert-Blumenthal varieties using the polylogarithm for families of topological tori and a decomposition with respect to the units in the center of $GL_2$. These classes are explicitly calculated in de Rham chomology and compared with Harder's Eisenstein classes. For non-trivial coefficient systems the whole Eisenstein cohomology in positive degrees is generated by these topological Eisenstein classes. This gives an alternative proof for the rationality of Harder's Eisenstein operator without using any multiplicity one arguments. The text constitutes my 2016 Regensburg PhD thesis.

\renewcommand{\baselinestretch}{1.4}\normalsize
\newcommand{\blsd}{\renewcommand{\baselinestretch}{1.4}} 
\newcommand{\asd}{\renewcommand{\arraystretch}{1.0}} 
\parskip1ex

\tableofcontents

\chapter{Introduction and overview}

\section{Eisenstein series and cohomology}\label{G_Eis}

Let $G$ be a reductive linear algebraic group over $\Q$ and $\A$ the adeles over $\Q$. Let us consider the double quotient space
\begin{equation*}
X_{G,K}:=K\backslash G(\A)/G(\Q)A_G(\R)^0.
\end{equation*}
Here $A_G$ is the maximal split torus in the center of $G$, $A_G(\R)^0\subset A_G(\R)$ is the connected component of the identity, $K=K_\infty K_f$ with $K_f\subset G(\A_f)$ compact open and $K_\infty \subset G(\R)$ maximal compact. The space $X_{G,K}$ has the Borel-Serre compactification $\overline {X_{G,K}}$. It is a manifold with corners, which has a stratification with respect to classes $\left\{P\right\}$ of associate parabolic $\Q$-subgroups of $G$.

A representation $(E,\rho)$ of $G$ in a finite dimensional $\Q$-vector space defines a local system $\tilde E$ on $X_{G,K}$ and one can consider the cohomology groups $H^\bullet(X_{G,K},\tilde E)$ and $H^\bullet(X_{G},\widetilde E):=\varinjlim_{K_f} H^\bullet(X_{G,K},\widetilde E)$. The last cohomology group can be considered as the colimit of the group cohomologies of all arithmetic subgroups $\Gamma\subset G(\Q)$ with values in $E$. From this point of view it is clear that these cohomology groups are of fundamental arithmetic interest.\\
When coefficients are extended to $\C$, \cite{Fr} showed that the cohomology classes in these groups may be represented by automorphic forms and that  one has a direct sum decomposition of the cohomology with respect to classes $\left\{P\right\}$ of associate parabolic $\Q$-subgroups of $G$ as $G(\A_f)$-module, whenever $A_G$ acts by a central character on $E$.
The summand corresponding to the parabolic $G$ itself is denoted by $H^\bullet _{cusp}(X_G,\widetilde{E\otimes\C})$ and is called the cuspidal cohomology. It is build up by cusp forms, in other words, automorphic forms whose constant terms at parabolic subgroups different from $G$ are zero. The cuspidal cohomology does not contribute to the cohomology of the boundary. The complement to the cuspidal cohomology in $H^\bullet(X_{G},\widetilde {E\otimes \C})$ is build up by Eisenstein series and residues of such.

As this decomposition of the cohomology is obtained by the theory of Eisenstein series and therefore by an analytic summation process, it is not clear whether it respects the natural $\Q$-structure of the cohomology. The $\Q$-rationality of the decomposition has been proven in the case $G=Res_{K/\Q}GL_n$ with $K/\Q$ a number field, see \cite{Fr} Theorem 20 and $\cite{FrSch}$ 4.3. Theorem, using a classification theorem of \cite{J-Sh} of cuspidal automorphic representations of $GL_n(\A_K)$. 

The $\Q$-rationality of the decomposition above is of great arithmetic interest, as it can be used to derive rationality results for special values of $L$-functions, which may occur as constant terms of Eisenstein series or as integrals of Eisenstein series over cycles, see \cite{Ha3} 3.1 or \cite{Ha1} (4.2.2) and V.

\section{Eisenstein cohomology for Hilbert-Blumenthal varieties}

As the main goal of this thesis is the construction of so called Eisenstein cohomology classes for $G=Res_{F/\Q}GL_2$ with $F$ a totally real number field, we recall Harder's definition of Eisenstein cohomology in this situation in more detail. \cite{Ha1} considers the space $X^\prime _{G,K} :=K\backslash G(\A)/G(\Q)Z(\R)^0$, where $Z\subset G$ is the center and $Z(\R)^0$ the connected component of the identity, and local systems $\widetilde E$ on $X^\prime _{G,K}$ associated to $G$-representations $(E,\rho)$. The space $X^\prime _{G,K}$ is a disjoint union of so called Hilbert-Blumenthal varieties and also has the Borel-Serre compactification whose boundary is homotopy equivalent to $\partial X^\prime_{G,K}:= K\backslash G(\A)/B(\Q)Z(\R)^0$ with $B\subset G$ the standard Borel subgroup of upper triangular matrices. The natural map $\partial X^\prime_{G,K}\to X^\prime_{G,K}$ induces by pullback a restriction map $\res_K:H^\bullet(X^\prime_{G,K},\widetilde E)\to H^\bullet(\partial X^\prime_{G,K},\widetilde E)$ on cohomology and by passing to the colimit a $G(\A_f)$-equivariant map $\res:H^\bullet(X^\prime_{G},\widetilde E)\to H^\bullet(\partial X^\prime_{G},\widetilde E)$.
The subspace $H^\bullet _{\Eis}(X^\prime_{G},\widetilde E)\subset H^\bullet(X^\prime_{G},\widetilde E)$, such that $\res:H^\bullet_{\Eis}(X^\prime_{G},\widetilde E)\to \image(\res)$ is an isomorphism, is called the Eisenstein cohomology. 

Harder determines the Eisenstein cohomology in two steps. First he describes the cohomology of the boundary. It may be understood as a sum of induced modules $\Ind^{G(\A_f)}_{B(\A_f)}\C\phi_f$, where the $\phi_f$ are the finite components of algebraic Hecke characters $\phi:T(\A)/T(\Q)\to \C^*$ and $T$ is the maximal torus in $G$, see \cite{Ha1} Theorem 1. More precisely, if a cohomology class on $\partial X^\prime_{G}$ is represented by a $B(\Q)$-invariant differential form $\omega$, then $\omega$ is in particular $U(\Q)$- invariant, where $U\subset B$ is the unipotent radical, and therefore $\omega$ may be developed into a Fourier series with respect to the group $U(\Q)$. The class $\omega$ is then already determined by the constant term of the differential form $\omega$, this means by its zeroth Fourier coefficient. 

As a second step Harder constructs a $G(\A_f)$-equivariant operator
\begin{equation*}
\Eis: \image(\res)\otimes \C\to H^\bullet_{\Eis}(X^\prime_{G},\widetilde E\otimes \C),
\end{equation*}
 which is a section for $\res$. Explicitly it may be described as follows: If $\omega \in\image(\res)\otimes \C$ is represented by a $B(\Q)$-invariant differential form, then $\Eis(\omega):=\sum_{\gamma\in G(\Q)/B(\Q)}\gamma^*\omega$ where the sum has to be defined by analytic continuation in general. If we have a trivialization $\tilde E\cong \C$, it turns out that $\res:H^{2\xi-1}(X^\prime_{G},\C)\to H^{2\xi-1}(\partial X^\prime_{G},\C)$ is not surjective, where we set $\xi:=[F:\Q]$. 

\section{The topological polylogarithm and associated Eisenstein classes}\label{pol_intro}

The connection between Eisenstein classes and special values of $L$-functions may be seen as motivation to construct $\Q$-rational Eisenstein classes geometrically. One way to do so is to specialize polylogarithms. \cite{Be-L} constructed polylogarithms for relative elliptic curves. In the case of the universal elliptic curve with level-$N$-structure $\Ec\to X_{G,K}$, $G=SL_{2}$, $K_f=\ker(SL_2(\hat\Z)\to SL_2(\Z/N\Z))$, Beilinson proved that the polylogarithm specialized along non-zero $N$-torsion sections actually yields Eisenstein classes for $SL_2$.\\  
Following  the ideas of Beilinson  and \cite{No} we can adapt this construction easily to our topological situation. Given a group $G$ as before with a finite dimensional representation $\psi:G\to Aut(V)$ we may consider the space
\begin{equation*}  
\pi:T_{G,K}:=V(\hat \Z)\rtimes KA_G(\R)^0\backslash V(\A)\rtimes G(\A)/V(\Q)\rtimes G(\Q)\to X_{G,K}
\end{equation*}
for $K_f\subset \ker(Aut(V(\hat\Z))\to Aut(V(\Z/N\Z)))$. It is a group object over $X_{G,K}$ and its fibers are topological tori, in other words, isomorphic to products of $S^1$. One has the pro local system $Log$ associated to the $V(\Q)\rtimes G(\Q)$-module $\prod_{k\geq 0}\Sym^kV(\Q)$, where $V(\Q)$ acts by multiplication with the exponential series and $G(\Q)$ via $\psi$. Furthermore, we consider $D\subset T_{G,K}$, which is the union of the images of non-zero $N$-torsion sections with open complement $U$, and the relative orientation bundle $\mu:=\widetilde{\det(V) }$. From the cohomological vanishing properties of $Log$ one easily derives short exact localization sequences
\begin{equation*}
 0\to H^{dim(V)-1}(U,Log\otimes \pi^{-1}\mu^{n+1})\to H^0(D, \pi_{|D}^{-1}\prod_{k\geq0 }\Sym^k\widetilde{V}\otimes \mu^n)\to H^{0}(X_{G,K},\mu^n)\to 0
\end{equation*}
for all $n\in \Z$.
Given a section $f\in H^0(D, \pi_{|D}^{-1}\prod_{k\geq0 }\Sym^k\widetilde{V}\otimes \mu^n)$ mapping to zero on the right-hand side we get a unique cohomology class called the polylogarithm associated to $f$
\begin{equation*}
\pol(f)\in H^{dim(V)-1}(U,Log\otimes \pi^{-1}\mu^{n+1}).
\end{equation*}
This class may be specialized along the zero section to obtain polylogarithmic Eisenstein classes
\begin{equation*}
(\Eis^k(f))_{k\geq 0}:=0^*pol(f)\in \prod_{k\geq 0}H^{dim(V)-1}(X_{G,K},\Sym^k\widetilde{V}\otimes\mu^{n+1}).
\end{equation*}
If we have $n=0$, typical examples for such $f$ may be given by functions $f:V(\Z/N\Z)\to \Q$ with $f(0)=\sum_{v\in V(\Z/N\Z)}f(v)=0$. By the naturality properties of the polylogarithm these polylogarithmic Eisenstein classes glue to $G(\A_f)$-equivariant operators
\begin{equation*}
\Eis^k:\Sc(V(\A_f),\Q)^0\otimes H^0(X_G,\mu^{n})\to H^{dim(V)-1}(X_{G},\Sym^k\widetilde{V}\otimes\mu^{n+1}),
\end{equation*}
where $\Sc(V(\A_f),\Q)^0$ are the $\Q$-valued Schwartz-Bruhat functions with $\int_{V(\A_f)}f(v)dv=f(0)=0$.
\cite{L} and \cite{No} represented the polylogarithm cohomology classes by explicit currents in the cohomology with $\C$-coefficients. The polylogarithmic Eisenstein classes are then represented by Eisenstein-Kronecker series. With this explicit description \cite{Bl2} and \cite{Ki1} proved in the situation $G=Res_{F/\Q}SL_2$ with $F/\Q$ a totally real number field, that again the polylogarithmic Eisenstein classes are non-trivial Eisenstein classes, as their constant terms turned out to be special values of partial $L$-functions associated to the field $F$. In this way they also proved that these special $L$-values have to be rational numbers, a result which already goes back to Siegel.

Here we only addressed the topological realization of the polylogarithm. Indeed, by \cite{Ki2} the polylogarithm for abelian schemes is of motivic origin and a powerful tool to tackle deep arithmetic problems and conjectures. See for example \cite{Ki3}, where the \'{e}tale elliptic polylogarithm is a decisive instrument for the proof of the Tamagawa number conjecture for CM elliptic curves.

\section{Decomposition of Eisenstein classes}\label{Eis_intro}

Even though the polylogarithmic Eisenstein classes are a mighty tool in arithmetic, there is one obvious flaw: They are all stuck in cohomological degree $dim(V)-1$. 

In this work we present a decomposition principle for polylogarithmic Eisenstein classes in the case $G=Res_{F/\Q}GL_2$, $V=Res_{F/K}\G_a^2$ and $F/\Q$ a totally real number field. \\
The idea is very easy. We may see $X_{G,K}$ as a fiber bundle
\begin{equation*}
\varphi:X_{G,K}\to X^\prime _{G,K}= K\backslash G(\A)/G(\Q)Z(\R)^0,
\end{equation*}
where again $X^\prime_{G,K}$ is a space considered by \cite{Ha1} in his landmark paper. 
The fiber is $\varphi^{-1}(1)=(A_G(\R)^0K_\infty \cap Z(\R)) \backslash Z(\R)/Z_K$
and $Z_K:=Z(\Q)\cap K_f\subset \Oc_F ^{\times}$ is a subgroup of finite index of the units of $\Oc_F\subset F$ the ring of integers. By Dirichlet's unit theorem we know that the fiber is compact and that we have for the cohomology 
\begin{equation*}
H^\bullet(\varphi^{-1}(1),\Q)=H^\bullet(Z_K,\Q)=H^\bullet(\Oc_F ^{\times},\Q)=\bigwedge^\bullet Hom(\Oc_F ^{\times},\Q).
\end{equation*}
Using the coordinate on $G(\R)$ coming from the determinant we see that $\varphi$ is up to a finite covering a trivial bundle. So one has cohomology classes in $H^\bullet (X_{G,K},\Q)$ restricting to a basis of the cohomology of all fibers of $\varphi$ and this gives a Leray-Hirsch isomorphism 
\begin{equation*}
H^\bullet(\varphi^{-1}(1),\Q)\otimes H^\bullet(X^\prime _{G,K},\varphi_*\Sym^k\widetilde{V}\otimes \mu^{n+1}) \to H^\bullet(X_{G,K},\Sym^k\widetilde{V}\otimes \mu^{n+1}). 
\end{equation*}
The polylogarithmic Eisenstein classes may then be decomposed with respect to this isomorphism.
We get by evaluation
\begin{equation*}
\Sc(V(\A_f),\Q)^0\otimes H^0(X_G,\mu^{n})\otimes H^q(\varphi^{-1}(1),\Q)^{*}\stackrel{\Eis^k_q}{\rightarrow} H^{dim(V)-1-q}(X_{G}^\prime,\varphi_*\Sym^k\widetilde{V}\otimes\mu^{n+1})
\end{equation*}
for $q=0,...,\xi-1=dim(\varphi^{-1}(1))$.

We want to show that this decomposition of polylogarithmic Eisenstein classes is as non-trivial as possible. In other words, we want to get as many of Harder's Eisenstein cohomology classes as possible. To do so we follow \cite{L} and \cite{No} to represent the polylogarithm by a current and hence the polylogarithmic Eisenstein classes by differential forms in de Rham cohomology. The decomposition isomorphism will be made explicit by fiber integration on the level of de Rham cohomology. This gives the decomposed polylogarithmic Eisenstein classes as Mellin transforms of theta series as considered by \cite{Wi}. 

To determine the image of our polylogarithmic Eisenstein classes we calculate the constant terms, in other words, the restriction to the cohomology of the boundary. 
We already have mentioned that we have thanks to \cite{Ha1} a complete understanding of the cohomology of the boundary of $\overline{X^\prime_{G,K}}$ as a sum of induced modules $\Ind^{G(\A_f)}_{B(\A_f)}\C\phi_f$.
The precise relation between $f$ and the constant term of the polylogarithmic Eisenstein class $\Eis^k(f)$ is then controlled by the horospherical map
\begin{equation*}
\rho:\Sc(V(\A_f),\Q)^0\otimes H^0(X_G,\mu^{n})\to \bigoplus \Ind^{G(\A_f)}_{B(\A_f)}\C\phi_f.
\end{equation*}
The horospherical map gives also the relation to special values of $L$-functions.

To understand the image of our polylogarithmic Eisenstein classes we have to determine which of the induced functions above may be realized by the horospherical map.

Our main result is as follows
\begin{thm}
The operators $\Eis^k_q$ factor through the Eisenstein cohomology. The image of $\Eis^k_q$ is the whole Eisenstein cohomology in degree $2\xi-1-q$, if $k>0$. If $k=0$, $\Eis^0_0$ generates the Eisenstein cohomology in degree $2\xi-1$, for $q>0$ we get all Eisenstein classes but not those associated to spherical functions (see \cref{spherical}).
\end{thm}
We want to remark a few things. 
\begin{enumerate}
\item The spherical functions are exactly those functions which are not in the image of the restriction map to the cohomology of the boundary in degree $2\xi-1$. As the $\Eis^k_q(f)$ have the same constant term for all $q$, we cannot obtain the Eisenstein classes associated to spherical functions in cohomological degrees $\xi,...,2\xi-1$.

\item The Eisenstein cohomology of $X^\prime _{G,K}$ is supported in cohomological degrees $0$ and $\xi,...,2\xi-1$, so that our Eisenstein operator actually generates most of the relevant part.

\item Even though Harder started with more general coefficient systems than we do, see \cite{Ha1} 1.4, we get all Eisenstein cohomology classes for non-trivial representations. The reason is that after extending coefficients to $\overline \Q$ his representations are direct summands in $\Sym^k V\otimes \det(V)^n$, $k\geq1$, and all weights $\phi_{|T(\R)^0}$ occurring in the cohomology of the boundary have to factor through the norm character, see \cite{Ha1} 2.8.   

\item We even get a much finer result. The sheaf $Log$ has an integral structure which allows us to define our polylogarithmic Eisenstein classes over the Ring $\Z[\frac{1}{N}]$, when we consider $X_{G,K}$ and $K_f$ defines the level-$N$-structure. From this we also deduce integrality results for special values of partial $L$-functions for totally real fields.

\item However, it seems as if one cannot expect that the polylogarithm always generates much of the Eisenstein cohomology. To demonstrate this problem look at $F=\Q(\sqrt{-1})$ and $G=Res_{F/\Q}GL_2$. One would naively choose $V=Res_{F/\Q}\G_a^2$. We have $dim(V)=4$ and this gives polylogarithmic Eisenstein classes in cohomological degree $3$. Now we have $X_{G,K}=X^\prime_{G,K}$, since $A_{G}(\R)^0K_\infty=\R_{>0}U(2)=\C^{\times}U(2)=Z(\R)^0K_\infty$ and the part of $Z(\R)$, which does not split over $\R$, already lies in the maximal compact subgroup. So there is no decomposition and $H^3(X_{G,K},\widetilde{E})=0$ by Poincar\'{e} duality. In this case the polylogarithm gives nothing. Nevertheless, we hope for applications of the polylogarithm to other groups. Especially the case of $G$ different from $GL_{n,F}$ with $F$ a number field would be interesting, as one does not know the $\Q$-rationality of the decomposition of the cohomology in this case, see \cref{G_Eis}. 

\end{enumerate}
\section{Outline of the thesis}

This thesis has three main parts. We want to discuss them here in more detail.

\subsection{The topological polylogarithm}

In the first part we define the general geometric setup. We begin with the definition of a family of topological tori over a manifold $S$ as a proper submersion $\pi:T\to S$ with connected fibers, which is a commutative group object over $S$. If $S$ is connected one identifies the category of families of topological tori with the category of finitely generated free abelian groups $L$ with $\pi_1(S)$-action (\Cref{tori}). This is done by associating to $L$ the quotient-torus $(L\otimes \R/L)\times_{\pi_1(S)}\widetilde S $. Here $\widetilde S$ is the universal cover of $S$ where $\pi_1(S)$ acts on by deck transformations and $\pi_1(S)$ acts on the left factor by the module structure. So the torus $T$ is determined by a representation $\pi_1(S)\to Aut(L)$, which we call the representation associated to $T$. 

This allows us to define locally constant sheaves on $T$ by using $L\rtimes \pi_1(S)$-modules. Let $A$ be any noetherian ring. $L\rtimes \pi_1(S)$ acts on $L$ by affine transformations. This makes the group ring $A[L]$ a $L\rtimes \pi_1(S)$-module. The augmentation $A[L]\stackrel{\aug}{\rightarrow} A$ is equivariant making the augmentation Ideal $\mathfrak a:=\ker(\aug)$ a  $L\rtimes \pi_1(S)$-module. Consequently, $A[[L]]:=\varprojlim_{n\in \N_0}A[L]/\mathfrak a ^{n+1}$ is a $L\rtimes \pi_1(S)$-module and the associated locally constant sheaf on $T$ is called the logarithm sheaf $Log$. 

Next we deduce a purity result for projective systems of local systems indexed over the natural numbers. If $i:D\to T$ is a closed submanifold of codimension $c$ of our torus we get $Ri^!Log=i^{-1}Log\otimes or_{D/T}[-c]$ (\Cref{lim_purity}), where $or_{D/T}$ is the relative orientation bundle. This gives a natural localization triangle 
\begin{equation*}
i_*i^{-1}Log\otimes or_{D/T}[-c]\to Log\to j_*j^{-1}Log\stackrel{+1}{\rightarrow},
\end{equation*}
where $j:U\to T$ is the inclusion of the open complement of $D$ (\Cref{localization_triangle}). We calculate the well known higher right derived images of $Log$
\begin{equation*}
R^p\pi_*(Log)=0, p\neq d,\ R^d\pi_*(Log)\stackrel{\aug \cong}{\rightarrow}R^d\pi_*(A),
\end{equation*}
where $d$ is the fiber dimension of $\pi$. From this we deduce the localization sequence 
\begin{equation*}
0\to H^{d-1}(U,Log\otimes or_{T/S}^{n+1})\to H^0(D,i^{-1}Log\otimes or_{T/S}^{n})\stackrel{\aug}{\rightarrow} H^0(S,or_{T/S}^{n}),
\end{equation*}
when $\pi:D\to S$ is a finite cover (\Cref{global_Log_localization}). 

To profit from the localization sequence, this means to get non-trivial polylogarithmic classes, we need a good understanding of $H^0(D,i^{-1}Log\otimes or_{T/S}^{n})$. Following \cite{BKL} we show that there is a section $A\to i^{-1}Log$ of $\aug$, if $D$ is the union of the images of torsion sections whose order is invertible in $A$ (\Cref{trivializationLog}). This gives us the definition of $\pol(f)\in H^{d-1}(U,Log\otimes or_{T/S}) $ for locally constant functions $f:D\to A$ with $\aug(f)=0$, as we have $H^0(D,A)\subset H^0(D,i^{-1}Log)$ (\Cref{pol}). 

For computational purposes we also need a trivialization of the pro vector bundle associated to $Log$. It comes from the isomorphism 
\begin{equation*}
\Q[[L]]\to \prod_{k\geq0}\Sym^k(L\otimes \Q),l\mapsto \exp(l),
\end{equation*}
which gives the nowhere vanishing $\Cc^\infty$-section $v\mapsto \exp(-v)$, $v\in L\otimes\R$, and the relation of $Log$ to the symmetric powers of the representation $\pi_1(S)\to Aut(L)$ associated to $T$: $0^{-1}Log=\prod_{k\geq0}\Sym^k\widetilde{(L\otimes \Q)}$ (\Cref{cont_triv}), where $0:S\to T$ is the zero section of the group object $T$.

Moreover, we show that the  localization sequence of $Log$ may be calculated by resolving the logarithm sheaf by non-continuous functionals. This gives a characterization of the polylogarithms by differential equations (\Cref{{differential_equ}}). 
    
\subsection{Polylogarithmic Eisenstein classes for Hilbert-Blumenthal varieties}

First we recall the geometric situation $\varphi:\Mc_K:=X_{G,K}\to \Sc_K:=X^\prime_{G,K}$, when we have $G=Res_{F/\Q}GL_2$ for a totally real field $F$. Let us keep the notation already fixed in \Cref{pol_intro} and \Cref{Eis_intro}. We construct a torus over $\Mc_K$ using the standard representation $G\to Aut(V)$, $V=Res_{F/\Q}\G_a^2$. This allows us to define polylogarithms and polylogarithmic Eisenstein classes for $\Mc_K$. We use the naturality properties of $Log$ and the localization sequence to show that we can glue our Eisenstein classes to $G(\A_f)$-equivariant operators as described in \Cref{pol_intro} (\Cref{Eis^k}).

As we have the polylogarithmic Eisenstein classes on $\Mc_K$, we want to decompose the cohomology. The proof is divided up into several parts. 

First we calculate the higher direct images of $\Sym^k\widetilde V\otimes \mu^{n+1}$. As $\varphi$ is a fiber bundle it suffices to understand the cohomology of the fibers, which we calculate using group cohomology $H^\bullet(Z_K,\Sym^k V\otimes \det(V)^{n+1})$. As the action of $Z_K$ is semi-simple we easily calculate these groups as $H^0(Z_K,\Sym^k V\otimes \det(V)^{n+1})\otimes H^\bullet(Z_K,\Q)$ (\Cref{{cohomology_abelian_group}}). 

This gives $R^\bullet\varphi_*(\Sym^k\widetilde V\otimes \mu^{n+1})=\varphi_*\Sym^k\widetilde V\otimes \mu^{n+1}\otimes R^\bullet\varphi_*(\Q)$ (\Cref{fiber_cohom}) and we trivialize $R^\bullet\varphi_*(\Q)$ by forms coming from global classes  $\mathfrak H^\bullet \subset H^\bullet (\Mc_K,\Q)$ using the fact that $\varphi$ is up to a finite cover a trivial fiber bundle (\Cref{global_classes}). 

These global classes define then by cup-product the decomposition isomorphism (\Cref{decomp})
$H^\bullet(\Sc_K,\varphi_*\Sym^k\widetilde V\otimes \mu^{n+1})\otimes \mathfrak H^\bullet\stackrel{\cup}{\rightarrow} H^\bullet(\Mc_K,\Sym^k\widetilde V\otimes \mu^{n+1})$,
which we finally discuss in the setting of de Rham cohomology with the theory of fiber integration (\Cref{trace_integration}).

\subsection{Comparison with Harder's Eisenstein classes}

In the last part we represent the polylogarithmic Eisenstein classes by differential forms and compare them with those of Harder.

Using the ideas of \cite{No} we represent the polylogarithms explicitly by currents describing the Eisenstein classes as differential forms. We give the description in adelic coordinates and calculate the decomposition of the polylogarithmic Eisenstein classes by fiber integration (\Cref{decomp_Eis}). We reinterpret this as a Mellin transform of a theta series as considered by \cite{Wi} (\Cref{Wi-series}).

As Harder's Eisenstein classes are determined by their restriction to the cohomology of the boundary, we recall Harder's calculation of the cohomology of the boundary of $\Sc_K$ (\Cref{boundary}). We restrict the polylogarithmic Eisenstein classes to the boundary where they are determined by their constant terms (\Cref{boundary_residueI}) and derive rationality and integrality results for these constant terms (\Cref{integral_L_value}). Then we define the horospherical map controlling the relation between $f$ and the constant term of $\Eis^k(f)$ and explain how the horospherical map determines the image of our Eisenstein operators (\Cref{hor}).

Next we translate our cohomology classes to $(\mathfrak g,K)$-cohomology and compare them there with Harder's Eisenstein operator. We see that the polylogarithmic Eisenstein classes are in the image of Harder's Eisenstein operator and therefore our operators $\Eis^k_q$ actually factor through the Eisenstein cohomology (\Cref{comparison}).

Finally, we determine the image of our operators $\Eis^k_q$ by studying the horospherical map. The main ingredient is to show that the horospherical map is surjective, when one allows general Schwartz-Bruhat functions, this means functions $f$ where we do not necessarily have $\int_{V(\A_f)}f(v)dv=f(0)=0$ (\Cref{hor_surjective}).

\section{Acknowledgments}

I would like to thank my advisor Guido Kings for giving me the opportunity to work on his beautiful idea to decompose the polylogarithmic Eisenstein classes for $Res_{F/\Q}GL_2$ by using the units in the center. His encouragement was decisive for the success of this work and I want to thank him for everything I have learned during my studies.

\chapter{The topological polylogarithm}

\section{The logarithm sheaf on families of topological tori}

We want to start with the definition and construction of the logarithm sheaf on families of topological tori. These will be our main geometric objects. As we want to calculate cohomology classes on them explicitly, we need practical and explicit descriptions for locally constant sheaves on them. This is achieved by the theory of representations of the fundamental group and equivariant sheaves on the universal cover.
\subsection{Families of topological tori}

For any site $C$ we denote by $Sh(C)$ the category of abelian sheaves on $C$.
For a $\Cc^\infty$-manifold $\Sc$ we denote by $Mfd/\Sc$ the category of manifolds over $\Sc$. 
\begin{definition}
Let $\pi: \Tc\rightarrow \Sc$ be a proper submersion of $\Cc^\infty$-manifolds. We call $\pi: \Tc\rightarrow \Sc$ a \textit{family of topological tori} or simply a \textit{torus over $\Sc$}, if it is a commutative group object in $Mfd/\Sc$ with connected fibers. Families of topological tori together with group homomorphisms form a subcategory of $Mfd/\Sc$, which we denote by $Tori/\Sc$.
\end{definition}
\begin{remark}
As the fibers of $\pi$ are compact commutative Lie groups, they are topological tori, in other words products of $S^1:=\left\{z\in \C:|z|=1\right\}$.
\end{remark}
\begin{remark}
We will always consider $Mfd/\Sc$ with the usual topology of open covers. It follows that we have a Yoneda embedding
\begin{equation*}
Tori/\Sc\rightarrow Sh(Mfd/\Sc),\ \Tc \mapsto \left\{X\mapsto \Tc(X)=Hom_{Mfd/\Sc}(X,\Tc)\right\}
\end{equation*}
\end{remark}
Let us fix a base point $s_0\in \Sc$ and let us take $t_0:=0(s_0)\in \Tc$ as base point for $\Tc$, where $0:\Sc\rightarrow \Tc$ is the zero section. If $\Sc$ is not connected, we choose a base point for any connected component and do all our constructions for each connected component separately.  
Any proper submersion of $\Cc^\infty$-manifolds is a fiber bundle, see \cite{E}. From this we deduce a short exact sequence
\begin{equation*}
0\rightarrow \pi_1(\Tc_{s_0},t_0)\stackrel{i_*}{\rightarrow} \pi_1(\Tc,t_0)\stackrel{\pi_*}{\rightarrow} \pi_1(\Sc,s_0)\rightarrow 0  ,
\end{equation*}  
where as usual $\Tc_{s_0}:=\pi^{-1}(s_0)$ denotes the fiber over $s_0$.
This sequence is split exact, since $\pi_*$ has the section $0_*$, and we have a natural left $\pi_1(\Sc,s_0)$-action on $\pi_1(\Tc_{s_0},t_0)=H_1(\pi^{-1}(s_0),\Z)$ given by conjugation, as the latter group is commutative. In other words, there is an isomorphism
\begin{equation*}
 \pi_1(\Tc_{s_0},t_0)\rtimes \pi_1(\Sc,s_0)\cong \pi_1(\Tc,t_0),\ (l,g)\mapsto i_*(l)0_*(g),
 \end{equation*} 
where the structure of the semidirect product is determined by $0$ and the by the action above. \\
We are going to construct locally constant sheaves on $\Tc$, which are induced by $\pi_1(\Tc,t_0)$-modules. The functor $\Fc\mapsto \Fc_{t_0}$, which assigns to a locally constant sheaf $\Fc$ its stalk at the base point, establishes an equivalence between the category of locally constant sheaves on a connected topological manifold $\Tc$ and the category of left-$\pi_1(\Tc,t_0)$-modules, see  \cite{Iv} IV Theorem 9.7. Note that we consider the universal cover $\widetilde{\Tc}$ of $\Tc$ equipped with a right $\pi_1(\Tc,t_0)$-action. 
\begin{example}
If $\pi:\Tc\to \Sc$ is a torus, we have the locally constant sheaf $\Hc:=\underline{Hom}_\Z(R^1\pi_*(\Z),\Z)$ and $\Hc_{s_0}=H_1(\Tc_{s_0},\Z)$ by Poincar\'{e} duality.
\end{example}
\begin{definition}
A \textit{lattice $L$} is a free $\Z$-module of finite rank. Let $A\to A^\prime $ be a ring homomorphism and $M$ a $A$-module. We set $M_{A^\prime}:=M\otimes_A A^\prime$.
\end{definition}

\begin{proposition}\label{tori}
Let $\Sc$ be a connected $\Cc^\infty$-manifold, say with base point $s_0\in \Sc$. Then we have an equivalence of categories
\begin{equation*}
Tori/\Sc\rightarrow \left\{\text{lattices $L$ with left-$\pi_1(\Sc,s_0)$-action and equivariant maps}\right\}
\end{equation*}
\begin{equation*}
\Tc\mapsto H_1(\Tc_{s_0},\Z)
\end{equation*}
\begin{proof}
Given a $\pi_1(\Sc,s_0)$-module $L$ we have the torus $T:=L_\R/L$. As $\pi_1(\Sc,s_0)$ acts on $L$, we get an induced action on $T$. We denote by $\tilde{\Sc}$ the universal cover of $\Sc$ and set $\Tc_{L}:=\tilde{\Sc}\times_{\pi_1(\Sc,s_0)}T$, where $\tilde{\Sc}\times_{\pi_1(\Sc,s_0)}T:= \tilde{\Sc}\times T/\pi_1(\Sc,s_0)$ and $\gamma\in \pi_1(\Sc,s_0)$ acts on $(s,t)\in\tilde{\Sc}\times T$ by $(s,t)\gamma=(s\gamma,\gamma^{-1}t)$. We have a structure map $\pi:\Tc_{L}\rightarrow \Sc$ coming from the first projection.
The action of $\pi_1(\Sc,s_0)$ on $\tilde{\Sc}$ is properly discontinuous and fixpoint free. So $\tilde{\Sc}$ may be covered by $\tilde{U}\subset \tilde{\Sc}$ open where $\tilde{U}\gamma\cap \tilde{U}\neq \emptyset$ implies $\gamma=\id$. If we denote by $U$ the image of $\tilde{U}$ in $\Sc$ under the canonical map, we get $\pi^{-1}(U)\cong U\times T$. So our quotient space is a fiber bundle with typical fiber $T$ and $\pi$ is proper. Since the action is by diffeomorphisms, $\Tc_{L}$ is even a manifold.\\ We endow $pr_1:U\times T\rightarrow U$ with the constant group structure as a group object in manifolds over $U$. By construction all transition maps of our fiber bundle come from the $\pi_1(\Sc,s_0)$-action on $T$. So all transition maps are group homomorphisms and this means that we may glue all group structures $pr_1:U\times T\rightarrow U$ to obtain a global group structure on $\Tc_{L}$.
Obviously, equivariant maps $f:L\rightarrow L^\prime$ induce maps on the corresponding spaces. This gives the desired quasi-inverse functor
\begin{equation*}
\left\{\text{lattices $L$ with left-$\pi_1(\Sc,s_0)$-action and equivariant maps}\right\}\to Tori/\Sc
\end{equation*}
\begin{equation*}
L\mapsto \tilde\Sc\times_{\pi_1(\Sc,s_0)} L_\R/L=:\Tc_L.
\end{equation*}
\end{proof}
\end{proposition}
\begin{remark}\label{cover_tori}
We deduce from the theorem above that our family of topological tori $\Tc$ has the universal cover $\tilde{\Sc}\times H_1(\Tc_{s_0},\R)$ where the fundamental group $\pi_1(\Tc,t_0)=H_1(\Tc_{s_0},\Z)\rtimes \pi_1(\Sc,s_0)$ acts by $(x,v)\cdot(l,\gamma)=(x\gamma,\gamma^{-1}(v+l))$ and the action of $\pi_1(\Sc,s_0)$ on $H_1(\Tc_{s_0},\Z)$ is the monodromy action coming from the locally constant sheaf $\Hc$.
\end{remark}
\begin{definition}
A group homomorphism $\phi:\Tc^\prime\rightarrow \Tc$ of tori over $\Sc$ is called an \textit{isogeny}, if $\phi$ is surjective and $\ker(\phi)$ is a finite covering of $\Sc$. Here 
\begin{equation*}
\ker(\phi)(X):=\ker(\phi:\Tc^\prime(X)\rightarrow \Tc(X)),\ X\in Ob(Mfd/\Sc)
\end{equation*}
\end{definition}
\begin{remark}
The kernel of a homomorphism $\phi:\Tc^\prime\rightarrow \Tc$ of tori is always representable as a base change with the zero section. More precisely, \Cref{tori} shows that $\ker(\phi)=\tilde{\Sc}\times_{\pi_1(\Sc,s_0)}\ker(\phi_{s_0})$ where $\phi_{s_0}:\Tc^\prime_{s_0}\rightarrow \Tc_{s_0}$ is the induced homomorphism on the fiber over $s_0$.
\end{remark}
\begin{remark}
Given a torus $\Tc$ over $\Sc$ we have for any $N\in \Z$ the \textit{multiplication by $N$ isogeny} $[N]:\Tc\rightarrow \Tc$.
On points $X\in Ob(Mfd/\Sc)$ it is determined by the functorial group homomorphism
\begin{equation*}
[N]:\Tc(X)\rightarrow \Tc(X), v\mapsto N\cdot v.
\end{equation*}
We call $\ker([N])=:\Tc[N]$ the \textit{$N$-torsion subgroup}, it is a finite covering of $\Sc$. 
\end{remark}
\begin{definition}
Let $\phi:\Tc^\prime\rightarrow \Tc$ be an isogeny over $\Sc$. The function 
\begin{equation*}
deg(\phi):\Sc\rightarrow \Z,\ s\mapsto |\ker(\phi_s)|
\end{equation*}
is locally constant and is called the \textit{degree} of $\phi$.
\end{definition}
\begin{remark}
We have $deg([N])=N^{dim(\Tc/\Sc)}$ where 
\begin{equation*}
dim(\Tc/\Sc):\Sc\rightarrow\Z,\ s\mapsto dim(\Tc_s).
\end{equation*}
\end{remark}
\begin{definition}
A family of topological tori $\pi:\Tc\rightarrow \Sc$ has a \textit{level-$N$-structure}, if the sheaf of sections associated to $\Tc[N]$ in $Sh(Mfd/\Sc)$ is constant. If the fiber dimension $dim(\Tc/\Sc)$ is constant, this just means that there is an isomorphism
$\Tc[N]\cong \Sc\times (\Z/N\Z)^{dim(\Tc/\Sc)}$.
\end{definition}
\begin{remark}
For any ring $A$ we set $GL(M):=Aut_{A-Mod}(M)$, if $M$ is a finitely generated free $A$-module.
If $\pi:\Tc\rightarrow \Sc$ has a level-$N$-structure, the associated representation $\rho:\pi_1(\Sc,s_0)\rightarrow GL(H_1(\Tc_{s_0},\Z))$ factors as $\rho:\pi_1(\Sc,s_0)\rightarrow GL(H_1(\Tc_{s_0},\Z))(N)$ where 
\begin{equation*}
1\rightarrow GL(H_1(\Tc_{s_0},\Z))(N)\rightarrow GL(H_1(\Tc_{s_0},\Z))\rightarrow GL(H_1(\Tc_{s_0},\Z/N\Z))\rightarrow 1
\end{equation*}
is exact with the natural projection on the right hand side.
\end{remark}

\subsection{Definition of the logarithm sheaf}

Consider the group algebra $A[\pi_1(\Tc_{s_0},t_0)]$ for any commutative ring $A$.
We denote basis elements by $l\in \pi_1(\Tc_{s_0},t_0)$ and write elements as $\sum_lf(l)l$, where $f(l)\in A$ and $f(l)=0$ for almost all $l$.
The group $\pi_1(\Tc,t_0)=\pi_1(\Tc_{s_0},t_0)\rtimes \pi_1(\Sc,s_0)$ acts naturally on 
$\pi_1(\Tc_{s_0},t_0)$ by
$(l,g)* v=l+g\cdot v$, with $l,v\in \pi_1(\Tc_{s_0},t_0)$ and $g\in \pi_1(\Sc,s_0)$. We can define this action directly in $\pi_1(\Tc,t_0)$ by
$x* v:=xv0_*\pi_*(x^{-1})$, $x\in \pi_1(\Tc,t_0)$ and $v\in \pi_1(\Tc_{s_0},t_0)$. This action induces a natural $\pi_1(\Tc,t_0)$-action on the group $A[\pi_1(\Tc_{s_0},t_0)]$. The group algebra comes with the natural augmentation map
\begin{equation*}
 \aug:A[\pi_1(\Tc_{s_0},t_0)]\rightarrow A,\ \sum_{l}f(l)l\mapsto \sum_l f(l).
\end{equation*}
The augmentation map is $\pi_1(\Tc,t_0)$-equivariant and the kernel $\mathfrak a=\ker(\aug)$ is again a $\pi_1(\Tc,t_0)$-module. In this manner we get for $n\in \N_0$ the $\pi_1(\Tc,t_0)$-modules $A[\pi_1(\Tc_{s_0},t_0)]/\mathfrak{a}^{n+1}$ and $A[[\pi_1(\Tc_{s_0},t_0)]]:=\varprojlim A[\pi_1(\Tc_{s_0},t_0)]/\mathfrak{a}^{n+1}$.
\begin{definition}\label{def_Log}
Let $\Tc$ be a torus over $\Sc$. We define the locally constant sheaf $Log_{\Tc/\Sc} ^{n}$ on $\Tc$ as the sheaf associated to the 
$\pi_1(\Tc,t_0)$-module
$A[\pi_1(\Tc_{s_0},t_0)]/\mathfrak{a}^{n+1}$, $n\in \N_0$, and call it the \textit{$n$-th logarithm sheaf}.
$\varprojlim Log_{\Tc/\Sc} ^{n}=:Log_{\Tc/\Sc}$ is called the \textit{logarithm sheaf}.
\end{definition}
\begin{remark}
$Log_{\Tc/\Sc}$ is also locally constant and is associated to the $\pi_1(\Tc,t_0)$- module $A[[\pi_1(\Tc_{s_0},t_0)]]$. Moreover, we have the augmentation map $\aug:Log_{\Tc/\Sc}\to A$ coming from the $\pi_1(\Tc,t_0)$-equivariant augmentation map $\aug:A[[\pi_1(\Tc_{s_0},t_0)]]\to A$.
\end{remark}
\begin{lemma}\label{basechange_Log}
$Log^n_{\Tc/\Sc}$ and $Log_{\Tc/\Sc}$ are natural in families of topological tori $\pi:\Tc\rightarrow \Sc$. In other words, given a commutative square of (pointed) families of topological tori
\begin{equation*}
\begin{xy}
\xymatrix{
\Tc^\prime\ar[r]^p\ar[d]^{\pi^\prime} & \Tc \ar[d]^\pi \\
\Sc^\prime\ar[r]^q & \Sc
}
\end{xy}
\end{equation*} 
with $p\circ 0^\prime=0\circ q$,
we have natural morphisms $Log^{n}_{\Tc^\prime/\Sc^\prime}\rightarrow p^{-1}Log^{n}_{\Tc/\Sc} $ and $Log_{\Tc^\prime/\Sc^\prime}\rightarrow p^{-1}Log_{\Tc/\Sc} $. If the diagram above is Cartesian, these maps are isomorphisms.
\begin{proof}
We choose $s_0 ^\prime \in \Sc^\prime$ with $q(s_0 ^\prime)=s_0$ and set $t_0 ^\prime =0^\prime(s_0 ^\prime)$.
Given a locally constant sheaf $\Fc$ on $\Tc$ associated to a $\pi_1(\Tc,t_0)$-module $M$ we recognize $p^{-1}\Fc$ as the locally constant sheaf associated to the $\pi_1(\Tc^\prime,t_0 ^\prime)$-module $M$ where the action is induced by the map $p_*:\pi_1(\Tc^\prime,t_0 ^\prime)\rightarrow \pi_1(\Tc,t_0)$. The morphism $p_*:\pi_1(\Tc^\prime _{s_0 ^\prime},t_0 ^\prime)\rightarrow \pi_1(\Tc_{s_0},t_0)$ is a map of $\pi_1(\Tc^\prime,t_0 ^\prime)$-sets. To see this choose $x\in \pi_1(\Tc^\prime,t_0 ^\prime)$ and $v\in\pi_1(\Tc^\prime _{s_0 ^\prime},t_0 ^\prime)$ arbitrarily.
We have
\begin{equation*}
 p_*(x* v)=p_*(xv0^\prime _*\pi^\prime _*(x^{-1}))=p_*(x)p_*(v) p_* 0^\prime _*\pi^\prime _*(x^{-1})=
\end{equation*}
\begin{equation*} 
p_*(x)p_*(v)0 _*\pi _*p_*(x^{-1})=p_*(x)* p_*(v)
\end{equation*} 
So $p_*:A[\pi_1(\Tc^\prime _{s^\prime _0},t^\prime _0)]\rightarrow A[ \pi_1(\Tc_{s_0},t_0)]$ is $\pi_1(\Tc^\prime, t_0 ^\prime)$-equivariant. This induces the desired maps on the level of sheaves. If our diagram is Cartesian the morphism $p_*:\pi_1(\Tc^\prime _{s^\prime _0},t^\prime _0)\rightarrow \pi_1(\Tc_{s_0},t_0)$ is an isomorphism inducing the desired isomorphisms on our sheaves. 
\end{proof}
\end{lemma}

\section{Limits of locally constant sheaves}

\subsection{Limits over the natural numbers}

In the case of sheaves there can be non-trivial $R^p\varprojlim$ for $p\geq2$ and there may be no such thing as a Mittag-Leffler condition on the system ensuring the vanishing of higher limits. We want to prove now that locally constant sheaves have, as expected, the same behavior with limits as abelian groups.
Let $\Ac$ be an abelian category. We denote by $\Ac^\N$ the category of all inverse systems in $\Ac$ indexed by the natural numbers. Objects of $\Ac^\N$ are families $(A_n,a_n)$ with objects $A_n$ in $\Ac$ and maps $a_n:A_{n+1}\rightarrow A_n$. Morphisms between $(A_n,a_n)$ and $(B_n,b_n)$ in $\Ac^\N$ are just commutative diagrams in $\Ac$ of the form
\begin{equation*}
\begin{xy}
\xymatrix{
\cdots & A_{n+1}\ar[r]\ar[d]^{f_{n+1}} & A_n \ar[r]\ar[d]^{f_n} & A_{n-1}\ar[r]\ar[d]^{f_{n-1}} & \cdots \\
\cdots & B_{n+1}\ar[r] & B_n \ar[r] & B_{n-1}\ar[r] & \cdots \\
}
\end{xy}
\end{equation*}
$\Ac^\N$ is an abelian category with kernels and cokernels defined componentwise.
\begin{proposition}
$\Ac^\N$ has enough injective objects, if and only if $\Ac$ has. We can characterize injectives in $\Ac^\N$ as follows:\\
$\Ic=(\Ic_n,d_n)$ is injective, if and only if each $\Ic_n$ is injective in $\Ac$ and all the $d_n$ are split epimorphisms.
\begin{proof}
\cite{Ja} (1.1) proposition.
\end{proof}
\end{proposition}
From now on we want to assume that $\Ac$ has enough injectives and that inverse limits indexed over $\N$ exist in $\Ac$. The functor $\varprojlim:\Ac^\N\rightarrow \Ac$ sending an inverse system to its limit is additive and left exact as right adjoint to the exact diagonal functor $\Delta$ sending $A$ to the system
$\cdots \rightarrow A \stackrel{\id}{\rightarrow}A \stackrel{\id}{\rightarrow}A \stackrel{\id}{\rightarrow}\cdots$. In particular $\varprojlim$ preserves injective objects.
Since $\Ac^\N$ has enough injectives and $\varprojlim$ is left exact, we form the higher right derived functors $R^p\varprojlim$. We will often omit the transition maps of our system and simply write $(A_n)$ for objects in $\Ac^\N$.
Let us consider $\Ac=Ab$ the category of abelian groups first.
\begin{definition}
$(A_n,d_n)$ in $Ab^\N$ satisfies the Mittag-Leffler condition (M-L), if for any $n\in \N$ the decreasing filtration of $A_n$ $F_m:=\image(A_{n+m}\rightarrow A_n)$ induced by the $d_{k}$, $k\geq n$, becomes stationary.
\end{definition}
\begin{remark}
If all $d_n$ are surjective, then $(A_n,d_n)$ obviously satisfies (M-L).
\end{remark}
\begin{proposition}
Let $(A_n)$ be in $Ab^\N$. Then $R^p\varprojlim(A_n)=0$ for $p\geq2$. If $(A_n)$ satisfies (M-L), then $R^1\varprojlim (A_n)=0$.
\begin{proof}
\cite{We} Corollary 3.5.4. and Proposition 3.5.7. 
\end{proof}
\end{proposition}

\subsection{(M-L) systems of locally constant sheaves}

Let $X$ be a topological manifold. The category $Sh(X)$ of sheaves on $X$ has enough injectives and all limits over the natural numbers exist. For $(\Fc_n)\in Sh(X)^\N$ the formula 
\begin{equation*}
(\varprojlim \Fc_n)(U)=\varprojlim \Fc_n(U),\ U\subset X\ \text{open}
\end{equation*}
can be used as definition.\\
\begin{definition}
We say $(\Fc_n)\in Sh(X)^\N$ is a \textit{system of locally constant sheaves on $X$}, if each $\Fc_n$ is a locally constant sheaf. We say that the system of locally constant sheaves $(\Fc_n)$ satisfies the \textit{Mittag-Leffler condition} (M-L), if for any $x\in X$ the induced system of abelian groups $(\Fc_{n,x})$ does so. 
\end{definition}
\begin{proposition}
Let $(\Fc_n)$ be a system of locally constant sheaves on a topological manifold $X$. Then $R^p\varprojlim(\Fc_n)=0$ for $p\geq 2$ and if $(\Fc_n)$ satisfies (M-L) $R^1\varprojlim (\Fc_n)=0$.
\begin{proof}
Consider an injective resolution $(\Fc_n)\hookrightarrow (\Ic_n)^\bullet$ in $Sh(X)^\N$. We have by definition $R^p{\varprojlim}(\Fc_n)=H^p(\varprojlim \Ic_n ^\bullet)$. $H^p(\varprojlim \Ic_n ^\bullet)$ is the sheaf associated to the presheaf
\begin{equation*}
U\mapsto H^p((\varprojlim \Ic_n )(U)^\bullet)=H^p(\varprojlim \Ic_n (U)^\bullet).
\end{equation*}
So it suffices to show that $H^p(\varprojlim \Ic_n (U)^\bullet)$ vanishes for any $U\subset X$ open and contractible for $p\geq2$ or $p\geq 1$, if $(\Fc_n)$ is (M-L).
The complexes $\Ic_n (U)^\bullet$ are injective resolutions of $M_n:=\Fc_n(U)$.
To see this recall that $\Ic^\bullet_{n|U}$ is an injective resolution of the \textit{constant} sheaf $\underline{M_n}=\Fc_{n|U}$. Here we use that locally constant sheaves are constant on simply connected manifolds (\cite{Iv} IV.9). So we may compute
\begin{equation*}
H^p(U,M_n)=H^p(U,\Fc_{n|U})=H^p(\Ic^\bullet_{n|U}(U))=H^p(\Ic^\bullet_n(U))
\end{equation*}
Now $U$ is contractible and by homotopy invariance of sheaf cohomology with constant coefficients (\cite{Iv} IV.1) we get $H^p(U,M_n)=H^p(\left\{pt.\right\},M_n)$ which is zero for $p\geq 1$ and $M_n$ for $p=0$. So $\Ic_n (U)^\bullet$ is a resolution of $M_n$ and it is an injective resolution, because the section functor $\Gamma(U,\ )$ is right adjoint to the constant sheaf construction which is exact. Since $(\Ic^p _n)$ is injective in $Sh(X)^\N$, all its transition functions are split epimorphisms. In particular, all transition functions of $(\Ic^p _n(U))$ are split epimorphisms. This implies that $(\Ic^p _n(U))$ is injective in $Ab^\N$. Therefore $(M_n)\hookrightarrow (\Ic_n (U))^\bullet$ is an injective resolution in $Ab^\N$ implying $R^p\varprojlim(M_n)=H^p(\varprojlim \Ic_n (U)^\bullet)$. Now we can use the results about limits in $Ab$ to conclude that $R^p\varprojlim(\Fc_n)=0$ for $p\geq 2$ and if $(\Fc_n)$ satisfies (M-L) $R^1\varprojlim(\Fc_n)=0$.
\end{proof}
\end{proposition}
\begin{remark}\label{lim_injective_resolution}
If $(\Fc_n)$ is a (M-L) system of locally constant sheaves and $(\Fc_n)\hookrightarrow (\Ic_n)^\bullet$ is an injective resolution in $Sh(X)^\N$, then $\varprojlim \Fc_n\hookrightarrow \varprojlim \Ic_n^\bullet$ is an injective resolution. This follows easily from $H^p(\varprojlim \Ic_n^\bullet)=R^p\varprojlim(\Fc_n)=0$, if $p\geq 1$.
\end{remark}

\section{Purity and localization}

In this section we want to prove a purity result for (M-L) systems of local systems on topological manifolds. This gives rise to localization sequences on cohomology, which finally will be used to derive a very easy localization sequence for the logarithm sheaf and to define the polylogarithm cohomology classes.

\subsection{Purity for local systems}

Let $A$ be a noetherian ring and $X$ a topological space.
\begin{definition}
A sheaf $\Fc$ on $X$ is called a \textit{local $A$-system}, or simply a \textit{local system}, if it is locally constant and if for all $x\in X$ the stalks $\Fc_x$ are free $A$-modules of finite rank.
The sheaf $\Fc^*:=\underline{Hom}_A(\Fc,A)$ is also a local $A$-system and is called the \textit{dual of $\Fc$}.
\end{definition}
\begin{lemma}\label{injectives_local_systems}
If $\Ic$ is an injective $A$-sheaf on a topological space $X$, then $\Fc\otimes_A\Ic$ is injective for any local system $\Fc$.
\begin{proof}
Given any inclusion $C\hookrightarrow D$ of $A$-sheaves we have a commutative diagram where the vertical arrows are isomorphisms:
\begin{equation*}
\begin{xy}
\xymatrix{
Hom(C,\Ic\otimes_A \Fc)\ar[r]\ar[d] & Hom(D,\Ic\otimes_A \Fc)\ar[d]\\
Hom(C\otimes_A \Fc^*,\Ic)\ar[r] & Hom(D\otimes_A \Fc^*, \Ic)
}
\end{xy}
\end{equation*}
To see this one uses the universal property of the tensor product, $(\Fc^*)^*=\Fc$ and $\Fc^*\otimes_A \Ic=\underline{Hom}_A(\Fc,\Ic)$.
Since $\Fc^*$ is a flat $A$-sheaf $C\otimes_A \Fc^*\hookrightarrow D\otimes_A \Fc^*$ is injective. Therefore the lower horizontal arrow is surjective, because $\Ic$ is injective. So $\Fc\otimes_A \Ic$ is injective.
\end{proof}
\end{lemma}
Let $Z\subset X$ be a closed subspace. The pushforward $i_*:Sh(Z,A)\rightarrow Sh(X,A)$ on the categories of sheaves of $A$-modules has a right adjoint $i^{!}:Sh(X,A)\rightarrow Sh(Z,A)$ (\cite{Iv} II. Proposition 6.6.).
\begin{lemma}\label{!and-1}
Let $\Fc$ be a local $A$-system on a topological space $X$ and let $\Gc$ be any sheaf on $X$. There is a functorial isomorphism
\begin{equation*}
i^{!}(\Fc\otimes_A \Gc)=i^{-1}\Fc\otimes_A i^{!}\Gc.
\end{equation*}
\begin{proof}
We use the Yoneda principle: Let $\Hc$ be any $A$-sheaf on $Z$. By adjunction and the properties of $^*$ and $\otimes_A$ we have isomorphisms 
\begin{equation*}
Hom(\Hc,i^!\Gc\otimes_A i^{-1}\Fc)\cong Hom(\Hc\otimes_A i^{-1}\Fc^*,i^!\Gc)\cong Hom(i_*(\Hc\otimes_A i^{-1}\Fc^*),\Gc)
\end{equation*}
Moreover,
$i_*(\Hc\otimes_A i^{-1}\Fc^*)\cong i_*\Hc\otimes_A i_*i^{-1}\Fc^*\cong i_*\Hc\otimes_A \Fc^*$ where the last arrow comes from the adjunction morphism $\Fc\rightarrow i_*i^{-1}\Fc^*$. Therefore we get natural in $\Hc$
\begin{equation*}
Hom(\Hc,i^!\Gc\otimes_A i^{-1}\Fc)\cong Hom(i_*\Hc\otimes_A \Fc^*,\Gc)\cong  Hom(i_*\Hc,\Gc\otimes_A \Fc)\cong 
\end{equation*}
\begin{equation*}
Hom(\Hc,i^!(\Gc\otimes_A \Fc))
\end{equation*}
\end{proof}
\end{lemma}
As $i_*$ is exact, $i^{!}$ is left exact. Purity is concerned with the calculation of the higher right derived functors of $i^{!}$. To do so we need some notation. \\
Let $X$ be a topological manifold and $A$ a noetherian ring.
\begin{itemize}
\item{ $D^+(X,A)$ denotes the derived category of bounded below complexes of sheaves of $A$-modules. We identify this category with the homotopy category of bounded below complexes of injective sheaves of $A$-modules.}
\item{ $D^+(A)$ is the derived category of bounded below complexes of $A$-modules. We identify this category with the homotopy category of bounded below complexes of injective $A$-modules.}
\item If $C=(C,d)$ is a complex of sheaves and $p\in \Z$, we denote by $C[p]$ the shifted complex, in other words, $C[p]^d=C^{p+d}$ with differential $d[p]:=(-1)^pd$. If $D$ is a sheaf, we identify $D$ with the complex which equals $D$ in degree zero and is zero elsewhere.
\end{itemize}
\begin{definition}
Let $X$ be a topological manifold of dimension $d$ and $A$ a noetherian ring.
The presheaf
\begin{equation*}
U\subset X \text{ open}\mapsto Hom_A(H^d _c(U,A),A)=:or_X(U)
\end{equation*}
is a sheaf on $X$. It is even a local $A$-system of rank 1 (\cite{Iv} III 8.14.) and we call it the \textit{orientation bundle of $X$}
\end{definition}
\begin{definition}
Let $X$ be a topological manifold of dimension $d$. The \textit{dualizing complex} $\Dc_X$ on $X$ is a complex of injective sheaves on $X$ characterized up to homotopy by the following property (\cite{Iv}VI.2):
There is an isomorphism natural in $A\in D^+(X,A)$: $Hom_{D^+(X,A)}(A,\Dc_X)=Hom_{D^+(A)}(\Gamma_c(X,A),A)$. Moreover, one has a quasi-isomorphism $or_X[d]\rightarrow \Dc_X$ (\cite{Iv} VI.3.2).
\end{definition}

\begin{lemma}\label{dualizing_purity}
Let $X$ be a topological manifold of dimension $d$ and $i:Z\hookrightarrow X$ a closed submanifold of codimension $c$. Then one has $i^{!}\Dc_X=\Dc_Z$.
\begin{proof}
\cite{Iv}VIII proposition 1.7.
\end{proof}
\end{lemma} 
\begin{definition}
Let $X$ be a topological manifold of dimension $d$ and $i:Z\hookrightarrow X$ a closed submanifold of codimension $c$. We define the \textit{orientation of $Z$ relative $X$} as the rank $1$ local system $or_{Z/X}:=or_Z\otimes_A i^{-1}or_X^*$. 
\end{definition}
\begin{proposition}\label{purity1}
Let $X$ be a topological manifold of dimension $d$, $i:Z\hookrightarrow X$ a closed submanifold of codimension $c$ and $\Fc$ a local system on $X$. Then we have $Ri^{!}\Fc=or_{Z/X}\otimes_A i^{-1}\Fc[-c]$ in $D(Z,A)$ functorial in $\Fc$.
\begin{proof}
We start with the quasi-isomorphism $or_X[d]\rightarrow \Dc_X$. $or_X ^*\otimes_A\Fc$ is flat and if we tensor the last quasi-isomorphism with it, we get a quasi-isomorphism $\Fc[d]\rightarrow \Dc_X\otimes_A or_X ^*\otimes_A\Fc$. The right-hand side is a complex of injectives by \Cref{injectives_local_systems}.
We get
\begin{equation*}
Ri^{!}\Fc[d]= i^{!}(\Dc_X\otimes_A or_X ^*\otimes_A\Fc)=i^{!}\Dc_X\otimes_A i^{-1}or_X^*\otimes_A i^{-1}\Fc 
\end{equation*}
by \Cref{!and-1} and this equals by \Cref{dualizing_purity}
\begin{equation*}
\Dc_Z\otimes_A i^{-1}or_X^*\otimes_A i^{-1}\Fc=or_Z [d-c]\otimes_A i^{-1}or_X^*\otimes_A i^{-1}\Fc=or_{Z/X}\otimes_A i^{-1}\Fc[d-c]
\end{equation*}
\end{proof}
\end{proposition}

\subsection{Purity for systems of local systems and naturality}

\begin{proposition}\label{lim_purity}
Let $X$ be a topological manifold of dimension $d$ and $i:Z\hookrightarrow X$ a closed submanifold of codimension $c$.
Let $(\Fc_n)$ be a system of local $A$-systems satisfying (M-L) on $X$. We have in $D^+(X,A)$
\begin{equation*}
R i^{!}(\varprojlim\Fc_n)=\varprojlim(i^{-1}\Fc_n\otimes_A or_{Z/X})[-c] =i^{-1}(\varprojlim\Fc_n)\otimes_A or_{Z/X}[-c]
\end{equation*}
functorial in $(\Fc_n)$.
\begin{proof}
Take an injective resolution $a:(\Fc_n)\hookrightarrow (\Ic_n)^\bullet$ in $Sh(X,A)^\N$. We already had the resolution $b:(\Fc_n)\hookrightarrow (\Dc_X\otimes_A or_X ^*\otimes_A\Fc_n)_n$ in $Sh(X,A)^\N$ functorial in systems $(\Fc_n)$. Since $(\Ic_n)^\bullet$ is a complex of injectives in $Sh(X,A)^\N$, we get a map of complexes $c$ unique up to homotopy making the following diagram commute
\begin{equation*}
\begin{xy}
\xymatrix{
 (\Fc_n)\ar[rd]^b\ar[rr]^{a} &  & (\Ic_n)^\bullet \\
                           & (\Dc_X\otimes_A or_X ^*\otimes_A\Fc_n)\ar[ur]^c &
                      }
\end{xy}
\end{equation*}
(\cite{Iv} I Theorem 6.2). Let us consider the functor $i^{!\N}:Sh(X,A)^\N\to Sh(Z,A)^\N$, $(\Gc_n)\mapsto (i^!\Gc_n)$. We have now
\begin{equation*}
Ri^{!\N}(\Fc_n)=(i^!\Ic_n)^\bullet\stackrel{i^{!\N}c}{\leftarrow}(i^!(\Dc_X\otimes_A or_X ^*\otimes_A\Fc_n))\stackrel{p\cong}{\rightarrow}(i^{-1}\Fc_n\otimes_A or_{Z/X})[-c],
\end{equation*}
where the map on the right-hand side is the functorial purity quasi-isomorphism $p$ from \cref{purity1}. The morphism $i^{!\N}c$ is also a quasi-isomorphism, as
\begin{equation*}
H^p(i^!\Ic_n)^\bullet=(H^pi^!\Ic_n^\bullet)=(R^pi^!\Fc_n)=H^p(i^!(\Dc_X\otimes_A or_X ^*\otimes_A\Fc_n)).
\end{equation*}
This gives a purity isomorphism 
\begin{equation*}
Ri^{!\N}(\Fc_n)=(i^{-1}\Fc_n\otimes_A or_{Z/X})[-c]\in D^+(Sh(X,A)^\N)
\end{equation*}
in the derived category of bounded below complexes of systems of sheaves of $A$-modules, which is functorial in systems of local systems $(\Fc_n)$. As $(\Fc_n)$ is (M-L), we have that $\varprojlim\Fc_n\hookrightarrow \varprojlim \Ic_n^\bullet$ is also an injective resolution by \Cref{lim_injective_resolution}. In other words, $R\varprojlim(\Fc_n)=\varprojlim\Fc_n\in D^+(X,A)$ and one gets in $D^+(X,A)$
\begin{equation*}
Ri^{!}\varprojlim\Fc_n=Ri^!R\varprojlim (\Fc_n)=R(i^!\circ \varprojlim)(\Fc_n)=R(\varprojlim\circ i^{!\N})(\Fc_n),
\end{equation*}
since $i^!$ commutes as right adjoint with the limit. We go on with
\begin{equation*}
R\varprojlim Ri^{!\N}(\Fc_n)=R\varprojlim (i^{-1}\Fc_n\otimes_A or_{Z/X})[-c]=\varprojlim (i^{-1}\Fc_n\otimes_A or_{Z/X})[-c],
\end{equation*}
where we used purity and the fact that $(i^{-1}\Fc_n\otimes_A or_{Z/X})$ is again (M-L). Finally, we have
\begin{equation*}
\varprojlim (i^{-1}\Fc_n\otimes_A or_{Z/X})=i^{-1}(\varprojlim\Fc_n)\otimes_A or_{Z/X}
\end{equation*}
which will follow from the following two lemmas and will complete the proof.
\end{proof}
\end{proposition}
\begin{lemma}
Let $X$ be a topological manifold and let $(\Fc_n)$ be a system of locally constant sheaves on $X$. For any $x\in X$ the canonical map
\begin{equation*}
(\varprojlim \Fc_n)_x\rightarrow\varprojlim(F_{n,x})
\end{equation*}
is an isomorphism.
More generally, if $Y$ is another topological manifold and $f:Y\rightarrow X$ a continuous map, the canonical map $f^{-1}(\varprojlim \Fc_n)\rightarrow \varprojlim f^{-1}\Fc_n$ is an isomorphism.
\begin{proof}
Let $\Gc$ be any sheaf on $X$ and $x\in X$. As $X$ is a topological manifold any open neighborhood $U$ of $x$ contains a contractible open neighborhood $B$ of $x$. It follows that the family $\mathfrak B _x$ of contractible open neighborhoods of $x$ is a cofinal inductive system in the family $\mathfrak U _x$ of all open neighborhoods of $x$. Therefore
$\Gc_x=\varinjlim_{U\in \mathfrak U _x}\Gc(U)=\varinjlim_{B\in \mathfrak B _x}\Gc(B)$.
If $\Fc$ is locally constant, $\Fc_{|B}$ is constant for $B\in \mathfrak B_x$, as contractible spaces are simply connected, and we have $\Fc_x=\varinjlim_{B\in \mathfrak B _x}\Fc(B)=\Fc(B)$ for any $B\in \mathfrak B _x$.
We calculate
\begin{equation*}
 (\varprojlim \Fc_n)_x=\varinjlim_{B\in \mathfrak B _x}(\varprojlim  \Fc_n)(B)=\varinjlim_{B\in \mathfrak B _x}(\varprojlim \Fc_n(B))
\end{equation*}
as the limit commutes with the section functor. The inductive limit becomes constant and therefore
\begin{equation*}
\varinjlim_{B\in \mathfrak B _x}(\varprojlim \Fc_n(B))=\varprojlim \Fc_n(B)=\varprojlim \Fc_{n,x}.
\end{equation*}
This proves the first claim. The second follows easily from the first by considering the induced maps on the stalks.
\end{proof}
\end{lemma}
\begin{lemma}
Let $(\Fc_n)$ be a system of locally constant $A$-sheaves on a topological manifold $X$ and $\Fc$ a local $A$-system on $X$. Then the natural map $(\varprojlim \Fc_n)\otimes_A \Fc\rightarrow \varprojlim(\Fc_n\otimes_A \Fc)$ is an isomorphism.
\begin{proof}
By the preceding lemma the problem is local so it suffices to proof the corresponding statement for a system of $A$-modules $(F_n)$ and a finitely generated free module $F$. For any $A$-module $M$ one has the functorial isomorphisms
\begin{equation*}
Hom(M,\varprojlim (F_n\otimes_A F)
)= \varprojlim Hom(M,F_n\otimes_A F)=\varprojlim Hom(M\otimes_A F^*,F_n)=
\end{equation*}	
\begin{equation*}
Hom(M\otimes_A F^*,\varprojlim F_n)=Hom(M,(\varprojlim F_n)\otimes_A F) 
\end{equation*}	
from which the claim follows by the Yoneda lemma.
\end{proof}
\end{lemma}
Now we come to our important application of the purity result
\begin{proposition}\label{localization_triangle}
Let $(\Fc_n)$ be a system of local systems on a topological manifold $X$ satisfying (M-L). Let $i:Z\hookrightarrow X$ be the inclusion of a closed submanifold of codimension $c$ and $j:U:=X\setminus Z\to X$ the inclusion of the open complement. One has an exact triangle in $D^+(X,A)$
\begin{equation*}
i_*(i^{-1}(\varprojlim\Fc_n)\otimes_A or_{Z/X})[-c]\rightarrow\varprojlim\Fc_n\rightarrow Rj_*j^{-1}\varprojlim\Fc_n\stackrel{+1}{\rightarrow} 
\end{equation*}
\begin{proof}
Let $\Fc$ be any sheaf on $X$. Let us recall 
\begin{equation*}
i_*i^!\Fc(V)=\left\{s\in \Fc(V):\text{supp}(s)\subset Z\right\}
\end{equation*}
(\cite{Iv} II.6.6). Therefore, we have an exact sequence
\begin{equation*}
0\rightarrow i_*i^!\Fc\rightarrow \Fc\rightarrow j_* j^{-1}\Fc,
\end{equation*}
where the arrows are the adjunction morphisms. If $\Fc$ is flabby, the last sequence is even right exact.
This means, if we replace $\Fc$ by an injective resolution $\Ic$, we get an exact sequence
\begin{equation*}
0\rightarrow i_*i^!\Ic\rightarrow \Ic\rightarrow j_* j^{-1}\Ic\rightarrow 0
\end{equation*}
giving rise to an exact triangle in the derived category $D^+(X,A)$
\begin{equation*}
i_*Ri^!\Fc\rightarrow \Fc\rightarrow Rj_* j^{-1}\Fc\stackrel{+1}{\rightarrow}.
\end{equation*}
If we apply now \Cref{lim_purity} to substitute $i_*Ri^!\Fc$, we get the exact triangle as claimed.
\end{proof}
\end{proposition}
\begin{remark}\label{natural_localization}
\Cref{localization_triangle} is functorial with respect to pairs $(X,Z)$, $(X^\prime,Z^\prime)$. By this we mean given a Cartesian square of topological spaces
\begin{equation*}
\begin{xy}
\xymatrix{
Z^\prime\ar[r]^{i^\prime}\ar[d]^{f^\prime} & X^\prime\ar[d]^f \\
Z\ar[r]^{i}& X
}
\end{xy}
\end{equation*}
with a topological manifold $X$, a closed submanifold $Z$ of codimension $c$ and topological submersion $f$ (see \cite{KS} Definition 3.3.1) we get by pullback a morphism of exact triangles (apply $\id\rightarrow Rf_*f^{-1}$) from
\begin{equation*}
i_*(i^{-1}(\varprojlim\Fc_n)\otimes_A or_{Z/X})[-c]\rightarrow \varprojlim\Fc_n\rightarrow Rj_*j^{-1}\varprojlim\Fc_n\stackrel{+1}{\rightarrow}
\end{equation*}
to
\begin{equation*}
Rf^{\prime}_*i^\prime_*(i^{\prime -1}(\varprojlim f^{-1}\Fc_n)\otimes_A or_{Z^\prime/X^\prime})[-c]\rightarrow Rf_*\varprojlim f^{-1}\Fc_n
\rightarrow Rf_*Rj^\prime _* j^{\prime -1}\varprojlim f^{-1}\Fc_n\stackrel{+1}{\rightarrow}
\end{equation*}
To see this we recall that $X^\prime$ and $Z^\prime$ have the structure of topological manifolds and that $Z^\prime$ has codimension $c$ in $X^\prime$. After exploiting the commutation rules of the functors involved the only thing that remains to be shown is that there is a natural pullback isomorphism $f^{\prime -1}or_{Z/X}\rightarrow or_{Z^\prime/X^\prime}$. As our diagram is Cartesian, we have a natural pullback morphism $i^{\prime -1}f^!A\rightarrow f^{\prime !} i^{-1}A$ by \cite{KS} proposition 3.1.9 (iii). If we apply $H^{dim(X) -dim(X^\prime)}$, we get the natural morphism $i^{\prime -1}or_{X^\prime/X}\rightarrow or_{Z^\prime/Z}$. Following the proof of \cite{KS}  proposition 3.3.2 we may see that it is actually an isomorphism. By the formalism provided by \cite{KS} remark 3.3.5 we may conclude that we also have a natural isomorphism $f^{\prime -1}or_{Z/X}\rightarrow or_{Z^\prime/X^\prime}$. 

\end{remark}
\begin{remark}\label{shrinking_localization}
If $(X,Z)$, $(X,Z^\prime)$ are two pairs of codimension $c$ as above with $Z^\prime\subset Z$, we have another compatibility. $Z^\prime \subset Z$ is also open, as both manifolds have the same dimension. From this we derive a natural extension by zero morphism
\begin{equation*}
i^{\prime }_*(i^{\prime -1}\varprojlim\Fc_n\otimes_A or_{Z^\prime/X})\rightarrow i_*(i^{ -1}\varprojlim\Fc_n\otimes_A or_{Z/X}) ,
\end{equation*}
with the inclusions $i:Z\rightarrow X$ and $i^\prime:Z^\prime\rightarrow X$. If we denote by $j:X\setminus Z\rightarrow X$ and $j^\prime:X\setminus Z^\prime\rightarrow X$ the inclusions of the open complements, the morphism above fits into a morphism of triangles
\begin{equation*}
\begin{xy}
\xymatrix{
i^\prime_*(i^{\prime-1}(\varprojlim\Fc_n)\otimes_A or_{Z^\prime/X})[-c] \ar[r]\ar[d] & \varprojlim\Fc_n\ar[r]\ar[d] & Rj^\prime_*j^{\prime-1}\varprojlim\Fc_n\ar[d]\ar[r]^>{+1} &\\
i_*(i^{-1}(\varprojlim\Fc_n)\otimes_A or_{Z/X})[-c]\ar[r] & \varprojlim\Fc_n\ar[r] & Rj_* j^{ -1}\varprojlim\Fc_n\ar[r]^>{+1}&
}
\end{xy} 
\end{equation*}
The right vertical arrow is just the restriction map.
\end{remark}
\begin{corollary}\label{localization}
Let $(\Fc_n)$ be a system of local systems on a topological manifold $X$ satisfying (M-L). Let $i:Z\hookrightarrow X$ be the inclusion of a closed submanifold of codimension $c$ with open complement $U:=X\setminus Z$ and $f:X\rightarrow Y$ a continuous map of topological spaces. There is a long exact cohomology sequence
\begin{equation*}
\rightarrow R^{p-1}f_{|U *}(\varprojlim\Fc_{n\ |U})\rightarrow R^{p-c}f_{|Z *}(\varprojlim i^{-1}\Fc_n\otimes_A or_{Z/X})\rightarrow R^pf_*(\varprojlim\Fc_n)
\end{equation*}
\begin{equation*}
\rightarrow R^p f_{|U*}(\varprojlim\Fc_{n\ |U})\rightarrow R^{p+1-c}f_{|Z *}(\varprojlim i^{-1}\Fc_n\otimes_A or_{Z/X})\rightarrow
\end{equation*}
called \textit{the localization sequence for the pair $(X,Z)$}. The sequence is natural in pairs $(X,Z)$ as described in \Cref{natural_localization} and \Cref{shrinking_localization}.
\begin{proof}
We take the exact triangle
\begin{equation*}
i_*(i^{-1}(\varprojlim\Fc_n)\otimes_A or_{Z/X})[-c]\rightarrow\varprojlim\Fc_n\rightarrow Rj_*j^{-1}\varprojlim\Fc_n\stackrel{+1}{\rightarrow}
\end{equation*}
from \Cref{localization_triangle} and apply $Rf_*$ to obtain an exact triangle
\begin{equation*}
Rf_*i_*(i^{-1}(\varprojlim\Fc_n)\otimes_A or_{Z/X})[-c]\rightarrow Rf_*\varprojlim\Fc_n\rightarrow Rf_*Rj_*j^{-1}\varprojlim\Fc_n \stackrel{+1}{\rightarrow},
\end{equation*}
which equals
\begin{equation*}
Rf_{|Z *}(i^{-1}(\varprojlim\Fc_n)\otimes_A or_{Z/X})[-c]\rightarrow Rf_*\varprojlim\Fc_n\rightarrow Rf_{|U*}j^{-1}\varprojlim\Fc_n\stackrel{+1}{\rightarrow}
\end{equation*}
by the composition rules of derived functors. If we apply the functor $H^0$, we obtain the statement. 
\end{proof}
\end{corollary}

\section{Construction of the topological polylogarithm }

\subsection{The localization sequence for the logarithm sheaf}

In this section we prove the important vanishing result for the cohomology of the logarithm sheaf. We will then derive a simple localization sequence for the logarithm sheaf, which enables us to define the polylogarithm cohomology classes.

From now on we will often drop the subscript $A$ at $\otimes$, when the coefficient ring has been fixed.
\begin{proposition}
Let $\pi:\Tc\rightarrow \Sc$ be torus of constant fiber dimension $d$ and $Log_{\Tc/\Sc}=Log$ the associated logarithm sheaf over the noetherian ring $A$. Then
\begin{equation*}
R^p\pi_* (Log)=0,\ p\neq d,\ \text{and}\ \aug:R^{d}\pi_* (Log)\rightarrow R^{d}\pi_* (A)
\end{equation*}
is an isomorphism.
\begin{proof}
As the problems are local on the base, we may check the corresponding statements for the stalks. Because $\pi$ is proper, we have for any sheaf $\Fc$ on $\Tc$ and $s\in \Sc$
\begin{equation*}
R^p\pi_* (\Fc)_s=H^p(\pi^{-1}(s),\Fc_{|\pi^{-1}(s)}).
\end{equation*}
Now $Log_{\Tc/\Sc,s}=Log_{\pi^{-1}(s)/\left\{s\right\}}$ by \Cref{basechange_Log} and we use \cite{BKL} Theorem 1.6.1 to conclude.
\end{proof}
\end{proposition}
\begin{remark}
If $\pi:\Tc\to \Sc$ is a topological submersion, there is a relative orientation bundle $or_{\Tc/\Sc}$ on $\Tc$, see \cite{KS} Definition 3.3.3.
If $\pi:\Tc\rightarrow \Sc$ is a torus, it is in particular a fiber bundle with orientable fibers and therefore $or_{\Tc/\Sc}=\pi^{-1}R^{d}\pi_*(A)^*$.
\end{remark}
Now let $D\subset \Tc$ be closed such that $\pi_{|D}:D\rightarrow \Sc$ is a (finite) covering. We set $U:=\Tc\setminus D$. One has $or_{D}=\pi_{|D}^{-1}or_{\Sc}$.
\begin{proposition}\label{relative_Log_localization}
Let $\Fc$ be a local $A$-system on $\Sc$, $\pi:\Tc\rightarrow \Sc$ a family of topological tori of constant fiber dimension $d$ and $D\subset \Tc$ be closed such that $\pi_{|D}$ is a covering. The localization sequence for the twisted logarithm sheaf $\pi^{-1}\Fc\otimes Log$ reduces to a short exact sequence
\begin{equation*}
0\rightarrow R^{d-1}\pi_{|U*}(\pi^{-1}\Fc\otimes Log\otimes or_{\Tc/\Sc})\rightarrow \pi_{|D*}(\pi^{-1}\Fc\otimes Log_{|D})\stackrel{\aug}{\rightarrow}\Fc\rightarrow 0
\end{equation*}
\begin{proof}
Because $\pi_{|D}:D\rightarrow \Sc$ is a covering, one has for any locally constant sheaf $\Gc$ on $D$ that $R^p\pi_{|D*}(\Gc)=0$, $p\neq 0$. To see this it suffices to show $H^p(\pi_{|D*} ^{-1}(U),\Gc)=0$, $p\neq 0$, for all small $U\subset \Sc$. Take $U$ contractible and uniformly covered, in other words $\pi_{|D*}^{-1}(U)=\coprod _{i\in I} U_i$ with $\pi_{|D|U_i}: U_i\rightarrow U$ a homeomorphism for all $i\in I$. The sheaf $\Gc_{|U_i}$ is constant, as locally constant on something simply connected, and we calculate
\begin{equation*}
H^p(\pi_{|D*}^{-1}(U),\Gc)=H^{p}(\coprod_{i\in I} U_i,\Gc)=\prod_{i\in I} H^{p}(U_i,\Gc)=0,\ p\neq0,
\end{equation*}
by homotopy invariance of sheaf cohomology with constant coefficients (\cite{Iv} IV.1).
The localization sequence for $(\pi^{-1}\Fc\otimes Log^n)_n$
\begin{equation*}
\rightarrow R^{p-1}\pi_{|U *}(\pi^{-1}\Fc \otimes Log\otimes or_{\Tc/\Sc})\rightarrow R^{p-d}\pi_{|D *}(\pi_{|D}^{-1}\Fc \otimes Log_{|D})
\end{equation*}
\begin{equation*}
 \rightarrow R^p\pi_*(\pi^{-1}\Fc \otimes Log\otimes or_{\Tc/\Sc})\rightarrow R^p \pi_{|U*}(\pi^{-1}\Fc \otimes Log\otimes or_{\Tc/\Sc})\rightarrow 
\end{equation*}
provided by \Cref{localization} gives immediately isomorphisms
\begin{equation*}
R^p \pi_{*}(\pi^{-1}\Fc \otimes Log\otimes or_{\Tc/\Sc})\rightarrow R^p \pi_{|U*}(\pi^{-1}\Fc \otimes Log\otimes or_{\Tc/\Sc})
\end{equation*}
for $p<d-1$ and $p>d$. We extract the exact sequence
\begin{equation*}
R^{d-1}\pi_{*}(\pi^{-1}\Fc \otimes Log\otimes or_{\Tc/\Sc})\rightarrow R^{d-1}\pi_{|U *}(\pi^{-1}\Fc \otimes Log\otimes or_{\Tc/\Sc})
\end{equation*}
\begin{equation*}
\rightarrow \pi_{|D *}(\pi^{-1}\Fc \otimes Log_{|D})\rightarrow R^{d}\pi_{*}(\pi^{-1}\Fc \otimes Log\otimes or_{\Tc/\Sc})
\end{equation*}
from our localization sequence.
We have
\begin{equation*}
R^p\pi_*(\pi^{-1}\Fc \otimes Log\otimes or_{\Tc/\Sc})=\Fc\otimes R^p\pi_*(Log)\otimes R^{d}\pi_* (A)^*=0
\end{equation*}
for $p\neq d$ and $R^{d}\pi_*(\pi^{-1}\Fc\otimes Log\otimes or_{\Tc/\Sc})=\Fc$ by the projection formula, since $\pi$ is proper and $\Fc$ is flat.
So we have the short exact sequence
\begin{equation*}
0 \rightarrow R^{d-1}\pi_{|U *}(\pi^{-1}\Fc\otimes Log\otimes or_{\Tc/\Sc})\rightarrow \pi_{|D *}(\pi^{-1}\Fc\otimes Log_{|D})\rightarrow \Fc
\end{equation*}
But the last arrow is obtained by applying $\Fc\otimes$ and the projection formula to the epimorphism  
$\pi_{|D *}(Log_{\Tc/\Sc|D})\stackrel{\aug}{\rightarrow}\pi_{|D *}A\stackrel{tr}{\rightarrow}A$, and therefore is itself an epimorphism (for $tr$ see \cite{Iv} VII.4.).
In particular
\begin{equation*}
R^p\pi_{|U*}(\pi^{-1}\Fc\otimes Log\otimes or_{\Tc/\Sc})=0,\ p\neq d-1.
\end{equation*}
\end{proof}
\end{proposition}

\begin{remark}\label{aug}
Let us recall the augmentation morphism $\pi_{|D *}(Log_{|D})\stackrel{\aug}{\rightarrow}\pi_{|D *}A\stackrel{tr}{\rightarrow}A$.
It is enough to do so on the fibers $\pi^{-1}(s)$. So we have a topological torus $T$ and $D\subset T$ is just a finite set of points. We have to examine the chain of maps
\begin{equation*}
H^0(D,Log_{T|D})\stackrel{purity\cong}{\rightarrow}H^{d}_D(T,Log_{T}\otimes or_T)\to H^{d}(T,Log_{T}\otimes or_T) \stackrel{\aug}{\rightarrow}
\end{equation*}
\begin{equation*}
 \to H^{d}(T,or_T)\stackrel{\cong}{\rightarrow} H^{d}(T,A)\otimes or_T(T)\stackrel{ev}{\rightarrow}A
\end{equation*}

We have $H^{d}_D(T,Log\otimes or_T)=\bigoplus_{x\in D} H^{d} _{\left\{x\right\}}(T,Log\otimes or_T)$ and isomorphisms
\begin{equation*}
H^0(\left\{x\right\},A)\stackrel{purity}{\rightarrow}H^{d}_{\left\{x\right\}}(T,A)\rightarrow H^{d}(T,A),
\end{equation*} 
as $T$ is orientable. We choose a generator $1(x)\in H^0(\left\{x\right\},A)$ and denote the corresponding images with $1_{\left\{x\right\}}\in H^{d}_{\left\{x\right\}}(T,A)$ and $\omega(x)\in H^{d}(T,A)$. We set $\omega(x)^*\in or_T(T)$ the element dual to $\omega(x)$.
Given $\sum_x f(x)1(x)\in H^0(D,Log_{|D})$, $f(x)\in Log_x$, purity maps it to $\sum_x f(x)1_{\left\{x\right\}}\otimes \omega(x)^* \in H^{d}_Z(T,Log\otimes or_T)$. This is mapped by $\aug$ to $\sum_x \aug(f(x))\omega(x)\otimes\omega(x)^*\in H^{d}(T,or_T)$ and then by $ev$ to $\sum_x \aug(f(x))$ as claimed.
\end{remark}
\begin{proposition}\label{global_Log_localization}
With notations from above we have an exact sequence
\begin{equation*}
0\rightarrow H^{d-1}(U,\pi^{-1}\Fc\otimes Log\otimes or_{\Tc/\Sc})\stackrel{\res}{\rightarrow} H^0(D,\pi^{-1}\Fc\otimes Log_{|D})
\stackrel{\aug}{\rightarrow} H^{0}(\Sc,\Fc)
\end{equation*}
\begin{proof}
Using the Grothendieck-Leray spectral sequence for the composition of functors $\Gamma(U,-)=\Gamma(\Sc,-)\circ \pi_{|U*}$ we easily see
\begin{equation*}
H^{p}(U,\pi^{-1}\Fc\otimes Log\otimes or_{\Tc/\Sc})=
H^{p-(d-1)}(\Sc,R^{d-1}\pi_{|U*}(\pi^{-1}\Fc\otimes Log\otimes or_{\Tc/\Sc})),
\end{equation*}
as $R^p\pi_{|U*}(\pi^{-1}\Fc\otimes Log\otimes or_{\Tc/\Sc})=0$, $p\neq d-1$. Applying $H^0(\Sc,-)$ to \Cref{relative_Log_localization} yields the result.
\end{proof}
\end{proposition}
\begin{definition}\label{def_pol}
We keep the notations from above.\\
Given $f\in H^0(D,\pi_{|D}^{-1}\Fc\otimes Log_{|D})$ with $\aug(f)=0$ we define the \textit{polylogarithm associated to $f$} as the unique cohomology class
\begin{equation*}
\pol(f)\in H^{d-1}(U,\pi^{-1}\Fc\otimes Log\otimes or_{\Tc/\Sc})\
\end{equation*}
with $\res(\pol(f))=f$. If $0(\Sc)\cap D=\emptyset$, we define the \textit{polylogarithmic Eisenstein class associated to $f$} to be the pullback
\begin{equation*}
\Eis(f):= 0^{*}\pol(f)\in H^{d-1}(\Sc,\Fc\otimes 0^{-1}(Log\otimes or_{\Tc/\Sc}) )
\end{equation*}
\end{definition}

\subsection{Trivializations of the logarithm sheaf}

Finding non-trivial polylogarithm classes is by definition equivalent to finding non-trivial $f\in H^0(D,\pi^{-1}\Fc\otimes Log_{|D})$ with $\aug(f)=0$. For which $D$ do we have sections or a trivialization for $Log_{|D}$?\\
We just write
$Log^n=Log_{\Tc/\Sc} ^n$, when the topological situation $\pi:\Tc\rightarrow \Sc$ has been fixed.
\begin{lemma}\label{0Log}
The augmentation map $0^{-1}Log\rightarrow A$ splits. The image of $1\in A$ is denoted by $1_0\in H^0(\Sc,0^{-1}Log)$. 
$0^{-1}Log=:\Rc$ is a sheaf of augmented algebras and $Log$ is an invertible $\pi^{-1}\Rc$-module.
\begin{proof}
$0^{-1}Log$ is the local system associated to the $\pi_1(\Sc,s_0)$-module $A[[\pi_1(\Tc_{s_0},t_0)]]$. The augmentation map corresponds to the augmentation map of this algebra. It has the obvious $\pi_1(\Sc,s_0)$-equivariant section providing the section of the corresponding sheaves. Obviously $\pi_1(\Sc,s_0)$ acts by algebra homomorphisms on $A[[\pi_1(\Tc_{s_0},t_0)]]$, therefore $\Rc$ is a sheaf of algebras. $\pi^{-1}\Rc$ is associated to the $\pi_1(\Tc,t_0)$-module $A[[\pi_1(\Tc_{s_0},t_0)]]$ with the action $(v,\gamma)\sum_l f(l)(l)=\sum_lf(l)(\gamma l)$ for $(v,\gamma)\in\pi_1(\Tc,t_0)$. Let us denote this module by $A[[\pi_1(\Tc_{s_0},t_0)]]_{\Rc}$ and the module associated to $Log$ by $A[[\pi_1(\Tc_{s_0},t_0)]]_{Log}$ to emphasize their different $\pi_1(\Tc,t_0)$-structures. Then we see that the multiplication map
\begin{equation*}
A[[\pi_1(\Tc_{s_0},t_0)]]_{\Rc}\otimes_A A[[\pi_1(\Tc_{s_0},t_0)]]_{Log}\rightarrow A[[\pi_1(\Tc_{s_0},t_0)]]_{Log}, f\otimes g\mapsto f\cdot g
\end{equation*}
is $\pi_1(\Tc,t_0)$-equivariant proving that $Log$ is an $\pi^{-1}\Rc$-module. 
\end{proof}
\end{lemma}
We set $\Rc^n:=0^{-1}Log^n$.
\begin{lemma}
Let $\phi:\Tc^\prime \rightarrow \Tc$ be an isogeny of tori over $\Sc$. We consider the logarithm sheaf over a ring $A$ and suppose that $deg(\phi)(s)\in A^{\times}$ for all $s\in \Sc$, where $A^{\times}$ always denotes the units in a ring $A$. 
Then the natural morphism $Log_{\Tc^\prime/\Sc}\rightarrow \phi^{-1}Log_{\Tc/\Sc} $
is an isomorphism.
\begin{proof}
It suffices to consider the map of underlying $\pi_1(\Tc^\prime,t_0 ^\prime )$-modules. It is given by the natural map $A[[\pi_1(\Tc^\prime _{s_0},t_0^\prime)]]\stackrel{\phi_*}{\rightarrow} A[[\pi_1(\Tc _{s_0},t_0)]]$. This map is an isomorphism by \cite{BKL} proposition 1.1.5.
\end{proof}
\end{lemma}
\begin{definition} 
Let $X$ be a topological space, $A\to A^\prime$ a ring homomorphism and $\Gc$ a sheaf of $A$-modules. We set $\Gc_{A^\prime}:=\Gc\otimes_A A^\prime$.
\end{definition}
Given a family of topological tori $\pi:\Tc\rightarrow \Sc$ with logarithm sheaf $Log$ we may consider $Log^{\times} :=\aug^{-1}(1)$ where $\aug:Log\rightarrow A$ is the augmentation map. $Log^{\times}$ is a sheaf of sets on $\Tc$.  We may realize this sheaf as a sheaf of sections of a topological space $Log^{\times sp}$ over $\Sc$. To do so we set
$R:=A[[\pi_1(\Tc_{s_0},t_0)]]$ and denote again by $\mathfrak a:=\aug^{-1}(0)$ its augmentation ideal. In $R$ we have the group of $1$-units $1+\mathfrak a=:R^{\times}_1$. The action of $\pi_1(\Tc,t_0)$ on $R$ restricts to $R^{\times}_1 $ making it a $\pi_1(\Tc,t_0)$-module. Considering $R^\times _1 $ with the discrete topology we get 
$Log^{\times sp}:=\tilde{\Tc}\times_{\pi_1(\Tc,t_0)} R^{\times} _1$. We had $\tilde{\Tc}=\tilde{\Sc}\times H_1(\Tc_{s_0},\R)$ and therefore we may write $Log^{\times sp}=\tilde \Sc\times_{\pi_1(\Sc,s_0)}(H_1(\Tc_{s_0},\R)\times_{\pi_1(\Tc_{s_0},t_0)}R^{\times} _1)$.
The second description tells us that $Log^{\times sp}$ may be considered as a group object over $\Sc$ with typical fiber $H_1(\Tc_{s_0},\R)\times_{\pi_1(\Tc_{s_0},t_0)}R^{\times} _1$. 
We also have the sheaf $\Rc^{\times}_1 $ of $1$-units of the augmented algebra $\Rc$. One has the inclusion $\delta:\pi_1(\Tc _{s_0},t_0)\subset R^{\times}_1$ giving rise to the inclusion $\delta: \Hc\rightarrow \Rc^{\times}_1 $. If we consider the fibers of $\Tc$ with the discrete topology, we may summarize this in the following pushout diagram of abelian sheaves on $\Sc$
\begin{equation*}
\begin{xy}
\xymatrix{
0\ar[r] & \Rc^{\times}_1 \ar[r] & Log^{\times sp}\ar[r]^p & \Tc\ar[r] & 0\\
0\ar[r] & \Hc_\Z\ar[u]^\delta \ar[r] & \Hc_\R  \ar[u]\ar[r] & \Tc\ar[r]\ar[u]^{\id} & 0
}
\end{xy}
\end{equation*}
where we consider $\Tc$ and $Log^{\times sp}$ as sheaves of sections.
Let $H\subset \Tc$ be a subgroup. By this we mean that $H$ is a sub local system of $\Tc$. We also may consider $H$ as a subspace of $\Tc$ covering $\Sc$. Now we follow \cite{BKL} Definition 1.5.1. 
\begin{definition}
A \textit{multiplicative trivialization of $Log$ along $H$} is a morphism $t:H\rightarrow Log^{\times sp}$, which is a section for $p:Log^{\times sp}\rightarrow \Tc$.
\end{definition}
\begin{remark}
A multiplicative trivialization of $Log$ gives a section $1_H\in \Gamma(H,Log_{|H})$, $\aug(1_H)=1$, trivializing $Log_{|H}$ as an invertible $\pi_{|H}^{-1}\Rc$-module. 
\end{remark}
Let $\Tc^{(A)}$ be the subsheaf of $\Tc$ consisting of all torsion sections whose order is invertible in $A$ and $\Tc^{tors}$ the subsheaf of all torsion sections. 
\begin{proposition}\label{trivializationLog}
There is a unique multiplicative trivialization $\rho_{can}:\Tc^{(A)}\rightarrow Log^{\times sp}$ for $Log$ along $\Tc^{(A)}$. It is compatible with isogenies and  for $t\in \Tc[N](\Sc)\subset \Tc^{A}(\Sc)$ it is explicitly given by the isomorphisms
\begin{equation*}
t^{-1}Log\cong t^{-1}[N]^{-1}Log\cong ([N]\circ t)^{-1}Log\cong 0^{-1}Log,
\end{equation*}
where one uses $Log\cong [N]^{-1}Log$, as $deg[N]\in A^{\times}$, and $[N]\circ t= Nt= 0$, as $t$ is $N$-torsion.
\begin{proof}
By \cite{BKL} proposition 1.5.3. we have the statement, if $\Sc=\left\{pt.\right\}$. But as $Log^{\times sp}$ and $\Tc$ are fiber bundles over $\Sc$, we may easily glue these local multiplicative trivializations to a global one by using unicity. 
\end{proof}
\end{proposition} 
This proposition tells us that torsion sections are a powerful tool to construct $D\subset \Tc$ for which we can write down sections of $Log_{|D}$. But before we can start calculations with $Log$, let us have a closer look at the situation over $A=\C$.
\begin{definition} 
Let $X$ be a topological space and $\Rc$ a sheaf of rings. For a projective system $(\Fc_n)$, $\Fc:=\varprojlim \Fc_n$, of $\Rc$-module sheaves and an $\Rc$-module sheaf $\Gc$ on $X$ we set $\Fc\hat\otimes_\Rc \Gc:=\varprojlim (\Fc_n\otimes_\Rc \Gc)$.
\end{definition}
Let $\Cc^\infty _\Tc$ be the sheaf of complex valued smooth functions on $\Tc$ and $\Omega_\Tc$ the de Rham complex of complex valued smooth forms on $\Tc$. We set
\begin{equation*}
Log^{\infty n}:=Log^n\otimes_\C \Cc^\infty_\Tc\ \text{and}\ Log^\infty:=\varprojlim Log^{\infty n}=Log\hat{\otimes}_\C \Cc^\infty _\Tc
\end{equation*}
to be the pro vector bundle associated to $Log$. As it comes from local systems, it is naturally endowed with an integrable connection
$\id\otimes d: Log^\infty\rightarrow Log^\infty\hat{\otimes}_{\Cc_\Tc ^\infty} \Omega_\Tc ^1$, whose kernel is exactly $Log$.
\begin{proposition}\label{R_alg}
Given a lattice $L$ one has the isomorphism of augmented algebras
\begin{equation*}
\Q[[L]]\rightarrow \prod_{k\geq 0} \Sym^k L_\Q, l\mapsto \exp(l):=\sum_{k\geq 0}\frac{l^{\otimes k}}{k!},
\end{equation*}
where the right hand side carries the Cauchy product and the augmentation is the projection onto the $\Sym^0L_\Q=\Q$- part.
In particular, one has an isomorphism $\Rc\stackrel{\exp}{\rightarrow}\prod_{k\geq 0}\Sym^k\Hc_\Q$ of augmented $\Q$-algebras and an isomorphism of systems 
\begin{equation*}
(\Rc^n)\stackrel{\exp}{\rightarrow}\left(\prod_{k=0}^nSym^k\Hc_\Q\right)_n
\end{equation*}
\begin{proof}
\cite{BKL} Corollary 1.1.10.
\end{proof}
\end{proposition}
\begin{remark}\label{Log_integral}
This isomorphism does not work integrally. Nevertheless, let $A$ be a torsion free ring. We have
\begin{equation*}
 A[[L]]\rightarrow \prod_{k\geq 0} \Sym^k L_{PD}, l\mapsto \exp(l):=\sum_{k\geq 0}\frac{l^{\otimes k}}{k!},
\end{equation*}
into a subring of $\prod_{k\geq 0} \Sym^k L_{A_\Q}$. Here we define the $A$-submodule $\Sym^k L_{PD}\subset \Sym^k L_{A_\Q}$ to be generated by $x_1^{[i_1]}\cdot..\cdot x_n^{[i_n]}$ with $x_1,...,x_n \in L$, $i_l\in \N_0$, $l=1,...,n$ and $\sum_{l=1}^ni_l=k$, where we define in any $\Q$-algebra $B$ the \textit{divided powers} $x^{[n]}:=\frac{x}{n!}$ for $x\in B$. If $e_1,...,e_d$ is a basis of $L$, consider $1\leq j_1<...< j_n\leq d$ and $i_l\in \N_0$, $l=1,...,n$, with $\sum_{l=1}^ni_l=k$, then $e_{j_1}^{[i_{1}]}\cdot..\cdot e_{j_n}^{[i_n]}$ is a basis of $\Sym^k L_{PD}$. The $\Sym^k L_{PD}$ give rise to local systems $\Sym^k \Hc_{PD}$ which are integral structures for $\Sym^k \Hc_{A_\Q}$. We have the projection map $p_k:\Rc\rightarrow \Sym^k \Hc_{PD}$ which is induced by
$A[[L]]\stackrel{\exp}{\rightarrow} \prod_{n\geq 0} \Sym^n L_{PD}\stackrel{pr_k}{\rightarrow}\Sym^k L_{PD}$. 
\end{remark}
\begin{remark}
Transporting the action of the fundamental group via $\exp$ we may understand the logarithm sheaf with complex coefficients via the completed symmetric algebra, what we will always do from now on. In particular, we have $\Rc=\prod_{n\geq 0}\Sym^n\Hc_\C$.
\end{remark}
Set
\begin{equation*}
Log^{sp}:=\tilde{\Tc}\times_{\pi_1(\Tc,t_0)}\prod_{n\geq 0} \Sym^n H_1(\Tc_{s_0},\C).
\end{equation*}
The sections of $Log^\infty$ may now be seen as the smooth sections $s:\Tc\rightarrow Log ^{sp}$. These are $\pi_1(\Tc,t_0)$-equivariant maps $f:\tilde{\Tc}\rightarrow \prod_{n\geq 0} \Sym^n H_1(\Tc_{s_0},\C)$ such that for all $n\geq 0$ the component maps $f_n:\tilde{\Tc}\rightarrow \Sym^n H_1(\Tc_{s_0},\C)$ are smooth. In this sense it is also clear what we mean by smooth or continuous sections of $Log^{\times sp}$.
\begin{proposition}\label{cont_triv}
There is a unique \textit{continuous multiplicative trivialization of $Log$}. By this we mean a continuous homomorphism $\rho_{cont}:\Tc\rightarrow Log^{\times sp}$, which is a section for $p:Log^{\times sp}\rightarrow \Tc$. It is compatible with $N$-multiplication and we have $\rho_{cont|\Tc^{tors}}=\rho_{can}$. Explicitly $\rho_{cont}$ is given by $\exp(-\underline v)=\sum_{n\geq 0}\frac{(-\underline v)^{\otimes n}}{n!}$, where  
\begin{equation*}
\underline v:\tilde{\Sc}\times H_1(\Tc_{s_0},\R), (s,v) \mapsto v\in \prod_{n\geq 0} \Sym^n H_1(\Tc_{s_0},\C)
\end{equation*}
is considered as a $\Cc^\infty$-section of $Log$ on the universal cover of $\Tc$.
\begin{proof}
$\pi_1(\Tc,t_0)$ acts on $ \prod_{n\geq 0}\Sym^n H_1(\Tc_{s_0},\C)$ by
\begin{equation*}
(l,\gamma)\sum_{n\geq 0}f_n=\exp(l)\cdot\sum_{n\geq 0}\gamma f_n.
\end{equation*}
From the equivariance of $\exp$
\begin{equation*}
\exp(-(v(l,\gamma))=\exp(-(\gamma^{-1}(v+l)))=\exp(-\gamma^{-1}l)\cdot \gamma^{-1} \exp(-v)=(l,\gamma)^{-1}\exp(-v)
\end{equation*}
it is clear that the section above descends to a multiplicative trivialization as claimed. The problem of unicity is local on the base $\Sc$ and as all spaces are fiber bundles over $\Sc$, we may reduce to the case where $\Sc$ is a point. This is \cite{BKL} Lemma 2.1.2. 
\end{proof}
\end{proposition}
We set $\Rc^{\infty n}:=\pi^{-1}\Rc^n\otimes\Cc_\Tc ^\infty\ \text{and}\ \Rc^\infty=\varprojlim\Rc^{\infty n}=\pi^{-1}\Rc\hat{\otimes}\Cc_\Tc ^\infty$ the complex pro vector bundle associated to $\pi^{-1}\Rc$.
$Log^\infty$ is now a module over this sheaf of algebras.
We also have the sheaf of differential graded algebras $\pi^{-1}\Rc\hat{\otimes}\Omega_\Tc=\Rc^\infty\hat{\otimes}_{\Cc^\infty_\Tc}\Omega_\Tc$
with the trivial differential $\id\otimes d$. One has $\Hc_\Z\subset \Rc$ and $\Hc_\C^*\subset \Omega^1_\Tc$ interpreting sections of $\Hc_\C^*$ as invariant differential forms along the fibers of $\pi$. The canonical isomorphism
\begin{equation*}
End(\Hc_\Z)\cong \Hc_\Z\otimes \Hc_\Z^*,\ \id \mapsto \kappa
\end{equation*}
gives the canonical global section $\kappa$ of $\pi^{-1}\Rc\hat{\otimes}\Omega_\Tc$. The section $\kappa$ is closed and multiplication with $\kappa$ gives an operator on our differential graded algebra. 
\begin{corollary}
One has the isomorphism of pro vector bundles $\Rc^\infty\rightarrow Log^\infty$, $f\mapsto f\cdot\rho_{cont}$.
If endow $\Rc^\infty$ with the integrable connection $\nabla:=d-\kappa$, this isomorphism commutes with the connections.
\begin{proof}
As $Log^\infty$ is an invertible $\Rc^\infty$-module the isomorphism is clear, since $\aug(\rho_{cont})=1$ and $\rho_{cont}$ is a global nowhere vanishing section of $Log^\infty$. To prove the compatibility of the connections we need to calculate $\id\otimes d(\rho_{cont})$. This can be done on the universal cover $\tilde{\Tc}$. We have
\begin{equation*}
\id\otimes d(\exp(-\underline v))=\id\otimes d\sum_{n\geq 0}\frac{(-\underline v)^{\otimes n}}{n!}=\sum_{n\geq0}\frac{\id\otimes d(-\underline v)^{\otimes n}}{n!}=
\end{equation*}
\begin{equation*}
\sum_{n\geq0}\frac{n(-\underline{v} )^{\otimes n-1}\cdot (\id\otimes d(-\underline v))}{n!}=\sum_{n\geq0}\frac{(-\underline v)^{\otimes n}}{n!}\id\otimes d(-\underline v)
\end{equation*}
Write $\underline v=\sum_{i=1} ^{d}e_i\otimes x_i$ for some $\R$-basis $(e_i)$ of $H_1(\Tc_{s_0},\R)$, where $(x_i)$ is the corresponding dual basis of $(e_i)$ considered as coordinate functions on $H_1(\Tc_{s_0},\R)$. Then
\begin{equation*}
\id\otimes d(\underline v)=\id\otimes d\sum_{i=1} ^{d}e_i\otimes x_i=\sum_{i=1} ^{d}e_i\otimes dx_i=\kappa
\end{equation*}
We conclude $\id\otimes d(\rho_{cont})=\rho_{cont}\cdot(-\kappa)$. This completes the proof.
\end{proof}
\end{corollary}
Now we come to our case of interest.
\begin{proposition}\label{pol}
Let $\pi:\Tc\rightarrow \Sc$ be family of topological tori with level-$N$-structure and constant fiber dimension $d$. We consider $D=\Tc[N]\setminus 0(\Sc)\cong ((\Z/N\Z)^{d}\setminus \left\{0\right\})\times \Sc$
and the logarithm sheaf $Log$ over a noetherian ring $A$ with $N\in A^{\times}$.
We have then a canonical trivialization
$H^0(D,Log_{|D})=H^0(((\Z/N\Z)^{d}\setminus \left\{0\right\})\times \Sc,\Rc)$. One has $A\subset \Rc$ and thus any locally constant function 
\begin{equation*}
f:((\Z/N\Z)^{d}\setminus \left\{0\right\})\times \Sc\rightarrow A,\ \aug(f)=\sum_{v\in(\Z/N\Z)^{d}\setminus \left\{0\right\}}f(v,\ )=0\in H^0(\Sc,A)
\end{equation*}
gives rise to a section $f\in H^0(D,Log_{|D})$ of augmentation zero and therefore to a unique polylogarithmic cohomology class $\pol(f)\in H^{d-1}(U,Log\otimes or_{\Tc/\Sc})$. 
Moreover, we get the polylogarithmic Eisenstein class
$0^{*}\pol(f)=:\Eis(f)\in H^{d-1}(\Sc,\Rc\otimes 0^{-1}or_{\Tc/\Sc})$, which induces polylogarithmic Eisenstein classes
\begin{equation*}
\Eis^k(f)\in H^{d-1}(\Sc,\Sym^k\Hc_{PD}\otimes 0^{-1}or_{\Tc/\Sc}),
\end{equation*}
whenever $A$ is torsion-free.
\begin{proof}
$A\subset \Rc$ is \Cref{0Log}. The trivialization is \Cref{trivializationLog}. Finally, we just need \Cref{aug} and \Cref{def_pol} to get $\pol(f)$. We set $\Eis^k(f):=p_k(\Eis(f))$, where $p_k:\Rc\rightarrow \Sym^k\Hc_{PD}$ is the projection defined in \Cref{Log_integral}.
\end{proof}
\end{proposition}
\begin{remark}
The restriction to torsion-free rings in the proposition above is not necessary. To get the appropriate substitute for our $\Sym^k\Hc_{PD}$
one has to use the universal divided power algebra $\Gamma_k(L)$, see \cite{BKL} Lemma 1.1.7.
\end{remark}

\subsection{A differential equation for the polylogarithm}

We want explicit cohomology classes representing our polylogarithms over $\C$. To do so we need an explicit resolution of $Log$ over $\C$, which allows us to calculate the localization sequence in \Cref{global_Log_localization}.
\begin{lemma}
Let $X$ be an oriented $n$-dimensional $\Cc^\infty$-manifold and $\Omega_X=(\Omega^\bullet _X ,d)$ the de Rham complex of $\C$-valued $\Cc^\infty$-differential forms. We may represent the dualizing complex of $\Dc_X$ by non-continuous functionals on the de Rham complex:
\begin{equation*}
\Dc_X ^p (U)=\Gamma_c(U,\Omega_X ^{-p})^*,\ U\subset X\ \text{open}
\end{equation*}
with differential $d$ defined by
\begin{equation*}
d(T)(\omega):=(-1)^{p+1}T(d(\omega)),\ \omega\in \Gamma_c(U,\Omega_X ^{-p})^*,\ T\in  \Gamma_c(U,\Omega_X ^{-p})^*
\end{equation*}
\begin{proof}
The Poincar\'{e} Lemma gives a quasi-isomorphism $\C\rightarrow \Omega_X$. So $\Omega_X$ is a bounded resolution of the constant sheaf $\C$ by soft $\C$-sheaves. We conclude with \cite{Iv} V.2.1.
\end{proof}
\end{lemma}
\begin{remark}
$\Dc_X$ is a complex of injective $\C$-sheaves which also have the structure of $\Cc_X ^\infty$-modules. One has a canonical quasi-isomorphism $or_X[n]=H^{-n}(\Dc_X)[n]\rightarrow \Dc_X$ (\cite{Iv}VI.3.). By choosing a volume form on $X$ we may define integration of $n$-forms. In other words, integration defines a functional $\Gamma_c(U,\Omega_X ^n)\rightarrow \C$, $\omega\mapsto \int_X \omega$ for any open $U\subset X$. This functional vanishes on exact forms by Stoke's theorem and therefore $\int_X$ defines an element in $or_X(X)=H^{-n}(\Dc_X)(X)$. It is nowhere vanishing and we get the isomorphism $\C\to or_X,\ 1\mapsto \int_X$.
In particular, we get the quasi-isomorphism $\C[0]\to \Dc_X[-n]$.

Now let $i:Z\hookrightarrow X$ be a closed oriented submanifold of codimension $c$.\\
We have the pullback morphism $i^*:\Omega_X\rightarrow i_*\Omega_Z$. As $i$ is proper we also get for any $U\subset X$ open $i^*:\Gamma_c(U,\Omega_X)\rightarrow \Gamma_c(U,i_*\Omega_Z)$ and by taking duals 
\begin{equation*}
i_*:\Gamma_c(U,i_*\Omega_Z)^*\rightarrow \Gamma_c(U,\Omega_X)^*
\end{equation*}
compatible with restrictions in $U$. So identifying $\Dc_X$ with the sheaf of non-continuous functionals on $\Omega_X$ we get $i_*\Dc_Z\rightarrow \Dc_X$  and by adjunction $i_*:\Dc_Z\rightarrow i^!\Dc_X$. This is exactly the arrow inducing \Cref{dualizing_purity}. 
We have then
\begin{equation*}
H^0(Z,\C)=H^0(Z,or_Z)=H^0(Z,\Dc_Z[c-n])=H^0(Z,i^!\Dc_X[c-n])=
\end{equation*}
\begin{equation*}
H^0_Z(X,\Dc_X[c-n])=H^c_Z(X,or_X)=H^c_Z(X,\C)
\end{equation*}
Here we have exactly 
\begin{equation*}
H^0(Z,\C)\rightarrow H^c_Z(X,\C),\ f\mapsto \int_f
\end{equation*}
with $\int_f\omega:=\int_Z f\cdot i^*\omega$ and integration is defined by the orientation defining $\C\cong or_Z$.
\end{remark}  
\begin{proposition}
Let $\pi:\Tc\rightarrow \Sc$ be a family of topological tori, $Log=\varprojlim Log^n$ the associated logarithm sheaf over $\C$ and $\Dc_\Tc$ the dualizing complex on $\Tc$. The system of complexes
$
(Log^n\otimes or_{\Tc}^{-1}  \otimes  \Dc_\Tc[-dim(\Tc)])_n
$
is an injective resolution of $(Log^n)$ in $Sh(\Tc,\C)^\N$. In particular,
$
\varprojlim   (Log^n\otimes or_{\Tc}^{-1} \otimes \Dc_\Tc)[-dim(\Tc)]
$
is an injective resolution of $Log$.
\begin{proof}
We already know that $\Dc_\Tc ^p$ is injective in $Sh(\Tc,\C)$. As $Log^n \otimes or_\Tc ^{-1}$ is a local system the sheaf $Log^n\otimes  or_{\Tc}^{-1}\otimes \Dc_\Tc ^p$ is injective in $Sh(\Tc,\C)$. In order to see that
$
(Log^n\otimes or_{\Tc}^{-1} \otimes \Dc_\Tc ^p)_n
$
is an injective object in $Sh(\Tc,\C)^\N$ we need to see that all transition maps are split epimorphisms. Now $\Dc_\Tc ^p$ is a $\Cc^\infty_\Tc$-module and therefore
$
( Log^n\otimes or_{\Tc}^{-1} \otimes \Dc_\Tc ^p)_n\cong ( Log^{\infty n}\otimes_{\Cc^\infty _\Tc} (or_{\Tc}^{-1}\otimes \Dc_\Tc ^p)).
$
We have an isomorphism of systems $(Log^{\infty n})_n=(\Rc^{\infty n})_n$ and the latter system has split epimorphic transition maps by \Cref{R_alg}. We conclude that $( Log^n\otimes or_{\Tc}^{-1} \otimes \Dc_\Tc ^p)_n$ has split epimorphic transition maps and therefore is injective in $Sh(\Tc,\C)^\N$. $(Log^n\otimes or_{\Tc}^{-1} \otimes \Dc_\Tc[-dim(\Tc)] )_n$ being a resolution of $( Log^n )_n$ may be checked for each $n\in \N_0$ separately by using the quasi-isomorphism $or_\Tc[dim(\Tc)] \rightarrow \Dc_\Tc$.
\end{proof}
\end{proposition}

\begin{proposition}\label{differential_equ}
Let $\Sc$ be an oriented manifold and $\pi:\Tc\rightarrow \Sc$ a family of topological tori with level-$N$-structure and constant fiber dimension $d$. We consider the logarithm sheaf $Log$ with $\C$-coefficients. 
The polylogarithms as defined in \Cref{pol} may be represented by non-continuous functionals 
\begin{equation*}
\pol(f)\in \Gamma(\Tc,\Rc^\infty\hat{\otimes}_{\Cc^\infty _\Tc} (\pi^{-1}or_\Sc ^{-1}\otimes \Dc_\Tc))\ \text{with}\ \nabla (\pol(f)) =\vol_\Sc\otimes\int_f
\end{equation*}
Here the differential $\nabla$ is induced by $(\Rc^\infty,\nabla)$ and integration on $\Sc$ is defined with respect to a chosen volume form $\vol_\Sc$ which we also use to trivialize $or_\Sc^{-1}$.
\begin{proof}
We set $dim(\Tc)=m$.
To calculate cohomology we can use the injective resolutions
$
\varprojlim (Log^n\otimes or_{\Tc/\Sc}\otimes or_\Tc ^{-1}\otimes \Dc_\Tc )[-m]
$
with differential $\id\otimes d$. We write these complexes as
$
Log^\infty\hat{\otimes}_{\Cc^\infty _\Tc} (\pi^{-1}or_\Sc ^{-1}\otimes \Dc_\Tc )[-m] 
$
using the $\Cc^\infty_\Tc$-module structure of $\Dc_\Tc$ and the rules for orientations. We have a chain of quasi-isomorphisms
\begin{equation*}
\pi_{|D}^{-1}\Rc\stackrel{\rho_{can}}{\rightarrow}Log_{|D}\rightarrow \varprojlim (Log_{|D}^n\otimes or_D ^{-1}\otimes \Dc_D)[-(m-d)]\stackrel{\id\otimes i_*}{\rightarrow}
\end{equation*}
\begin{equation*}
\varprojlim(Log_{|D}^n\otimes or_D ^{-1}\otimes i^!\Dc_\Tc)[-(m-d)]\stackrel{\cong}{\rightarrow}i^!\varprojlim (Log^n\otimes \pi^{-1} or_\Sc ^{-1}\otimes \Dc_\Tc)[-(m-d)].
\end{equation*}
The last isomorphism is \Cref{!and-1}. If we identify
\begin{equation*}
i^!\varprojlim (Log^n\otimes \pi^{-1} or_\Sc ^{-1}\otimes \Dc_\Tc)[-(m-d)]=i^!(Log^\infty\hat{\otimes}_{\Cc_\Tc ^\infty}(\pi^{-1} or_\Sc ^{-1}\otimes  \Dc_\Tc))[-(m-d)],
\end{equation*}
we may give a global section of this sheaf namely $\rho_{cont}\otimes \vol_\Sc \otimes \int_f$. If we go through the previous morphisms remembering $\rho_{can}=\rho_{cont|\Tc^{tors}}$, we see that $f$ maps exactly to this section. So we just have proved
\begin{equation*}
H^0(D,\Rc)\stackrel{\rho_{can}}{\rightarrow}H^0(D,Log_{|D})\stackrel{purity}{\rightarrow}H^d_D(\Tc,Log\otimes or_{\Tc/\Sc}),\ f\mapsto \rho_{cont}\otimes \vol_\Sc \otimes \int_f.
\end{equation*}
Now we use the continuous trivialization $\rho_{cont}$ of $Log$ to identify our injective resolutions from above with
$\Rc^\infty\hat{\otimes}_{\Cc_\Tc ^\infty} (\pi^{-1}or_\Sc ^{-1}\otimes \Dc_\Tc) [-m]$. 
We call the differential of the last complex $\nabla$, as it is induced by $(\Rc^\infty,\nabla)$.
Clearly, $f$ corresponds now to
\begin{equation*}
\vol_\Sc\otimes \int_f\in H^d(D,i^!(\Rc^\infty\hat{\otimes}_{\Cc_\Tc ^\infty} (\pi^{-1}or_\Sc ^{-1}\otimes \Dc_\Tc) [-m])).
\end{equation*}
The residue map $\res$ of our localization sequence is nothing but a connecting homomorphism of a long exact cohomology sequence. So being in the image of $\res$ just means that there is a 
\begin{equation*}
\pol(f)\in \Gamma(\Tc,\Rc^\infty\hat{\otimes}_{\Cc_\Tc ^\infty} (\pi^{-1}or_\Sc ^{-1}\otimes \Dc_\Tc )[-m]),\ \nabla(\pol(f))=\vol_\Sc\otimes \int_f.
\end{equation*}
\end{proof}
\end{proposition}

\subsection{Invariant functionals and quotient manifolds}

To represent polylogarithms with functionals it seems promising to pullback the situation to the universal cover of $\Tc$, because the logarithm sheaf is constant there. We would like to write down sections invariant under the $\pi_1(\Tc,t_0)$-action descending to sections representing polylogarithms. This would be exactly the right thing, if we knew that polylogarithms may be represented by differential forms. The problem is that we cannot pullback functionals in a naive way, as the universal covering $p:\tilde{\Tc}\rightarrow \Tc$ is usually not proper (finite). But things still work the other way round, invariant functionals on $\tilde{\Tc}$ descend to functionals on $\Tc$. This is what we explain now.

Let $X$ be a $\Cc^\infty$-manifold and $G$ a group acting on $X$ by diffeomorphisms from the right.
So for $g\in G$ we have the diffeomorphism $r_g:X\to X,\ x\mapsto xg $.
We assume this action properly discontinuous and fixpoint free. The quotient space with the quotient topology has a unique $\Cc^\infty$-structure making the natural cover $\pi:X\rightarrow X/G$ a local diffeomorphism. We are interested in currents on $X/G$. For a definition of currents see \cite{DR} chapitre III $8$. Denote the currents on a manifold $V$ by $\mathscr{D}(V)$ and the sheaf of currents on $V$ by $\mathscr{D}_V$.
\begin{proposition}\label{invariant_functionals}
Let $X$ be a $G$-manifold as described above. There is an isomorphism
\begin{equation*}
(\pi_*\mathscr{D}_X)^G\rightarrow \mathscr{D}_{X/G}, T\mapsto \tilde T
\end{equation*}
\begin{proof}
Let us first construct the map. Take $U\subset X/G$ open and $T\in \mathscr{D}(\pi^{-1}(U))^G$ and a $\omega\in \Omega_c(U)$. We have to define $\tilde T(\omega)$. To do so take an open cover $\bigcup_i U_i= X/G$ together with $U_i ^\prime \subset X$ open, such that $\pi_{|U_i ^\prime}:U_i ^\prime\rightarrow U_i$ is a diffeomorphism. Choose a partition of unity $\sum_i \epsilon_i=1$ subordinated $U_i$. The support of $\epsilon_i \omega$ is compact and contained in $U_i\cap U$, so the the support of $\pi_{|U_i ^\prime} ^* \epsilon_i \omega$ is compact and contained in $U_i^\prime\subset \pi^{-1}(U)$. We extend the latter differential form by zero outside $U_i ^\prime$ and set 
\begin{equation*}
\tilde T(\omega):=\sum_i T(\pi_{|U_i ^\prime} ^* \epsilon_i \omega)=T(\sum_i \pi_{|U_i ^\prime} ^* \epsilon_i \omega)
\end{equation*}
The sums are finite (\cite{DR} Th\'{e}or\`{e}me 1) and $\tilde T$ is linear. 

Consider a sequence of forms $(\omega_n)$ with $\text{supp}(\omega_n)\subset K\subset V$, where $K$ is compact, $V$ a coordinate patch and $\omega_n\rightarrow 0$ uniformly together with all its derivatives. By the Leibniz-rule $\epsilon_i \omega_n$ converges in the same way as $\omega$. We have that $\pi_{|U_i ^\prime} ^{-1}(V)$ is again a coordinate patch on $X$ and, since $\pi_{|U_i ^\prime}$ is a diffeomorphism, $\pi_{|U_i ^\prime} ^* \epsilon_i \omega_n$ is a sequence of forms with $\text{supp}(\pi_{|U_i ^\prime} ^* \epsilon_i \omega_n)\subset K^\prime\subset V^\prime$, where $K^\prime$ is compact, $V^\prime$ a coordinate patch and $\pi_{|U_i ^\prime} ^* \epsilon_i \omega_n\rightarrow 0$ uniformly together with all its derivatives. Since $T$ is continuous, we get
$\tilde T(\omega_n)=\sum_i T(\pi_{|U_i ^\prime} ^* \epsilon_i \omega_n)\rightarrow 0$
as $n$ tends to infinity.
Therefore $\tilde T$ is continuous.

Up to now $\tilde T$ seems to depend on the partition of unity $(U_i,\epsilon_i)$ and on the choice of the $U_i ^\prime$. We will show that this is not the case.
We start with the independence of choice of the $U_i ^\prime$. Given $U_i ^\prime $ any other set mapping via $\pi$ diffeomorphic to $U_i$ is of the form $r_g(U_i ^\prime) $. We get:
\begin{equation*}
 T(\pi_{|r_g(U_i ^\prime)}^*\epsilon_i \omega)= T((\pi_{|U_i ^\prime}\circ r_g ^{-1})^*\epsilon_i \omega)=T(r_g ^{-1*}\pi_{|U_i ^\prime}^*\epsilon_i \omega)=T(\pi_{|U_i ^\prime}^*\epsilon_i \omega),
 \end{equation*}
since $T$ is $G$-invariant.
Next we show the independence of the partition of unity.
We assume now that we have two partitions of unity $(U_i,\epsilon_i )$ and $(V_j,\eta_j )$. By the first step we may assume that whenever $U_i\cap V_j\neq \emptyset$ we have $U_i ^\prime\cap V_j ^\prime\neq \emptyset$. We may consider $(U_i \cap V_j ,U_i ^\prime \cap V_j ^\prime,\epsilon_i \eta_j)$. It is a partition of unity as the others. We get
\begin{equation*}
\sum_iT(\pi_{|U_i ^\prime}^*\epsilon_i \omega)=\sum_iT(\pi_{|U_i ^\prime}^*\epsilon_i \sum_j\eta_j\omega)=
\sum_{i,j}T(\pi_{|U_i ^\prime\cap V_j ^\prime}^*\epsilon_i \eta_j\omega),
\end{equation*}
from which the independence of the partition follows. The construction is of course compatible with restrictions in $U$. This means that we have the desired morphism of sheaves $(\pi_*\mathscr{D}_X)^G\rightarrow \mathscr{D}_{X/G}$, $T\mapsto \tilde T$. It remains to show that it is an isomorphism. This is a local problem. Choose a coordinate patch $U\subset X/G$ such that $\pi^{-1}(U)=\coprod_{g\in G} r_g(U^\prime)$ with all $r_g(U^\prime)$ coordinate patches. One has the isomorphism
\begin{equation*}
(\pi_*\mathscr{D}_X)^G(U)=(\prod_{g\in G} \mathscr{D}(r_g(U^\prime)))^G\rightarrow \mathscr{D}(U^\prime),\ (T_g)\mapsto T_{\id}
\end{equation*}
with inverse $T\mapsto (g^{-1}T)$. Now $\pi_{^|U^\prime}: U^\prime\rightarrow U$ is a diffeomorphism and we see by pushforward $\pi_{|U^\prime *}:\mathscr{D}(U^\prime)\cong \mathscr{D}(U)$ that the composite map
\begin{equation*}
\mathscr{D}(U^\prime)\cong(\pi_*\mathscr{D}_X)^G(U)\rightarrow \mathscr{D}_{X/G}(U), T\mapsto \pi_{|U^\prime *}T=\tilde T.
\end{equation*}
is an isomorphism.
This completes the proof.
\end{proof}
\end{proposition} 
\begin{remark}
\begin{enumerate}
\item Differentiation of currents is a $G$-equivariant morphism.
\item 
$(\pi_*\mathscr{D}_X)^G\rightarrow \mathscr{D}_{X/G}$, $T\mapsto \tilde T$
is a morphism of complexes with respect to the differentiation of currents. Let $\omega$ be a form with compact support on $X/G$ and $T$ a invariant current on $X$. As the problem is local, we may suppose that $\text{supp}(\omega)\subset U$, with $\pi_{|U^\prime}:U^\prime\rightarrow U$ a diffeomorphism. Then 
\begin{equation*}
\widetilde{dT}(\omega)=(dT)(\pi_{|U ^\prime}^*\omega)=(-1)^{deg(T)}T(d\pi_{|U ^\prime}^*\omega)=
(-1)^{deg(T)}T(\pi_{|U ^\prime}^* d\omega)=d\tilde T(\omega)
\end{equation*}
\item The whole theory of also works for non-continuous functionals. 
\end{enumerate}
\end{remark}

\chapter{Polylogarithmic Eisenstein classes for Hilbert-Blumenthal varieties}

We will recall the Hilbert-Blumenthal varieties $\Sc_K$ also considered by \cite{Ha1}. Our goal is to give a direct $\Q$-rational construction of Harder's Eisenstein cohomology classes. To do so we construct a fiber bundle $\Mc_K$ over $\Sc_K$. Over $\Mc_K$ there is a family of topological tori and we can apply our construction of polylogarithmic Eisenstein classes. In order to get classes on $\Sc_K$ in different cohomological degrees we have to decompose the cohomology of $\Mc_K$.

\section{The topological situation}

\subsection{Notation}

Let $F$ be a totally real number field of degree $[F:\Q]=\xi$ and $\Oc\subset F$ its ring of integers. We denote places of $F$ by $\nu$. If $\nu$ is an archimedean place, we write $\nu|\infty$ and call these places also the infinite places. Non-archimedean places correspond to maximal ideals $\mathfrak p\subset \Oc$ and we call these places the finite places. If we need to emphasize that $\nu$ corresponds to the maximal ideal $\mathfrak p$, we write $\nu_{\mathfrak p}$ instead of $\nu$. We have $\mathfrak p\cap \Z=(p)$ for a prime number $p\in\N$ and write $\nu| p$ in this case.
Infinite places may be identified with field embeddings $\sigma:F\rightarrow \R$. If $A$ is any $\R$-algebra, we have an isomorphism
\begin{equation*}
F_A:=F\otimes_\Q A= \prod_{\sigma:F\rightarrow \R} A=\prod_{\nu|\infty}A,\ x\otimes a\mapsto (\sigma(x)a)_\sigma.
\end{equation*}
In this manner we may speak of the $\nu$-component of $z\in F\otimes_\Q A$.
The completion of $F$ with respect to a place $\nu$ is denoted by $F_\nu$. If $\nu|p$, $\Oc_\nu$ is the integral closure of $\Z_p$ inside $F_\nu$ and $\mathfrak p_\nu=(\pi_\nu)\subset\Oc_\nu$ is the maximal ideal with residue field $\kappa(\mathfrak p _\nu):=\Oc_\nu/\mathfrak p_\nu=\Oc/\mathfrak p=:\kappa(\mathfrak p)$. 
For any place $\nu$ we have a normalized absolute value:
If $\nu|\infty$, $|x|_\nu$ is the usual absolute value of $x\in F_\nu=\R$. If $\nu$ is finite, we have $|x|_\nu:=|\kappa(\mathfrak p_\nu)|^{-\nu(x)}=|\Oc_\nu/(x)|^{-1}$ provided $x\in \Oc_\nu\setminus \left\{0\right\}$. $\nu$ is considered here as a discrete valuation on $F_\nu$ and may also be defined by $x=u\pi_\nu^{\nu(x)}$,  $u\in \Oc_\nu ^{\times}$.

We have the ring of adeles $\A_F$ and the group of ideles $\I_F$ over $F$. We set $\A=\A_\Q$ and $\I_\Q=\I$.
We have a decomposition into a finite and an infinite part
$\A_F=\A_{F,\infty}\times \A_{F,f}$, $\I_F=\I_{F,\infty}\times \I_{F,f}$.
We define $\hat\Z:=\varprojlim_{n\in\N} \Z/n\Z$ to be the profinite completion of $\Z$ and recall $\A_{F,f}=\hat\Z\otimes_\Z F$. We also set $\hat \Oc:=\hat\Z\otimes_\Z \Oc$.

We consider $G_0:=GL_{2}$ and $V_0:=\mathbb \G_{a} ^2$ as algebraic groups over $spec(\Oc)$, where the first group acts naturally on the second by matrix multiplication from the left. For any algebraic subgroup $H_0 \subset V_0\rtimes G_0$ we write $H:=Res_{\Oc/\Z}H_0$ for the restriction of scalars. We have for example $B_0\subset G_0$ the standard Borel of upper triangular matrices, $U_0\subset B_0$ its unipotent radical, $T_0\subset B_0$ the maximal torus and $Z_0\subset G_0$ the center. For any ring $R$ we write $H(R)$ for the group of $R$ valued points of $H$.\\
We have $H(\A)\subset \prod_\nu H(F_\nu)$, so $h=(h_\nu)\in H(\A)$ is determined by its $\nu$-components $h_\nu\in H(F_\nu)$. 
$h\in H(\A)$ may also be uniquely written as $h=h_\infty h_f$ with $h_f\in H(\A_f)$ and $h_\infty\in H(\R)$. We will deal with closed subgroups $K\subset H(\A)$ of the form $K=K_\infty\times K_f$ with closed subgroups $K_f\subset H(\A_f)$ and $K_\infty \subset H(\R)$. We consider $K_\infty$ as a Lie group in a natural way and for any Lie group $L$ we denote the connected component of the identity by $L^0$.
We will identify $Z(\R)=(F\otimes_\Q\R)^{\times}=\prod_{\nu |\infty}\R^{\times}$ and $\R^{\times}$ embedded diagonally.
We have $G(\R)=G_0(F\otimes_\Z \R)=\prod_{\nu|\infty} GL_2(\R)$
and subgroups
\begin{equation*}
K_\infty:=\prod_{\nu|\infty}K_\nu,\ K_\infty ^1:=\R_{>0}\cdot\prod_{\nu|\infty}K_\nu ^1,
\end{equation*}
where $K_\nu:=\R_{>0} SO(2)$, $K_\nu ^1:=SO(2)$ and $\R_{>0}$ the positive real numbers. We identify the group of connected components of $GL_2(\R)$ with $\left\{\begin{pmatrix}1 & 0\\0 & \epsilon\end{pmatrix}:\ \epsilon=\pm 1\right\}$. In this way we may view $\left\{\pm 1\right\}^\xi\cong\pi_0(G(\R))\subset G(\R)$ as a subgroup normalizing $K_\infty$.
\subsection{The spaces}
Take a compact open subgroup $K_f\subset G(\A_f)$. Our main examples will be
\begin{equation*}
\ker\left(G(\hat\Z)\rightarrow G(\Z/N\Z)\right),\ N\in \N,
\end{equation*}
and its conjugates.
Set $K:=K_\infty K_f\subset G(\A)$, $K^1:=K_\infty ^1 K_f\subset G(\A)$ and 
\begin{equation*}
K_N:=\ker\left(G(\hat\Z)\rightarrow G(\Z/N\Z)\right)K_\infty.
\end{equation*}
We define spaces 
\begin{equation*}
\Sc_K:=K\backslash G(\A)/G(\Q)\ \text{and }\Mc_K:=K^1\backslash G(\A)/G(\Q)
\end{equation*}
by taking double quotients. These two spaces are related by the obvious projection map $\phi_K:\Mc_K\rightarrow \Sc_K$.

\begin{remark}\label{spaces_facts}
Let us recall some facts which can be found mainly in \cite{Ha1} 1.0 and 1.1.
\begin{itemize}
\item 
$\Sc_K$ parametrizes abelian varieties $A$ over $spec(\C)$ of dimension $\xi$ with real multiplication by $\Oc$ and a level-$K_f$-structure (\cite{Mi}, Theorem 1.2). $\Sc_K$ can be realized as the $\C$-points of a scheme defined over $\Q$. If $\xi>1$, there is no universal abelian scheme above $\Sc_K$, as there are always non-trivial automorphisms of our abelian varieties respecting the $\Oc$-structure and the level-$K_f$-structure . To see this consider
\begin{equation*}
Z_K:=Z(\Q)\cap K_f\subset Z(\Q)\cap G(\hat\Z)=\Oc^{\times}
\end{equation*}   
a subgroup of finite index. If $\xi>1$, $Z_K$ is never trivial, because the rank of this abelian group is positive by Dirichlet's unit theorem. 
\item $K_f\backslash G(\A_f)$ is discrete and countable.
\item Set $\H^\xi_{\pm}:=(F\otimes_\Q \C)^{\times}\setminus (F\otimes_\Q \R)^{\times}$ and denote the the connected component of $1\otimes i$ by $\H^\xi$. We identify $1\otimes i$ with the imaginary unit $i$ and get the diffeomorphism
\begin{equation*}
K_\infty\backslash G(\R) \rightarrow \H^{\xi}_{\pm},\ g=\begin{pmatrix} a&b\\c&d\end{pmatrix}\mapsto \frac{b+id}{a+ic}=:\tau=\tau(g)
\end{equation*}
So $G(\R)$ acts on $\H^{\xi}_{\pm}$ from the right by $\tau\cdot g:=\frac{b+d\tau }{a+c\tau}$. We also have the diffeomorphism
\begin{equation*}
K_\infty ^1\backslash G(\R) \rightarrow \H^{\xi}_{\pm}\times (\R_{>0}\backslash Z(\R)^0),\ g\mapsto (i\cdot g, |\det(g)|).
\end{equation*}
The inverse map is
\begin{equation*}
(\tau=x+iy,r)\mapsto \sqrt{r|y|^{-1}}\cdot \begin{pmatrix}1&x\\0&y\end{pmatrix}.
\end{equation*}
The absolute value and the square root are taken componentwise. 

\item The stabilizer of a connected component of $K\backslash G(\A)$ or $K^1\backslash G(\A)$ in $G(\Q)$ is given by $G(g_f):=\left(g_f ^{-1}K_f g_f \cdot G(\R)^0\right)\cap G(\Q)$ for a $g_f\in G(\A_f)$. Therefore, the connected components of $\Sc_K$ are of the form $\H^\xi/G(g_f)$ and this means that they are \textit{Hilbert-Blumenthal varieties}, in other words, quotients of $\H^\xi$ by arithmetic subgroups of $GL_2(F)$.
\item $G(\Q)/Z_K$ acts properly discontinuously on $K\backslash G(\A)$ and fixpoint free, if $K$ is small enough, say $K=K_N$ and $N\geq 3$. Using the diffeomorphism above we see that $G(\Q)$ acts properly discontinuously and without fixpoints on $K^1\backslash G(\A)$, if we suppose $K_f$ to be small enough. 
\item 
Consider the ray class group $Cl_F ^{K}:=\det(K_f)\prod_{\nu|\infty}\R_{>0}\backslash \I_F/F^{\times} $. The determinant induces maps
$\det:\Sc_K\rightarrow Cl_F ^{K}$ and $\det:\Mc_K\rightarrow  Cl_F ^{K}$.
Since we have strong approximation for $SL_2$, the fibers of these maps are the connected components of $\Sc_K$ and $\Mc_K$. So our spaces have a finite number of connected components, since the class number of $F$ is finite.
\item If $K_f$ is small enough, $\Mc_K$ has a structure of a $\Cc^\infty$-manifold making the canonical projection $K^1\backslash G(\A)\rightarrow \Mc_K$ a local diffeomorphism.  
\end{itemize} 
We will always assume that we have chosen $K_f$ small enough so that the quotient spaces $\Sc_K$ and $\Mc_K$ have the natural structure of $\Cc^\infty$-manifolds. Moreover, we consider
\begin{equation*}
\left\{\pm 1\right\}^\xi=\prod_{\nu|\infty }\left\{\pm1\right\}\subset Z(\R)=\prod_{\nu|\infty}\R^{\times}.
\end{equation*}
\end{remark}
\begin{lemma}
$\phi_K:\Mc_K\rightarrow \Sc_K$ is a principal $\R_{>0}\left\{\pm 1\right\}^\xi\backslash Z(\R)/Z_K$-bundle.
\begin{proof}
$Z(\R)$ acts clearly by matrix multiplication on $\Mc_K$ from the right, as it lies in the center of $G(\A)$. This action factors over the group $\R_{>0}\left\{\pm 1\right\}^\xi\backslash Z(\R)/Z_K$, is by diffeomorphisms and preserves the fibers of $\phi_K$. What is left to be shown is that $\phi_K$ is a fiber bundle and that the action of $\R_{>0}\left\{\pm 1\right\}^\xi\backslash Z(\R)/Z_K$ on the fibers is simply transitive.
This problem is local on the base.  
\begin{equation*}
\H^\xi _{\pm}/G(g_f)\subset \Sc_K,\ \tau G(g_f)\mapsto  K\begin{pmatrix}1&x\\0&y\end{pmatrix}g_f G(\Q)
\end{equation*}
is the inclusion of some connected components. We have
\begin{equation*}
\phi_K ^{-1}(\H^\xi _{\pm}/G(g_f))=\left(\H^\xi _{\pm}\times (\R_{>0}\backslash Z(\R)^0)\right)/G(g_f),
\end{equation*}
where $G(g_f)$ acts on the second factor by multiplication with the absolute value of the determinant. Here
\begin{equation*}
\phi_K:\left(\H^\xi _{\pm}\times (\R_{>0}\backslash Z(\R)^0)\right)/G(g_f)\rightarrow \H^\xi _{\pm}/G(g_f)
\end{equation*}
is induced by the first projection. The cosets $\gamma Z_K\in G(g_f)/Z_K$ act properly discontinuously and fixpoint free on $\H^{\xi}_{\pm}$, so we find for any point in $\H^{\xi}_{\pm}$ an open neighborhood $U$ such that $U\gamma \cap U\neq \emptyset$ implies $\gamma\in Z_K $. If we denote the image of $U$ in $\Sc_K$ by $U^\prime$, we get a diffeomorphism $\phi_K^{-1}(U^\prime)\cong U^\prime \times  (\R_{>0}\backslash Z(\R)^0/\det(Z_K))$. So $\phi_K$ is a fiber bundle. But
\begin{equation*}
\det:\R_{>0}\left\{\pm 1\right\}^\xi\backslash Z(\R)/Z_K\rightarrow \R_{>0}\backslash Z(\R)^0/\det(Z_K)
\end{equation*}
is a diffeomorphism showing that the action of $\R_{>0}\left\{\pm 1\right\}^\xi\backslash Z(\R)/Z_K$ on the fibers of $\phi_K$ is simply transitive. This completes the proof.
\end{proof}
\end{lemma}
\begin{remark}
If $K=K_N$, $N\geq 3$, we have $-1\notin Z_K$ and Dirichlet's unit theorem tells us that $Z_K$ acts properly discontinuously and fixpoint free on $\R_{>0}\left\{\pm 1\right\}^\xi\backslash Z(\R)$ and that the quotient $\R_{>0}\left\{\pm 1\right\}^\xi\backslash Z(\R)/Z_K$ is compact.
\end{remark}
Let us now define families of topological tori above $\Mc_K$. Take $W_f\subset V(\A_f)\rtimes G(\A_f)$ a compact open subgroup such that the projection of $W_f$ onto $G(\A_f)$ gives $K_f$. Set $W:=W_f\cdot K_\infty ^1 \subset G(\A)$ where we always consider $G(\A)\subset V(\A)\rtimes G(\A) $ via $g\mapsto (0,g)$. Set
\begin{equation*}
\pi_W:\Tc_W:=W\backslash V(\A)\rtimes G(\A)/V(\Q)\rtimes G(\Q)\rightarrow \Mc_K,\ (v,g)\mapsto g 
\end{equation*}

\begin{lemma}
$\pi_W:\Tc_W\rightarrow \Mc_K$ is a family of topological tori.
\begin{proof}
Write $W_f=V_f\rtimes K_f$, where $V_f\subset V(\A_f)$ is a compact open subgroup.
We have the bijection
\begin{equation*}
V(\A)\rtimes  G(\A)\rightarrow V(\A)\times  G(\A),\ (v,g)\mapsto (g^{-1}v,g)=:(w,g).
\end{equation*}
In the new coordinates $(w,g)$ the group actions look like 
\begin{equation*}
(w,g)(q,\gamma)=(\gamma^{-1}(w+q),g\gamma),\ (q,\gamma)\in V(\Q)\rtimes G(\Q),
\end{equation*}
\begin{equation*}
(u,k)(w,g)=((kg)^{-1}u +w,kg),\  (u,k)\in W. 
\end{equation*}
When we write $V(\A)\times G(\A)$ we always work in the coordinates $(w,g)$ with induced group actions.
The pullback of $\Tc_W$ to $K^1\backslash G(\A)$ is $W\backslash V(\A)\times G(\A)/V(\Q)$. The action of $V_f$ on $V(\A)$ depends on $g_f\in G(\A_f)$. For $g_f\in G(\A_f)$ consider $K^1\backslash K^1g_f G(\R)$. This is just a collection of connected components in $K^1\backslash G(\A)$. 
We may now identify
\begin{equation*}
pr^{-1}(K^1\backslash K^1g_f G(\R))\cong (g_f^{-1}V_f\backslash V(\A)/V(\Q))\times (K^1\backslash K^1g_f G(\R))
\end{equation*}
\begin{equation*}
W(w,h)V(\Q)\mapsto (g_f ^{-1}V_f+w+V(\Q),K^1h)
\end{equation*}
and the fiber
$
g_f ^{-1}V_f\backslash V(\A)/V(\Q)=V(\R)/(g_f ^{-1}V_f\cap V(\Q))
$
over each point
is a topological torus. This shows that $W\backslash V(\A)\times G(\A)/V(\Q)\rightarrow K\backslash G(\A)$ is a family of topological tori. Because $G(\Q)$ acts on the latter by homomorphisms and properly discontinuously without fixpoints on the base, $\Tc_W$ is a family of topological tori.
\end{proof}
\begin{remark}\label{connected_component}
We keep the notation $(w,g):=(g^{-1}v,g)$, if $(v,g)\in V(\A)\rtimes G(\A)$ denotes a general element.
Moreover, consider $K_f\subset G(\hat\Z)$ and $W_f=V(\hat\Z)\rtimes K_f$.
Over 
\begin{equation*}
\left(\H^{\xi}_{\pm}\times (\R_{>0}\backslash Z(\R) ^0)\right)/G(g_f)
\end{equation*}
our torus looks like
\begin{equation*}
\left(F\otimes_\Q\C\times\H^{\xi}_{\pm}\times (\R_{>0}\backslash Z(\R)^0)\right)/V(g_f)\rtimes G(g_f), 
\end{equation*}
where $V(g_f):=g_f ^{-1}V(\hat\Z)\cap V(\Q)$ and $V(g_f)\rtimes G(g_f)$ acts in the following way:
\begin{equation*}
(z,\tau,r)\cdot\left(\begin{pmatrix}l_1\\l_2\end{pmatrix},\gamma=\begin{pmatrix}a& b\\c&d\end{pmatrix}\right)=\left(\frac{z+l_1+\tau l_2}{a+c\tau},\tau\cdot\gamma,r\cdot \det(\gamma)\right).
\end{equation*}
The natural inclusion is then given by
\begin{equation*}
 \left(F\otimes_\Q\C\times\H^{\xi}_{\pm}\times (\R_{>0}\backslash Z(\R)^0)\right)/V(g_f)\rtimes G(g_f)\rightarrow \Tc_W,\ (z,\tau,r)\mapsto (w,g)
\end{equation*}
\begin{equation*}
z=u^1 +u^2\tau,\ (w,g)_\infty=\left(\begin{pmatrix}u^1\\u^2 \end{pmatrix},\sqrt{r|y|^{-1}}\begin{pmatrix}1& x\\0&y\end{pmatrix}\right)\text{ and } (w,g)_f=(0,g_f).
\end{equation*}
\end{remark}
\end{lemma}

\section{Construction of adelic polylogarithmic Eisenstein classes}

Next we will discuss how we can associate not only to formal linear combinations of torsion sections of tori but also to Schwartz-Bruhat functions $f$ on $V(\A_f)$ polylogarithmic Eisenstein classes $\Eis(f)$. The naturality of the polylogarithmic construction guarantees the $G(\A_f)$-equivariance of the operator $\Eis$.

\subsection{Coefficient sheaves and integral structures}

We introduce some locally constant sheaves on our spaces. Given a topological space $X$ on which a group $G$ acts continuously from the right we associate to any $G$-(left)-module $V$ a sheaf $\Vc$ on the quotient space $X/G$ in a functorial way. We make the constant sheaf $V$ on $X$ into a $G$-sheaf by defining for $g\in G$ morphisms
\begin{equation*}
g:V\rightarrow r_{g*} V,\ f\mapsto \left\{x\mapsto gf(xg)\right\},
\end{equation*} 
if $f:U\rightarrow V$ is a locally constant map on $U\subset X$. If $p:X\rightarrow X/G$ is the canonical map, these morphisms make $p_*V$ into a sheaf of $G$-modules and we can consider the fixed sheaf $p_* ^GV $ which is defined by $p_* ^GV (U):=p_*V(U)^G$. This fixed sheaf is our $\Vc$. We will often just speak about the sheaf induced by a $G$-module or associated to a $G$-module.

On $\Tc_W$ we have the logarithm sheaf $Log_W$ with $\Q$-coefficients. It is induced by the $V(\Q)\rtimes G(\Q)$-left-module $\prod_{k\geq0}\Sym^k V(\Q)$, where $V(\Q)$ acts by multiplication with the exponential series and $G(\Q)$ by its standard action, see \Cref{R_alg}. \\
The sheaf $Log_W$ carries an integral structure. Suppose we are on a connected component cut out by $g_f\in G(\A_f)$ and $K_f\subset G(\hat\Z)$, $W_f=V(\hat\Z)\rtimes K_f$, compare \Cref{connected_component}. The fundamental group of the family of topological tori above this connected component is $V(g_f)\rtimes G(g_f)$. $Log_W$ is the locally constant sheaf associated to the $V(g_f)\rtimes G(g_f)$- module $\Z[[V(g_f)]]$ as considered in \Cref{def_Log}.\\
Moreover, we have the locally constant sheaves $\Sym^k\Hc_K$ on $\Mc_K$ associated to the $G(\Q)$-modules $\Sym^kV(\Q)$ and the sheaf $\Sym^k\Hc^\prime _K$ on $\Sc_K$, which is associated to the same $G(\Q)$-module. These sheaves also have integral structures coming from $\Sym^k V(g_f)_{PD}$, which we denote by $\Sym^k\Hc_{PD}$ and $\Sym^k\Hc^{\prime}_{PD}$. 
\begin{remark}
One has the identification $\phi_{K*}\Sym^k\Hc_K=\Sym^k\Hc^\prime _K$.
This is best understood by looking at the stalks of the sheaves. We have $\Sym^k\Hc_{K,x}=\Sym^kV(\Q)$, since the action of $G(\Q)$ on $K^1\backslash G(\A)$ is without fixpoints. But we have $\Sym^k\Hc_{K,x}^\prime=\Sym^kV(\Q)^{Z_K}$, as $Z_K$ acts trivially on $K\backslash G(\A)$. 
\end{remark}
If we consider sheaves on $\Mc_K$ arising from $G(\Q)$-modules, the sheaves on $\Sc_K$ corresponding to the same module will always be denoted by a prime. \\
The sheaf $\mu_K$ on $\Mc_K$ is induced by the $G(\Q)$-module $\Q$, where $G(\Q)$ acts by multiplication with the norm of the determinant $\Norm\circ \det$.  As $\Norm\circ \det(Z_K)=1$, there is no need to distinguish between $\mu_K$ on $\Mc_K$ and $\mu_K ^\prime$ on $\Sc_K$, as $\mu_K$ would be the pullback of $\mu_K ^\prime$. We have a trivialization of $\mu_K$ given by the global section 
\begin{equation*}
s_K: K \backslash G(\A) \rightarrow \Q,\ g\mapsto \left\|\det\ g_f\right\|_f\cdot \sgn(\Norm(\det\ g_\infty))^{-1},
\end{equation*}
where $\left\|\right\|_f:=\prod_{\nu\nshortmid \infty}|\ |_\nu$ is the product of the finite local norms of $F$. Clearly, $or_{\Tc_W/\Mc_K}=\mu_K=\bigwedge^{2\xi}\Hc_K$.

As we vary the level on $\Mc_K$ or $\Tc_W$, we see that our sheaves $Log_W$ are obtained by pullback. For example, given $W_f ^\prime\subset W_f$ we have the canonical projection $p_W ^{W^\prime}:\Tc_{W^\prime}\rightarrow\Tc_W $ and $p_W ^{W^\prime  *}Log_W=Log_{W^\prime}$ by \Cref{basechange_Log}. The same is true for $\Sym^k\Hc_K$ and $\mu_K$. We will drop the subscript $K$, when it does not cause any confusion.

\subsection{Torsion sections}\label{torsion sections}

We are looking for a good class of compact open subgroups $W_f\subset V(\A_f)\rtimes G(\A_f)$ such that we can give explicit sections for $\pi_W$.
 
From now on we will always deal with $W_f$ of the kind: $W_f=V_K\rtimes K_f$, where $K_f\subset G(\A_f)$, $V_K\subset V(\A_f)$ are compact open subgroups and $V_K$ is an $\hat{\Oc}$-module on which $K_f$ acts as usual by matrix-multiplication. Moreover, we suppose that we have a compact open $\hat{\Oc}$-module $V_K ^\prime\supset V_K$ such that $K_f$ acts trivially on $V_K ^\prime /V_K$, $|V_K ^\prime /V_K|>1$. Our main example will be $K=K_{N^2}$, $V_K=NV(\hat{\Z})$ and $V_K ^\prime=N^{-1}V(\hat{\Z})$. Given such a $W$ we have a natural injective arrow
$V_K ^\prime/V_K\times \Mc_K\rightarrow \Tc_W$, which is induced by the inclusion $V_K ^\prime\times G(\A)\rightarrow V(\A)\rtimes G(\A)$.\\
We define $D_W$ as the image of this map without the image of the zero section. $D_W$ is a closed submanifold of the torus $\Tc_W$, which is a finite covering of $\Mc_K$. We set $U_W$ to be the open complement of $D_W$ in $\Tc_W$.

Let us apply the construction of the polylogarithm to this situation. We have the natural short exact sequence
\begin{equation*}
0\rightarrow H^{2\xi-1}(U_W,Log_W\otimes_\Q \pi_W ^{-1}\mu_K)\rightarrow H^0(D_W,Log_{W\ |D_W})\rightarrow H^0(\Mc_K,\Q)
\end{equation*}
induced by the localization sequence for the pair $(D_W,\Tc_W)$. 
We have the inclusion $H^0(D_W,\Q)\subset H^0(D_W,Log_{W\ |D_W})$ and on this submodule the last arrow in the short exact sequence is just given by the trace map $H^0(D_W,\Q)\rightarrow H^0(\Mc_K,\Q)$. Exactness of the sequence ensures that given a section $f\in H^0(D_W,\Q)$ of trace $0$ we have a unique element $\pol_W(f)\in H^{2\xi-1}(U_W,Log_W\otimes_\Q \pi_W^{-1}\mu_K)$
mapped to $f$. As $D_W\cap \image(0)=\emptyset$, we may specialize $\pol_W(f)$ along the zero section of $\Tc_W$ to get a class 
\begin{equation*}
\Eis_W(f)\in H^{2\xi-1}(\Mc_K,\prod_{k\geq0}\Sym^k\Hc_K\otimes \mu_K)= H^{2\xi-1}(\Mc_K,0^{-1}Log_W\otimes \mu_K),
\end{equation*}
in other words $(\Eis_W ^k(f))_k\in \prod_{k\geq0}H^{2\xi-1}(\Mc_K,\Sym^k\Hc_K\otimes \mu_K)$.

\subsection{Definition of adelic Eisenstein classes}

We want to extend the idea of polylogarithmic Eisenstein classes to more general functions $f$.
Let us denote the space of Schwartz-Bruhat functions on $V(\A_f)$ with values in a subfield $k\subset \C$ by $\Sc(V(\A_f),k)$, compare for example \cite{Ge} p.267. Set
\begin{equation*}
\Sc(V(\A_f),k)^0=\left\{f\in \Sc(V(\A_f),k):f(0)=0,\int_{V(\A_f)}f dv_f=0 \right\}
\end{equation*}
where $dv_f$ is any Haar measure on $V(\A_f)$.
$\Sc(V(\A_f),k)^0$ is a right $G(\A_f)$-module, when we set as usual $f\cdot g_f(v)=f(g_fv)$. This follows from the integral transformation formula
\begin{equation*}
\int_{V(\A_f)}f(g_f v) dv_f(v)= \left\|\det(g_f)\right\|_f ^{-1}\int_{V(\A_f)}f( v) dv_f(v).
\end{equation*}
We will need a more general space of functions. For any $n\in \Z$ denote by \\$\Sc(V(\A_f),\mu^{\otimes n}\otimes k)$ the space of functions $f:V(\A_f)\rtimes G(\A_f)\times \pi_0(G(\R))\rightarrow k$ such that 
\begin{equation*}
 \forall x\in G(\A_f)\times \pi_0(G(\R)):\ f(\ ,x )\in \Sc(V(\A_f),k),
\end{equation*}
there is some compact open subgroup $K_f\subset G(\A_f)$ such that $f$ factors as
\begin{equation*}
f(\ ,\ ):K_f\backslash G(\A_f)\times \pi_0(G(\R))  \rightarrow \Sc(V(\A_f),k),\ x\mapsto f(\ ,x )
\end{equation*}
and for all $v\in V(\A_f)$, $x\in G(\A_f)\times \pi_0(G(\R))$ and $\gamma\in G(\Q)$
\begin{equation*}
f(v,x\gamma)=N(\det(\gamma))^{-n}f(v,x).
\end{equation*}
We define $\Sc(V(\A_f),\mu^{\otimes n}\otimes k)^0$ as the subspace of $\Sc(V(\A_f),\mu^{\otimes n}\otimes k)$ consisting of those functions $f$ with $f(\ ,x)\in \Sc(V(\A_f),k)^0$ for all $x\in G(\A_f)\times \pi_0(G(\R))$. Again $\Sc(V(\A_f),\mu^{\otimes n}\otimes k)^0$ is a right $G(\A_f)\times \pi_0(G(\R))$-module, when we take the action induced by the canonical left $G(\A_f)\times \pi_0(G(\R))$ action on $V(\A_f)\rtimes G(\A_f)\times \pi_0(G(\R))$. If $k=\Q$ we simply set $\Sc(V(\A_f),\mu^{\otimes n}\otimes k)=:\Sc(V(\A_f),\mu^{\otimes n})$ and  $\Sc(V(\A_f),\mu^{\otimes n}\otimes k)^0 =:\Sc(V(\A_f),\mu^{\otimes n})^0 $. 
\begin{lemma}\label{varphi}
$\Sc(V(\A_f),\mu^{\otimes n}\otimes\C)^0$ is generated by $\varphi$ of the form
\begin{equation*}
\varphi(v,g)=f(v)\cdot\eta(\det(g))\left(\left\|\det(g)\right\|_f\sgn(N(\det(g)))\right)^{n},
\end{equation*}
for $f\in\Sc(V(\A_f),\C)^0$ and $\eta:\prod_{\nu|\infty}\R_{>0}\backslash\I_{F}/F^{\times}\rightarrow \C^{\times}$ a Dirichlet character.
\begin{proof}
Consider $\varphi\in\Sc(V(\A_f),\mu^{\otimes n}\otimes\C)^0$ arbitrarily. Write
\begin{equation*}
\varphi(v,g)=\phi(v,g)\left(\left\|\det(g)\right\|_f\sgn(N(\det(g)))\right)^{n},
\end{equation*}
for $\phi\in\Sc(V(\A_f),\mu^{\otimes 0}\otimes \C)^0$. $\phi$ factors in the second argument over some $K_f\subset G(\A_f)$ compact open, in other words
\begin{equation*}
\phi:K_f\backslash G(\A_f)\times \pi_0(G(\R))/G(\Q)\rightarrow \Sc(V(\A_f),\C)^0.
\end{equation*}
We have strong approximation for $SL_{2,F}$ and therefore
\begin{equation*}
\det: K_f\backslash G(\A_f)\times \pi_0(G(\R))/G(\Q)\rightarrow \det(K_f)\prod_{\nu|\infty}\R_{>0}\backslash\I_{F}/F^{\times}=Cl_F ^{K}
\end{equation*}
is an isomorphism. In this manner we may interpret $\phi$ in the second argument as a function on $Cl_F ^{K}$. The ray class group $Cl_F ^{K}$ is a finite group of cardinality $h_K$ and we can describe functions $Cl_F ^{K}\rightarrow \C$ as finite linear combinations of characters by the Fourier inversion theorem.
We get explicitly
\begin{equation*}
\phi(v,g)=h_{K}^{-1}\sum_{\eta\in \widehat{Cl_F ^{K}}}\hat{\phi}(v,\eta)\cdot\eta(\det(g)),\ \hat{\phi}(v,\eta)=\sum_{x\in Cl_F ^{K}}\phi(v,x)\overline{\eta(x)}.
\end{equation*}
\end{proof}
\end{lemma}
Finally, let us set
\begin{equation*}
\varinjlim_{K_f} H^\bullet(\Mc_K,\Sym^k\Hc_K\otimes\mu_K ^{\otimes n})=:H^\bullet(\Mc,\Sym^k\Hc\otimes\mu^{\otimes n}).
\end{equation*}
It is a $G(\A_f)\times \pi_0(G(\R))$-module as explained in \cite{Ha1} 1.2.
\begin{remark}
$\Sc(V(\A_f),\mu^{\otimes n})^0=\Sc(V(\A_f),\Q)^0\otimes_\Q H^0(\Mc,\mu^{\otimes n})$.
\end{remark}
\begin{proposition}\label{Eis^k}
We have a $G(\A_f)\times\pi_0(G(\R))$-equivariant operator 
\begin{equation*}
\Eis^k:\Sc(V(\A_f),\mu^{\otimes n})^0 \rightarrow H^{2\xi-1}(\Mc,\Sym^k\Hc\otimes\mu^{\otimes n+1}),
\end{equation*}
which is induced by the operators $\Eis^k _W$.
\begin{proof}
Given a $\phi\in \Sc(V(\A_f),\Q)^0$ there is an $M\in \N$ such that $\text{supp}(\phi)\subset M^{-1}V(\hat\Z)$ and an $M^\prime \in \N$ such that $\phi$ is well-defined on $V(\A_f)/M^\prime V(\hat\Z)$. So considering $N=MM^\prime$, we can identify $\phi$ with a function on $N^{-1}V(\hat\Z)/NV(\hat\Z)$.
Let us investigate $f\in \Sc(V(\A_f),\mu^{\otimes n})^0$. As $f$ factors over some compact open subgroup $K_f$ of $G(\A_f)$ and $K_f\backslash G(\A_f)\times \pi_0(G(\R))/G(\Q)$ is finite, we may find an $N\in \N$ such that $f$ is well-defined modulo $K_{N^2}$ in the second argument and $f$ behaves like $\phi$ above in the first argument. In other words, $f$ induces a function
\begin{equation*}
(N^{-1}V(\hat\Z)/NV(\hat\Z)\setminus \left\{0\right\})\times K_{N^2}\backslash G(\A_f)\times \pi_0(G(\R))\rightarrow \Q, 
\end{equation*}
which we also denote by $f$. 
We are now in the situation of our groups $W$ of \Cref{torsion sections}, namely $K_f=K_{N^2}$, $NV(\hat\Z)=V_K$, $V^\prime _K =N^{-1}V(\hat \Z)$ and $W_f=NV(\hat\Z)\rtimes K_{N^2}$. The function $f$ may be identified with a section $f\in H^0(D_W,\mu_K ^{\otimes n})$ which actually has trace zero:
We have for all $x\in G(\A_f)\times\pi_0(G(\R))$:
\begin{equation*}
\sum_{v\in N^{-1}V(\hat\Z)/NV(\hat\Z)\setminus \left\{0\right\}}f(v,x)=\sum_{v\in N^{-1}V(\hat\Z)/NV(\hat\Z)}f(v,x)=0 
\end{equation*}
if and only if
\begin{equation*}
0=\sum_{v\in N^{-1}V(\hat\Z)/NV(\hat\Z)}f(v,x)\int_{NV(\hat\Z)}dv_f=
\end{equation*}
\begin{equation*}
=\sum_{v\in N^{-1}V(\hat\Z)/NV(\hat\Z)}\int_{NV(\hat\Z)}f(v+u,x)dv_f(u)=\int_{V(\A_f)}f(v ,x)dv_f .
\end{equation*}
Consequently, we get $\Eis_W ^k(f)\in H^{2\xi-1}(\Mc_K,\Sym^k\Hc_K\otimes \mu_K ^{\otimes n +1})$. We want to define
 \begin{equation*}
\Eis^k(f)=\text{image of } \Eis_W ^k(f)\text{ in }H^{2\xi-1}(\Mc,\Sym^k\Hc\otimes \mu^{\otimes n+1}).
\end{equation*}
We have to show that this construction is independent of the choices made.
Consider again $f\in \Sc(V(\A_f),\mu^{\otimes n}) ^0$. Suppose that we have $W_1$ and $W_2$ both fulfilling the properties above. In other words, we have for $i=1,2$ compact open subgroups $K_{f,i}\subset G(\A_f)$ acting on the compact open $\hat\Oc$-modules $V_{K_i}^\prime\supset V_{K_i}$ and acting trivially on $V_{K_i}^\prime/ V_{K_i}$, such that
$f$ is well-defined modulo $K_i$ in the second argument, $f$ is well-defined modulo $V_{K_i}$ in the first argument and for all $x\in G(\A_f)\times \pi_0(G(\R))$ $\text{supp}(f(\ ,x))\subset V_{K_i}^\prime$. \\
If we set
\begin{equation*} 
K=K_1\cap K_2,\ V_K=V_{K_1}\cap V_{K_1},\ V_K ^\prime=V_{K_1} ^\prime\cap V_{K_1} ^\prime,\text{ and }W=W_1\cap W_2, 
\end{equation*}
all the properties are also fulfilled by $W$ so that we have $\Eis^k_W(f),\ \Eis^k_{W_1}(f)$ and $\Eis^k _{W_2}(f)$. To see that $\Eis^k$ is well-defined we have to show that $\Eis^k_W(f)$ is obtained by pullback along the canonical maps $q^W _{W_i}:\Mc_K\rightarrow \Mc_{K_i}$ from $\Eis_{W_i} ^k(f)$ for $i=1,2$. It is enough to consider $i=1$. We have a commutative square 
\begin{equation*} 
\begin{xy}
\xymatrix{
\Tc_W\ar[r]^{p^W _{W_1}} \ar[d] & \Tc_{W_1}\ar[d] \\
\Mc_K\ar[r]^{q^W _{W_1}} & \Mc_{K_1}
}
\end{xy}
\end{equation*}
Set $D_{W_1}^\prime :=p^{W -1 }_{W_1}(D_{W_1})\supset D_W$, it is a finite cover of $\Mc_K$. As the localization-sequence, $\mu$ and $Log$ are natural in pullbacks, see \Cref{natural_localization}, we get the following commutative diagram 
\begin{equation*} 
\begin{xy}
\xymatrix{
 H^{2\xi-1}(\Tc_{W_1}\setminus D_{W_1},Log\otimes\mu^{\otimes n+1})\ar[d]\ar[r] & H^0(D_{W_1},Log_{|D_{W_1}}\otimes \mu^{\otimes n})\ar[d] \\
H^{2\xi-1}(\Tc_{W}\setminus D_{W_1}^\prime ,Log\otimes\mu^{\otimes n+1})\ar[r] & H^0(D_{W_1}^\prime,Log_{|D_{W_1}^\prime}\otimes \mu^{\otimes n})
}
\end{xy}
\end{equation*}
The diagram shows that $p^{W * }_{W_1}\pol_{W_1}(f)$ corresponds to the section $p^{W * }_{W_1}f\in H^0(D_{W_1}^\prime,\mu^{\otimes n})$. By naturality of the localization sequence for the inclusion $i:D_W\subset D_{W_1}^\prime$, see \Cref{shrinking_localization}, we also get the commutative diagram
\begin{equation*} 
\begin{xy}
\xymatrix{
H^{2\xi-1}(\Tc_{W}\setminus D_{W},Log\otimes\mu^{\otimes n+1})\ar[d]^{\res}\ar[r] & H^0(D_{W},Log_{|D_{W_1}}\otimes \mu^{\otimes n})\ar[d]^{i_*} \\
H^{2\xi-1}(\Tc_{W}\setminus D_{W_1}^\prime ,Log\otimes\mu^{\otimes n+1})\ar[r] & H^0(D_{W_1}^\prime,Log_{|D_{W_1}^\prime}\otimes\mu^{\otimes n})
}
\end{xy}
\end{equation*}
where $i_*:H^0(D_{W},\mu^{\otimes n})\rightarrow H^0(D_{W_1}^\prime,\mu^{\otimes n})$ is given by extension by zero. $i_*f$ and $p^{W * }_{W_1}f\in H^0(D_{W_1}^\prime,\mu^{\otimes n})$ coincide and consequently $\res(\pol_W(f))=p^{W * }_{W_1}\pol_{W_1}(f)$.
As $0:\Mc_K\rightarrow\Tc$ meets neither $D_W$ nor $D_{W_1}^\prime$, we have 
\begin{equation*}
\Eis_W(f)= 0^*\pol_W(f)=0^* \res(\pol_W)=
\end{equation*}
\begin{equation*}
0^*p^{W * }_{W_1}\pol_{W_1}(f)= 
q^{W *} _{W_1}0^*\pol_{W_1}(f)=q^{W *} _{W_1}\Eis_{W_1}(f)
\end{equation*}
It follows that our operator is well-defined.\\
$G(\A_f)\times \pi_0(G(\R))$-equivariance follows easily from our construction. Given a section $f$ with corresponding groups $K$, $V_K$, $V_K ^\prime$ and $W$ as above and given $g\in G(\A_f)\times \pi_0(G(\R))$ the corresponding groups for $f\cdot g$ can be chosen $g^{-1}Kg$, $g^{-1}V_K$, $g^{-1}V_K ^\prime$ and $g^{-1}Wg$. We have the morphism $\Tc_{g^{-1}Wg}\to\Tc_W$,\ $(v,h)\mapsto (gv,gh) =g(v,h)$
and use \Cref{natural_localization} to get the commutative diagram
\begin{equation*} 
\begin{xy}
\xymatrix{
H^{2\xi-1}(\Tc_{W}\setminus D_{W},Log\otimes\mu^{\otimes n+1})\ar[d]\ar[r] & H^0(D_{W},Log_{|D_{W_1}}\otimes \mu^{\otimes n})\ar[d] \\
H^{2\xi-1}(\Tc_{g^{-1}Wg}\setminus D_{g^{-1}Wg},Log\otimes\mu^{\otimes n+1})\ar[r] & H^0(D_{g^{-1}Wg},Log_{|D_{g^{-1}Wg}}\otimes\mu^{\otimes n})
}
\end{xy}
\end{equation*}
As the pullback $g^*f$ is exactly the section corresponding to $f\cdot g$, we have $g^*\pol_W(f)=\pol_{g^{-1}Wg}(f\cdot g)$. It follows
\begin{equation*}
\Eis(f\cdot g)=0^*\pol_{g^{-1}Wg}(f\cdot g)=0^*g^*\pol_W(f)=g^*0^*\pol_W(f)=g^*\Eis(f) 
\end{equation*}
and we are done.
\end{proof}
\end{proposition}

\section{Pushforward of polylogarithmic Eisenstein classes}

We are going to decompose the cohomology of $\Mc_K$ by a Leray-Hirsch theorem into a product of the cohomology of $\Sc_K$ and the cohomology of the fiber of $\phi_K$. This will yield a decomposition of the polylogarithmic Eisenstein classes on $\Mc_K$ and a way to integrate these classes to gain classes on $\Sc_K$.
\begin{lemma}
Let $G$ and $H$ be groups and $R$ a Dedekind domain such that $R$ has a resolution $P\rightarrow R$ by finitely generated free $R[G]$-modules (for example $G\cong \Z^r$). Let $M$ be an $R[G]$-module finitely generated free as $R$-module and let $N$ be any $R[H]$-module. Then we have the Künneth sequence
\begin{equation*}
0\rightarrow \bigoplus_{p+q=n} H^p(G,M)\otimes_R H^q(H,N)\rightarrow H^n(G\times H,M\otimes_R N)\rightarrow  
\end{equation*} 
\begin{equation*}
\rightarrow\bigoplus_{p+q=n+1} Tor_1 ^{R}(H^p(G,M),H^q(H,N))\rightarrow 0
\end{equation*} 
\begin{proof}
We choose a resolution $P^\bullet\to R$ by finitely generated free $R[G]$-modules and a resolution $Q^\bullet\to R$ by free $R[H]$-modules. Since $R$ is a Dedekind domain the torsion-free $R$-modules are exactly the flat $R$-modules and therefore any submodule of a flat $R$-module is flat. With this at hand we may use \cite{We} Theorem 3.6.3 to conclude that  $P^\bullet\otimes_RQ^\bullet\to R$ is a resolution by free $R[G\times H]$-modules. So $H^\bullet(G\times H,M\otimes_R N)=H^\bullet(Hom_{R[G\times H]}(P^\bullet\otimes_RQ^\bullet,M\otimes_RN))$. Let us consider the canonical map
\begin{equation*}
\mu:Hom_{R[G]}(P,M)\otimes_R Hom_{R[H]} (Q,N)\rightarrow Hom_{R[G\times H]}(P\otimes_R Q,M\otimes_R N) 
\end{equation*}
$\mu(f\otimes g)(p\otimes q)=f(p)\otimes g(q)$ for $p\in P$ and $q\in Q$ for a finitely generated free $R[G]$-module $P$ and a free $R[H]$-module $Q$. We have $Hom_{R[G\times H]}(P\otimes_R Q,M\otimes_R N)$ equals
\begin{equation*}
Hom_{R[G\times H]}(R[G]^n\otimes_R Q,M\otimes_R N)=
Hom_{R[H]}(Q,M\otimes_R N)^n
\end{equation*}
Now $M=R^m$ is free, so we continue by
\begin{equation*}
Hom_{R[H]}(Q,N^m)^n=Hom_{R[H]}(Q,N)^{mn}=Hom_{R[H]}(Q,N)\otimes_R M^n=
\end{equation*}
\begin{equation*}
=Hom_{R[H]}(Q,N)\otimes_R Hom_{R[G]}(R[G]^n,M)=Hom_{R[H]}(Q,N)\otimes_R Hom_{R[G]}(P,M)
\end{equation*}
So $\mu$ is an isomorphism. Moreover, $M^n=Hom_{R[G]}(P,M)$ is a flat $R$-module and we can apply \cite{We} Theorem 3.6.3 to the complex $Hom_{R[G]}(P,M)\otimes_R Hom_{R[H]} (Q,N)$ to get the exact sequence
\begin{equation*}
0\to \bigoplus_{p+q=n}H^p(Hom_{R[G]}(P,M))\otimes_RH^q(Hom_{R[H]} (Q,N))\to 
\end{equation*}
\begin{equation*}
\to H^n(Hom_{R[G]}(P,M)\otimes_A Hom_{R[H]} (Q,N))\to 
\end{equation*}
\begin{equation*}
\to\bigoplus_{p+q=n+1}Tor^R_1(H^p(Hom_{R[G]}(P,M)),H^q(Hom_{R[H]} (Q,N))\to 0
\end{equation*}
Using $\mu$ we compute this sequence as
\begin{equation*}
0\rightarrow \bigoplus_{p+q=n} H^p(G,M)\otimes_R H^q(H,N)\rightarrow H^n(G\times H,M\otimes_R N)\rightarrow  
\end{equation*} 
\begin{equation*}
\rightarrow\bigoplus_{p+q=n+1} Tor_1 ^{R}(H^p(G,M),H^q(H,N))\rightarrow 0
\end{equation*} 
\end{proof}
\end{lemma}
\begin{proposition}\label{cohomology_abelian_group}
Let $G$ be a finitely generated free abelian group, $k$ a field and $M$ an $k[G]$-module finitely generated as $k$-module. We suppose that we have a field extension $k\subset K$ such that there is a decomposition $M_K=\bigoplus_{\chi}K(\chi)^{m(\chi)}$ with $\chi:G\rightarrow K^{\times}$ running through all characters and $K(\chi)$ the $K[G]$-module $K$ with action by $\chi$. Then cup-product induces an isomorphism $H^0(G,M)\otimes_A H^\bullet(G,k)=H^\bullet(G,M)$.
\begin{proof}
We begin with $H^\bullet(G,K(\chi))=0$, if $\chi\neq 1$. Write $G\cong \Z^r$ and for $j=1,...,r$ define $\iota_j:\Z\to G$ to be the inclusion of the $j$-component. As $\chi\neq 1$, there is a $j$ such that $\chi\circ \iota_j\neq 1$. So we may assume $G=\Z\times G^\prime$ with $\chi_{|\Z}\neq 1$ and $G^\prime$ a subgroup. One has $K(\chi)=K(\chi_{|\Z})\otimes_K K(\chi_{|G^\prime})$, as well as $H^\bullet (\Z,K(\chi_{|\Z}))=0$ by \cite{We} Example 6.1.4. and we use the Künneth sequence
\begin{equation*}
0\rightarrow \bigoplus_{p+q=n} H^p(\Z,K(\chi_{|\Z}))\otimes_A H^q(G^\prime,K(\chi_{|G^\prime}))\rightarrow H^n(G,K(\chi))\rightarrow  
\end{equation*} 
\begin{equation*}
\rightarrow\bigoplus_{p+q=n+1} Tor_1 ^{K}(H^p(\Z,K(\chi_{|\Z})),H^q(G^\prime,K(\chi_{|G^\prime})))\rightarrow 0
\end{equation*} 
to conclude. Now we may calculate
\begin{equation*}
H^\bullet(G,M_K)=H^\bullet(G,\bigoplus_{\chi}K(\chi)^{m(\chi)})=\bigoplus_{\chi}H^\bullet(G,K(\chi))^{m(\chi)}.
\end{equation*}
But all summands with $\chi\neq 1$ are zero, so this equals
\begin{equation*}
H^\bullet(G,K(1))^{m(1)}=H^0(G,M_K)\otimes_K H^\bullet(G,K) 
\end{equation*}
and this isomorphism is induced by cup-product. We want to see that this map comes from extension of scalars from the cup-product map $\cup:H^0(G,M)\otimes_k H^\bullet(G,k)\rightarrow H^\bullet (G,M)$. This is the case, because $H^\bullet(G, N_K)=H^\bullet(G, N)\otimes_k K$ for any $G$-module $N$ compatible with cup-product. To see this take a resolution $P\rightarrow k$ by finitely generated free $k[G]$-modules. Then $H^\bullet(G, N_K)=H^\bullet(Hom_{K[G]}(P_K,N_K))$. As $k\rightarrow K$ is flat and $P$ consists of finitely generated free $k[G]$-modules one identifies
\begin{equation*}
H^\bullet(Hom_{K[G]}(P_K,N_K))=H^\bullet(Hom_{k[G]}(P,N))\otimes_k K=H^\bullet(G, N)\otimes_k K
\end{equation*}
proving the claim.
Finally, the extension of scalars from $k$ to $K$ of the morphism $\cup$ is an isomorphism, therefore $\cup$ is already an isomorphism, since $k\rightarrow K$ is faithfully flat.
\end{proof}
\end{proposition}
\begin{corollary}\label{fiber_cohom}
We consider $\phi_K:\Mc_K\rightarrow \Sc_K$ and $\Sym^n\Hc$. Then cup-product induces an isomorphism $\cup:\phi_{K*}\Sym^n\Hc\otimes_\Q R^p\phi_{K*}(\Q)\rightarrow R^p\phi_{K*}(\Sym^n\Hc)$.
\begin{proof}
It is clear that we have a global cup-product map and being an isomorphism is a local problem on $\Sc_K$. As $\phi_K$ is a fiber bundle, it suffices to show that for all $s\in \Sc_K$ the induced maps on the fibers $\cup:H^0(\phi_{K}^{-1}(s),\Sym^n\Hc)\otimes_\Q H^p(\phi_{K}^{-1}(s),\Q)\rightarrow H^p(\phi_{K}^{-1}(s),\Sym^n\Hc)$ are isomorphisms. Following \cite{GrT} Chapitre V we may prove this by showing that the corresponding map on the level of group cohomology $\cup:H^0(Z_K,\Sym^nV(\Q))\otimes_\Q H^p(Z_K,\Q)\rightarrow H^p(Z_K,\Sym^nV(\Q))$ is an isomorphism. This is the case, as our situation fulfills the requirements of \Cref{cohomology_abelian_group}. $Z_K$ is a finitely generated free abelian group by Dirichlet's unit theorem and we have the field extension
$\Q\rightarrow \overline\Q$ such that $\Sym^nV(\Q)$ decomposes: Take the isomorphism $\Oc\otimes_\Z \overline\Q\rightarrow \prod_{\sigma:\Oc\rightarrow \overline \Q} \overline\Q,$ $x\otimes r\mapsto (\sigma(x)r)_\sigma$	and obtain 
\begin{equation*}
\Sym^n _\Q V (\Q)\otimes_\Q \overline\Q=\bigoplus_{\sum_\sigma n_\sigma =n}\bigotimes_\sigma \Sym^{n_\sigma} _{\overline\Q}\overline\Q^2.
\end{equation*}
If we take the standard basis $e^1=:X$, $e^2=:Y$ of $\overline\Q^2$, we get the $\overline\Q$-basis $\prod_\sigma X_\sigma ^{k_\sigma}Y_\sigma ^{n_\sigma -k_\sigma }$, $k_\sigma=0,...,n_\sigma$, of $\Sym^n _\Q V (\Q)\otimes_\Q \overline\Q$, which also yields the decomposition into one-dimensional $Z_K$-representations, as $\epsilon\cdot\prod_\sigma X_\sigma ^{k_\sigma} Y_\sigma ^{n_\sigma -k_\sigma }=\prod_\sigma \sigma(\epsilon)^{n_\sigma}\prod_\sigma X_\sigma ^{k_\sigma}Y_\sigma ^{n_\sigma -k_\sigma }$
for $\epsilon \in Z_K\subset \Oc^{\times}$.
\end{proof}
\end{corollary}
\begin{remark}\label{invariants}
From now on we set $A:=\Z[\frac{1}{d_FN}]$ where $d_F>0$ is the discriminant of $F$. We have the Hilbert class field $F^{\text{Hil}}$ of $F$ with ring of integers $\Oc^{\text{Hil}}$. Let us set $R:=\Oc^{\text{Hil}}[\frac{1}{Nd_F}]$. Then $A\to R$ is a faithfully flat ring extension. From now on we will always consider the sheaf $\Sym^k \Hc_{PD}$ over the ring $A$. 
The natural adjunction morphism 
\begin{equation*}
\phi_K ^{-1}\phi_{K*}\Sym^k \Hc_{PD}\rightarrow \Sym^k \Hc_{PD}
\end{equation*}
is a split monomorphism. Indeed, the corresponding map of $G(g_f)$-modules is a split mono. To see this note first that 
\begin{equation*}
\Sym^k V(g_f)_{PD} ^{Z_K}\otimes_A R=\left(\Sym^k V(g_f)_{PD}\otimes_A R  \right)^{Z_K},
\end{equation*}
 as $R$ is a free $A$-module. By faithful flatness it suffices now to find a retraction for the map
\begin{equation*}
\left(\Sym^k V(g_f)_{PD}\otimes_A R  \right)^{Z_K} \rightarrow \Sym^k V(g_f)_{PD}\otimes_A R. 
\end{equation*}
There are two fractional ideals $\mathfrak a$ and $\mathfrak b$ of $\Oc$ such that $V (g_f)=\mathfrak a\oplus \mathfrak b$. As we have $d_F\in R^{\times}$, there are again isomorphisms
\begin{equation*}
\mathfrak a\otimes_A R \rightarrow \prod_{\sigma:\Oc\rightarrow \overline \Q} \sigma(\mathfrak a)R,\ x\otimes r\mapsto (\sigma(x)r)_\sigma,\	
\mathfrak b\otimes_A R \rightarrow \prod_{\sigma:\Oc\rightarrow \overline \Q} \sigma(\mathfrak b)R,\ x\otimes r\mapsto (\sigma(x)r)_\sigma	
\end{equation*}
and we get
\begin{equation*}
\Sym^k V(g_f)_{PD}\otimes_A R=\bigoplus_{\sum_\sigma k_\sigma =k}\bigotimes_\sigma \Sym^{k_\sigma} _{R}\left(\sigma(\mathfrak a)R\oplus \sigma(\mathfrak b)R\right)_{PD}.
\end{equation*}
In the Hilbert class field any Ideal of $\Oc$ becomes principal. Therefore we get $\alpha_\sigma,\beta_\sigma\in F^{\text{Hil}}$ with $\sigma(\mathfrak a)R=\alpha_\sigma R$ and $\sigma(\mathfrak b)R=\beta_\sigma R$. If again $e^1=:X$, $e^2=:Y$ is the standard basis of $R^2$, we may split $\Sym^k V(g_f)_{PD}\otimes_A R$ as in \Cref{fiber_cohom} into one dimensional $Z_K$-representations 
$R\prod_\sigma (\alpha_\sigma X_\sigma) ^{[n_\sigma]}(\beta_\sigma Y_\sigma) ^{[k_\sigma -n_\sigma ]}$. 
In $\left(\Sym^k V(g_f)_{PD}\otimes_A R  \right)^{Z_K}$ those one dimensional representations occur, on which $Z_K$ acts trivially. From this description it is clear that the inclusion of the invariants is a split mono. But we can say more. The one dimensional subrepresentations occurring in the invariants induce characters $Res_{F/\Q}\G_m(\overline \Q)\rightarrow \G_m (\overline \Q)$ which have to be trivial on $Z_K$. By \cite{Se} II.3.3 we know that these have to factor through the norm character, as our field $F$ is totally real. In other words, $R\prod_\sigma (\alpha_\sigma X_\sigma) ^{[n_\sigma]}(\beta_\sigma Y_\sigma) ^{[k_\sigma -n_\sigma ]}$ occurs in the invariants, if and only if there is an $m\in \N_0$ such that $m=k_\sigma$ for all $\sigma :F\rightarrow\overline \Q$. It follows that the module $\left(\Sym^k V(g_f)_{PD}\otimes_A R\right) ^{Z_K}$ does not depend on the level $K$, as soon as $Z_K\subset \Oc^{\times}$ does not contain an element of negative norm. By faithful flatness $\Sym^k V(g_f)_{PD} ^{Z_K} $ does not depend on the level $K$ in this case. As we have usually $K_f=K_N$, $N\geq 3$, we do not have elements of negative norm.
\end{remark}
\begin{corollary}
We have $\Sym^k \Hc_{PD} ^\prime=0$, if $\xi\nmid k$. If $k=\xi m$ and $Z_K$ contains no element of negative norm then $\Sym^k \Hc_{PD} ^\prime\neq0$.
\end{corollary}
\begin{remark}\label{invariants_integral}
We want to describe the image
\begin{equation*}
\Sym^k_ A V(g_f)_{PD}^{Z_K} \rightarrow \left(\Sym^k_A V(g_f)_{PD}\otimes_A R \right)^{Z_K}=\Sym^k_R (V(g_f)\otimes_A R)_{PD}^{Z_K}
\end{equation*}
explicitly. $\text{Gal}(F^{\text{Hil}}/\Q)$ acts on $R$ and we get $R^{\text{Gal}(F^{\text{Hil}}/\Q)}=A$ by Galois-theory. This implies
\begin{equation*}
\Sym^k _AV(g_f)_{PD} = \left(\Sym^k_A V(g_f)_{PD}\otimes_A R\right)^{\text{Gal}(F^{\text{Hil}}/\Q)}=\Sym^k_R (V(g_f)\otimes_A R)_{PD}^{Gal(F^{\text{Hil}}/\Q)},
\end{equation*}
as $\Sym^k _AV(g_f)_{PD}$ is a free $A$-module. We identify
\begin{equation*}
\Sym^k _AV(g_f)_{PD} ^{Z_K}=\left(\Sym^k_R (V(g_f)\otimes_A R)_{PD}^{Z_K}\right)^{\text{Gal}(F^{\text{Hil}}/\Q)},
\end{equation*}
since the Galois-action and the $Z_K$-action commute. Let us suppose $k=\xi m$, otherwise $\Sym^k _AV(g_f)^{Z_K}=0$ by the considerations above. Write $V(g_f)=\mathfrak a \oplus \mathfrak b$ as in \Cref{invariants}. In $\Sym^k_R (V(g_f)\otimes_A R)_{PD}^{Z_K}$ we have $R\prod_\sigma (\beta_\sigma Y_\sigma) ^{[m]}$ as direct summand. $\text{Gal}(F^{\text{Hil}}/\Q)$ acts trivially on $\prod_\sigma  Y_\sigma ^{[m]}$, because $\tau\cdot Y_{\sigma}=Y_{\tau^{-1}\circ\sigma}$ for $\tau\in \text{Gal}(F^{\text{Hil}}/\Q)$.
Therefore
\begin{equation*}
\prod_\sigma \sigma(\mathfrak b)^mR\stackrel{\cong}{\rightarrow}R\prod_\sigma  (\beta_\sigma Y_\sigma) ^{[m]}=R\prod_\sigma  \beta_\sigma^m Y_\sigma ^{[m]},\ x\mapsto x\cdot \prod_\sigma Y_\sigma ^{[m]}   
\end{equation*} 
as Galois-modules, where on the left-hand side we have the product of ideals.
\begin{lemma}
\begin{equation*}
\left(\prod_\sigma \sigma(\mathfrak b)R\right)^{\text{Gal}(F^{\text{Hil}}/\Q)}=\Norm(\mathfrak b)A
\end{equation*} 
where $\Norm(\mathfrak b):=\prod_{\text{$\nu$ finite place of $F$}}|\kappa(\mathfrak p_\nu)|^{\nu(\mathfrak b)}$ and $\nu(\mathfrak b)$ is defined by $\mathfrak b\Oc_\nu=(\pi_\nu ^{\nu(\mathfrak b)})$.
\begin{proof}
We may assume that $\mathfrak b$ is integral, otherwise we multiply the ideal with a number $z\in \Z\setminus \left\{0\right\}$. As a localization of a Dedekind ring $A$ is itself a Dedekind ring. So we have unique prime factorization of ideals in $A$ and therefore our problem is local at each maximal ideal $(p)\subset A$. First we have
\begin{equation*}
\left(\left(\prod_\sigma \sigma(\mathfrak b)R\right)^{\text{Gal}(F^{\text{Hil}}/\Q)}\right)_{(p)}=\left(\left(\prod_\sigma \sigma(\mathfrak b)R\right)\cap \Q\right)_{(p)}=\left(\prod_\sigma \sigma(\mathfrak b)R_{(p)}\right)\cap \Q
\end{equation*} 
by Galois-theory. As $(p)\subset \Q$ and $\sigma$ is $\Q$-linear, we identify the last expression with $\left(\prod_\sigma \sigma(\mathfrak b_{(p)})R_{(p)}\right)\cap \Q$. Now $\Oc_{(p)}$ is a Dedekind ring with only finitely many primes and hence a principal ideal domain. So $\mathfrak b_{(p)}=b\Oc_{(p)}$ for $b\in \Oc$ and we identify our ideal with 
\begin{equation*}
\left(\prod_\sigma \sigma(b)R_{(p)}\right)\cap \Q=\Norm(b)R_{(p)}\cap \Q.
\end{equation*} 
$R_{(p)}$ is the integral closure of $A_{(p)}$ inside $F^{\text{Hil}}$ and as $A_{(p)}$ is integrally closed in $\Q$ we get $\Norm(b)R_{(p)}\cap \Q=\Norm(b)A_{(p)}$.
We calculate
\begin{equation*}
\Norm(\mathfrak b)A_{(p)}=\prod_{\text{$\nu$ finite place of $F$}}|\kappa(\mathfrak p_\nu)|^{\nu(\mathfrak b)}A_{(p)}=\prod_{\nu|p}|\kappa(\mathfrak p_\nu)|^{\nu(\mathfrak b)}A_{(p)}.
\end{equation*}
In $\Oc_{(p)}$ we have the unique prime factorization
\begin{equation*}
 b \Oc_{(p)}=\mathfrak b_{(p)}=\prod_{\text{$\mathfrak p\subset \Oc$ max., $\mathfrak p|p$}}\mathfrak p^{\nu_{\mathfrak p}(\mathfrak b)}\Oc_{(p)}
\end{equation*}
and we see $\nu_{\mathfrak p}(b)=\nu_{\mathfrak p}(\mathfrak b)$ for $\mathfrak p| p$. We get
\begin{equation*}
\left(\left(\prod_\sigma \sigma(\mathfrak b)R\right)^{\text{Gal}(F^{\text{Hil}}/F)}\right)_{(p)}=|\Norm( b)|A_{(p)}=\prod_{\text{$\mathfrak p\subset \Oc$ max.}}|\kappa(\mathfrak p)|^{\nu_{\mathfrak p}(b)}A_{(p)}=
\end{equation*}
\begin{equation*}
\prod_{\text{$\mathfrak p\subset \Oc$ max., $\mathfrak p|p$}}|\kappa(\mathfrak p)|^{\nu_{\mathfrak p}(b)}A_{(p)}=\prod_{\text{$\mathfrak p\subset \Oc$ max., $\mathfrak p|p$}}|\kappa(\mathfrak p)|^{\nu_{\mathfrak p}(\mathfrak b)}A_{(p)}=\Norm(\mathfrak b)A_{(p)}
\end{equation*}
and we are done.
\end{proof}
\end{lemma}
\begin{proposition}\label{Y_integral}
We consider $A:=\Z[\frac{1}{Nd_F}]$. $V(g_f)=g_f^{-1}V(\hat\Z)\cap V(\Q)=\mathfrak a\oplus \mathfrak b$. Then $\Norm(\mathfrak b)^m\prod_\sigma Y_\sigma ^{[m]}$ is part of an $A$-basis of $\Sym^k _AV(g_f)_{PD} ^{Z_K}$ and $\Norm(\mathfrak b)=\left\|t_2\right\|_f $, if we write
\begin{equation*}
g_f=xb=x\begin{pmatrix}t_1 & t_1 u\\ 0 & t_2 \end{pmatrix}, x\in G(\hat\Z).
\end{equation*}
\begin{proof}
The first part of the proposition is clear by the preceding lemma. We have $V(g_f)=V(b)$, since $x^{-1}V(\hat\Z)=V(\hat\Z)$. From this description it is clear that $\mathfrak b = t_2 ^{-1} \hat\Oc\cap F$ and therefore $\Norm(\mathfrak b)=\left\|t_2\right\|_f $.
\end{proof}
\end{proposition}
\end{remark}
\begin{proposition}\label{global_classes}
Suppose $\det(K_f)=K_f\cap Z(\A_f)$, what we always will assume from now on. 
We have cohomology classes in $H^\bullet (\Mc_K,\Z[\frac{1}{2}])$ restricting for all $s\in \Sc_K$ to a basis of $H^\bullet(\phi_K ^{-1}(s),\Z[\frac{1}{2}])$. In particular, the local system $R^\bullet\phi_{K*}(\Z[\frac{1}{2}])$ is trivial.
\begin{proof}
Suppose that we had the cohomology classes with the desired properties.
Let $H^\bullet\subset H^\bullet (\Mc_K,\Z[\frac{1}{2}])$ be the $\Z[\frac{1}{2}]$-submodule generated by the elements restricting to bases of the cohomology of the fibers. The natural morphism $H^\bullet (\Mc_K,\Z[\frac{1}{2}])\rightarrow H^0(\Sc_K,R^\bullet\phi_{K*}(\Z[\frac{1}{2}]))$ induces $H^\bullet \rightarrow H^0(\Sc_K,R^\bullet\phi_{K*}(\Z[\frac{1}{2}]))$ and therefore an isomorphism of sheaves $\underline H^\bullet \rightarrow R^\bullet\phi_{K*}(\Z[\frac{1}{2}])$, which would be our trivialization. So let us find these cohomology classes. We have determinant maps 
\begin{equation*}
\Sc_K\rightarrow \det(K_f)\prod_{\nu|\infty}\R^\times\backslash \I_F/F^{\times}=:Cl_F ^{\overline{K}},\ g\mapsto \det(g),
\end{equation*}
\begin{equation*}
\Mc_K\rightarrow \det(K_f)\R_{>0}\left\{\pm 1\right\}^\xi\backslash \I_F/F^{\times}=:T_K,\ g\mapsto \det(g)
\end{equation*}
We consider the fiber product
$\Sc_K\times_{Cl_F ^{\overline{K}}}T_K$, where on the right hand side the structure map is the natural projection. We have the natural maps
\begin{equation*}
p:\Mc_K\rightarrow \Sc_K\times_{Cl_F ^{\overline{K}}}T_K,\ K^1gG(\Q)\mapsto (KgG(\Q), \det(g)),
\end{equation*}
\begin{equation*}
pr_2:\Sc_K\times_{Cl_F ^{\overline{K}}}T_K\rightarrow T_K.
\end{equation*}
Everything is fibered over $Cl_F ^{\overline{K}}$ so the cohomology groups decompose into direct sums corresponding to components which are cut out by elements of $Cl_F ^{\overline{K}}$, so we can treat these components separately. We consider $G(g_f)$ for some $g_f\in G(\A_f)$. We get
\begin{equation*}
\phi_K:\left(\H^{\xi}_{\pm}\times (\R_{>0}\backslash Z(\R) ^0)\right)/G(g_f)\rightarrow \H^{\xi}_{\pm}/G(g_f),\ (\tau,r)\mapsto \tau
\end{equation*} 
and
\begin{equation*}
p:\left(\H^{\xi}_{\pm}\times (\R_{>0}\backslash Z(\R)^0)\right)/G(g_f)\rightarrow(\H^{\xi}_{\pm}/G(g_f))\times (\R_{>0}\left\{\pm 1\right\}^\xi\backslash Z(\R)/\det(G(g_f)))
\end{equation*} 
is a map of fiber bundles over $\H^{\xi}_{\pm}/G(g_f)$. We have the projection
\begin{equation*}
pr_2:(\H^{\xi}_{\pm}/G(g_f))\times (\R_{>0}\left\{\pm 1\right\}^\xi\backslash Z(\R)/\det(G(g_f)))\rightarrow \R_{>0}\left\{\pm 1\right\}^\xi\backslash Z(\R)/\det(G(g_f))
\end{equation*}
and given a basis 
$(\eta_i)\subset H^\bullet(\R_{>0}\left\{\pm 1\right\}^\xi\backslash Z(\R)/\det(G(g_f)),\Z[\frac{1}{2}])$
the classes 
\begin{equation*}
(pr_2 ^*\eta_i)\subset H^\bullet((\H^\xi/G(g_f))\times (\R_{>0}\left\{\pm 1\right\}^\xi\backslash Z(\R)/\det(G(g_f))),\Z[\frac{1}{2}])
\end{equation*} 
restrict back to a basis of the fibers. $p$ respects the fibers and therefore the same is true for $(p^*pr_2^*\eta_i)=:(\overline\eta_i)$. To see this restrict $p$ to the fiber
\begin{equation*}
 p:\R_{>0}\backslash Z(\R) ^0/Z_K\rightarrow \R_{>0}\left\{\pm 1\right\}^\xi\backslash Z(\R)/\det(G(g_f)), x\cdot Z_K\mapsto x\cdot \det(G(g_f))
 \end{equation*}
We may suppose $Z_K=\det(G(g_f))$. On the left hand side $\epsilon \in Z_K$ acts by multiplication with $\epsilon ^2$, whereas $\epsilon \in \det(G(g_f))$ acts on the right hand side by multiplication with $\epsilon$. Therefore the map above is an isogeny of real tori of degree $(Z_K:Z_K^2)$, which is a power of $2$. As $2\in \Z[\frac{1}{2}] ^{\times}$ the map induces an isomorphism on cohomology.

Let us now denote the set of connected components of $\Mc_K$ by $J$. We have constructed cohomology classes $(\overline\eta_i ^j)$ for each $j\in J$ restricting to basis of the cohomology of the fibers. The family of cohomology classes $(\sum_{j\in J} \overline \eta_i ^j)$ on $\Mc_K$ has the desired property: If $s\in \Sc_K$ is given, then $\phi_K ^{-1}(s)$ will be contained in a connected component corresponding to some $j_0\in J$. We get $(\sum_{j\in J} \overline \eta_{i\ |\phi^{-1}(s)} ^j)=(\overline \eta_{i\ |\phi^{-1}(s)} ^{j_0})$ and the right-hand side is a basis of the cohomology of the fiber $\phi_K ^{-1}(s)$ by the construction above.
\end{proof}
\end{proposition}
\begin{remark}\label{invariant_Z_cohom}
The cohomology classes described above are pullback by $\det:\Mc_K\rightarrow \det(K_f)\R_{>0}\left\{\pm1\right\}^\xi\backslash \I_F/F^{\times}=:T_K$. We have a short exact sequence of groups
\begin{equation*}
0\to \R_{>0}\left\{\pm 1\right\}^\xi\backslash Z(\R)/\det(G(g_f))\to T_K \to Cl_F ^{\overline {K}}\to 0
\end{equation*}

The group $T_K$ acts on itself by translation and therefore also on $H^\bullet(T_K,\Z)$ its cohomology. This action factors through the quotient $T_K\to Cl_F ^{\overline{K}}$ and we get 
\begin{equation*}
H^\bullet (T_K,\Z)^{T_K}=H^\bullet(T_K,\Z)^{Cl_F ^{\overline{K}}}\cong H^\bullet(\R_{>0}\backslash Z(\R) ^0/\det(G(g_f)),\Z).
\end{equation*}
The pullback of these forms by $\det$ have the desired property of the proposition above and we denote the corresponding subspace by $\mathfrak H_K ^\bullet\subset H^\bullet(\Mc_K,\Z)$. We have the invariant $1$-forms $\frac{dr_\nu}{r_\nu}$ on $F_\R^{\times}$, which we may also interpret as invariant $1$-forms on $T_K$. The invariant forms $\mathfrak H _K ^\bullet$
are $\R$-linear combinations of wedge-products these forms. Moreover, we set $\mathfrak H ^\bullet:=\mathfrak H_K ^\bullet\otimes\Q$, as this space does not depend on the level $K$.
\end{remark}

\begin{theorem}\label{decomp}
We have an isomorphism
\begin{equation*}
H^\bullet(\Sc_K,\Sym^k\Hc ^\prime)\otimes_\Q \mathfrak H ^\bullet\rightarrow H^\bullet(\Mc_K,\Sym^k\Hc),\ \omega\otimes\eta\mapsto \omega\cup\eta 
\end{equation*}
\begin{proof}
Let us choose a resolution $I^\bullet$ of $\Q$ on $\Mc_K$ by injective sheaves of $\Q$-modules. We take a basis $\eta_1,...,\eta_l$ of $\mathfrak{ H} ^\bullet$ in other words cohomology classes restricting to a basis of the cohomology of the fibers of $\phi_K$. Let us denote the cohomological degree of the class $\eta_i$ by $d_i\in \N_0$ ($i\geq j$, then $d_i\geq d_j$).
We may interpret $\eta_i$ as a morphism $\eta_i: \Q[-d_i]\rightarrow I^\bullet$. These morphisms induce
\begin{equation*}
\id\otimes\eta_i=\eta_i: \Sym^k\Hc[-d_i]\rightarrow \Sym^k\Hc\otimes_\Q I^\bullet
\end{equation*}
By \Cref{injectives_local_systems} we may use $\Sym^k\Hc\otimes_\Q I^\bullet$ as an injective resolution of $\Sym^k\Hc$. Applying $\phi_{K *}$ and summing over all $i$ yields
\begin{equation*}
\sum_i \eta_i:\bigoplus_i \phi_{K *}\Sym^k\Hc [-d_i]\rightarrow R\phi_{K *}\Sym^k\Hc.
\end{equation*}
Now apply $H^p$ and get a map
\begin{equation*}
\sum_{\text{$i$ with $d_i=p$}} \eta_i:\bigoplus_{\text{$i$ with $d_i=p$}} \phi_{K *}\Sym^k\Hc[-d_i]\rightarrow R^p\phi_{K *}(\Sym^k\Hc).
\end{equation*} 
\begin{equation*}
R^p\phi_{K *}(\Sym^k\Hc)=\phi_{K *}\Sym^k\Hc\otimes_\Q R^p\phi_{K *}(\Q)
\end{equation*}
and the properties of $\eta_i$ show that the last map is an isomorphism. In particular,
\begin{equation*}
\sum_i \eta_i:\bigoplus_i \phi_{K *}\Sym^k\Hc[-d_i]\rightarrow R\phi_{K *}\Sym^k\Hc
\end{equation*}
is a quasi-isomorphism. We conclude
\begin{equation*}
R\Gamma(\Sc_K,\phi_{K *}\Sym^k\Hc)\otimes_\Q \mathfrak H ^\bullet:=\bigoplus_i R\Gamma(\Sc_K,\phi_{K *}\Sym^k\Hc)[-d_i]=
\end{equation*}
\begin{equation*}
 =R\Gamma(\Sc_K,R\phi_{K *}\Sym^k\Hc)=R\Gamma(\Mc_K,\Sym^k\Hc).
\end{equation*}
Moreover, we have $\Sym^k\Hc^\prime=\phi_{K*}\Sym^k\Hc $ completing the proof.
\end{proof}
\end{theorem}
Let us consider $\mathfrak H ^\bullet$ as trivial $G(\A_f)\times \pi_0(G(\R))$-module. Set 
\begin{equation*}
\varinjlim H^\bullet(\Sc_K,\Sym^k\Hc ^\prime):=H^\bullet(\Sc,\Sym^k\Hc ^\prime).
\end{equation*}
 \begin{corollary}\label{equiv_decomp}
We have an isomorphism of $G(\A_f)\times \pi_0(G(\R))$-modules
\begin{equation*}
H^\bullet(\Sc,\Sym^k\Hc ^\prime)\otimes_\Q \mathfrak H ^\bullet\rightarrow H^\bullet(\Mc,\Sym^k\Hc)
\end{equation*}
\begin{proof}
The isomorphism is \Cref{decomp}. $G(\A_f)\times \pi_0(G(\R))$-equivariance follows from the fact that $\mathfrak H ^\bullet$ consists of invariant classes.
\end{proof}   
\end{corollary}
\begin{corollary}\label{Eis^k _q}
We have an $G(\A_f)\times \pi_0(G(\R))$-equivariant operator
\begin{equation*}
\Eis^k _q:\Sc(V(\A_f),\mu^{\otimes n})^0 \otimes_\Q \mathfrak H ^{q*}\rightarrow H^{2\xi-1-q}(\Sc,\Sym^k\Hc^\prime\otimes\mu^{\otimes n+1})
\end{equation*}
$f\otimes \eta^*\mapsto \Eis^k(f)(\eta^*)$ for $q=0,...,\xi-1$.
\begin{proof}
This is just \Cref{Eis^k} and \Cref{equiv_decomp}
\end{proof}
\end{corollary}

\subsection{Explicit description of the decomposition of the cohomology}\label{trace_integration}

We want to calculate our polylogarithmic Eisenstein classes by differential forms in order to compare them with those of Harder. Therefore, we need to make the decomposition of the cohomology of $\Mc_K$ explicit. This can be done by the theory of fiber integration on de Rham cohomology.

Let us work over the coefficient ring $\Q$ first. We want to give the inverse map of the natural decomposition isomorphism provided by \Cref{decomp}. In other words, given a basis $\eta_1,...,\eta_l$ of $\mathfrak H ^\bullet$ and $\omega\in H^\bullet(\Mc_K,\Sym^k\Hc)$ we want to know explicitly $\omega_1,...,\omega_l\in H^\bullet(\Sc_K,\Sym^k\Hc ^\prime)$ such that $\omega=\sum_{i=1} ^l\omega_i\cup \eta_i$ holds.\\
We choose the retraction $\Sym^k\Hc\rightarrow \phi_K ^{-1}\Sym^k\Hc ^\prime$ of the canonical monomorphism $\phi_K ^{-1}\Sym^k\Hc ^\prime\rightarrow \Sym^k\Hc$ provided by \Cref{invariants}. 
\begin{lemma}\label{Z_proj1}
The retraction induces an isomorphism 
\begin{equation*} 
H^\bullet(\Mc_K,\Sym^k\Hc)\rightarrow H^\bullet(\Mc_K,\phi_{K} ^{-1}\Sym^k\Hc ^\prime)
\end{equation*}
on cohomology.
\begin{proof}
We know by \Cref{fiber_cohom} and the projection formula
\begin{equation*} 
R^p\phi_{K*}(\phi_K ^{-1}\Sym^k\Hc ^\prime)=\Sym^k\Hc ^\prime\otimes_\Q R^p \phi_{K*}(\Q)=R^p\phi_{K*}(\Sym^k\Hc)
\end{equation*}
so that our projection induces an isomorphism
\begin{equation*} 
R\phi_{K*}(\Sym^k\Hc ^\prime)=R\phi_{K*}(\Sym^k\Hc )\in D^+(\Sc_K,\Q)
\end{equation*}
in the derived category of bounded below complexes of sheaves of $\Q$-modules.
Applying $R^p\Gamma(\Sc,\ )$ yields the isomorphism in question.
\end{proof}
\end{lemma}
\begin{remark}
We may modify the proof of \Cref{decomp} to obtain an isomorphism
\begin{equation*}
H^\bullet(\Sc_K,\Sym^k\Hc^\prime)\otimes_\Q \mathfrak H ^\bullet\stackrel{\cup}{\rightarrow} H^\bullet(\Mc_K,\phi_K^{-1}\Sym^k\Hc ^\prime)
\end{equation*}
compatible with the isomorphism of \Cref{Z_proj1}.
\end{remark}
\begin{definition}\label{def_tr}
Let $B$ and $F$ be topological manifolds and $f:E\rightarrow B$ a fiber bundle with typical fiber $F$. We assume $F$ to be compact of dimension $r$. So $f$ is proper. Moreover, for a ring $A$ we consider a $A$-local system $\Vc$ on $B$. There is the natural \textit{edge morphism}
\begin{equation*}
 e:H^p(E,f^{-1}\Vc)\rightarrow H^{p-r}(B,\Vc\otimes_kR^rf_*(A)).
\end{equation*}	
Let us recall how it is constructed. $R^pf_*(f^{-1}\Vc)=0$, $p>r$, since
\begin{equation*}
R^pf_*(f^{-1}\Vc)_b=H^p(f^{-1}(b),f^{-1}\Vc)=0,\ p>r,
\end{equation*}
by \cite{Iv}VII 1.4, III 9.6 and 9.10.
Therefore the inclusion of the truncation induces a quasi-isomorphism $\tau_{\leq r}Rf_{*}(f^{-1}\Vc)\rightarrow Rf_{*}(f^{-1}\Vc)$. We have the canonical projection
\begin{equation*}
\tau_{\leq r}Rf_{*}(f ^{-1}\Vc)\rightarrow R^{r}f_{*}(f ^{-1}\Vc)[-r]=\Vc\otimes_k R^{r}f_{*}(A)[-r],
\end{equation*}
where we used the projection formula on the right hand side (\cite{Iv} VII 2.4). Both maps together give us the morphism
\begin{equation*}
e:Rf_{*}(f ^{-1}\Vc)\rightarrow \Vc\otimes_k R^{r}f_{*}(A)[-r]
\end{equation*}
in $D^+(B,k)$. Applying $H^p(B,\ )$ yields the edge morphism.
\end{definition}
\begin{remark}\label{tr_proj}
Since $\phi_K:\Mc_K\rightarrow \Sc_K$ is a fiber bundle with compact fibers, we get
\begin{equation*}
e:H^p(\Mc_K,\phi_K ^{-1}\Sym^k\Hc ^\prime)\rightarrow H^{p-(\xi-1)}(\Sc_K,\Sym^k\Hc ^\prime)\otimes_\Q \mathfrak H ^{\xi-1}.
\end{equation*}
Going through the proof of \Cref{decomp} we see that the composition of
\begin{equation*}
\bigoplus_{q} H^{p-q}(\Sc_K,\Sym^k\Hc ^\prime)\otimes_\Q \mathfrak H ^{q}\stackrel{\cup}{\rightarrow}H^p(\Mc_K,\phi_K ^{-1}\Sym^k\Hc ^\prime)
\end{equation*}
with the edge morphism is just the projection onto $H^{p-(\xi-1)}(\Sc_K,\Sym^k\Hc ^\prime)\otimes_\Q \mathfrak H ^{\xi-1}$.
\end{remark}
Poincaré duality $H^\bullet(\phi_K ^{-1}(s),\Z)\cong H^{\xi-1-\bullet}(\phi_K ^{-1}(s),\Z)^*$ induces a Poincaré duality $\mathfrak H ^\bullet\cong \mathfrak H ^{\xi-1-\bullet\ *}$. In particular, we have the isomorphism $\int:\mathfrak H ^{\xi-1}\rightarrow \Q=\mathfrak H ^{0}$, which is called the \textit{integration or trace morphism}, and the morphism
\begin{equation*}
(\id\otimes\int )\circ e:H^p(\Mc_K,\phi_K ^{-1}\Sym^k\Hc ^\prime)\rightarrow H^{p-(\xi-1)}(\Sc_K,\Sym^k\Hc ^\prime),
\end{equation*}
which we also denote by $\int$. 
We extend $\int$ to $\mathfrak H ^\bullet$ by setting it zero on $\mathfrak H ^p$, $p\neq \xi-1$. By Poincaré duality we get a perfect pairing $\mathfrak H ^\bullet\times \mathfrak H ^\bullet\to \Q$, $(\omega,\eta)\mapsto \int\omega\cup\eta$ and given our basis $\eta_1,...,\eta_l$, with cohomological degrees $deg(\eta_1)\leq...\leq deg(\eta_l)$, we choose a dual basis $\eta^* _1,...,\eta_l ^*$ with $\int\eta_i \cup\eta_j ^*=\delta_{ij}$.
\begin{lemma}\label{cup_decomp}
We consider $H^\bullet(\Mc_K,\phi_K ^{-1}\Sym^k\Hc ^\prime)$ as right $\mathfrak H ^{\bullet}$-module via cup-product.
If $\omega\in H^p(\Mc_K,\phi_K ^{-1}\Sym^k\Hc^\prime)$ is given, we have $\omega=\sum_{i=1}^l(\int\omega\cup\eta_i^*)\cup \eta_i$.
\begin{proof}
To proof this we may suppose by \Cref{decomp} $\omega=\sum_{i=1}^l\omega_i\cup \eta_i$, with $\omega_i\in H^p(\Sc_K,\Sym^k\Hc ^\prime)$. But we see with \Cref{tr_proj}
\begin{equation*}
\int \omega\cup\eta_j ^*=\sum_{i=1}^l\int(\omega_i\cup \eta_i)\cup\eta_j ^*=\sum_{i=1}^l\int\omega_i\cup (\eta_i\cup\eta_j ^*)=\sum_{i=1}^l\omega_i\cup \int(\eta_i\cup\eta_j ^*)=\omega_j
\end{equation*}
and the claim follows.
\end{proof}
\end{lemma}

As we really want to calculate cohomology classes, we finally have to make $e$ explicit in de Rham cohomology. So let us consider complex coefficients and the situation $f:E\rightarrow B$ from above, where we assume now $B,E,F$ to be $\Cc^\infty$-manifolds and $\Vc$ a $\C$-local system on $B$. We have $Rf_{*}(f^{-1}\Vc)=\Vc\otimes f_{*}(\Omega_{E} ^\bullet)$ and the quasi-isomorphism $\Vc\otimes R^{r}f_{*}(\C)\rightarrow\Vc\otimes R^{r}f_{*}(\C)\otimes\Omega_{B} ^\bullet$. We need a chain map $e^\bullet$ such that 
\begin{equation*}
\begin{xy}
\xymatrix{
\Vc\otimes f_{*}(\Omega_{E} ^\bullet)\ar[d]^{e} \ar[r]^ -{e^\bullet}& \Vc\otimes R^{r}f_{*}(\C)\otimes\Omega_{B} ^\bullet[-r]\\
\Vc\otimes R^{r}f_{*}(\C)[-d]\ar[r]^{=} & \Vc\otimes R^{r}f_{*}(\C)[-r]\ar[u]
}
\end{xy}
\end{equation*}
commutes. It suffices to find $e^\bullet$ such that
\begin{equation*}
\begin{xy}
\xymatrix{
f_{*}(\Omega_{E} ^\bullet)\ar[d]^{e} \ar[r]^ -{e^\bullet}& R^{r}f_{*}(\C)\otimes\Omega_{B} ^\bullet[-r]\\
R^{r}\phi_{*}(\C)[-r]\ar[r]^{=} & R^{r}\phi_{*}(\C)[-r]\ar[u]
}
\end{xy}
\end{equation*}
commutes. This is the theory of fiber integration. We assume the fiber bundle $f:E\rightarrow B$ to be orientated in the sense of \cite{HalGr}. This means that we have a $r$-form $\eta$ on $E$ restricting to a volume form on each fiber of $f$. We demand additionally $\eta$ to be closed. Therefore $\eta$ defines a trivialization of $R^{r}f_{*}(\C)$.
\begin{definition}\label{fiber_integration}
We define the \textit{fiber integration operator}
\begin{equation*}
e^\bullet:f_{*}\Omega_{E}^\bullet\rightarrow R^{r}f_{*}(\C)\otimes\Omega_{B}^\bullet[-r],\ \omega\mapsto (-1)^{r(\bullet-r)}\eta\otimes \vol(F)^{-1}\int_F \omega,
\end{equation*}
where $\int_F \omega$ is defined as in \cite{HalGr} with respect to $\eta$ and $\vol(F):=\int_F\eta$.
\end{definition}
\begin{remark}
\begin{enumerate}
\item $\vol(F)$ is a well-defined locally constant function, as $\vol(F)(b)=\int_{f^{-1}(b)}\eta$ is finite and $\eta$ is closed. 
\item The factor of $-1$ is just for the proper sign conventions, as $d[p]=(-1)^p d$ for the differential of a shifted complex. So $e^\bullet$ is a chain map. 
\item We have $e^p=0$ for $p<r$.
\end{enumerate}

\end{remark}
\begin{proposition}
The map $e^\bullet$ has the desired property, in other words
\begin{equation*}
e=H^p(B,e^\bullet):H^p(E,f^{-1}\Vc)\rightarrow H^{p-r}(B,\Vc\otimes R^df_*(\C))
\end{equation*}
is the map from \Cref{def_tr}. 
\begin{proof}
By \Cref{def_tr} we have to show that the following diagram commutes
\begin{equation*}
\begin{xy}
\xymatrix{
\tau_{\leq r}f_{*}(\Omega_{E} ^\bullet)\ar[d]^{e} \ar[r]^ -{e^\bullet}& R^{r}f_{*}(\C)\otimes\Omega_{B} ^\bullet[-r]\\
R^{r}\phi_{*}(\C)[-r]\ar[r]^{=} & R^{r}\phi_{*}(\C)[-r]\ar[u]
}
\end{xy}
\end{equation*}
where $e$ just comes from the canonical map 
$\ker(f_*d:f_*\Omega_E ^r\to f_*\Omega_E^{r+1})\to R^{r}\phi_{*}(\C)$. As $e^\bullet $ is a chain map and zero in degrees $p<r$, the only thing that has to be proven is that
$e=e^r:\ker(f_*d:f_*\Omega_E ^r\to f_*\Omega_E^{r+1})\to R^{r}\phi_{*}(\C)$
and since both maps factor over exact forms we just have to prove $e^r=\id:R^{r}\phi_{*}(\C)\to R^{r}\phi_{*}(\C)$. The sheaf $R^{r}\phi_{*}(\C)$ is trivialized by the section $\eta$, so it suffices to prove $e^r(\eta)=\eta\otimes 1\in \Gamma(B,R^{r}f_{*}(\C)\otimes\Omega_{B} ^\bullet[-r])$. But we easily calculate
\begin{equation*}
e^r(\eta)=(-1)^{r(r-r)}\eta\otimes \vol(F)^{-1}\int_F \eta=\eta\otimes \vol(F)^{-1} \vol(F)=\eta\otimes 1.
\end{equation*}
\end{proof}
\end{proposition}
\begin{definition}
Set $\left\langle n\right\rangle:=\left\{1,...,n\right\}$ considered as an ordered set. If $A$ is any $\R$-algebra (non-commutative), $I\subset\left\langle n\right\rangle$ an ordered subset and $x\in F_A=\prod_{i=1}^n A$, set $x_I:=\prod_{i\in I}x_i$ and $\Norm(x):=x_{\left\langle n\right\rangle}$. If $i\in \left\langle n\right\rangle$ we also write $I\setminus i$ for $I\setminus \left\{i\right\}$, $I\cup i$ for $I\cup\left\{i\right\}$ and $I^c=\left\langle n\right\rangle\setminus I$ considering them as ordered sets.
Moreover, if $I\subset \left\langle n\right\rangle$ is an subset we denote by $\sgn(I)=\pm1$ the sign of the permutation describing $\left\langle n\right\rangle\to\left\{i_1,...,i_k,j_1,...,j_{n-k}\right\}$ with $I=\left\{i_1,...,i_k\right\}$ and $\left\langle n\right\rangle\setminus I=\left\{j_1,...,j_{n-k}\right\}$.
\end{definition}
\begin{remark}\label{int_decomp}
Let us come back to $\Mc_K$ and summarize what it means to do fiber integration to decompose the cohomology. First of all it is enough to know how things work on connected components of $\Mc_K$, so let us consider $(\H^{\xi}_{\pm}\times (\R_{>0}\left\{\pm1\right\}^\xi\backslash Z(\R)))/G(g_f)$. If we consider $F_\R ^1:=\left\{t\in F_\R ^{\times}:|\Norm(t)|=1 \right\}$, we may also parametrize the space above by $(\H^{\xi}_{\pm}\times (\left\{\pm1\right\}^\xi\backslash F_\R ^1))/G(g_f)$. We denote the typical coordinate on $F_\R ^1$ by $\tilde{r}:=r_{|F_\R ^1}$. The cohomology of the fiber $\left\{\pm1\right\}^\xi\backslash F_\R ^1/\det(Z_K)$ may be generated as a ring by $\frac{d\tilde{r}}{\tilde{r}}_i$, $i=1,...,\xi$. The $\frac{d\tilde{r}}{\tilde{r}}_i$ are not linear independent, because $0=\frac{d(|\Norm(\tilde{r})|)}{|\Norm(\tilde{r})|}=\sum_{i=1}^\xi\frac{d\tilde{r}}{\tilde{r}}_i$. But $\frac{d\tilde{r}}{\tilde{r}}_i$, $i=1,....,\xi-1$, are linear independent and therefore the elements $\frac{d\tilde{r}}{\tilde{r}}_I$, $I\subset \left\langle \xi-1\right\rangle$, form a basis of $\mathfrak H_\C^\bullet=H^\bullet(\left\{\pm1\right\}^\xi\backslash F_\R ^1/\det(Z_K),\C)$. We get the volume form 
$\vol_{F_\R ^1}:=\frac{d\tilde{r}}{\tilde{r}}_{\left\langle \xi-1\right\rangle}$ satisfying the formula $\vol_{F_\R ^1 }\wedge\frac{dt}{t}=\frac{dr}{r}_{\left\langle \xi\right\rangle}=:\vol_{F_\R ^\times}$ for $t:=|\Norm(r)|$. Note that the volume form $\vol_{F_\R ^1 }$ allows integration on $F_\R^1$ and induces a Haar measure $d^\times\tilde{r}$ on $F_\R^1$. We fix the measure $\frac{dt}{t}$ on $\R^\times$ and the product measure on $\prod_{\nu|\infty}\R^\times=\F_\R^\times$, which we denote by $d^\times r$. The Haar measure $d^\times\tilde{r}$ is now uniquely determined by the Fubini formula $d^\times r=d^\times \tilde{r}\frac{dt}{t}$. 
The fiber $\left\{\pm1\right\}^\xi\backslash F_\R ^1/\det(Z_K)$ is also oriented with respect to the volume form $\vol_{F_\R ^1}$ and we define integration on this space with respect to this volume form. Set $R_K:=\vol(\left\{\pm1\right\}^\xi\backslash F_\R ^1/\det(Z_K))$. The trace morphism $\int:\mathfrak H_\C ^{\xi-1}\to \C$ is given by integration of top forms and the basis which under Poincaré duality is dual to $(\frac{d\tilde{r}}{\tilde{r}}_I)_{I\subset \left\langle \xi-1\right\rangle}$ is $(\frac{(-1)^{\sgn(I)}}{R_K}\frac{d\tilde{r}}{\tilde{r}}_{\left\langle \xi-1\right\rangle\setminus I})_{I\subset \left\langle \xi-1\right\rangle}$. 

Now we want to make \cref{cup_decomp} explicit. Suppose we are given a differential form $\omega=\sum_{J\subset \left\langle \xi\right\rangle}\omega_J(r)\wedge\frac{dr}{r}_J$, such that $\omega_J(r)$ is a form of degree zero with respect to $r$, representing a cohomology class in
$H^\bullet(\Mc_K,\phi_K ^{-1}\Sym^k\Hc ^\prime)$. First we may pullback the form to $(\H^{\xi}_{\pm}\times (\left\{\pm1\right\}^\xi\backslash F_\R ^1))/G(g_f)$ giving the form $\sum_{J\subset \left\langle \xi\right\rangle}\omega_J(\tilde{r})\wedge\frac{d\tilde{r}}{\tilde{r}}_J$. 
We rewrite this form as 
\begin{equation*}
\sum_{J\subset \left\langle \xi\right\rangle}\omega_J(\tilde{r})\wedge\frac{d\tilde{r}}{\tilde{r}}_J=\sum_{I\subset\left\langle \xi-1\right\rangle}\tilde{\omega}(\tilde{r})_I\wedge \frac{d\tilde{r}}{\tilde{r}}_{I}
\end{equation*}
where the $\tilde{\omega}_I$ are just linear combinations of the $\omega_J$. So we get by \cref{cup_decomp} 
\begin{equation*}
\omega =\sum_I \int(\tilde{\omega}(\tilde{r})_I\wedge \frac{d\tilde{r}}{\tilde{r}}_{I}\wedge \frac{(-1)^{\sgn (I)}}{R_K}\frac{d\tilde{r}}{\tilde{r}}_{\left\langle \xi-1\right\rangle\setminus I})\wedge \frac{d\tilde{r}}{\tilde{r}}_{I}=\sum_I \frac{1}{R_K}\int(\tilde{\omega}(\tilde{r})_I\wedge \vol_{F_\R ^1})\wedge \frac{d\tilde{r}}{\tilde{r}}_{I}
\end{equation*}
But we have $\int(\tilde{\omega}(\tilde{r})_I\wedge \vol_{F_\R ^1})=\int_{\left\{\pm1\right\}^\xi\backslash F_\R ^1/\det(Z_K)}\tilde{\omega}(\tilde{r})_Id^\times\tilde{r}$ by the definition of fiber integration, see \cite{HalGr}. So we get
\begin{equation*}
\omega =\sum_{I} \frac{1}{R_K}\int_{\left\{\pm1\right\}^\xi\backslash F_\R ^1/\det(Z_K)}\tilde{\omega}(\tilde{r})_Id^\times\tilde{r}\wedge \frac{d\tilde{r}}{\tilde{r}}_{I}=\frac{1}{R_K}\sum_{J}\int_{\left\{\pm1\right\}^\xi\backslash F_\R ^1/\det(Z_K)}\omega_J(\tilde{r})d^\times\tilde{r}\wedge\frac{d\tilde{r}}{\tilde{r}}_J.
\end{equation*}

\end{remark}

\chapter{Comparison with Harder's Eisenstein classes}

\cite{Ki1} and \cite{Bl1} have shown that polylogarithmic Eisenstein classes may restrict non-trivially to cohomology classes of the boundary of the Borel-Serre compactification of Hilbert-Blumenthal varieties. This is done by evaluating the polylogarithmic Eisenstein classes on the boundary.\\
Harder on the other hand started with cohomology classes on the boundary and constructed an operator from the cohomology of the boundary into the cohomology of the Hilbert-Blumenthal variety, which is a section for the restriction map. The image of this operator is called the Eisenstein cohomology. \\
Now we want to represent the polylogarithmic Eisenstein classes by differential forms so that we can compare them with Harder's Eisenstein cohomology classes finally.

\section{Nori's calculation of the polylogarithm as current}

One of the big advantages of the polylogarithm is that it can actually be calculated. For example, this has been done by \cite{No} and A. Levin in \cite{L} with $\C$-coefficients by using currents. We will follow Nori's approach. We recall a well known lemma which justifies the extension of coefficients to $\C$.
\subsection{Integral cohomology  classes}
\begin{lemma}\label{ext_scalars}
Let $X$ be a topological manifold which has a finite good cover $\mathfrak U=(U_i)_{i=1,...,m}$, in other words a cover by finitely many open sets each homeomorphic to $\R^n$ such that all finite intersections $U_{i_0}\cap...\cap U_{i_p}$ are again homeomorphic to $\R^n$. Moreover, we consider $\Vc$ a locally constant sheaf of $A$-modules on $X$ and $A\rightarrow R$ a flat ring extension. 
Then the natural morphism $H^\bullet(X,\Vc)\otimes_A R\rightarrow H^\bullet(X,\Vc\otimes_A R)$ is an isomorphism.
\begin{proof}
Let $\mathfrak U=(U_i)$ be a finite good cover of $X$. 
We have the the $\check C$ech to derived functor spectral sequence \cite{Go} Théorème 5.4.1. $E_2 ^{p,q} =\check{H}^p(\mathfrak U, \underline H^q(\Vc))\Rightarrow H^{p+q}(X,\Vc)$, with the presheaf $\underline H^q(\Vc)(V)=H^q(V,\Vc)$, $V\subset X$ open. Consider $U_{i_0},...,U_{i_p}\in \mathfrak U$. As $\mathfrak U$ is good, we get that $U_{i_0}\cap...\cap U_{i_p}$ is homeomorphic to $\R^n$ and in particular contractible. It follows that $\Vc_{|U_{i_0}\cap...\cap U_{i_p}}$ is constant and therefore
\begin{equation*}
\underline H^q(\Vc)(U_{i_0}\cap...\cap U_{i_p})=H^q(U_{i_0}\cap...\cap U_{i_p},\Vc)=0,\ q>0,
\end{equation*}
by homotopy invariance of sheaf cohomology with constant coefficients \cite{Iv} IV. Theorem 1.1. We conclude
$C^p(\mathfrak U, \underline H^q(\Vc))=0$, if $q>0$, and therefore $E_2 ^{p,q}=0$, $q\neq 0$. The spectral sequence degenerates and the edge morphisms yield isomorphisms
\begin{equation*}
H^p(X,\Vc)=\check{H}^p(\mathfrak U, \underline H^0(\Vc))=\check{H}^p(\mathfrak U, \Vc).
\end{equation*}
If $A\rightarrow R$ is a flat ring extension, we see
\begin{equation*}
H^p(X,\Vc\otimes_A R)=\check{H}^p(\mathfrak U,\Vc\otimes_A R)=\check{H}^p(\mathfrak U,\Vc)\otimes_A R=H^p(X,\Vc)\otimes_A R,
\end{equation*}
as the covering $\mathfrak U$ is finite and $\Gamma(U_{i_0}\cap...\cap U_{i_p},\Vc\otimes_A R)=\Gamma(U_{i_0}\cap...\cap U_{i_p},\Vc)\otimes_A R$, since $\Vc$ is constant on $U_{i_0}\cap...\cap U_{i_p}$.
\end{proof}
\end{lemma}
\begin{remark}
Let $X$ a topological manifold with a finite good cover and $\Vc$ an $A$-local system. By induction on the cardinality of good covers one easily shows using the Mayer-Vietoris sequence that $H^\bullet(X,\Vc)$ is a finitely generated $A$-module. Let now $A\subset \C$ be a principal ideal domain. By \Cref{ext_scalars} we may embed our cohomology groups of local $A$-systems - up to torsion - in the cohomology groups with $\C$-coefficients. More precisely, $H^\bullet(X,\Vc)$ decomposes into a direct sum of a free $A$-module $H^\bullet(X,\Vc)_{\text{free}}$ and a torsion $A$-module $H^\bullet(X,\Vc)_{\text{tor}}$. We have isomorphisms 
\begin{equation*}
H^\bullet(X,\Vc)_{\text{free}}\otimes_A\C\cong H^\bullet(X,\Vc)\otimes_A\C \cong H^\bullet(X,\Vc\otimes_A\C)
\end{equation*}
and $H^\bullet(X,\Vc)_{\text{free}}$ is isomorphic to $\image(H^\bullet(X,\Vc)\stackrel{can}{\rightarrow}H^\bullet(X,\Vc\otimes_A\C))$.
We call the latter group the $A$-integral or simply the \textit{integral classes} inside $H^\bullet(X,\Vc\otimes_A\C)$.
If $X$ is a manifold, integral classes may be calculated by differential forms or currents with values in these local systems. In particular, we will represent the integral part of the polylogarithm by a current as proposed in \Cref{differential_equ}.
\end{remark}

\subsection{Solving the differential equation}

Let us consider $\pi_W:\Tc_W\rightarrow \Mc_K$. Let us suppose $K=K_N$ and $W=V(\hat\Z)\rtimes K_N$. To find a representative for the polylogarithm we may treat each connected component of $\Mc_{K_N}$ separately. So our topological situation is of the form
\begin{equation*}
\prod_{\nu|\infty}SO(2) \backslash G(\R) \stackrel{\cong}{\rightarrow}M:=\left\{(\tau,r)\in\H_{\pm}^\xi\times Z(\R):\text{Im}(r\tau)\in Z(\R)^0\right\},\ g\mapsto (i\cdot g, \det(g)),
\end{equation*}
\begin{equation*}
\pi:\left(V(\R)/V(g_f)\times M\right)/G(g_f) \rightarrow M/G(g_f) 
\end{equation*} 
and $\gamma=\begin{pmatrix}a&b\\c&d\end{pmatrix}\in G(\R)$ acts on $M$ by $(\tau,r)\gamma:=\left(\frac{b+d\tau}{a+c\tau},r\cdot \det(\gamma)\right)$. We write as before $V(g_f)=\mathfrak a\oplus \mathfrak b$ as a sum of two fractional ideals. Set $L:=V(g_f)$, $\Gamma:=G(g_f)$ and $\tilde{\pi}: \tilde T:=V(\R)/L\times M\rightarrow M$ the projection. Moreover, we consider $\R_{>0}\subset Z(\R)$ acting on $M$ by multiplication. \\
By \Cref{invariant_functionals} we have 
\begin{equation*}
\Rc^\infty\otimes_{\Cc^\infty_{\tilde T/\Gamma}}\otimes \mathscr{D}_{\tilde T/\Gamma}=q_*\left(\prod_{k\geq 0}\Sym^kV(\C)\otimes \mathscr{D}_{\tilde T}\right)^\Gamma,
\end{equation*}
where $q:\tilde T\rightarrow \tilde T/\Gamma$ is the canonical projection.
We will construct a $\Gamma$-invariant $2\xi-1$-current on $\tilde T$, which descends to the right current $\pol(f)$ on $\tilde T/\Gamma$. Let us fix coordinates first. We have $F_\R \stackrel{\cong}{\rightarrow}\prod_\sigma \R$, $x\otimes r\rightarrow (\sigma(x)r)_\sigma$ and after fixing an ordering of the set of embeddings $\sigma:F\rightarrow \R$ we get an isomorphism $\prod_\sigma \R \rightarrow \R^\xi$. Since $V(\R)=F_\R^2$, we take the coordinates $w=(w^1,w^2)$, where $w^{i}=(w^{i}_\sigma)_\sigma$, $i=1,2$. $\H_{\pm}^\xi$ is considered as an open subspace of the complex manifold $F_\C$ and we take the induced coordinates $\tau=(\tau_\sigma)_\sigma$. The corresponding real coordinates are denoted by $(x,y)$, i.e. $\tau=x+iy$. The ($F_\C$-valued) differential forms $dw^1,dw^2$ are globally defined on $\tilde T$ making the algebra of global differentials on $\tilde T$ into a free $\Cc^\infty(\tilde T)$-module of finite rank. 
We trivialize $or^{-1}_{M/\Gamma}\otimes \C$ and $or^{-1}_{\R_{>0}\backslash M/\Gamma}\otimes \C$ by choosing volume forms 
$\vol_{M/\Gamma}$ and $\vol_\Mc:=\vol_{\R_{>0}\backslash M/\Gamma}$. We also need an orientation of the fibers of our family of tori. We take $\theta:=\frac{1}{\vol(T)}\bigwedge_\sigma dw^1_\sigma\wedge dw^2_\sigma$, where the volume is calculated with respect to $\vol_V:=\bigwedge_\sigma dw^1_\sigma\wedge dw^2_\sigma$, but the current we are going to construct will not depend on the particular choice of $\theta$. 

We can identify distributions $\phi$ on $T$ with $0$-currents, when we use our volume form $\theta$: Given a form $\eta$ with compact support of top degree on $T$ we may write it uniquely as $f\theta$ with a compactly supported $\Cc^\infty$-function $f$ on $T$. Then we define $\phi(\eta)=:\phi(f)$. For example, we have the distribution $\delta_f=\sum_{t\in T}f(t)\delta_t$, where $f:T\rightarrow \C$ is a map of finite support and $\delta_t$ is the Delta distribution in $t$.\\
Currents $\phi$ on $T$ can be considered as linear functionals $\phi:\Omega_c(\tilde T)\rightarrow \Omega_c(M)$; compare \cite{DR} Th\'{e}or\`{e}me 9. In this manner we can consider $\delta_f \theta$ as a $2\xi$-current on $\tilde T$ with values in forms on $M$. Given a $\eta\in \Gamma_c(\tilde T, \Cc^\infty _{\tilde T}\otimes_{\tilde\pi ^{-1}\Cc^\infty _M}\tilde\pi ^{-1}\Omega_M)$ we have explicitly
\begin{equation*}
\delta_f \theta(\eta)=\delta_f(w)(\theta\wedge\eta)=\sum_{t\in T}f(t)\eta(t,-)
\end{equation*}
where $\delta_f(w)$ means that $\delta_f$ just acts on the coordinates of $T$. Moreover, we have $1\otimes\delta_f \theta=\sum_{t\in T} f(t)\int_{\left\{t\right\}\times M}$ where $1\otimes\delta_f \theta$ means the tensor product of currents on $M$ and $T$ and integration is defined with respect to our fixed orientations.
In general, all currents on $T$ can be identified with currents on $\tilde T$ with values in forms on $M$ and these can be identified via integration over $M$ with currents on $\tilde T$. So we may identify $1\otimes\delta_f \theta$ with $\delta_f \theta$.

With this overview at hand we can say in which space we want to look for $\pol(f)$. We consider the space generated by forms $\phi \wedge\omega$, where $\phi$ is a $0$-current on $T$ and $\omega$ a smooth form on $\tilde T$. Then $\phi\wedge \omega(\eta)=\phi(w)(\omega\wedge \eta)$ as currents on $\tilde T$ with values in forms on $M$. Consider the the pullback of $Log^\infty$ to $\tilde T$. We identify it with the trivial pro vector bundle $\prod_{k\geq0} \Sym^k V(\C)\otimes \Cc_{\tilde T} ^\infty$. We have the connection $\nabla=d-\kappa$ and are interested in solving the equation $\nabla(\phi)=\delta_{f} \theta$ of currents on $\tilde T$ with values in forms on $M$, where the support of $f$ consists of points of order $N$ on the torus $T$ and $\aug(f)=\sum_{t\in T}f(t)=f(0)=0$.
\begin{remark}\label{delta_int}
Our arithmetic group $\Gamma$ guarantees that $t:(\tau,r)\mapsto (t,\tau,r)$ defines a torsion section of $\tilde T/\Gamma\stackrel{\pi}{\rightarrow}M/\Gamma$ so that $1\otimes\delta_{f} \theta$ descends to a current on $\tilde T/\Gamma$, which can be identified with the functional
$\sum_{t\in T} f(t)\int_{M/\Gamma} t^*=\int_f$
which does not depend on the choice of $\theta$.
\end{remark}
Currents $\phi$ on $T$ have a Fourier expansion. Since each current $\phi$ on $T$ can be uniquely written as differential form with distributional coefficients in the basis $dw^1,dw^2$, we may define Fourier expansion coefficientwise. Therefore it suffices to settle the case, when $\phi$ itself is a distribution. But then we have $\phi=\sum_{\lambda\in L^*}\phi_\lambda \exp(2\pi i \lambda(w))$, where $Hom(L,\Z)=:L^*$ and $\phi_\lambda=\phi(\exp(-2\pi i\lambda(w)))$. Now we may evaluate $\phi$ as 
\begin{equation*}
\phi(g)= \sum_{\lambda\in L^*}\phi_\lambda g_{-\lambda},\ \text{with }g=\sum_{\lambda\in L^*}g_\lambda \exp(2\pi i \lambda(w))\in \Cc^\infty (T).
\end{equation*}
For this see \cite{S} VII 1.
\begin{remark}\label{eval_distr}
If $\phi\in \Cc^\infty(T)$ and $T_\phi$ denotes the distribution $T_\phi(f):=\int_Tf\phi dx$, we have $\phi=\vol(T)^{-1}\sum_\lambda T_{\phi,\lambda}\exp(2\pi i \lambda(w))$.
\end{remark}
So we have a Fourier expansion
\begin{equation*}
\delta_{f} \theta=\sum_{\lambda\in L^*}\sum_{t\in T}f(t)\exp(2\pi i \lambda(w-t))\theta=\sum_{0\neq\lambda\in L^*}\sum_{t\in T}f(t)\exp(2\pi i \lambda(w-t))\theta,
\end{equation*}
as $\aug(f)=0$.
In this spirit let us write $\phi=\sum_{\lambda\in L^*} \phi_{\lambda}\exp(2\pi i \lambda(w))$ as Fourier series where $\phi_{\lambda}$ are smooth forms on $M$ with values in the completed symmetric algebra tensored with closed invariant forms on $T$. Applying $\nabla$ we get the differential equation
\begin{equation*}
d(\phi_\lambda)+A_\lambda \wedge \phi_\lambda=\sum_{t\in T}f(t)\exp(2\pi i \lambda(-t))\theta,\ \lambda\in L^*,\ A_\lambda := 2\pi i\, d\lambda - \kappa
\end{equation*}
for each Fourier coefficient. We try to solve this equation in
\begin{equation*}
R=\prod_k \Sym^k V(\C)\otimes\Omega(M)\otimes_\C \bigwedge ^\bullet V(\C)^*\subset \prod_k \Sym^k V(\C)\otimes\Omega_{\tilde T}(\tilde T)
\end{equation*}
$R=\bigoplus R^{p,q}$ is bigraded: $R^{p,q}:=\prod_k \Sym^k V(\C)\otimes\Omega^p(M)\otimes_\C \bigwedge ^q V(\C)^*$. Let us first consider the case $\lambda= 0$. We get the equation $d(\phi_0)-\kappa \wedge \phi_0=0$, as $\delta_f$ has $0$ zeroth Fourier coefficient. So we simply set $\phi_0=0$. We endow $V(\R)$ with a symplectic structure:
\begin{equation*}
\left\langle .,.\right\rangle:V(\R)\wedge V(\R)\rightarrow \R,\ w_1\wedge w_2 \rightarrow  \Tr( \det(w_1,w_2))
\end{equation*}
with $\Tr(x\otimes y):=\Tr_{F/\Q}(x\otimes y)=\sum_\sigma \sigma(x)y$ for $x\otimes y\in F_\R$. Next we fix global ($F_\C$-valued) vector fields $u:=\frac{\overline{\tau}e^1-e^2}{\overline{\tau}-\tau}$, $\overline{u} :=-\frac{\tau e^1-e^2}{\overline{\tau}-\tau}$ on $\tilde T$, where the tangent bundle of $T$ is trivialized by $e^1\mapsto \partial_{w_1}$, $e^2\mapsto \partial_{w_2}$. These vector fields reflect the fact that we have a decomposition into types $\Hc_\Z\otimes\Cc^\infty _{M/\Gamma}=\Hc^{0,-1}\oplus \Hc^{-1,0}$ coming from the Hodge theory of the fibers. The two subspaces are Lagrangian with respect to $\left\langle .,.\right\rangle_\C$ and we have $\det( u,\overline u)= \frac{1}{\overline \tau-\tau}$. Whenever $s$ is a local section of $\Hc_\Z\otimes\Cc^\infty _{M/\Gamma}$ denote by $s^{0,-1}$, $s^{-1,0}$ the decomposition into types. In this spirit $u$ defines a frame for $\Hc^{-1,0}$ and $\overline{u}$ for $\Hc^{0,-1}$ and, since we have
\begin{equation*}
w=w^1e^1+w^2e^2=(w^1+\tau w^2)u+(w^1+\overline{\tau} w^2)\overline{u}
\end{equation*}
for $ w\in V(\R)$, we see $w^{-1,0}=(w^1+w^2\tau)u$. Moreover, the decomposition into types is $\Gamma$-equivariant and we get
\begin{equation*}
(\gamma w)^{-1,0}=\gamma w^{-1,0},\ \gamma u =(a+c\tau)u,\ \gamma=\begin{pmatrix}a&b\\c&d\end{pmatrix}.
\end{equation*}
Since $\left\langle .,.\right\rangle$ is non-degenerate, we may express $L^*$ by it. If we set $\mathfrak d\subset \Oc$ to be the different-ideal of $\Oc$, we have an isomorphism
\begin{equation*}
\mathfrak b ^{-1} \mathfrak d^{-1}\oplus \mathfrak a ^{-1}\mathfrak d^{-1} \rightarrow L^*,\ l\mapsto \left\langle l,.\right\rangle:=\lambda
\end{equation*}
of $\Oc$-modules. We will stick to the vector fields $\frac{l^{-1,0}}{r}$, $0\neq l =(l_1,l_2)\in L^*\subset V(\R)$. We have $\gamma\frac{l^{-1,0}}{r}=\frac{\left(\det\gamma^{-1}\gamma l\right)^{-1,0}}{r}$ compatible with the natural $\Gamma$ action on $L^*$
\begin{equation*}
\gamma\lambda(w):=\lambda(\gamma^{-1}w)=\left\langle l,\gamma^{-1}w\right\rangle=\left\langle \det\gamma^{-1} \gamma l,w\right\rangle,\ \lambda=\left\langle l,.\right\rangle.
\end{equation*}
Now we follow Nori's construction of $\pol(f)$. The vector field $\frac{l^{-1,0}}{r}$ on $\tilde T$ defines a linear operator on the differential algebra $R$ onto itself by contraction. Let us denote this operator by $i_l$. Furthermore, define the operator $C_l:=d+A_l\wedge$ on $R$, where $A_l$ is the $A_\lambda$ from above. Then $C_l\circ C_l=0$, $i_l\circ i_l=0$ and $C_l\circ i_l+i_l\circ C_l$ is an isomorphism. The first assertion is immediately clear. For the second write $F_l=(d i_l +i_l d) +(A_l i_l+i_l A_l)$. The first summand is the Lie derivative $\Lc_{\frac{l^{-1,0}}{r}}$ with respect to the vector field $\frac{l^{-1,0}}{r}$. It is a sum of maps $R^{p,q}\rightarrow R^{p+1,q-1}$ and therefore $\Lc_{\frac{l^{-1,0}}{r}}^{2\xi+1}=0$ and $\Lc_{\frac{l^{-1,0}}{r}}$ is nilpotent on $R$. 
\begin{equation*}
(A_l i_l+ i_l A_l)(\omega)=A_l(i_l\omega)+i_l(A_l\wedge\omega)=A_l(i_l\omega)+i_l(A_l)\wedge\omega-A_l\wedge i_l(\omega)
\end{equation*}
by the rules of interior multiplication and therefore
\begin{equation*}
(A_l i_l+ i_l A_l)(\omega)=i_l(A_l)\wedge\omega=i_l(A_l)\omega
\end{equation*}
is just multiplication by
\begin{equation*}
i_l(A_l)=2\pi i\, d\lambda(\frac{l^{-1,0}}{r})-\kappa(\frac{l^{-1,0}}{r})=2\pi i\left\langle l,\frac{l^{-1,0}}{r}\right\rangle - \frac{l^{-1,0}}{r}
\end{equation*}
a function with values in the completed symmetric algebra of $V(\C)$. We calculate
\begin{equation*}
\left\langle l,\frac{l^{-1,0}}{r}\right\rangle=\left\langle l^{0,-1},\frac{l^{-1,0}}{r}\right\rangle,
\end{equation*}
since the Hodge decomposition is Lagrangian, and endow $F_\C$ with the natural Hermitian metric induced by $\Tr$ to get
\begin{equation*}
2\pi i\left\langle l^{0,-1},\frac{l^{-1,0}}{r}\right\rangle=\pi \left\|\frac{l_1+l_2\tau}{\sqrt{r\text{Im}(\tau)}}\right\|^2=F_l (\tau,r).
\end{equation*}
Note $\gamma F_l(\tau,r)=F_l(\tau\gamma,r\det\gamma)=F_{\det\gamma^{-1}\gamma l}(\tau,r)$. Since $l\neq 0$, we have $F_l(\tau,r)\in \C^{\times}$ and therefore $i_l(A_l)$ is certainly invertible in $\prod_k \Sym^k V(\C)\otimes \Cc^\infty_M $, because the leading coefficient is invertible. It follows that $A_l i_l+ i_l A_l$ is invertible.\\
Now $\Lc_{\frac{l^{-1,0}}{r}}$ and $i_l(A_l)$ commute, since the vector field $\frac{l^{-1,0}}{r}$ takes its values in the tangent space of $T$ and therefore its Lie derivative is $\Cc_M ^\infty$-linear. We conclude that the sum $\Lc_{\frac{l^{-1,0}}{r}}+i_l(A_l)$ is also invertible. Explicitly,
\begin{equation*}
(i_lC_l+C_li_l)^{-1}=(\Lc_{\frac{l^{-1,0}}{r}}+i_l(A_l))^{-1}=\sum_{k\geq0}(-1)^k i_l(A_l)^{-(k+1)}\Lc_{\frac{l^{-1,0}}{r}}^k.
\end{equation*}
This is a finite sum.
Set $\omega_{t,l}:=\exp(-2\pi i\lambda(t))\theta$. Then $i_l(i_lC_l+C_li_l)^{-1}(\omega_{t,l})$ is a solution for $d(\phi_\lambda)+A_\lambda \wedge \phi_\lambda=C_l(\phi_\lambda)=\exp(2\pi i \lambda(-t))\theta$. To prove this note that $C_l$ and $i_lC_l+C_li_l$ commute. Therefore $C_l$ and $(i_lC_l+C_li_l)^{-1}$ commute and
\begin{equation*}
\omega_{t,l}=(i_lC_l+C_li_l)(i_lC_l+C_li_l)^{-1}(\omega_{t,l})=C_li_l(i_lC_l+C_li_l)^{-1}(\omega_{t,l}),
\end{equation*}
since $i_lC_l(i_lC_l+C_li_l)^{-1}(\omega_{t,l})=i_l(i_lC_l+C_li_l)^{-1}C_l(\omega_{t,l})=0$, as $C_l\omega_{t,l}=0$. Because $i_l\Lc_{\frac{l^{-1,0}}{r}}=i_ldi_l$ we get 
\begin{equation*}
\phi_l=\phi_\lambda=\sum_{t\in T}f(t)\sum_{k\geq0}(-1)^k i_l(A_l)^{-(k+1)}i_l(di_l)^k\omega_{t,l}
\end{equation*}
So we should have
\begin{equation*}
\phi=  \sum_{k\geq0} \sum_{0\neq l\in L^*}(-1)^k \frac{\sum_{t\in T}f(t)\exp(2\pi i\left\langle l,-t\right\rangle)i_l(d i_l)^k\theta}{\left(\pi \left\|\frac{l_1+l_2\tau}{\sqrt{r\text{Im}(\tau)}}\right\|^2-\frac{l^{-1,0}}{r}\right)^{k+1}}\exp(2\pi i\left\langle l,w\right\rangle)
\end{equation*}
The sum over $k$ only has non-trivial summands between zero and $2\xi-1$. We rewrite
\begin{equation*}
\left(\pi \left\|\frac{l_1+l_2\tau}{\sqrt{r\text{Im}(\tau)}}\right\|^2-\frac{l^{-1,0}}{r}\right)^{-(k+1)}
\end{equation*}
as binomial series
\begin{equation*}
\sum_{n\geq0} \frac{(k+n)!}{k!n!} \left\|\frac{\sqrt{\pi}(l_1+l_2\tau)}{\sqrt{r\text{Im}(\tau)}}\right\|^{-2(k+n+1)}\left(\frac{l^{-1,0}}{r}\right)^{\otimes n}.
\end{equation*}
So in total $\phi$ should be given by
\begin{equation*}
\sum_{n,k\geq0}\sum_{0\neq l\in L^*}  \frac{(k+n)!(-1)^{k}\left(\frac{l^{-1,0}}{r}\right)^{\otimes n}\sum_{t\in T}f(t)\exp(2\pi i\left\langle l,w-t\right\rangle)i_l(d i_l)^k\theta}{k!n!\left\|\frac{\sqrt{\pi}(l_1+l_2\tau)}{\sqrt{r\text{Im}(\tau)}}\right\|^{2(k+n+1)}} 
\end{equation*}
We want to see that this expression defines a current with values in the completed symmetric algebra. To do so we have to control the coefficients of the Fourier series.
The first part is to control 
\begin{equation*}
\frac{1}{\left\|\frac{\sqrt{\pi}(l_1+l_2\tau)}{\sqrt{r\text{Im}(\tau)}}\right\|^{2(k+n+1)}}.
\end{equation*}
Consider $\frac{\sqrt{\pi}(l_1+l_2\tau)}{\sqrt{r\text{Im}(\tau)}}$ as a continuous map $S:V(\R)\times M\rightarrow F_\C$, which is $F_\R$-linear in the first argument.
Take any compact set $K$ in $M$ and $K^\prime$ the one-sphere in $V(\R)$. The product is again compact and $\left\|S\right\|$ has a minimal value on $K^\prime\times K$. Since $S(.,\tau,r)$ is injective for each $(\tau,r)$, this minimum has to be bigger than zero. Therefore, there is a constant $C>0$ such that $\left\|S(v,\tau,r)\right\|>C\left\|v\right\|$ for all $0\neq v\in V(\R)$ and $(\tau,r)\in K$. It follows that 
\begin{equation*}
\frac{1}{\left\|\frac{\sqrt{\pi}(l_1+l_2\tau)}{\sqrt{r\text{Im}(\tau)}}\right\|^{2(k+n+1)}}\leq \frac{1}{C\left\|l\right\|^{2(k+n+1)}}
\end{equation*}
for all $l$ and for any $(\tau,r)$ in a compact set $K$ with $C$ only depending on $K$.
Next we have to consider $\left(\frac{l^{-1,0}}{r}\right)^{\otimes n}i_l(d i_l)^k\theta$. We have to estimate each coefficient in front of an $\Cc^\infty (M)$-basis of $\Sym^k V(\C)\hat{\otimes}\Omega^p(M)\otimes_\C \bigwedge ^q V(\C)^*$. Each coefficient can be controlled by $C^\prime\left\|l\right\|^{n+k+1}$ for $(\tau,r)$ in a compact set.\\
Now we want to see that each coefficient of $\phi$, which is given as a formal Fourier series, actually converges to a distribution on $\tilde T$. But we have just shown that the coefficients can be controlled by $C\left\|l\right\|^{n+k+1}$ for all $(\tau,r)\in K$ compact. \cite{S} VII.1 tells us that $\phi$ is a current, since our estimates are locally uniform in $(\tau,r)$. 
The next step is to prove that the current $\phi$ is $\Gamma$-invariant. First we have $\gamma^*\phi_l= \phi_{\det\gamma^{-1}\gamma l}$, as
\begin{equation*}
\gamma^*i_l(d i_l)^k\theta=i_{\det\gamma^{-1}\gamma l}d(i_{\det\gamma^{-1}\gamma l})^k \theta,
\end{equation*}
since $\gamma^*\circ i_l=i_{\det\gamma^{-1}\gamma l}\circ \gamma^*$ and $\theta$ is $\Gamma$-invariant. Therefore
\begin{equation*}
\gamma^*\sum_l \phi_l \exp(2\pi i\left\langle l,w\right\rangle)=\sum_l(\gamma^*\phi_l) \exp(2\pi i\left\langle l,\gamma^{-1}w\right\rangle)=
\end{equation*}
\begin{equation*}
\sum_l\phi_{\det\gamma^{-1}\gamma l} \exp(2\pi i\left\langle \det\gamma^{-1}\gamma  l,w\right\rangle).
\end{equation*}
Evaluating on a test form and using locally uniform convergence we conclude that this equals $\sum_l\phi_{ l} \exp(2\pi i\left\langle l,w\right\rangle)$ again. Finally, we have thanks to \cite{No}.
\begin{theorem}
\begin{equation*}
\phi= \sum_{n,k\geq0}\sum_{0\neq l\in L^*}  \frac{(k+n)!(-1)^{k}\left(\frac{l^{-1,0}}{r}\right)^{\otimes n}\sum_{t\in T}f(t)\exp(2\pi i\left\langle l,w-t\right\rangle)i_l(d i_l)^k\theta}{k!n!\left\|\frac{\sqrt{\pi}(l_1+l_2\tau)}{\sqrt{r\text{Im}(\tau)}}\right\|^{2(k+n+1)}} \otimes \vol_{\Mc}
\end{equation*}
defines a $\R_{>0}\times \Gamma$-invariant current with values in the completed symmetric algebra on $\tilde T$ that satisfies $\nabla(\phi)=\delta_f\theta$. In particular, $\phi$ induces a current on $\R_{>0}\backslash\tilde T/\Gamma$ with values in $Log^\infty\otimes \pi^{-1}or^{-1}_{\R_{>0}\backslash M/\Gamma}$, which equals $\pol(f)$ as cohomology class. 
\begin{proof}

See \Cref{differential_equ}, \Cref{invariant_functionals} and \Cref{delta_int}.
 \end{proof}
\end{theorem}

\section{Global description of the polylogarithm in adelic coordinates}

We already have an explicit description of the polylogarithm as a current on each connected component of $\Tc_W$. However, the choice of a connected component of $\Mc_K$ always means a choice of a specific lattice $L$ inside $V(\R)$, which represents this connected component. We want a description of the polylogarithm, which does not depend on a particular choice of a lattice. To achieve this we use adelic coordinates. We start by recalling Fourier analysis on $V(\A)$.

\subsection{Fourier analysis}\label{Fourier}

Let us take the skew-symmetric non-degenerate bilinear Form
\begin{equation*}
V(\A)\times V(\A)\rightarrow \A,\ \left\langle x,y\right\rangle:=\Tr_{F/\Q}(\det(x,y))
\end{equation*}
Let us fix a character $\psi_0:\A\rightarrow\C^{\times}$ trivial on $\Q$. We define $\psi_{0}=\exp(2\pi i \lambda) $ where we take $\lambda: \A\rightarrow \R/\Z$, $\lambda =\sum_{\nu\text{ place of }\Q} \lambda_\nu$, with
 \begin{equation*}
  \lambda_p: \Q_p\stackrel{can.}{\rightarrow}\Q_p/\Z_p \subset \Q/\Z\subset \R/\Z \text{ and }\lambda_\infty: \R\rightarrow \R/\Z,\ x\mapsto -x\ \text{mod}\ \Z.
  \end{equation*}
We identify $V(\A)$ with its Pontryagin dual by $V(\A)\rightarrow \widehat{V(\A)}$, $x\rightarrow \psi_{0}\left\langle x,\ \right\rangle:=\psi_{0}(\left\langle x,\ \right\rangle)$. Next we need to fix a Haar measure $dv$ on $V(\A)=\A_F^2$. To do so we fix a Haar measure $dx$ on $\A_F$ and take $dv$ to be the product measure. We take Haar measures $dx_\nu$ on $F_\nu$ for each place $\nu$. If $\nu|\infty$, we have $F_\nu=\R$ and we take the Lebesgue measure. If $\nu|p$ we fix $dx_\nu$ by $dx_\nu(\Oc_\nu)=N(\mathfrak d_\nu)^{-\frac{1}{2}}$, where $\mathfrak d_\nu\subset \Oc_\nu$ is the different of $\Oc_\nu/\Z_p$. We take then $dx=\prod_{\nu\text{ place of }F} dx_\nu$ in the sense of Tate, see \cite{T} XV 3.3. In a similar manner we get the Haar measure $dv$ on $V(\A_f)$ and if we restrict $\psi_{0}$ to $V(\A_f)$, we get the topological isomorphism $V(\A_f)\rightarrow \widehat{V(\A_f)}$, $x\rightarrow \psi_{0}\left\langle x,\ \right\rangle$.

For any integrable function $\varphi:V(\A)\rightarrow\C$ and any $g\in G(\A)$ we have the transformation formula
\begin{equation*}
\int_{v\in V(\A)}\varphi(gv)dv=\left\|\det(g)\right\|^{-1}\int_{v\in V(\A)}\varphi(v)dv,
 \end{equation*}
with $\left\|\ \right\|:=\prod_{\nu\text{ place of }F} |\ |_\nu$ the Tate-character and similarly 
\begin{equation*}
\int_{v\in V(\A_f)}\varphi(g_fv)dv=\left\|\det(g_f)\right\|_f ^{-1}\int_{v\in V(\A_f)}\varphi(v)dv
 \end{equation*}
for integrable $\varphi:V(\A_f)\rightarrow\C$. We define the Fourier transform of a Schwartz-Bruhat function $\varphi:V(\A)\rightarrow \C$ by
\begin{equation*}
\hat\varphi(x):=\int_{V(\A)}\varphi(v)\overline{\psi_{0}\left\langle x,v\right\rangle} dv.
\end{equation*}
The measure $dv$ is self dual, in other words, the Fourier inversion formula holds: $\hat{\hat\varphi}=\varphi$. There is no twist by $-1^*$. We have the formula $\widehat{(\varphi\cdot g)}(x)=\left\|\det(g)\right\|^{-1}\hat\varphi (\hat g^{-1}x)$, for $\varphi \in \Sc(V(\A))$ and $g\in G(\A)$. Here $\hat g$ denotes the adjoint operator of $g$ with respect to $\left\langle \ ,\ \right\rangle$. It is explicitly given by $\hat {g}=\det g \cdot g^{-1}$. Similarly, we get $\widehat{(\varphi\cdot g_f)}(x)=\left\|\det(g_f)\right\|_f^{-1}\hat\varphi (\hat g_f^{-1}x)$, for $\varphi \in \Sc(V(\A_f),\C)$ and $g_f\in G(\A_f)$. 

We want to examine $\hat \varphi$ of a Schwartz-Bruhat function $\varphi:V(\A_f)\rightarrow \C$. Let us assume that $\text{supp}(\varphi)\subset N^{-1}V(\hat\Z)$ and $\varphi$ is well-defined modulo $NV(\hat \Z)$. We may write
\begin{equation*}
\varphi=\sum_{\overline u \in N^{-1}V(\hat\Z)/NV(\hat\Z) } \varphi(u)\chi_{\overline u},
\end{equation*}
where $\overline u=u+NV(\hat \Z)$ and $\chi_M$ always denotes the characteristic function of a subset $M$. We calculate
\begin{equation*}
\hat \varphi(x)=\sum_{\overline u \in N^{-1}V(\hat\Z)/NV(\hat\Z) } \varphi(u)\int_{V(\A_f)}\chi_{\overline u}(v)\overline{\psi_{0}\left\langle x,v\right\rangle} dv
\end{equation*}
The latter integral is
\begin{equation*}
\int_{V(\A_f)}\chi_{\overline u}(v)\overline{\psi_{0}\left\langle x,v\right\rangle} dv=\int_{V(\A_f)}\chi_{\overline u}(v+u)\overline{\psi_{0}\left\langle x,v+u\right\rangle} dv=\int_{V(\A_f)}\chi_{\overline u-u}(v)\overline{\psi_{0}\left\langle x,v+u\right\rangle} dv
\end{equation*}
\begin{equation*}
=\int_{V(\A_f)}\chi_{NV(\hat\Z)}(v)\overline{\psi_{0}\left\langle x,v+u\right\rangle} dv=\overline{\psi_{0}\left\langle x,u\right\rangle} \int_{NV(\hat \Z)}\overline{\psi_{0}\left\langle x,v\right\rangle} dv.
\end{equation*}
The last integral is now easy to compute. $\int_{NV(\hat \Z)}\overline{\psi_{0}\left\langle x,v\right\rangle} dv$ is zero, whenever $\overline{\psi_{0}\left\langle x,\ \right\rangle} $ is not the trivial character on $NV(\hat\Z)$. This is exactly the case, when $x\notin N^{-1}\mathfrak d^{-1} V(\hat\Z)$. Otherwise the integral is just $\int_{NV(\hat \Z)}dv=(N^{2\xi}d_F)^{-1}$. This gives the final formula
\begin{equation*}
\hat \varphi(x)=(N^{2\xi}d_F)^{-1}\sum_{\overline u \in N^{-1}V(\hat\Z)/NV(\hat\Z) } \varphi(u)\overline{\psi_{0}\left\langle x,u\right\rangle} \chi_{N^{-1}\mathfrak d^{-1} V(\hat\Z)}(x).
\end{equation*}

\subsection{The adelic polylogarithm}

Let us give the description of the polylogarithm current on $\Tc_W$ in global adelic coordinates.
First we trivialize the orientation bundle by the section 
\begin{equation*}
\vol_{\Tc/\Mc}^*:=\left\| \det g\right\|_f \sgn(\Norm(\det g))^{-1}\bigwedge_\sigma X_\sigma\wedge Y_\sigma\in H^0(\Mc,\mu)=H^0(\Mc,\bigwedge^{2\xi}\Hc)
\end{equation*}  
which is the element dual to the form
\begin{equation*}
\vol_{\Tc/\Mc}:=\frac{\bigwedge_\sigma dw^1_\sigma\wedge dw^2_\sigma}{\left\| \det g\right\|_f \sgn(\Norm(\det g))^{-1}}=\frac{\vol_V}{ \left\| \det g\right\|_f \sgn(\Norm(\det g))^{-1}}. 
\end{equation*}  
Here we extended the norm character $\Norm$ and $\left\|\ \right\|_f$ to the whole of $\I_F$ by setting $\Norm=1$ on the finite ideles and $\left\|.\right\|_f=1$ on $F_\R ^{\times}$. 

Suppose we are given a function in $\Sc(V(\A_f),\mu^{\otimes n})^0 $. Of course, any such function may be written as
\begin{equation*} 
(v,g)\mapsto \left(\left\|\det(g)\right\|_f \sgn(\Norm(\det(g)))^{-1}\right)^n f(v,g)=\varphi(v,g),
\end{equation*}
with $f\in \Sc(V(\A_f),\mu^{\otimes 0})^0 $ such that $f$ factors in the second argument over $K_f=K_{N^2}$, $\text{supp}(f(\ ,g))\subset N^{-1}V(\hat \Z)$ and $f(\ ,g)$ is well-defined modulo $NV(\hat\Z)$ for all $g\in G(\A)$. Let us recall the main steps in the calculation of $\pol(f)=\pol_{W}(f)$, $W_f=NV(\hat\Z)\rtimes K_{N^2}$.
First of all we did our calculations in the coordinates $(w,g)$, see \Cref{connected_component}. This means that we associated to $f$ on each fiber
\begin{equation*} 
V(\R)/V(\Q)\cap (g^{-1}NV(\hat\Z)+V(\R))\cong g^{-1}NV(\hat\Z)\backslash V(\A)/V(\Q)
\end{equation*}
the delta distribution
\begin{equation*} 
\sum_{w\in g^{-1}N^{-1}V(\hat \Z)/g^{-1}NV(\hat \Z)}f(g_fw,g)\delta_{w}=\sum_{v\in N^{-1}V(\hat \Z)/NV(\hat \Z)}f(v,g)\delta_{g_f^{-1}v} .
\end{equation*}
We fixed the lattice $L$ given by $ V(\Q)\cap(Ng^{-1}V(\hat \Z)+V(\R))$ and the dual lattice $L^*=V(\Q)\cap (N^{-1}\mathfrak d^{-1}\hat{g}V(\hat{\Z})+V(\R))$ in $V(\R)$ with respect to $\left\langle \ ,\ \right\rangle$.
The canonical map
\begin{equation*} 
N^{-2}L/L=g^{-1}N^{-1}V(\hat \Z)\cap V(\Q)/g^{-1}NV(\hat \Z)\cap V(\Q)\rightarrow g^{-1}N^{-1}V(\hat \Z)/g^{-1}NV(\hat \Z)
\end{equation*}
is an isomorphism. We constructed the polylogarithm $\pol(f)$ with respect to \\
$\sum_{w\in N^{-2}L/L}f(g_fw,g)\delta_{w}$ and $\theta=\frac{\vol_{\Tc/\Mc}}{d_F N^{2\xi}}$ as the current
\begin{equation*}
\sum_{j,k\geq0} \sum_{0\neq l\in L^*}  \frac{(j+k)! (-1)^{k}\tilde f(l,g)\left(\frac{l^{-1,0}}{r}\right)^{\otimes j}i_l(d i_l)^k\frac{\vol_{\Tc/\Mc}}{d_F N^{2\xi}}}{j!k!\left\|\frac{\sqrt{\pi}(l_1+l_2\tau)}{\sqrt{r\text{Im}(\tau)}}\right\|^{2(j+k+1)}} \exp(2\pi i\left\langle l,w_\infty\right\rangle)\otimes \vol_\Mc,
\end{equation*}
where $\tilde f(l,g)$ was actually given by
\begin{equation*}
\sum_{w\in N^{-2}L/L}f(g_fw,g)\exp(2\pi i\left\langle l,-w\right\rangle)=\sum_{v\in N^{-1}V(\hat \Z)/NV(\hat \Z)}f(v,g)\psi_{0}\left\langle l,-g_f^{-1}v\right\rangle,\ l\in L^*.
\end{equation*}
On the left-hand side of the last equation we interpret $w\in N^{-2}L \subset V(\R)$ in the infinite component, whereas on the right-hand side we have $v\in N^{-1}V(\hat \Z)\subset V(\A_f)$.
If we have Fourier analysis for Schwartz-Bruhat functions on $V(\A_f)$ in mind, we see $(N^{2\xi}d_F)^{-1}\tilde f(l,g)=\hat f(\hat g^{-1} l,g)$.
Here $\hat f$ means that we do Fourier transformation in the first argument. More precisely, $\hat{g}^{-1\ *}\chi_{N^{-1}\df^{-1}V(\hat\Z)}=\chi_{N^{-1}\df^{-1}\hat{g}V(\hat\Z)}$, so $\hat{g}^{-1\ *}\chi_{N^{-1}\df^{-1}V(\hat\Z)}(l)$ is not equal to zero for $l\in V(\Q)$ if and only if $l\in L^*$.
With this notation at hand we may write $\pol(f)$ as
\begin{equation*}
\sum_{j,k\geq0} \sum_{l\in V(\Q)\setminus \left\{0\right\}}  \frac{(j+k)!(-1)^{k}\hat f(\hat g^{-1}l,g)\left(\frac{l^{-1,0}}{r}\right)^{\otimes j}i_l(d i_l)^k\vol_{\Tc/\Mc}}{j!k!\left\|\frac{\sqrt{\pi}(l_1+l_2\tau)}{\sqrt{r\text{Im}(\tau)}}\right\|^{2(j+k+1)}} \overline{\psi_0\left\langle l,w_\infty\right\rangle}\otimes \vol_\Mc.
\end{equation*}
If we multiply with our trivialization $\vol_{\Tc/\Mc}^{*n}$, we may write $\pol(\varphi)$ as
\begin{equation*}
\sum_{j,k\geq0}\sum_{l\in V(\Q)\setminus \left\{0\right\}} \frac{(j+k)!(-1)^{k}\hat f(\hat g^{-1}l,g)\left(\frac{l^{-1,0}}{r}\right)^{\otimes j}i_l(d i_l)^k\vol_{\Tc/\Mc}\otimes\vol_{\Tc/\Mc}^{*n}}{j!k!\left\|\frac{\sqrt{\pi}(l_1+l_2\tau)}{\sqrt{r\text{Im}(\tau)}}\right\|^{2(j+k+1)}} \overline{\psi_0\left\langle l,w_\infty\right\rangle}\otimes \vol_\Mc.
\end{equation*}

\subsection{Analytic continuation of the polylogarithm}

If we want to specialize the polylogarithm along the zero section, we need to know, whether this current may be represented by a smooth differential form on $U$. This is what Levin proved in \cite{Bl1}. Let us quickly recall how things work. We have identifications $\Oc\otimes\R^2\cong \Oc\otimes\C$, $e^1\mapsto 1$, $e^2\mapsto i$,
\begin{equation*}
(\prod_{\nu|\infty}SO(2))\backslash G(\R)\rightarrow M,\ g_\infty\mapsto (ig_\infty,\det g_\infty)=(\tau,r),
\end{equation*}
\begin{equation*}
\left\|\frac{\sqrt{\pi}(l_1+l_2\tau)}{\sqrt{r\text{Im}(\tau)}}\right\|=\pi\left\|\hat g_\infty ^{-1}l\right\|,\ \text{and }
\frac{\Gamma(s)}{x^s}=\int_{\R_{>0}}\exp(-ux){u^s}\frac{d u}{u},\ x\in\R_{>0}.
\end{equation*}
This yields
\begin{equation*}
\left\|\frac{\sqrt{\pi}(l_1+l_2\tau)}{\sqrt{r\text{Im}(\tau)}}\right\|^{-2(j+k+1)}=\int_{\R_{>0}}\exp(-u\pi\left\|\hat g_\infty ^{-1}l\right\|^2)\frac{u^{j+k+1}}{(j+k)!}\frac{du}{u}.
\end{equation*}
Given $\varphi_f\in \Sc(V(\A_f),\C)$ and $\varphi_\infty\in \Sc(V(\R),\C)$ we define $\varphi_f\otimes \varphi_\infty \in\Sc(V(\A),\C)$ by 
\begin{equation*}
\varphi_f\otimes \varphi_\infty(v):=\varphi_f(v_f)\varphi_\infty(v_\infty),\ v\in V(\A).
\end{equation*}
$w\mapsto \exp(-\pi \left\|w\right\|^2)$ is fixed by Fourier transformation, so we have a $(V(\A_f)\rtimes G(\A_f))\times\pi_0(G_\infty)$-equivariant linear map
\begin{equation*}
\phi:\Sc (V(\A_f))\rightarrow \Sc(V(\A)),\ f\mapsto f\otimes \exp(-\pi\left\|x\right\|^2),
\end{equation*}
which is compatible with Fourier transformation, in other words, $\widehat{\phi(f)}=\phi(\hat f)$.
We can write the polylogarithm as 
\begin{equation*}
\sum_{j,k\geq0}\frac{(j+k)!}{j!k!}(-1)^{k}\phi^{(j,k)}\otimes \vol_{\Tc/\Mc}^{*n}\otimes\vol_\Mc  ,
\end{equation*}
\begin{equation*}
\phi^{(j,k)}:= \sum_{l\in V(\Q)\setminus \left\{0\right\}}  \int_{u\in \R_{>0}}\frac{\phi(\hat{f})(\sqrt{u}\hat g^{-1}l,g)\left(\frac{l^{-1,0}}{r}\right)^{\otimes j}u^{j+k+1}\frac{du}{u}i_l(d i_l)^k \vol_{\Tc/\Mc}}{(j+k)!}
\overline{\psi_0\left\langle l,w_\infty\right\rangle}
\end{equation*}
The first observation is that for $k+j+1>2\xi$ we already have honest differential forms, because the sum over $l$ converges locally uniformly. For these big $k+j+1$ we may also put the integral over $u\in \R_{>0}$ in front of all sums by uniform convergence. Let us consider the differential form
\begin{equation*}
\phi_s ^{(j,k)}:=\int_{u\in \R_{>0}}\sum_{l\in V(\Q)\setminus \left\{0\right\}}\frac{\phi(\hat{f})(\sqrt{u}\hat g^{-1}l,g)\left(\frac{l^{-1,0}}{r}\right)^{\otimes j}u^{s}\frac{du}{u}i_l(d i_l)^k \vol_{\Tc/\Mc}}{(j+k)!}\cdot \overline{\psi_0\left\langle l,w_\infty\right\rangle}
\end{equation*}
for complex $s\in \C$, $\Real(s)>2\xi$. Levin proves by using Poisson summation that $\phi^{(j,k)}_s$ has an analytic continuation as a holomorphic function in $s\in\C$ on $U$. This analytic continuation yields the representation of the polylogarithmic current as a differential form on $U$, see also \cite{BKL} Theorem 2.3.6.: $\pol(\varphi)=\sum_{j,k\geq 0}\phi^{(j,k)}_{j+k+1}\otimes \vol_{\Tc/\Mc}^*$.
The additional $\vol_{\Tc/\Mc}^*$ and the missing $\vol_\Mc$ are due to the fact that differential forms resolve the constant sheaf $\C$, whereas the currents resolve the orientation bundle of $\Tc_W$, and due to the definition of the Fourier series of our currents, see \cref{eval_distr}.
Now we may specialize our polylogarithm to get the representation of $\Eis^k(\varphi)$ as
\begin{equation*}
\lim{s\rightarrow 2\xi+k}\int_{u\in\R_{>0}}\sum_{l\in V(\Q)\setminus \left\{0\right\}}\frac{-\phi(\hat{f})(\sqrt{u}\hat g^{-1}l,g)}{k!(2\xi-1)!}\left(\frac{l^{-1,0}}{r}\right)^{\otimes k}u^{s}\frac{du}{u}i_l(d i_l)^{2\xi-1} \vol_{\Tc/\Mc}\otimes\vol_{\Tc/\Mc}^{*n+1}
\end{equation*}
for $k\geq0$ and the limit exists by analytic continuation.

\section{The decomposition isomorphism as Mellin transform}

Next we are going to apply \Cref{trace_integration} to $\Eis^k(\varphi)$. The fiber integral of $\Eis^k(\varphi)$ turns out to be a Mellin transform of a theta series. Given this description of $\Eis^k(\varphi)$ and having \cite{Wi} in mind the connection to Eisenstein series and Harder's Eisenstein classes is immediately apparent. But before we can apply the fiber integral, we need to express $\Eis^k(\varphi)$ by $\mathfrak H^\bullet _\C$ which is generated by $\frac{dr_i}{r_i}$, $i=1,...,\xi$.

\begin{lemma}\label{calculating_vol}
\begin{equation*}
\frac{i_l(d i_l)^{2\xi-1}\vol_V}{(2\xi-1)!}=
\Norm\left(\frac{l_1+\tau l_2}{r(\overline{\tau}-\tau)}\right)\sum_{i=1} ^\xi\prod_{i\neq j}\left(d\left(\frac{l_1+\tau l_2}{r(\overline{\tau}-\tau)}\right)_j
\wedge d\overline{\tau}_j\right)\wedge\frac{l_1+\tau l_2}{r(\overline{\tau}-\tau)}_i d\overline{\tau}_i
\end{equation*}
\begin{proof}
We follow the computation of \cite{Bl2} 4.7.
Set $dw^1+\tau dw^2=:\eta^1$, $\frac{\overline\eta^1}{\overline\tau-\tau}=:\eta^2$. We have
\begin{equation*}
\eta^1(u)=1,\ \overline\eta^1(u)=0,\ d\eta^1=d\tau\wedge dw^2=d\tau\wedge\frac{\overline \eta^1-\eta^1}{\overline\tau-\tau},\ \vol_V=\prod_{i=1}^\xi\eta^1_i\wedge\eta^2_i.
\end{equation*}
\begin{equation*}
i_l(d i_l)^{2\xi-1}\vol_V=i_l\Lc_{\frac{l^{-1,0}}{r}}^{2\xi-1}\vol_V
\end{equation*}
and we calculate
\begin{equation*}
\Lc_{\frac{l^{-1,0}}{r}}^{2\xi-1}\vol_V=\Lc_{\frac{l^{-1,0}}{r}}^{2\xi-1}\prod_{i=1}^\xi\eta^1_i\wedge\eta^2_i=
\end{equation*}
\begin{equation*}
(2\xi-1)!\sum_{i=1}^\xi\prod_{j\neq i}\left(\Lc_{\frac{l^{-1,0}}{r}}\eta^1_j\wedge\Lc_{\frac{l^{-1,0}}{r}}\eta^2_j\right)\wedge\left(\eta^1_i\wedge\Lc_{\frac{l^{-1,0}}{r}}\eta^2_i+\Lc_{\frac{l^{-1,0}}{r}}\eta^1_i\wedge \eta^2 _i\right),
\end{equation*}
as $\Lc_{\frac{l^{-1,0}}{r}}^2\eta^1_i=\Lc_{\frac{l^{-1,0}}{r}}^2\eta^2_i=0$ for $i=1,...,\xi$ and $\Lc_{\frac{l^{-1,0}}{r}}$ satisfies the Leibniz rule. Moreover, 
\begin{equation*}
\Lc_{\frac{l^{-1,0}}{r}}\eta^1=(di_l+i_l d)(\eta^1)=d(\frac{l_1+\tau l_2}{r})+\frac{l_1+\tau l_2}{r(\overline{\tau}-\tau)}d\tau,
\end{equation*}
\begin{equation*}
\Lc_{\frac{l^{-1,0}}{r}}\eta^2=(di_l+i_l d)(\eta^2)=i_ld(\eta^2)=\frac{l_1+\tau l_2}{r(\overline{\tau}-\tau)}\frac{d\overline\tau}{\overline{\tau}-\tau}
\end{equation*}
and $i_l\Lc_{\frac{l^{-1,0}}{r}}\eta^1_i=i_l\Lc_{\frac{l^{-1,0}}{r}}\eta^2_i=i_l\eta^2_i=0$, $i=1,...,\xi$. Now we get
\begin{equation*}
\frac{i_l(d i_l)^{2\xi-1}\vol_V}{(2\xi-1)!}=\sum_{i=1}^\xi\prod_{j\neq i}\left(\Lc_{\frac{l^{-1,0}}{r}}\eta^1_j\wedge\Lc_{\frac{l^{-1,0}}{r}}\eta^2_j\right)\wedge\left(i_l\eta^1_i\wedge\Lc_{\frac{l^{-1,0}}{r}}\eta^2_i\right)=
\end{equation*}
\begin{equation*}
\sum_{i=1}^\xi\prod_{j\neq i}\left(\left(d(\frac{l_1+\tau l_2}{r})+\frac{l_1+\tau l_2}{r(\overline{\tau}-\tau)}d\tau\right)\wedge\left(\frac{l_1+\tau l_2}{r(\overline{\tau}-\tau)}\frac{d\overline\tau}{\overline{\tau}-\tau}\right)\right)_j\wedge\left(\frac{l_1+\tau l_2}{r(\overline{\tau}-\tau)}\right)_i^2d\overline\tau_i=
\end{equation*}
\begin{equation*}
\sum_{i=1}^\xi\prod_{j\neq i}\left(\left(\frac{d(\frac{l_1+\tau l_2}{r})}{\overline{\tau}-\tau}+\frac{l_1+\tau l_2}{r(\overline{\tau}-\tau)}\frac{d\tau}{\overline{\tau}-\tau}\right)\wedge\left(\frac{l_1+\tau l_2}{r(\overline{\tau}-\tau)}d\overline\tau\right)\right)_j\wedge\left(\frac{l_1+\tau l_2}{r(\overline{\tau}-\tau)}\right)_i^2d\overline\tau_i=
\end{equation*}
\begin{equation*}
\sum_{i=1}^\xi\prod_{j\neq i}\left(d\left(\frac{l_1+\tau l_2}{r(\overline{\tau}-\tau)}\right)\wedge\left(\frac{l_1+\tau l_2}{r(\overline{\tau}-\tau)}d\overline\tau\right)\right)_j\wedge\left(\frac{l_1+\tau l_2}{r(\overline{\tau}-\tau)}\right)_i^2d\overline\tau_i
\end{equation*}
\end{proof}
\end{lemma}

\begin{lemma}\label{calc_vol}
\begin{equation*}
\Norm\left(\frac{l_1+\tau l_2}{r(\overline{\tau}-\tau)}\right)^{-1}\frac{i_l(d i_l)^{2\xi-1}\vol_V}{(2\xi-1)!}=
\end{equation*}
\begin{equation*}
\sum_{i=1} ^\xi\sum_{I\subset \left\langle \xi\right\rangle\setminus i}(-1)^{|I|}\frac{l_1+\overline\tau l_2}{r(\overline{\tau}-\tau)}_{(i\cup I)^c}\frac{l_1+\tau l_2}{r(\overline{\tau}-\tau)}_{I\cup i}
\frac{d\tau\wedge d\overline\tau}{\overline\tau-\tau}_{(i\cup I)^c}\wedge \frac{dr\wedge d \overline\tau}{r}_I\wedge d\overline\tau_i
\end{equation*}
\begin{proof}
We have
\begin{equation*}
d\left(\frac{l_1+\tau l_2}{r(\overline{\tau}-\tau)}\right)
\wedge d\overline{\tau}=\left(-\frac{l_1+\tau l_2}{r(\overline{\tau}-\tau)}\frac{dr}{r}+d\left(\frac{l_1+\tau l_2}{\overline{\tau}-\tau}\right)r^{-1}\right)\wedge d\overline{\tau}
\end{equation*}
and
\begin{equation*}
d\left(\frac{l_1+\tau l_2}{\overline{\tau}-\tau}\right)\wedge d\overline\tau=d\left(\frac{l_1+\overline\tau l_2 -\overline\tau l_2+\tau l_2}{\overline{\tau}-\tau}\right)\wedge d\overline\tau=
\end{equation*}
\begin{equation*}
d\left(\frac{l_1+\overline\tau l_2}{\overline{\tau}-\tau}-l_2\right)\wedge d\overline\tau=d\left(\frac{l_1+\overline\tau l_2}{\overline{\tau}-\tau}\right)\wedge d\overline\tau=\frac{l_1+\overline\tau l_2}{(\overline{\tau}-\tau)}\frac{d\tau\wedge d\overline\tau}{\overline{\tau}-\tau}
\end{equation*}
If we plug these formulas in \Cref{calculating_vol}, we get the desired formula.
\end{proof}
\end{lemma}
\begin{corollary}\label{dr_invariant}
The coefficients in front of the invariant forms $\frac{dr}{r}_I$, $I\subset\left\langle \xi\right\rangle$, of \\
$\Norm(r)^2i_l(di_l)^{2\xi-1}\vol_V$ do not depend on $r$.
\end{corollary} 
\begin{lemma}\label{Z_proj}
We set $\omega:=(\overline\tau-\tau)u$ and $\Norm(\omega):=\prod_{i=1}^\xi\omega_i$. The projection of $\frac{1}{k!}\left(\frac{l^{-1,0}}{r}\right)^{\otimes k}$ onto the $\Sym^k V(\C)^{Z_K}$-part is 
\begin{equation*}
\Norm\left(\frac{l_1+\tau l_2}{r(\overline\tau-\tau)}\right)^m\frac{\Norm(\omega)^m}{(m!)^\xi},\ \text{if}\ k=\xi m, 
\end{equation*}
and zero otherwise.
\begin{proof}
We set $x^{[k]}:=\frac{x^k}{k!}$ for $x$ in any $\Q$-Algebra $R$. We have $(x+y)^{[k]}=\sum_{l=0}^kx^{[l]}y^{[k-l]}$. Doing induction on $n\in \N$ we may easily proof the formula
\begin{equation*}
\left(\sum_{i=1}^nx_i\right)^{[k]}=\sum_{(k_i)\in\N_0 ^n:\sum k_i=k}\prod_{i=1}^nx^{[k_i]}
\end{equation*}
Write $\frac{l^{-1,0}}{r}=\sum_{i=1}^\xi\frac{l^{-1,0}}{r}_i$ and get
\begin{equation*}
\frac{1}{k!}\left(\frac{l^{-1,0}}{r}\right)^{\otimes k}=\left(\frac{l^{-1,0}}{r}\right)^{\otimes [k]}=\left(\sum_{i=1}^\xi\frac{l^{-1,0}}{r}_i\right)^{\otimes [k]}=
\sum_{(k_i)\in \N_0^{\xi}:\sum k_i=k} \prod_{i=1}^\xi\left(\frac{l^{-1,0}}{r}\right)_i ^{\otimes [k_i]}=
\end{equation*}
\begin{equation*}
\sum_{(k_i)\in \N_0^{\xi}:\sum k_i=k} \prod_{i=1}^\xi \frac{1}{k_i !}\left(\frac{l^{-1,0}}{r}\right)_i ^{\otimes k_i}=\sum_{(k_i)\in \N_0^{\xi}:\sum k_i=k} \prod_{i=1}^\xi \frac{1}{k_i !}\left(\frac{l_1+\tau l_2}{r(\overline\tau -\tau)}\omega\right)_i ^{\otimes k_i}
\end{equation*}
The projection of this element onto the $\Sym^k V(\C)^{Z_K}$-part are those summands where the $Z_K$ action factors through the norm character. Thus the only summand living in $\Sym^k V(\C)^{Z_K}$ is the one corresponding to $k_i =m$ for $i=1,...,\xi$, see \Cref{invariants}.  
\end{proof}
\end{lemma}
Before we can go on, we need to fix measures again. Let us consider the ideles $\I_F$. If $\nu\nmid\infty$, we take $d^{\times}x_\nu:=\frac{|\kappa(\mathfrak p_\nu)|}{|\kappa(\mathfrak p_\nu)|-1}\frac{dx_\nu}{|x_\nu|_\nu}$, as Haar measure of the multiplicative group $F_\nu ^{\times}$, where $dx_\nu$ is the Haar measure on the additive group $F_\nu$, which we already have fixed. If $\nu|\infty$, we have $F_\nu=\R$ and we take $d^{\times}x_\nu:=\frac{dx_\nu}{|x_\nu|}$, where on the right-hand side $dx_\nu$ is the usual Lebesgue measure on $\R$. These local measures induce global measures $d^{\times}x$ on $\I_F$ and $\I_{F,f}$ and $F_\R ^{\times}$ as in \Cref{Fourier}.
\begin{remark}
\begin{itemize}
\item Let $G$ be a locally compact group and $H\subset G$ a normal and closed subgroup. Suppose we have fixed (left) Haar measures $dg$ on $G$ and $dh$ on $H$. Then there is a uniquely determined (left) Haar measure $dgH$ on $G/H$ such that
\begin{equation*}
\int_Gf(g)dg=\int_{G/H}\left(\int_{H}f(gh)dh\right)dgH
\end{equation*}
holds for all integrable functions $f$ on $G$.
\item If $H\subset G$ is an open subgroup and we have fixed a Haar measure $dg$ on $G$, we always take on $H$ the restriction of $dg$ as Haar measure on $H$. 
\item 
We have $d^{\times}x(\R_{>0}\backslash F_\R ^{\times,0}/\det(Z_K))=R_K$ (see \cref{int_decomp}), where the measure $d^{\times}x$ is induced by the measure $\frac{dt}{t}$ on $\R_{>0}$, $d^{\times}x$ on $F_\R ^{\times,0}$ and the counting measure on $\det(Z_K)$.
\item In the following calculations the choice of the Haar measures will be clear from the context.
\end{itemize} 
\end{remark}
Set
\begin{equation}\label{Theta_1}
\Theta_l:=\Norm\left(\frac{l_1+\tau l_2}{\overline\tau-\tau}\right)^m \Norm(\omega)^m\Norm(r)^2i_l(di_l)^{2\xi-1}\vol_{\Tc/\Mc}\otimes \vol_{\Tc/\Mc}^{*n+1}  
\end{equation}
a form constant in $r$ and $C(k):=\frac{-1}{(m!)^\xi(2\xi-1)!\xi}$.

\begin{lemma}\label{decomp_Eis}
We have 
\begin{equation*}
\Eis^k(\varphi)=
\lim{s\rightarrow 0}\frac{C(k)2^\xi}{R_K}\int_{t\in\left\{\pm1\right\}^\xi\backslash Z(\R)/Z_K}\sum_{l\in V(\Q)\setminus \left\{0\right\}}\frac{\phi(\hat{f})(t\hat g ^{-1}l,g) |\Norm(t)|^{2(m+2+s)}}{\Norm(\det(g))^{m+2}} d^{\times}t\cdot \Theta_l
\end{equation*}
for $k=\xi m\geq0$. 
\begin{proof}
From \cref{Z_proj1}, \cref{Z_proj} and \cref{Theta_1} we already know that we may identify the cohomology class $\Eis^k(\varphi)$ with
\begin{equation*}
\lim{s\rightarrow 2\xi+k}\xi C(k)\int_{u\in \R_{>0}}\sum_{l\in V(\Q)\setminus \left\{0\right\}}\frac{\phi(\hat{f})\left(\sqrt{u}\hat g^{-1}l,g\right) u^{s}}{\Norm(r)^{m+2}}\frac{du}{u}\cdot\Theta_l
\end{equation*}
We write everything in base and fiber coordinates and have to work with $\tilde{r}\in F_\R ^1$: 
\begin{equation*}
\begin{pmatrix}1 & x\\ 0 & y\end{pmatrix}=:(1,\tau)\ \text{and } \sqrt{u}\hat g^{-1}=\frac{\sqrt{u}\hat g_f ^{-1}(1,\tau)}{\sqrt{\sqrt[\xi]{|N(r)|}|\tilde{r}\text{Im}(\tau)|}}.
\end{equation*}
and we get for $\Eis^k(\varphi)$
\begin{equation*}
\lim{s\rightarrow 2\xi+k}\xi C(k)\int_{u\in \R_{>0}}\sum_{l\in V(\Q)\setminus \left\{0\right\}}\frac{\phi(\hat{f})\left(\frac{\sqrt{u}\hat g_f ^{-1}(1,\tau)}{\sqrt{\tilde{r}\text{Im}(\tau)}}l  ,g\right) u^{s}}{\sgn\Norm(\det(g))^{m+2}}\frac{du}{u}\cdot\Theta_l
\end{equation*}
Note that $\Theta_l$ is just built up by the invariant forms $\frac{dr}{r}_i$, $i=1,...,\xi$, and is constant in $r$, therefore we may apply \cref{int_decomp} to identify $\Eis^k(\varphi)$ with
\begin{equation*}
\lim{s\rightarrow 2\xi+k}\frac{\xi C(k)}{R_K}\int_{\tilde{r}\in \left\{\pm1\right\}^\xi\backslash F_\R ^1/\det(Z_K)}\int_{u\in \R_{>0}}\sum_{l\in V(\Q)\setminus \left\{0\right\}}\frac{\phi(\hat{f})\left(\frac{\sqrt{u}\hat g_f ^{-1}(1,\tau)}{\sqrt{\tilde{r}\text{Im}(\tau)}}l,g\right) u^{s}}{ \sgn \Norm((\det(g)))^{m+2}}\frac{du}{u}d^\times\tilde{r}\cdot\Theta_l
\end{equation*}
Here we already have interchanged the limit and $\int_{\tilde{r}\in \left\{\pm1\right\}^\xi\backslash F_\R ^1/\det(Z_K)}$. This is possible by the dominated convergence theorem, as the limit in $s$ gives a continuous function in $r$, which is integrable over the compact space $\left\{\pm1\right\}^\xi\backslash F_\R ^1/\det(Z_K)$. We continue our calculation by changing variables $\tilde{r}\mapsto \tilde{r}^{-1}$ and changing the order of integration with Fubini to get for $\Eis^k(\varphi)$
\begin{equation*}
\lim{s\rightarrow 2\xi+k}\frac{\xi C(k)}{R_K}\int_{u\in \R_{>0}}\int_{\tilde{r}\in \left\{\pm1\right\}^\xi\backslash F_\R ^1/\det(Z_K)}\sum_{l\in V(\Q)\setminus \left\{0\right\}}\frac{\phi(\hat{f})\left(\frac{\sqrt{|u\tilde{r}|}\hat g_f ^{-1}(1,\tau)}{\sqrt{|\text{Im}(\tau)|}}l,g\right) u^{s}}{ \sgn \Norm((\det(g)))^{m+2}}d^\times\tilde{r}\frac{du}{u}\cdot\Theta_l=
\end{equation*}
\begin{equation*}
\lim{s\rightarrow 2\xi+k}\frac{C(k)}{R_K}\int_{v\in \R_{>0}}\int_{\tilde{r}\in \left\{\pm1\right\}^\xi\backslash F_\R ^1/\det(Z_K)}\sum_{l\in V(\Q)\setminus \left\{0\right\}}\frac{\phi(\hat{f})\left(\frac{\sqrt{|\sqrt[\xi]{v}\tilde{r}|}\hat g_f ^{-1}(1,\tau)}{\sqrt{|\text{Im}(\tau)|}}l,g\right) v^{\frac{s}{\xi}}}{ \sgn \Norm((\det(g)))^{m+2}}d^\times\tilde{r}\frac{dv}{v}\cdot\Theta_l
\end{equation*}
with $v=u^\xi$ and $\frac{dv}{v}=\xi\frac{du}{u}$. We have $d^\times r=d^\times\tilde{r}\frac{dv}{v}$ using the isomorphism (the trivialization of the fiber bundle $F_\R ^\times$)
\begin{equation*}
F_\R ^1\times \R_{>0}\to F_\R ^\times,\ (\tilde{r},v)\mapsto \sqrt[\xi]{v}\tilde{r}
\end{equation*}
as remarked in \cref{int_decomp} and therefore $\Eis^k(\varphi)$ equals
\begin{equation*}
\lim{s\rightarrow 2\xi+k}\frac{C(k)}{R_K}\int_{r\in \left\{\pm1\right\}^\xi\backslash F_\R ^\times/\det(Z_K)}\sum_{l\in V(\Q)\setminus \left\{0\right\}}\frac{\phi(\hat{f})\left(\frac{\sqrt{|r|}\hat g_f ^{-1}(1,\tau)}{\sqrt{|\text{Im}(\tau)|}}l,g\right) |\Norm(r)|^{\frac{s}{\xi}}}{ \sgn \Norm((\det(g)))^{m+2}}d^\times r\cdot\Theta_l
\end{equation*}
Using the isomorphism
\begin{equation*}
\varphi:\left\{\pm1\right\}^\xi\backslash Z(\R)/Z_K\rightarrow \left\{\pm1\right\}^\xi\backslash F_\R ^{\times}/\det(Z_K),\ t\mapsto t^2=r,\ \varphi^*d^{\times}r=2^{\xi} d^{\times}t.
\end{equation*}
we get for $\Eis^k(\varphi)$
\begin{equation*}
\lim{s\rightarrow 2\xi+k}\frac{2^{\xi}C(k)}{R_K}\int_{t\in \left\{\pm1\right\}^\xi\backslash Z(\R)/Z_K}\sum_{l\in V(\Q)\setminus \left\{0\right\}}\frac{\phi(\hat{f})\left(\frac{t\hat g_f ^{-1}(1,\tau)}{\sqrt{|\text{Im}(\tau)|}}l,g\right) |\Norm(t)|^{\frac{2s}{\xi}}}{ \sgn \Norm((\det(g)))^{m+2}}d^\times t\cdot\Theta_l=
\end{equation*}
\begin{equation*}
=\lim{s\rightarrow 0}\frac{2^{\xi}C(k)}{R_K}\int_{t\in \left\{\pm1\right\}^\xi\backslash Z(\R)/Z_K}\sum_{l\in V(\Q)\setminus \left\{0\right\}}\frac{\phi(\hat{f})\left(t\hat g ^{-1}l,g\right) |\Norm(t)|^{2(m+2+s)}}{ \Norm((\det(g)))^{m+2}}d^\times t\cdot\Theta_l
\end{equation*}
\end{proof}
\end{lemma}
\begin{remark}
From now on we assume $K_f$ is so small that $Z_K\subset Z(\R)^0$. Then we may write
\begin{equation*}
\Eis^k(\varphi)= \frac{C(k)2^\xi}{R_K}\lim{s\rightarrow 0}\int_{t\in Z(\R)^0/Z_K}\sum_{l\in V(\Q)\setminus \left\{0\right\}}\frac{\phi(\hat{f})(t\hat g ^{-1}l,g)\Norm(t)^{2(m+2+s)} d^{\times}t}{\Norm(\det(g_\infty))^{m+2}}\cdot\Theta_l.
\end{equation*}
\end{remark}
The next step is to interpret this integral as a global Tate integral. To do this we decompose $\Eis^k(\varphi)$ into several summands.\\
There is a short exact sequence
\begin{equation*}
1\rightarrow Z(\R)^0/Z_K\rightarrow K_f ^Z\backslash Z(\A)/Z(\Q)\rightarrow Cl_F ^K\rightarrow 1,
\end{equation*}
with $K ^Z:=K\cap Z(\A)$ and $Cl_F ^K=K_f ^Z Z(\R)^0\backslash Z(\A)/Z(\Q)$ (we always suppose $K_f ^Z=\det(K_f)$). Denote by $\widehat{Cl_F ^K}$ the Pontryagin dual of the last group and let $h_{K}$ be its cardinality.
Moreover, we define
\begin{equation*}
\widehat{Cl_F ^K}(m):=\left\{\chi\in \widehat{Cl_F ^K}:\chi(t)=\sgn(\Norm(t))^m,\ t\in Z(\R) \right\}.
\end{equation*}
\begin{remark}\label{Iwasawa}
Recall the Iwasawa decomposition on $G_0(\R)=G_0(F_\nu)$ for $\nu|\infty$:
\begin{equation*}
g=\begin{pmatrix}a&b\\c&d\end{pmatrix}=k(g)a(g)n(g):
\end{equation*}
\begin{equation*}
k(g)=\frac{1}{\sqrt{a^2+c^2}}\begin{pmatrix}a&-c\\c&a\end{pmatrix}=\begin{pmatrix}\cos(\theta)&-\sin(\theta)\\ \sin(\theta)&\cos(\theta)\end{pmatrix}\in SO(2),
\end{equation*}
\begin{equation*}
a(g)=\begin{pmatrix}\sqrt{a^2+c^2}&0\\0&\frac{\det(g)}{\sqrt{a^2+c^2}}\end{pmatrix}=\begin{pmatrix}t_1&0\\0& t_2\end{pmatrix}\in T_0(F_\nu),
\end{equation*}
\begin{equation*}
n(g)=\begin{pmatrix}1&\frac{ab+dc}{a^2+c^2}\\0&1\end{pmatrix}=\begin{pmatrix}1&x\\0& 1\end{pmatrix}\in U_0(F_\nu).
\end{equation*}
We can consider $k,a,n$ as continuous functions $G(F_\nu)\rightarrow \C$. Moreover, we set $b:=an$. 

One also has a Iwasawa decomposition for $g\in G_0(F_\nu)$ and $\nu$ a finite place: 
\begin{equation*}
g=kb \text{ with } k\in G_0(\Oc_\nu) \text{ and } b=\begin{pmatrix}t_1&t_1x\\0&t_2\end{pmatrix}\in B_0(F_\nu).
\end{equation*} 
This decomposition is not unique, as $b$ is just well-defined up to elements from $B_0(\Oc_\nu)$. Nevertheless, the Iwasawa decomposition often suffices to define functions on $G_0(F_\nu)$. For example, if $\chi:F_\nu^{\times}\to\C^{\times}$ is an unramified quasi-character, say $\chi(x)=|x|_\nu ^s$ for $s\in \C$, we have well-defined continuous functions $G_0(F_\nu)\to \C^{\times}$, $g\mapsto \chi(t_1)$ or $\chi(t_2)$. Finally, we get the Iwasawa decomposition for $g\in G(\A)$. Of course, it is defined place by place and is given by
\begin{equation*}
\begin{pmatrix}a&b\\c&d\end{pmatrix}=g=kb \text{ with } k\in G(\hat\Z)\cdot\prod_{\nu|\infty}SO(2) \text{ and } b=\begin{pmatrix}t_1&t_1x\\0&t_2\end{pmatrix}\in B(\A).
\end{equation*}
Again, this decomposition is just well-defined up to elements from $B(\hat\Z)$. For the infinite places $\nu|\infty$ we also have the functions $\theta_\nu:G_0(F_\nu)\to\R/2\pi\Z$ defined by $e^{i\theta}_\nu:=e^{i\theta_\nu}=\frac{(a_\nu+ic_\nu)}{\sqrt{a_\nu^2+c_\nu^2}}$. 
\end{remark}
For $i\in \left\langle \xi\right\rangle$ and $I\subset \left\langle \xi\right\rangle\setminus i$ we define the form
\begin{equation*}
\Theta_{I,i}:=\frac{(-1)^{|I|}(2\xi-1)! \Norm(e^{-(m+1)i\theta})e^{-i\theta}_{I\cup i}e^{i\theta}_{(I\cup i)^c}}{(-2i)^{\xi(m+2)} \left\| \det g\right\|_f \sgn(\Norm(\det g))^{-1}{\Norm(t_{2}) ^{m+2}}}\cdot
\end{equation*}
\begin{equation}\label{Theta_2}
\Norm(\omega)^m\frac{d\tau\wedge d\overline\tau}{\overline\tau-\tau}_{\left\langle \xi\right\rangle\setminus (i\cup I)}\wedge \frac{dr\wedge d \overline\tau}{r}_I\wedge d\overline\tau_i\otimes \vol_{\Tc/\Mc}^{*n+1}
\end{equation}
We have Schwartz-Bruhat-functions 
\begin{equation*}
\phi(\hat{f})^m_{I\cup i}(v,g):=\phi(\hat{f})(v,g)\Norm(v_\infty ^1+iv_\infty ^2)^{m+1}(v_\infty ^1+iv_\infty ^2)_{I\cup i}(v_\infty ^1-iv_\infty ^2)_{(I\cup i)^c}
\end{equation*}
on $V(\A)$. Write $\phi(\hat{f})^m_{I\cup i}=\phi(\hat{f})^m_{I\cup i,\infty}\otimes \phi(\hat{f})^m_{I\cup i,f}$ with $\phi(\hat{f})^m_{I\cup i,f}=\hat f$.
For $\chi\in \widehat{Cl_F ^K}$ and $s\in\C$ we may define differential forms
\begin{equation*}
\Eis^k _{I,i}(\varphi,\chi,s)=
\sum_{\gamma\in G(\Q)/B(\Q))}\int_{t\in Z(\A)}
\phi(\hat{f})^m _{I\cup i}(tg\gamma e^1,g) \chi(\det(g)t)\left\|\det(g)t\right\|^{s} d^{\times}t\cdot \Theta_{I,i}.
\end{equation*}
\begin{lemma}
If $Re(s)>2$, $\Eis^k _{I,i}(\varphi,\chi,s)$ is an honest differential form and a holomorphic function in $s$. $\Eis^k _{I,i}(\varphi,\chi,s)$ has an analytic continuation as meromorphic function in $s$ to the whole complex plane and does not have a pole at $s=2$. In particular, the limit
\begin{equation*}
\lim {s\rightarrow 0}\Eis^k _{I,i}(\varphi,\chi,m+2+2s)=:\Eis^k _{I,i}(\varphi,\chi)
\end{equation*}
exists and gives a well-defined differential form on $\Mc_K$.
\begin{proof}
The first statement is \cite{Wi} III Proposition 6. Wielonsky considers functions in the variable $\det(t)$, whereas we have the variable $t$. Therefore Wielonsky has the statement for $\sigma>1$ and we have our statement for $Re(s)>2$. The second statement is
\cite{Wi} III Proposition 9, as we have $\widehat{\phi(\hat{f})^m _{I\cup i}}(0)=0$.
\end{proof}
\end{lemma}
Let us set $\kappa_F:=d^{\times}t(\R_{>0}\backslash Z(\A)/Z(\Q))$.
\begin{theorem}\label{Wi-series}
Suppose we are given a function in
$\varphi\in\Sc(V(\A_f),\mu^{\otimes n})^0$. Write
\begin{equation*} 
(v,g)\mapsto \left(\left\|\det(g)\right\|_f \sgn(\Norm(\det(g)))^{-1}\right)^n f(v,g)=\varphi(v,g),
\end{equation*}
with $f\in \Sc(V(\A_f),\mu^{\otimes 0})^0$ such that $f$ factors in the second argument over $K_f$ with $Z_K\subset Z(\R)^0$. The polylogarithmic Eisenstein class associated to $\varphi$
\begin{equation*}
\Eis^k(\varphi)\in H^{2\xi-1}(\Mc,\Sym^k\Hc\otimes\mu^{\otimes n+1})=\bigoplus_{p+q=2\xi-1}H^p(\Sc,\Sym^k\Hc^\prime\otimes\mu^{\otimes n+1})\otimes \mathfrak H^q
\end{equation*}
may be represented by the differential form
\begin{equation*}
\frac{ 2C(k)}{\kappa_F}\sum_{i=1}^\xi\sum_{I\subset\left\langle \xi\right\rangle\setminus i}\sum_{\chi\in \widehat{Cl_F ^K}(m)}\Eis^k_{I,i}(\varphi,\chi)
\end{equation*}
\begin{proof}
Define $\sgn(r):=(\sgn(r)_\nu)_{\nu| \infty}:=(\sgn(r_\nu))_{\nu|\infty}$, $r\in F_\R ^{\times}$ and $S:V(\R)\to F_\C$, $v=(v^1,v^2)\mapsto v^1+iv^2$. For $g\in G(\R)$ we have $S(k(g)v)=e^{i\theta(g)}S(v)$ and
\begin{equation*}
\sgn(\det(g))k(\hat g^{-1})=\widehat{k(g)}^{-1},\ b(\hat g^{-1})=\sgn(\det(g))\widehat{b(g)}^{-1}=\frac{1}{\sqrt{ry}}\begin{pmatrix}1& x\\ 0& y\end{pmatrix}.
\end{equation*}
Therefore
\begin{equation*}
\Norm(e^{-(m+1)i\theta})e^{-i\theta}_{I\cup i}e^{i\theta}_{(I\cup i)^c}\frac{\sgn(\Norm(r))^{m+2}}{\Norm(\sqrt{ry})^{m+2}}\phi(\hat{f})^m_{I\cup i,\infty}(\hat{g}^{-1}v,g)=
\end{equation*}
\begin{equation*}
\phi(\hat{f})(\hat{g}^{-1}v,g)\Norm\left(\frac{v_\infty ^1+\tau v_\infty ^2}{ry}\right)^{m+1}\left(\frac{v_\infty ^1+\tau v_\infty ^2}{ry}\right)_{I\cup i}\left(\frac{v_\infty ^1+\overline\tau v_\infty ^2}{ry}\right)_{(I\cup i)^c}.
\end{equation*}
We see $\sgn(r)\sqrt{ry}=\sgn(r)|t_2|=t_2$ and with \Cref{calc_vol}, the definition of $\Theta_l$ (\Cref{Theta_1}) and $\Theta_{I,i}$ (\Cref{Theta_2}) we conclude for $t\in Z(\R)^0$
\begin{equation*}
\sum_{i=1}^\xi\sum_{I\subset\left\langle \xi\right\rangle\setminus i}\phi(\hat{f})^m _{I\cup i}(t\hat{g}^{-1}l,g) |\Norm(t)|^{m+2+2s}  \Theta_{I,i}=
\phi(\hat{f})(t\hat g ^{-1}l,g)\Norm(t)^{2(m+2+s)}\Norm(r)^{-(m+2)}\Theta_l.
\end{equation*}
We express the characteristic function of $Z(\R)^0/Z_K\subset {K_f ^Z}\backslash Z(\A)/Z(\Q)$ by $h_K^{-1}\sum_{\chi\in \widehat{Cl_F ^K}}\chi$ by considering a character $\chi$ as a function on $Z(\A)$. We get 
\begin{equation*}
\frac{ d^{\times}t(K_f ^Z)R_K h_K}{2^{\xi}C(k)}\Eis^k(\varphi)=
\end{equation*}
\begin{equation*}
\sum_{i=1}^\xi\sum_{I\subset\left\langle \xi\right\rangle\setminus i}\sum_{\chi\in \widehat{Cl_F ^K}}
\lim{s\rightarrow0}\int_{t\in Z(\A)/Z(\Q)}\sum_{l\in V(\Q)\setminus \left\{0\right\}}\phi(\hat{f})^m_{I\cup i}(t\hat g ^{-1}l,g) \chi(t)\left\|t\right\|^{m+2+2s} d^{\times}t\cdot \Theta_{I,i}.  
\end{equation*}
Using short exact sequences one calculates $2^{-(\xi-1)}d^{\times}t(K_f ^Z)R_K h_K=\kappa_F$. For $Re(s)>0$ the inner sum converges absolutely and locally uniformly. Therefore we may rearrange summation $\sum_{l\in V(\Q)\setminus \left\{0\right\}}=\sum_{l\in V(\Q)\setminus \left\{0\right\}/Z(\Q)}\sum_{z\in Z(\Q)}$, use the bijection $G(\Q)/B(\Q)\rightarrow V(\Q)\setminus \left\{0\right\}/Z(\Q)$, $\gamma\mapsto \gamma e^1$, and interchange $\int_{Z(\A)/\Z(\Q)}$ and $\sum_{\gamma\in G(\Q)/B(\Q)}$. We get in total
\begin{equation*}
\lim{s\rightarrow0}\int_{t\in Z(\A)/Z(\Q)}\sum_{l\in V(\Q)\setminus \left\{0\right\}}\phi(\hat{f})^m_{I\cup i}(t\hat g ^{-1}l,g) \chi(t)\left\|t\right\|^{m+2+2s} d^{\times}t\cdot \Theta_{I,i}  =
\end{equation*}
\begin{equation*}
\lim{s\rightarrow0}\sum_{\gamma\in G(\Q)/B(\Q)}\int_{t\in Z(\A)}\phi(\hat{f})^m_{I\cup i}(t\hat g ^{-1}\gamma e^1,g) \chi(t)\left\|t\right\|^{m+2+2s} d^{\times}t\cdot \Theta_{I,i}=  
\end{equation*}
\begin{equation*}
\lim{s\rightarrow0}\sum_{\gamma\in G(\Q)/B(\Q)}\int_{t\in Z(\A)}\phi(\hat{f})^m_{I\cup i}(t g \gamma e^1,g) \chi(\det(g)t)\left\|\det(g)t\right\|^{m+2+2s} d^{\times}t\cdot \Theta_{I,i}.
\end{equation*}
For this see also \cite{Wi} III Proposition 7.
So we already have
\begin{equation*}
\Eis^k(\varphi)=\frac{2 C(k)}{\kappa_F}\sum_{i=1}^\xi\sum_{I\subset\left\langle \xi\right\rangle\setminus i}\sum_{\chi\in \widehat{Cl_F ^K}}\Eis^k_{I,i}(\varphi,\chi).
\end{equation*}
To prove the theorem it suffices to show that $\Eis^k_{I,i}(\varphi,\chi)=0$, whenever $\chi\notin \widehat{Cl_F ^K}(m)$. So take a $\chi\notin \widehat{Cl_F ^K}(m)$ and an $\epsilon\in \left\{\pm 1\right\}^\xi\subset Z(\R)$ with $\chi(\epsilon)\sgn(\Norm(\epsilon))^m\neq 1$. Then
\begin{equation*}
\int_{t\in Z(\A)}
\phi(\hat{f})^m _{I\cup i}(t g \gamma e^1,g) \chi(\det(g)t)\left\|\det(g)t\right\|^{s} d^{\times}t=
\end{equation*}
\begin{equation*}
\int_{t\in Z(\A)}
\phi(\hat{f})^m _{I\cup i}(\epsilon t g \gamma e^1,g) \chi(\epsilon \det(g) t)\left\|\epsilon \det(g) t\right\|^{s} d^{\times}t=
\end{equation*}
\begin{equation*}
\chi(\epsilon)\sgn(\Norm(\epsilon))^m\int_{t\in Z(\A)}\phi(\hat{f})^m _{I\cup i}(t g \gamma e^1,g) \chi(\det(g)t)\left\|\det(g)t\right\|^{s} d^{\times}t,
\end{equation*}
since 
\begin{equation*}
\phi(\hat{f})^m _{I\cup i}(\epsilon t g \gamma e^1,g)=\sgn(\Norm(\epsilon))^m\phi(\hat{f})^m _{I\cup i}(tg\gamma e^1,g).
\end{equation*}
So we conclude
\begin{equation*}
\int_{t\in Z(\A)}\phi(\hat{f})^m _{I\cup i}(t g \gamma e^1,g) \chi(\det(g)t)\left\|\det(g)t\right\|^{s} d^{\times}t=0
\end{equation*}
and therefore $\Eis^k _{I,i}(\varphi,\chi)=0$.
\end{proof}
\end{theorem}

\section{Cohomology of the boundary}\label{boundary}

As Harder's Eisenstein classes are by definition determined by their restriction to the cohomology of the boundary of the Borel-Serre compactification of $\Sc_K$, we need a thorough understanding of the cohomology of the boundary to compare the polylogarithmic Eisenstein classes with those of Harder. In this section we quickly recall Harder's description of the cohomology of the boundary as a module induced from $B(\A_f)$ to $G(\A_f)$. 

From now on we consider $\Q$-coefficients, but we will emphasize, when things can be defined integrally. We introduce the spaces $\partial \Sc_K:=K\backslash G(\A)/B(\Q)$ and $\partial \Mc_K:=K^1\backslash G(\A)/B(\Q)$. Denote by $\overline{\Sc_K}$ the Borel-Serre compactification of the space $\Sc_K$ . It is a manifold with corners and we have that the boundary $\overline{\Sc_K}\setminus\Sc_K$ is homotopy equivalent to $\partial\Sc_K$, see \cite{Ha1} 2.1. Therefore, we refer to $H^\bullet(\partial \Sc_K,\Sym^k\Hc ^\prime\otimes \mu^{\otimes n})$ as the \textit{cohomology of the boundary}. \\ 
The natural map $j_{\Sc_K}:\partial \Sc_K\rightarrow\Sc_K$ induces by pullback a restriction map on cohomology
\begin{equation*}
\res_{\Sc_K}:H^\bullet(\Sc_K,\Sym^k\Hc^\prime\otimes \mu^{\otimes n})\rightarrow H^\bullet(\partial\Sc_K,\Sym^k\Hc^\prime\otimes \mu^{\otimes n})
\end{equation*}
and performing the colimits over $K$ these maps glue to $G(\A_f)\times \pi_0(G_\infty)$-equivariant maps
$\res_{\Sc}:H^\bullet(\Sc,\Sym^k\Hc^\prime\otimes \mu^{\otimes n})\rightarrow H^\bullet(\partial\Sc,\Sym^k\Hc^\prime\otimes \mu^{\otimes n})$ of the colimits of the groups above. \\
Analogously we define $\res_{\Mc}:H^\bullet(\Mc,\Sym^k\Hc^\prime\otimes \mu^{\otimes n})\rightarrow H^\bullet(\partial\Mc,\Sym^k\Hc^\prime\otimes \mu^{\otimes n})$. The cohomology groups of the boundary are induced $G(\A_f)\times \pi_0(G_\infty)$-modules. More precisely, set $K^B:=K\cap B(\A)$, $K^{1,B}:=K^1\cap B(\A)$ and get the natural inclusions $\Sc_{K}^B:=K ^B\backslash B(\A)/B(\Q)\rightarrow \partial \Sc_K$ and $\Mc_{K}^B:=K^{1,B}\backslash B(\A)/B(\Q)\rightarrow \partial \Mc_K$, which are actually inclusions of several connected components. As usual we set
\begin{equation*}
\varinjlim_K H^\bullet(\Sc_{K}^B,\Sym^k\Hc^\prime\otimes \mu^{\otimes n})=:H^\bullet(\Sc^B,\Sym^k\Hc^\prime\otimes \mu^{\otimes n}).
\end{equation*}
These groups are $B(\A_f)\times \pi_0(B(\R))$-modules and we have
\begin{equation*}
\Ind_{B(\A_f)\times \pi_0(B(\R))} ^{G(\A_f)\times \pi_0(G_\infty)}H^\bullet(\Sc^B,\Sym^k\Hc^\prime \otimes \mu^{\otimes n})=H^\bullet(\partial\Sc,\Sym^k\Hc^\prime\otimes \mu^{\otimes n})
\end{equation*}
as $G(\A_f)\times \pi_0(G_\infty)$-modules (\cite{Ha2}, p. 117.).

Next consider $W_0=\mathbb G_a\subset V_0=\mathbb G_a^2$ embedded in the first component. If $g_f\in G(\A_f)$ is given, we define $W(g_f):=W(\Q)\cap V(g_f)$. If we write $V(g_f)=\mathfrak a\oplus \mathfrak b$, we have $W(g_f)=\mathfrak a$. We recognize $W(\Q)$ as a $B(\Q)$-invariant submodule of $V(\Q)$ and $W(g_f)$ as a $B(g_f):=G(g_f)\cap B(\Q)$-invariant submodule.
Remember
\begin{equation*}
\Sym^kV(g_f)_{PD}=\Sym^k(\mathfrak a\oplus \mathfrak b)_{PD}=\bigoplus_{l=0}^k\Sym^{k-l}(\mathfrak a)_{PD}\otimes \Sym^l(\mathfrak b)_{PD}.
\end{equation*}
So we get the $B(g_f)$-invariant submodule
\begin{equation*}
W_k(g_f):=\bigoplus_{l=1}^k\Sym^{l}(\mathfrak a)_{PD}\otimes \Sym^{k-l}(\mathfrak b)_{PD}\subset \Sym^kV(g_f)_{PD}.
\end{equation*}
We set again $A:=\Z[\frac{1}{Nd_F}]$ and $R:=\Oc_{F^{\text{Hil}}}[\frac{1}{Nd_F}]$. When coefficients are extended to $A$ the submodule $W_k(g_f)_A^{Z_K}\subset \Sym^k _AV(g_f)_{PD}^{Z_K}$
is a direct summand. To see this we may do the faithfully flat ring extension $A\to R$. If we are in the non-trivial case $k=\xi m$, we have by \Cref{invariants} the elements $\prod_\sigma (\alpha_\sigma X_\sigma )^{[n_\sigma] }(\beta_\sigma Y_\sigma)^{[m-n_\sigma]}$, $n_\sigma=0,...,m$, as $R$-basis of $\Sym^k _AV(g_f)_{PD}^{Z_K}$ and we see immediately that the direct summand generated by those basis elements where $n_\sigma\geq1$ equals $W_k(g_f)_R^{Z_K}$.
In particular, we conclude using \Cref{Y_integral} that 
$\Sym^{k} _A V(g_f)_{PD}^{Z_K}/W_k(g_f)_A^{Z_K}$
is a free $A$-module of rank one generated by the residue class of
\begin{equation*}
\left\|t_2\right\|_f^m \prod_\sigma Y^{[m]} _\sigma,\ g_f=x\begin{pmatrix}t_1& t_1u\\0&t_2 \end{pmatrix},\ x\in G(\hat\Z),
\end{equation*}
if $k=\xi m$ and zero otherwise. Note that the $B(g_f)$-action factors over $T(g_f)$ the image of $B(g_f)$ in $T(\Q)$. 

Let us consider the locally constant sheaf $\underline\omega^k$ on $\partial\Sc_K$, which is associated to the $T(\Q)$-module $\Sym^{k}V(\Q)^{Z_K}/(\Sym^{k-1}V(\Q)\cdot W(\Q))^{Z_K}$. The sheaf associated to the $T(g_f)$-module $\Sym^{k} _A V(g_f)_{PD}^{Z_K}/W_k(g_f)_A^{Z_K}$ defines an integral structure $\underline\omega^k_{PD}$ for the sheaf $\underline\omega^k$. We have the natural projection map $\Sym^k\Hc^\prime\rightarrow \underline\omega^k$ of sheaves on $\Sc_{K}^B$. Set $K^T:=\text{image of $K^B$ in $T(\A)$}$. 

The natural map $B\rightarrow T$ induces a fiber bundle 
\begin{equation*}
p:\Sc_{K}^B\rightarrow K ^T\backslash T(\A)/T(\Q)=:\Sc_{K}^T
\end{equation*}
with fiber $p^{-1}(1)=K\cap U(\A)\backslash U(\A)/U(\Q)$, which is a topological torus of dimension $\xi$. As described in \Cref{trace_integration} we have the trace morphism
\begin{equation*}
Rp_*(\Sym^k\Hc^\prime\otimes \mu^{\otimes n})\rightarrow Rp_*(\underline\omega^k\otimes \mu^{\otimes n})\stackrel{tr}{\rightarrow} R^{\xi}p_*(\Q)\otimes\underline\omega^k \otimes \mu^{\otimes n}[-\xi].
\end{equation*}
Here we used the projection formula and that $\underline\omega ^k$ is pullback by $p$.
Of course, the trace morphism is defined integrally. Now we have the following result due to Harder.
\begin{theorem}\label{boundary_fiber}
The trace morphism
\begin{equation*}
H^p(\Sc_{K}^B,\Sym^k\Hc^\prime\otimes \mu^{\otimes n})\rightarrow H^{p-\xi}(\Sc_{K}^T,R^{\xi}p_*(\Q)\otimes\underline\omega^k \otimes \mu^{\otimes n}),\ p\geq \xi,
\end{equation*}
is an isomorphism.
\begin{proof}
By faithfully flatness and \Cref{ext_scalars} we may extend coefficients from $\Q$ to $\overline\Q$.
In this situation see \cite{Ha1} II. Important is \cite{Ha1} 2.8 for the occurring weights.
\end{proof} 
\end{theorem}
Harder gives a complete description of the cohomology group on the right-hand side, that we will recall now.
First we trivialize the sheaf $R^\xi p_*(\Q)$. We fix an ordering $\Sigma=\left\{\sigma_1,...,\sigma_\xi\right\}$.
We write elements
\begin{equation*}
b=\begin{pmatrix}t_1& t_1 x\\0 & t_2\end{pmatrix}\in B(\R)=B_0(F\otimes_\Q \R)
\end{equation*}
and $x=(x_1,...,x_\xi)\in F\otimes_\Q\R=\R^\xi$. Set
\begin{equation*}
\vol_{\Sc^B/\Sc^T}(t):=\frac{dx_1\wedge...\wedge dx_\xi}{\sqrt{d_F}\left\|t_2 ^{-1}t_1 \right\|_f \sgn(\Norm(t_2 ^{-1}t_1))^{-1}}.
\end{equation*}
This form is $T(\Q)$ invariant, closed and gives each fiber $p^{-1}(p(b))$ a rational volume, therefore we may interpret $\vol_{\Sc^B/\Sc^T}$ as a cohomology class in $H^0(\Sc^T _K,R^\xi p_*(\Q))$. If we have $K=K_N$, then $N^{-\xi}\vol_{\Sc^B/\Sc^T}$ gives each fiber exactly volume one and defines an integral trivialization. For $k=\xi m$ we have
\begin{equation*}
\left\|t_2\right\|_f ^m \sgn(t_2)^{-m}\prod_\sigma (-Y_\sigma)^{[m]} \in H^0(\Sc^T _K,\underline\omega^k_{PD})
\end{equation*}
in the cohomology with $A$-coefficients. It defines an $A$-integral trivialization of $\underline\omega^k_{PD}$. 
Finally,
\begin{equation*}
\left\|\det(g)\right\|_f \sgn(\Norm(\det(g)))^{-1}\vol_V^*=\vol_{\Tc/\Mc}^*,\ \vol_V^*:=\bigwedge_\sigma X_\sigma\wedge Y_\sigma
\end{equation*}
is a trivialization of $\mu$. If we have $K=K_N$, then it even defines a trivialization over the ring $A$, as one sees by considering the volumes of the fibers of $\Tc_W\rightarrow \Mc_{K_N}$, $W_f=V(\hat\Z)\rtimes K_N$. We suppose $K^T=Z(\A)\cap K  \times Z(\A)\cap K$ . We have a short exact sequence of groups
\begin{equation*}
1\rightarrow Z(\R)^0/ Z_K\stackrel{t\mapsto (s,t)}{\rightarrow}\Sc_K ^T\rightarrow \pi_0(\Sc_K ^T)=K^T T(\R)^0\backslash T(\A)/T(\Q)\rightarrow 1.
\end{equation*}
The cohomology of a connected component of $\Sc_K ^T$ is therefore isomorphic to the cohomology of  $Z(\R)^0/Z_K$. Because $Z_K$ acts properly discontinuously and fixpoint free on $Z(\R)^0$, we get $H^\bullet( Z(\R)^0/Z_K,\Z)=H^\bullet(Z_K,\Z)$. We may interpret these classes as invariant classes on $\Sc_K ^T$ as described in \Cref{invariant_Z_cohom} and as such they span a free $\Z$-submodule $\Hc^\bullet(T/Z)_K\subset H^\bullet(\Sc_K ^T,\Z)$ isomorphic to $H^\bullet(Z_K,\Z)$. We set $\Hc^\bullet(T/Z):=\Hc^\bullet(T/Z)_K\otimes_\Z\Q$ as this space does not depend on the level $K$.

For $m\in \N_0$ and $n\in \Z$ we fix characters
\begin{equation*}
\gamma_{m,n}:T(\R)\rightarrow \R^{\times},\ (t_1,t_2)\mapsto \Norm(t_2)^{-(m+n+1)}\Norm(t_1)^{-(n-1)}.
\end{equation*}
Harder considers algebraic Hecke characters of type $\gamma_{m,n}$ (\cite{Ha1} (2.5.2)), in other words, continuous homomorphisms
\begin{equation*}
\phi:T(\A)/T(\Q)\rightarrow \C^{\times}, \ \phi_{|T(\R)^0}=\gamma_{m,n|T(\R)^0}^{-1},
\end{equation*}
to decompose the cohomology of the boundary. We have $\phi_{|T(\R)}=\gamma_{m,n}\sgn(\phi)$ for a certain character $\sgn(\phi):T(\R)^0\backslash T(\R)\rightarrow \left\{\pm 1\right\}$. We set $\tilde{\phi}_f:=\phi_f\cdot \sgn(\phi)$, where as usual $\phi_f=\phi_{|T(\A_f)}$. Then $\tilde{\phi}_f(xa)=\tilde{\phi}_f(x)\gamma_{m,n}(a)$ for all $a\in T(\Q)$ and $x\in T(\A)$.
Moreover, $\tilde{\phi}_f$ takes its values in $\overline \Q$ and therefore $\text{Gal}(\overline\Q/\Q)$ acts on the set of $\tilde\phi_f$ by $ \sigma\cdot\tilde\phi_f(x):=\sigma(\tilde\phi_f(x))$, when we have $\sigma \in \text{Gal}(\overline\Q/\Q)$ and $x\in T(\A)$.

\begin{theorem}\label{cohom_boundary}
If $k=\xi m$ and $p\in \N_0$, we have an isomorphism of $\text{Gal}(\overline\Q/\Q)$-modules
\begin{equation*}
\bigoplus_{\phi:\ type(\phi)=\gamma_{m,n},\tilde{\phi}_{f|Z(\R)=1}}\overline\Q \cdot\tilde\phi_f\otimes \Hc^p(T/Z)\cong H^{p+\xi}(\Sc^B,\Sym^k\Hc^\prime\otimes \mu^{\otimes n}\otimes \overline\Q)
\end{equation*}
\begin{equation*}
\psi\otimes \eta\mapsto \frac{\psi \cdot \prod_\sigma (-Y_\sigma)^{[m]}\otimes \vol_V^{*n} dx_1\wedge...\wedge dx_{\xi }}{\sqrt{d_F}}\cup\eta
\end{equation*}
\begin{proof}
\cite{Ha1} Theorem 1.
\end{proof}
\end{theorem}
\begin{remark}
Harder's $\Q$-structure on cohomology, compare \cite{Ha1} 1.3, 2.4 and (2.7.1), agrees with our natural $\Q$-structure on $H^\bullet(\Sc_K,\Sym^k\Hc^\prime \otimes\mu^{\otimes n})$. The point is that Harder's $\Q$-structure on cohomology is $H^\bullet(\Sc_K,\Sym^k\Hc^\prime \otimes\mu^{\otimes n}\otimes \overline\Q)^{\text{Gal}(\overline\Q/\Q)}$, where the Galois action comes from functoriality, compare \Cref{invariants_integral} and \cite{Ha1} 1.3, 1.4. We have by \Cref{ext_scalars}
\begin{equation*}
H^\bullet(\Sc_K,\Sym^k\Hc^\prime \otimes\mu^{\otimes n}\otimes \overline\Q)^{\text{Gal}(\overline\Q/\Q)}=\left(H^\bullet(\Sc_K,\Sym^k\Hc^\prime \otimes\mu^{\otimes n})\otimes \overline\Q\right)^{\text{Gal}(\overline\Q/\Q)}=
\end{equation*}
\begin{equation*}
H^\bullet(\Sc_K,\Sym^k\Hc^\prime \otimes\mu^{\otimes n}).
\end{equation*}
\end{remark}
\begin{remark}
Even though Harder started with more general local systems than $\Sym^k\Hc^\prime\otimes \mu^{\otimes n}$, see \cite{Ha1} 1.4, he does not obtain more cohomology groups in degree $p=\xi,...,2\xi-1$ on the boundary than we do. \cite{Ha1} (2.8.2) tells us that we get all weights occurring in cohomological degrees $p=\xi,...,2\xi-1$.
\end{remark}

\section{Restriction of polylogarithmic Eisenstein classes to the boundary}

The next step is to consider the restriction of our polylogarithmic Eisenstein classes to the boundary.
We have the isomorphism of groups 
\begin{equation*}
\A_F\rightarrow U(\A),\ x\mapsto \begin{pmatrix}1&x\\0&1 \end{pmatrix}=u
\end{equation*}
inducing on $U(\A)$ a Haar measure $du$ corresponding to our Haar measure $dx$ on $\A_F$.
\begin{proposition}\label{boundary_residueI}
The cohomology class
\begin{equation*}
\res_{\Sc}(\Eis^k(\varphi))\in \bigoplus_{p+q=2\xi-1}H^{p}(\partial\Sc,\Sym^k\Hc^\prime\otimes \mu^{\otimes n+1})\otimes \mathfrak H^q
\end{equation*}
for $\varphi\in\Sc(V(\A_f),\mu^{\otimes n})^0 $ may be represented by the differential form
\begin{equation*}
\frac{ 2C(k)}{\kappa_F}\sum_{i=1}^\xi\sum_{I\subset\left\langle \xi\right\rangle\setminus i}\sum_{\chi\in \widehat{Cl_F ^K}(m)}
\int_{t\in Z(\A)}\phi(\hat{f})^m _{I\cup i}(tge^1,g) \chi(\det(g)t)\left\|\det(g)t\right\|^{m+2} d^{\times}t\cdot \Theta_{I,i}
\end{equation*}
In other words, $\res_{\Sc}(\Eis^k(\varphi))$ is determined by the zero Fourier coefficient of the $U(\Q)$-invariant form $\Eis^k(\varphi)$.
\begin{proof}
Let us set for a moment
\begin{equation*}
F_{I \cup i}(\chi,s,g):=\int_{t\in Z(\A)}
\phi(\hat{f})^m _{I\cup i}(tge^1,g) \chi(\det(g)t)\left\|\det(g)t\right\|^{m+2+2s} d^{\times}t.
\end{equation*}
By \Cref{cohom_boundary} and \Cref{boundary_fiber} we know that we may perform fiber integration, see \Cref{trace_integration}, to simplify our cohomology classes. The connected components of $\partial\Sc_K$ are fiber bundles and the fiber passing through $g\in G(\A)$ is 
\begin{equation*}
U(\R)/g_fK_fg_f^{-1}\cap U(\Q)=gKg^{-1}\cap U(\A)\backslash U(\A)/U(\Q)\hookrightarrow K\backslash G(\A)/B(\Q)=\partial \Sc_K,
\end{equation*}
$u\mapsto gu$. As we have
\begin{equation*}
\Eis^k(\varphi)=\frac{2 C(k)}{\kappa_F}\sum_{i=1}^\xi\sum_{I\subset\left\langle \xi\right\rangle\setminus i}\sum_{\chi\in \widehat{Cl_F ^K}(m)}\Eis^k_{I,i}(\varphi,\chi),
\end{equation*}
we have to do fiber integration for each $\Eis^k_{I,i}(\varphi,\chi)$ and for those we have
\begin{equation*}
du(g_fK_fg_f^{-1})^{-1}\int_{u\in U(\R)/g_fK_fg_f^{-1}\cap U(\Q)}\lim{s\rightarrow 0}\sum_{\gamma\in G(\Q)/B(\Q))}
F_{I \cup i}(\chi,s,tg_f k_\infty a_\infty u \gamma)du\cdot \Theta_{I,i},
\end{equation*}
where we have used the Iwasawa decomposition of $g_\infty=k_\infty a_\infty n_\infty\in G(\R)$ and that the forms $\Theta_{I,i}$ are constant in $u$ direction. If we use translation invariance of $du$, we may write the last integral as
\begin{equation*}
\int_{u\in U(\A)/U(\Q)}\lim{s\rightarrow 0}\sum_{\gamma\in G(\Q)/B(\Q))}
F_{I \cup i}(\chi,s,tgu\gamma)du\cdot \Theta_{I,i}.
\end{equation*}

By dominated convergence we may interchange the limit and integration, as the limit is a continuous function in $u$ and therefore integrable over the compact space $U(\A)/U(\Q)$. We obtain
\begin{equation*}
\lim{s\rightarrow 0}\int_{u\in U(\A)/U(\Q)}\sum_{\gamma\in G(\Q)/B(\Q))}
F_{I \cup i}(\chi,s,gu\gamma)du\cdot \Theta_{I,i}.
\end{equation*}
We write our integral as a sum
\begin{equation*}
\lim{s\rightarrow 0}\int_{u\in U(\A)/U(\Q)}F_{I \cup i}(\chi,s,gu)du\cdot \Theta_{I,i}+
\end{equation*}
\begin{equation*}
 \lim{s\rightarrow 0}\int_{u\in U(\A)/U(\Q)}\sum_{1\neq\gamma\in G(\Q)/B(\Q))}
 F_{I \cup i}(\chi,s,gu\gamma) du\cdot \Theta_{I,i}.
\end{equation*}
The first summand does not depend on $u$, as $ue^1=e^1$, so it equals
\begin{equation*}
\int_{t\in Z(\A)}\phi(\hat{f})^m _{I\cup i}(tge^1,g) \chi(\det(g)t)\left\|\det(g)t\right\|^{m+2} d^{\times}t\cdot \Theta_{I,i}
\end{equation*}
and this Tate integral exists by \cite{Wi} III Proposition 5.
To treat the second summand we consider the bijection
\begin{equation*}
U(\Q)\rightarrow \left\{\gamma B(\Q)\in G(\Q)/B(\Q): \gamma\notin B(\Q)\right\}
\end{equation*}
\begin{equation*}
u=\begin{pmatrix} 1 & \alpha\\ 0& 1\end{pmatrix}\mapsto uJ=\begin{pmatrix} \alpha & -1\\ 1& 0\end{pmatrix},\ J=\begin{pmatrix} 0 & -1\\ 1& 0\end{pmatrix}
\end{equation*} 
and obtain the integral
\begin{equation*}
\lim{s\rightarrow 0}\int_{u\in U(\A)}\int_{t\in Z(\A)}\phi(\hat{f})^m _{I\cup i}(tguJe^1,g) \chi(\det(g)t)\left\|\det(g)t\right\|^{m+2+2s} d^{\times}t du\cdot \Theta_{I,i}.
\end{equation*}
This integral is zero by the following lemma and our proposition is proved.
\end{proof}
\end{proposition}
\begin{lemma}
\begin{equation*}
\lim{s\rightarrow 0}\int_{u\in U(\A)}\int_{t\in Z(\A)}\phi(\hat{f})^m _{I\cup i}(tguJe^1,g) \chi(\det(g)t)\left\|\det(g)t\right\|^{m+2+2s} d^{\times}t du=0
\end{equation*}
\begin{proof}
By Fubini's theorem we may switch the order of integration and obtain the integral
\begin{equation*}
\lim{s\rightarrow 0}\int_{t\in Z(\A)}\int_{u\in U(\A)}\phi(\hat{f})^m _{I\cup i}(tguJe^1,g) \chi(\det(g)t)\left\|\det(g)t\right\|^{m+2+2s}  du d^{\times}t=
\end{equation*}
\begin{equation*}
\lim{s\rightarrow 0}\int_{t\in Z(\A)}\int_{u\in U(\A)}\phi(\hat{f})^m _{I\cup i|g}(tuJe^1,g) \chi(\det(g)t)\left\|\det(g)t\right\|^{m+2+2s}  du d^{\times}t,
\end{equation*}
where as usual $\psi_{|g}(v)=\psi(gv)$ for a Schwartz-Bruhat-function $\psi\in \Sc(V(\A))$.
The advantage of this integral is that we may interpret the inner integral as a Fourier transform, whereas the exterior integral may be examined with the theory of Tate integrals. 
We had $u=\begin{pmatrix}1&x\\0&1 \end{pmatrix}$ and get
\begin{equation*}
\lim{s\rightarrow 0}\int_{t\in Z(\A)}\int_{x\in \A_F}\phi(\hat{f})^m _{I\cup i|g}((tx,t),g) \chi(\det(g)t)\left\|\det(g)t\right\|^{m+2+2s}  dxd^{\times}t=
\end{equation*}
\begin{equation*}
\lim{s\rightarrow 0}\int_{t\in \I_F}\left\|t\right\|^{-1}\int_{x\in \A_F}\phi(\hat{f})^m _{I\cup i|g}((x,t),g) \chi(\det(g)t)\left\|\det(g)t\right\|^{m+2+2s}  dx d^{\times}t=
\end{equation*}
\begin{equation*}
\chi(\det(g))\left\|\det(g)\right\|^{m+2} \lim{s\rightarrow 0}\int_{t\in \I_F}\int_{x\in \A_F}\phi(\hat{f})^m _{I\cup i|g}((x,t),g) \chi(t)\left\|t\right\|^{m+1+2s}  dx  d^{\times}t
\end{equation*}
To examine the integral
\begin{equation*}
\lim{s\rightarrow 0}\int_{t\in \I_F}\int_{x\in \A_F}\phi(\hat{f})^m _{I\cup i|g}((x,t),g) \chi(t)\left\|t\right\|^{m+1+2s}  dx d^{\times}t
\end{equation*}
we treat the finite and the infinite places separately. Let us start with the finite part. 
It is of the form
\begin{equation*}
\int_{t\in \I_{F,f}}\int_{x\in \A_{F,f}}\hat{f}_{|g}((x,t),g)dx\left\|t\right\|^{1+m+2s}\chi(t) d^{\times}t.
\end{equation*}
We define the Schwartz-Bruhat function
\begin{equation*}
\varphi:\A_{F,f}\rightarrow \C,\ \varphi(u):=\int_{x\in \A_{F,f}}\hat{f}_{|g}((x,u),g)dx.
\end{equation*}
We denote the Fourier transform of Schwartz-Bruhat functions $\varphi:\A_F\rightarrow \C$ by
\begin{equation*}
P\varphi(v):=\int_{x\in \A_F}\phi(x)\overline{\psi_0(\Tr_{F/\Q}(xv))}dx, 
\end{equation*}
where $\psi_0:\A/\Q\rightarrow \C^{\times}$ is our fixed non-trivial character.
We see
\begin{equation*}
P\varphi(0)=\int_{x\in \A_{F,f}}\int_{u\in \A_{F,f}}\hat{f}_{|g}((u,x),g)dudx=\int_{v\in V(\A_f)}\hat{f}_{|g}(v,g)dv=
\end{equation*}
\begin{equation*}
\left\|\det(g)\right\|_f ^{-1}\hat{\hat{f}}(0,g)=\left\|\det(g)\right\|_f ^{-1}f(0,g)=0.
\end{equation*}
We add an appropriate Euler factor at infinity. For a place $\nu|\infty$ we set $\varphi^{p(m)}_\nu(u):=\exp(-\pi u^2)$, if $2|m$, and $\varphi^{p(m)}_\nu(u):=\exp(-\pi u^2)u$, if $2\nmid m$.
We define 
\begin{equation*}
\varphi^{p(m)}_\infty:F\otimes_\Q \R=\prod_{\nu|\infty} \R\rightarrow \C,\ (u_\nu)\mapsto \prod_{\nu|\infty}\varphi^{p(m)}_\nu(u_\nu),
\end{equation*}
\begin{equation*}
\varphi^{p(m)}:\A_F \rightarrow \C, u=(u_f,u_\infty)\mapsto \varphi^{p(m)}_\infty(u_\infty)\varphi(u_f),
\end{equation*}
and consider the Tate integral $\lim{s\rightarrow 0}\int_{t\in \I_F}\varphi^{p(m)}(t)\chi(t)\left\|t\right\|^{m+1+2s}d^{\times}t$. The integral exists in any case, since $P\varphi^{p(m)}(0)=0$ and we do not have a pole at $m+1+2s=1$, see \cite{T} Main Theorem 4.4.1. Moreover, the integral over the infinite places exists and is non-zero as computed in \cite{T} 2.5. We conclude that the finite integral $\lim{s\rightarrow 0}\int_{t\in \I_{F,f}}\varphi(t)\chi(t)\left\|t\right\|^{m+1+2s}d^{\times}t$ exists.
 
So it suffices to to show that 
\begin{equation*}
\lim{s\rightarrow 0}\int_{t\in F_\R^{\times}}\int_{x\in F_\R}\phi(\hat{f})^m _{I\cup i|g,\infty}((x,t),g) \chi(t)\left\|t\right\|^{m+1+2s}  dx  d^{\times}t=0.
\end{equation*}
We transform the integral back and and have to show the vanishing of
\begin{equation*}
\lim{s\rightarrow 0}\int_{t\in F_\R^{\times}}\int_{x\in F_\R}\phi(\hat{f})^m _{I\cup i|g,\infty}((tx,t),g) \chi(t)\left\|t\right\|^{m+2+2s}  dx d^{\times}t=
\end{equation*}
\begin{equation*}
2^\xi\lim{s\rightarrow 0}\int_{x\in F_\R}\int_{t\in F_\R^{\times,0}}\phi(\hat{f})^m _{I\cup i|g,\infty}((tx,t),g) \left\|t\right\|^{m+2+2s}   d^{\times}t dx.
\end{equation*}
The last equality holds, since we have $\chi\in \widehat{Cl_F ^K}(m)$. We use again the Iwasawa decomposition for $g\in G(\R)$ and calculate the last integral as
\begin{equation*}
\frac{2^\xi\Norm(e^{(m+1)i\theta})e^{i\theta}_{I\cup i}e^{-i\theta}_{(I\cup i)^c}}{\Norm(t_1)^{m+2}}\lim{s\rightarrow 0}\int_{x\in F_\R}\int_{t\in F_\R^{\times,0}}\phi(\hat{f})^m _{I\cup i,\infty}(t(x,y),g) \left\|t\right\|^{m+2+2s}   d^{\times}t dx,
\end{equation*}
with $y=\frac{t_2}{t_1}$. To end the proof we show   
\begin{equation*}
\lim{s\rightarrow 0}\int_{x\in F_\R}\int_{t\in F_\R^{\times,0}}\phi(\hat{f})^m _{I\cup i,\infty}(t(x,y),g) \left\|t\right\|^{m+2+2s}   d^{\times}t dx=0.
\end{equation*}
We have 
\begin{equation*}
\phi(\hat{f})^m_{I\cup i,\infty}(t(x,y),g)=
\exp(-\pi\left\|t(x,y)\right\|^2)\Norm(t(x+iy))^{m+1}(t(x+iy))_{I\cup i}(t(x-iy))_{(I\cup i)^c},
\end{equation*}
so
\begin{equation*}
\int_{t\in F_\R^{\times,0}}\phi(\hat{f})^m _{I\cup i,\infty}(t(x,y),g) \left\|t\right\|^{m+2+2s}   d^{\times}t=
\end{equation*}
\begin{equation*}
\frac{\Gamma(m+2+s)^\xi}{2^\xi\pi^{\xi(m+2+s)}}\frac{\Norm(x+iy)^{m+1}(x+iy)_{I\cup i}(x-iy)_{(I\cup i)^c}}{|\Norm(x+iy)|^{2(m+2+s)}}
\end{equation*}
and the integral in question is therefore
\begin{equation*}
\frac{\Gamma(m+2)^\xi}{2^\xi\pi^{\xi(m+2)}}\lim{s\rightarrow 0}\int_{x\in F_\R}\frac{dx}{\Norm(x-iy)^{(m+1)}(x-iy)_{I\cup i}(x+iy)_{(I\cup i)^c}|\Norm(x+iy)|^{2s}}.
\end{equation*}
This can be calculated as 
\begin{equation*}
\frac{\Gamma(m+2)^\xi}{\pi^{\xi(m+2)}}\int_{x\in F_\R}\frac{dx}{\Norm(x-iy)^{(m+1)}(x-iy)_{I\cup i}(x+iy)_{(I\cup i)^c}}
\end{equation*}
by dominated convergence, because the integral is dominated by
\begin{equation*}
c\int_{x\in F_\R}\frac{dx}{|\Norm(x+iy)|^{(m+2)}},\ c\in \R_{>0},
\end{equation*}
which exists by Fubini, as the corresponding one dimensional integrals $\int_{x\in\R}\frac{dx}{|x+iy|^{(m+2)}}$ exist for $y\in \R^{\times}$. This is seen as follows. The integral equals $\int_{x\in\R}\frac{dx}{(x^2+y^2)^{(\frac{m}{2}+1)}}$. We split it up as 
\begin{equation*}
\int_{x\in\R}\frac{dx}{(x^2+y^2)^{(\frac{m}{2}+1)}}=2\int_0 ^1 \frac{dx}{(x^2+y^2)^{(\frac{m}{2}+1)}} +2\int_1 ^\infty\frac{dx}{(x^2+y^2)^{(\frac{m}{2}+1)}}.
\end{equation*}
Because $\frac{1}{(x^2+y^2)^{(\frac{m}{2}+1)}}$ is as function in $x$ continuous on $[0,1]$, there is a constant $C>0$ such that $\frac{1}{(x^2+y^2)^{(\frac{m}{2}+1)}}<C$, $x\in [0,1]$. Therefore we may estimate our integral by $2C +2\int_1 ^\infty\frac{dx}{(x^2+y^2)^{(\frac{m}{2}+1)}}$. Moreover, we have $\frac{1}{(x^2+y^2)^{(\frac{m}{2}+1)}}\leq \frac{1}{x^{m+2}}$, $x\in [1,\infty[$, and our integral may be estimated by $2C +2\int_1 ^\infty\frac{dx}{x^{m+2}}$, which is finite as $m\geq 0$.
Let us come back to
\begin{equation*}
\frac{\Gamma(m+2)^\xi}{2^\xi\pi^{\xi(m+2)}}\int_{x\in F_\R}\frac{dx}{\Norm(x-iy)^{(m+1)}(x-iy)_{I\cup i}(x+iy)_{(I\cup i)^c}},
\end{equation*}
which is a product of integrals. If we want to show that it vanishes, it suffices to show the vanishing of a single factor, for example the $i$th-factor. This factor is essentially the integral
\begin{equation*}
\int_{x\in\R}\frac{1}{(x-iy)^{(m+2)}}dx=\left[\frac{1}{-(m+1)(x-iy)^{m+1}}\right]_{-\infty} ^{\ \infty}=0.
\end{equation*}
\end{proof}
\end{lemma}
Consider the Iwasawa decomposition 
\begin{equation*}
g =kb\in G(\A) \text{ with } b=\begin{pmatrix}t_1 & t_1x\\ 0 & t_2 \end{pmatrix} \text{ and } k\in G(\hat\Z)\cdot\prod_{\nu|\infty} SO(2).
\end{equation*}
We define a differential form on $K^1\backslash G(\A)$:
\begin{equation*}
\vol_P:=\sum_{i=1} ^\xi \frac{(-2)^{\xi}}{\xi\kappa_F}\prod_{j\neq i}\left(\frac{dt_{2,j}}{t_{2,j}}\wedge dx_j\right) \wedge dx_i.
\end{equation*}
\begin{lemma}\label{Theta_vol_P}
We have
\begin{equation*}
\frac{\xi\kappa_F\vol_P}{\left\|t_2 ^{-1}t_1\right\|_f\sgn(\Norm(t_2 ^{-1}t_1))} =
\sum_{i=1} ^\xi\sum_{I\subset \left\langle \xi\right\rangle\setminus i}\frac{-2(-1)^{|I|}}{\left\|t_2 ^{-1}t_1\right\|_f\sgn(\Norm(t_2 ^{-1}t_1))} \frac{d\tau\wedge d\overline\tau}{\overline\tau-\tau}_{(i\cup I)^c}\wedge \frac{dr\wedge d \overline\tau}{r}_I\wedge d\overline\tau_i
\end{equation*}
as cohomology classes in $H^{2\xi-1}(\partial\Mc_K,\C)=\bigoplus_{p+q=2\xi-1}H^{p}(\partial\Sc_K,\C)\otimes \mathfrak H^q$. More precisely, if $A=\Z[\frac{1}{Nd_F}]$, $K=K_N$ and $N\geq 3$
\begin{equation*}
\frac{d^{\times}t(K_{N,f} ^Z)h_{K_N}\vol_P}{(-2)^\xi\sqrt{d_F}\left\|t_2 ^{-1}t_1\right\|_f\sgn(\Norm(t_2 ^{-1}t_1))}
\end{equation*}
is in the image of the cohomology with $A$-coefficients and part of an $A$-basis of the $A$-integral classes.
\begin{proof}
We have $\frac{d\tau\wedge d\overline\tau}{\overline\tau-\tau}=\frac{dx\wedge dy}{y}$ and we see that the top $x$-component of
\begin{equation*}
\sum_{i=1} ^\xi\sum_{I\subset \left\langle \xi\right\rangle\setminus i}(-1)^{|I|}\frac{d\tau\wedge d\overline\tau}{\overline\tau-\tau}_{(i\cup I)^c}\wedge \frac{dr\wedge d \overline\tau}{r}_I\wedge d\overline\tau_i
\end{equation*}
equals
\begin{equation*}
\sum_{i=1} ^\xi\sum_{I\subset \left\langle \xi\right\rangle\setminus i}(-1)^{|I|}\frac{dx\wedge dy}{y}_{(i\cup I)^c}\wedge \frac{dr\wedge dx}{r}_I\wedge dx_i.
\end{equation*}
By \Cref{boundary_fiber} we already know that these two classes have to be cohomologous. We have the formulas
\begin{equation*}
ry\, d\left(\frac{1}{ry}\right)\wedge dx=-\frac{d(ry)}{ry}\wedge dx=\frac{dx\wedge dy}{y}-\frac{dr}{r}\wedge dx,\ ry\, d\left(\frac{1}{ry}\right)\wedge dx=(-2)\frac{dt_{2}}{t_{2}}\wedge dx
\end{equation*}
from which the claimed equality of forms follows. What remains to be shown is the integrality statement.
Consider components $\left\{\pm1\right\}^\xi\R_{>0}\backslash B(\R)/B(g_f)\subset\partial\Mc_{K_N}$. Set $P_0:=\left\{g\in G_0:ge^1=e^1\right\}$ and get $P:=Res_{\Oc/\Z}P_0$. Define
\begin{equation*}
P(\R)^1:=\left\{x\in P(\R):\ x\in P(\R)^0,\ \Norm(\det(x))=1\right\}.
\end{equation*}
This gives the submanifold
\begin{equation*}
i:P(\R)^1/P(g_f)\rightarrow \left\{\pm1\right\}^\xi\R_{>0}\backslash B(\R)/B(g_f)\text{ with }P(g_f):=B(g_f)\cap P(\R) ^0.
\end{equation*}
The form
\begin{equation*}
\frac{d^{\times}t(K_{N,f} ^Z)h_{K_N}}{(-2)^\xi\sqrt{d_F}\left\|t_2 ^{-1}t_1\right\|_f\sgn(\Norm(t_2 ^{-1}t_1))}\vol_P
\end{equation*}
is a volume form of the compact space on the left-hand side giving it volume $N^\xi\in A^{\times}$. So the form is (up to a factor of $N^\xi$) dual to the fundamental class under Poincaré duality and hence defines an $A$-integral class on $P(\R)^1/P(g_f)$. As we also have a right-inverse 
\begin{equation*}
p:\left\{\pm1\right\}^\xi\R_{>0}\backslash B(\R)/B(g_f)\rightarrow P(\R)^1/P(g_f),\ \begin{pmatrix}t_1&t_1x\\0&t_2 \end{pmatrix}\mapsto\begin{pmatrix}1 & x\\0&\frac{|t_2t_1^{-1}|}{\sqrt[n]{|\Norm(t_2t_1^{-1})|}} \end{pmatrix} 
\end{equation*}
for the map $i$, we may conclude that our form is part of an $A$-basis of the cohomology of $\partial \Mc_{K_N}$.
\end{proof}
\end{lemma}

\begin{proposition}\label{boundary_residueII}
For $\varphi\in\Sc(V(\A_f),\mu^{\otimes n})^0 $ the $K_f$-invariant cohomology class 
\begin{equation*}
\res_{\Sc}(\Eis^k(\varphi))\in \bigoplus_{p+q=2\xi-1}H^{p}(\partial\Sc,\Sym^k\Hc^\prime\otimes \mu^{\otimes n+1})\otimes \mathfrak H^q
\end{equation*}
equals
\begin{equation*}
\sum_{\chi\in \widehat{Cl_F ^K}(m)}\frac{\Gamma(m+2)^{\xi}}{(-2\pi i)^{\xi(m+2)}}\int_{t\in Z(\A_f)}\hat{f} (tg_fe^1,g) \chi(\det(g_f)t)\left\|\det(g_f)t\right\|^{m+2} d^{\times}t\cdot
\end{equation*}
\begin{equation*}
 \frac{\Norm(\omega)^m \vol_P\otimes \vol_V^{*n+1}}{(m!)^\xi\left(\left\|\det(g)\right\|_f\sgn(\Norm(\det(g)))^{-1}\right)^{-n}}
\end{equation*}
In particular,
\begin{equation*}
\sum_{\chi\in \widehat{Cl_F ^K}(m)}\frac{\Gamma(m+2)^{\xi}}{(-2\pi i)^{\xi(m+2)}}\int_{t\in Z(\A_f)}\hat{f} (tg_fe^1,g) \chi(\det(g_f)t)\left\|\det(g_f)t\right\|^{m+2} d^{\times}t
\end{equation*}
is a rational number.

\begin{proof}
We calculate
\begin{equation*}
\int_{t\in Z(\R)}\phi(\hat{f})^m _{I\cup i,\infty}(tg_\infty e^1,g_\infty) \chi(\det(g_\infty)t)\left\|\det(g_\infty)t\right\|^{m+2} d^{\times}t=
\end{equation*}
\begin{equation*}
\Norm(e^{(m+1)i\theta})e^{i\theta}_{I\cup i}e^{-i\theta}_{(I\cup i)^c}\frac{\Norm(t_{2}) ^{m+2}\Gamma(m+2)^{\xi}}{\pi^{\xi(m+2)}}
\end{equation*}
Remembering the definition of the $\Theta_{I,i}$ (\Cref{Theta_2}) the formula for $\res_{\Sc}(\Eis^k(\varphi))$ follows from \Cref{boundary_residueI} and \Cref{Theta_vol_P}. The rationality statement follows, since
\begin{equation*}
\left\|t_2\right\|_f ^{m+2}\sgn(\Norm(t_2))^{m+2}\left\|\det(g)\right\|_f ^{n}\sgn(\Norm(\det g))^{n}\Norm(\omega)^m \vol_P\otimes \vol_V^{*n+1}=
\end{equation*}
\begin{equation*}
\left\|t_2\right\|_f ^{m+2}\sgn(\Norm(t_2))^{m+2}\left\|\det(g)\right\|_f ^{n}\sgn(\Norm(\det g))^{n}\prod_\sigma(-Y_\sigma)^m \vol_P\otimes \vol_V^{*n+1}
\end{equation*}
as cohomology classes in $ H^\bullet(\partial\Sc,\underline\omega^k\otimes \mu^{\otimes n+1})\otimes \mathfrak H^\bullet$. The latter defines a rational class, compare \Cref{boundary}, \Cref{Theta_vol_P} and recall $d^{\times}t(K_f ^Z)\in \Q^{\times}\frac{1}{\sqrt{d_F}}$, since $K_f ^Z\subset \hat\Oc ^{\times}$ is a subgroup of finite index and $d^{\times}t(\hat\Oc^{\times})=\frac{1}{\sqrt{d_F}}$.
\end{proof}
\end{proposition}
Of course, we can refine this statement, if we use the integral structures on cohomology, which we have chosen.
\begin{proposition}\label{integral_L_value}
Let $N\geq 3$ be an integer, $A:=\Z[\frac{1}{d_FN}]$ and $f\in \text{map}(V(\Z/N\Z),A)$ with $f(0)=\sum_{v\in V(Z/N\Z)}f(v)=0$. We denote by $\Oc^{\times}(N)^+\subset \Oc^{\times}$ the group of totally positive units which are congruent to one modulo $N$. We write
\begin{equation*}
g=k\begin{pmatrix}t_1& t_1x\\ 0& t_2 \end{pmatrix}\in G(\A_f)\times \pi_0(G(\R)),\ k\in G(\hat\Z),
\end{equation*}
and get a well-defined function in $map(K_N\backslash G(\A_f)\times \pi_0(G(\R))/B(\Q),A)$
\begin{equation*}
g\mapsto \frac{(-1)^\xi\sqrt{d_F}\Gamma(m+2)^\xi }{(-2\pi i)^{\xi(m+2)}\sgn(\Norm(t_2))^{m+2}\left\|t_2\right\|_f^{m+2}}\sum_{l\in F^{\times}/\Oc^{\times}(N)^+}\frac{\hat{f}(\hat{g}_fle^1)}{ \Norm(l)^{m+2}}\in A
\end{equation*}
\begin{proof}
We apply the construction of polylogarithmic Eisenstein classes to $f$ and obtain by \Cref{pol}
\begin{equation*}
\Eis^k(f)\in H^{2\xi-1}(\Mc_{K_N}, \Sym^k\Hc_{PD}\otimes \mu)
\end{equation*}
with $A$-coefficients. Since $\Sym^k\Hc^\prime _{PD}\subset \Sym^k\Hc_{PD}$ is a direct summand, we obtain by restricting $\Eis^k(f)$ to the boundary $\partial \Mc_{K_N}$
\begin{equation*}
\res_{\Mc}(\Eis^k(f))\in H^{2\xi-1}(\partial\Mc_{K_N}, \underline{\omega}^k _{PD}\otimes \mu)
\end{equation*}
with $A$-coefficients. We have calculated $\res_{\Mc}(\Eis^k(f))$ with $\C$-coefficients as
\begin{equation*}
\sum_{\chi\in \widehat{Cl_F ^{K_N}}(m)}\frac{\Gamma(m+2)^{\xi}}{(-2\pi i)^{\xi(m+2)}}\int_{t\in Z(\A_f)}\hat{f} (tg_fe^1) \chi(\det(g_f)t)\left\|\det(g_f)t\right\|^{m+2} d^{\times}t \cdot
\end{equation*}
\begin{equation*}
 \frac{\Norm(\omega)^m \vol_P\otimes \vol_V^{*}}{(m!)^\xi}
\end{equation*}
Consider $\varphi_\infty\in \Sc(F_\R)$ defined by $\varphi_\infty(v):=\exp(-\left\|v\right\|^2)\Norm(v)^{m+2}$, $v\in F_\R$. We have for $\chi\in \widehat{Cl_F ^{K_N}}(m)$
\begin{equation*}
\int_{t\in F_\R^{\times}}\varphi_\infty(t)\chi_\infty(t)|\Norm(t)|^{m+2}d^{\times}t=\Gamma(m+2)^\xi.
\end{equation*}
Define $\varphi\in \Sc(\A_F)$ by $\varphi(v):=\varphi_\infty(v_\infty)\hat{f}(g_fv_fe^1)$ for $v\in \A_F$. We calculate the integral
\begin{equation*}
\sum_{\chi\in \widehat{Cl_F ^{K_N}}(m)}\frac{\Gamma(m+2)^{\xi}}{(-2\pi i)^{\xi(m+2)}}\int_{t\in Z(\A_f)}\hat{f} (tg_fe^1) \chi(\det(g_f)t)\left\|\det(g_f)t\right\|^{m+2} d^{\times}t=
\end{equation*}
\begin{equation*}
\sum_{\chi\in \widehat{Cl_F ^{K_N}}(m)}\frac{1}{(-2\pi i)^{\xi(m+2)}}\int_{t\in Z(\A)}\varphi(t)\chi(\det(g_f)t)\left\|\det(g_f)t\right\|^{m+2} d^{\times}t=
\end{equation*}
\begin{equation*}
\sum_{\chi\in \widehat{Cl_F ^{K_N}}}\frac{1}{(-2\pi i)^{\xi(m+2)}}\int_{t\in Z(\A)}\varphi(t)\chi(\det(g_f)t)\left\|\det(g_f)t\right\|^{m+2} d^{\times}t=
\end{equation*}
\begin{equation*}
\sum_{\chi\in \widehat{Cl_F ^{K_N}}}\frac{d^{\times}t(K_{N,f} ^Z)}{(-2\pi i)^{\xi(m+2)}}\int_{t\in K_{N,f} ^Z\backslash Z(\A)/F^{\times}}\sum_{l\in F^{\times}}\varphi(tl)\chi(\det(g_f)t)\left\|\det(g_f)t\right\|^{m+2} d^{\times}t=
\end{equation*}
\begin{equation*}
\frac{d^{\times}t(K_{N,f} ^Z)h_{K_N}}{(-2\pi i)^{\xi(m+2)}}\int_{t\in Z(\R)^0/\Oc^{\times}(N)^+}\sum_{l\in F^{\times}}\varphi(\det(g_f)^{-1}tl)\left\|t\right\|^{m+2} d^{\times}t=
\end{equation*}
\begin{equation*}
\frac{d^{\times}t(K_{N,f}^Z)h_{K_N}}{(-2\pi i)^{\xi(m+2)}}\sum_{l\in F^{\times}/\Oc^{\times}(N)^+}\int_{t\in Z(\R)^0}\varphi(\det(g_f)^{-1}tl)\left\|t\right\|^{m+2} d^{\times}t=
\end{equation*}
\begin{equation*}
\frac{d^{\times}t(K_{N,f}^Z)h_{K_N}\Gamma(m+2)^\xi}{2^\xi(-2\pi i)^{\xi(m+2)}}\sum_{l\in F^{\times}/\Oc^{\times}(N)^+}\frac{\hat{f}(\hat{g}_fle_1)}{\Norm(l)^{m+2}}.
\end{equation*}
The form
\begin{equation*}
\frac{d^{\times}t(K_{N,f} ^Z)h_{K_N}}{\sqrt{d_F}(-2)^\xi\left\|\frac{t_1}{t_2}\right\|_f\sgn(\Norm(\frac{t_1}{t_2}))}\vol_P
\end{equation*}
is part of an $A$-basis of the cohomology of $\partial\Mc_{K_N}$ as seen in \Cref{Theta_vol_P}. We already have mentioned in \Cref{boundary} that
\begin{equation*}
\left\|t_2\right\|_f^{m}\sgn(\Norm(t_2))^{m}\prod_{\sigma}(-Y_\sigma)^{[m]}=\frac{\left\|t_2\right\|_f^{m}\sgn(\Norm(t_2))^{m}\Norm(\omega)^m}{(m!)^\xi}\in H^0(\partial\Mc_{K_N},\underline\omega^k_{PD})
\end{equation*}
and 
\begin{equation*}
\left\|\det(g)\right\|_f\sgn(\Norm(\det(g)))\vol_V^*\in H^0(\partial\Sc_{K_N},\mu)
\end{equation*}
define $A$-trivializations and therefore 
\begin{equation*}
\frac{d^{\times}t(K_{N,f} ^Z)h_{K_N}\left\|t_2\right\|_f^{m}\sgn(\Norm(t_2))^{m}\left\|\det(g)\right\|_f\sgn(\Norm(\det(g)))\Norm(\omega)^m\vol_P\otimes \vol_V^* }{\sqrt{d_F}(m!)^\xi(-2)^\xi\left\|\frac{t_1}{t_2}\right\|_f\sgn(\Norm(\frac{t_1}{t_2}))}
\end{equation*}
is part of an $A$-basis of 
\begin{equation*}
\image(H^{2\xi-1}(\partial\Mc_{K_N}, \underline{\omega}^k _{PD}\otimes \mu)\rightarrow  H^{2\xi-1}(\partial\Mc_{K_N}, \underline{\omega}^k _{PD}\otimes \mu\otimes \C)).
\end{equation*}
Now we may conclude that
\begin{equation*}
\frac{(-1)^\xi\sqrt{d_F}\Gamma(m+2)^\xi }{(-2\pi i)^{\xi(m+2)}\sgn(\Norm(t_2))^{m+2}\left\|t_2\right\|_f^{m+2}}\sum_{l\in F^{\times}/\Oc^{\times}(N)^+}\frac{\hat{f}(\hat{g}_fle^1)}{ \Norm(l)^{m+2}}\in A,
\end{equation*}
as $\res_{\Mc}(\Eis^k(f))$ is $A$-integral.
\end{proof}
\end{proposition}
\begin{definition}
Let us define the \textit{horospherical map}
\begin{equation*}
\rho_{m,n}^0:\Sc(V(\A),\mu^{\otimes n}\otimes \overline\Q)^0\rightarrow \Ind_{B(\A_f)\times \pi_0(B(\R))}^{G(\A_f)\times \pi_0(G(\R))}\bigoplus_{\phi:\ type(\phi)=\gamma_{m,n+1},\tilde{\phi}_{f|Z(\R)=1}}\overline\Q \cdot\phi,
\end{equation*}
by
\begin{equation*}
\rho_{m,n}^0(\varphi)(g):= \sum_{\chi\in \widehat{Cl_F ^{K}}(m)}\frac{\Gamma(m+2)^{\xi}}{(-2\pi i)^{\xi(m+2)}}\frac{\int_{t\in Z(\A_f)}\hat{f} (tg_fe^1,g) \chi(\det(g_f)t)\left\|\det(g_f)t\right\|^{m+2} d^{\times}t}{\left(\left\|\det(g)\right\|_f\sgn(\Norm(\det(g)))^{-1}\right)^{-n}},
\end{equation*}
where we have written
\begin{equation*}
\varphi(v,g)=f(v,g)\left(\left\|\det(g)\right\|_f\sgn(\Norm(\det(g)))^{-1}\right)^{n}
\end{equation*}
with $f\in \Sc(V(\A_f),\mu^{\otimes 0}\otimes \overline\Q)^0$ and $K_f$-invariant.
\end{definition}
\begin{remark}
Let us quickly recall why $\rho_{m,n}^0$ is well-defined. Take $K^\prime_f \subset K_f$, $\chi\in Cl_F ^{K^\prime}(m)$ and $\epsilon \in K_f$ with $\chi(\epsilon)\neq 1$. Then
\begin{equation*}
\int_{t\in Z(\A_f)}\hat{f} (tg_fe^1,g) \chi(\det(g_f)t)\left\|\det(g_f)t\right\|^{m+2} d^{\times}t=
\end{equation*}
\begin{equation*}
\int_{t\in Z(\A_f)}\hat{f} (\epsilon tg_fe^1,g) \chi(\det(g_f)\epsilon t)\left\|\det(g_f)\epsilon t\right\|^{m+2} d^{\times}t=
\end{equation*}
\begin{equation*}
\chi(\epsilon)\int_{t\in Z(\A_f)}\hat{f} ( tg_fe^1,g) \chi(\det(g_f) t)\left\|\det(g_f) t\right\|^{m+2} d^{\times}t,
\end{equation*}
so we conclude for such $\chi$
\begin{equation*}
0=\int_{t\in Z(\A_f)}\hat{f} ( tg_fe^1,g) \chi(\det(g_f) t)\left\|\det(g_f) t\right\|^{m+2} d^{\times}t
\end{equation*}
and the definition of $\rho_{m,n}^0$ is independent of the chosen $K_f$. To see that $\rho_{m,n}^0$ takes its values in $\overline{\Q}$-valued functions we can actually imitate the proof of \cref{boundary_residueII}. Given a $\varphi\in\Sc(V(\A),\mu^{\otimes n}\otimes \overline\Q)^0$ we get that
\begin{equation*}
\res_{\Sc}(\Eis^k(\varphi))\in \bigoplus_{p+q=2\xi-1}H^{p}(\partial\Sc,\Sym^k\Hc^\prime\otimes \mu^{\otimes n+1})\otimes \mathfrak H^q\otimes \overline\Q
\end{equation*}
\begin{equation*}
\cong\bigoplus_{\phi:\ type(\phi)=\gamma_{m,n},\tilde{\phi}_{f|Z(\R)=1}}\Ind_{B(\A_f)\times \pi_0(B(\R))}^{G(\A_f)\times \pi_0(G(\R))}\overline\Q \cdot\tilde\phi_f\otimes \Hc^p(T/Z)
\end{equation*}
equals
\begin{equation*}
\sum_{\chi\in \widehat{Cl_F ^K}(m)}\frac{\Gamma(m+2)^{\xi}}{(-2\pi i)^{\xi(m+2)}}\int_{t\in Z(\A_f)}\hat{f} (tg_fe^1,g) \chi(\det(g_f)t)\left\|\det(g_f)t\right\|^{m+2} d^{\times}t\cdot
\end{equation*}
\begin{equation*}
 \frac{\Norm(\omega)^m \vol_P\otimes \vol_V^{*n+1}}{(m!)^\xi\left(\left\|\det(g)\right\|_f\sgn(\Norm(\det(g)))^{-1}\right)^{-n}}
\end{equation*}
and therefore
\begin{equation*}
\sum_{\chi\in \widehat{Cl_F ^K}(m)}\frac{\Gamma(m+2)^{\xi}}{(-2\pi i)^{\xi(m+2)}}\int_{t\in Z(\A_f)}\hat{f} (tg_fe^1,g) \chi(\det(g_f)t)\left\|\det(g_f)t\right\|^{m+2} d^{\times}t
\end{equation*}
is an algebraic number, since the cohomology class
\begin{equation*}
 \frac{\Norm(\omega)^m \vol_P\otimes \vol_V^{*n+1}}{(m!)^\xi\left(\left\|\det(g)\right\|_f\sgn(\Norm(\det(g)))^{-1}\right)^{-n}}
\end{equation*} 
is rational. We may even conclude that $\rho_{m,n}^0$ is $\text{Gal}(\overline\Q/\Q)$-equivariant.
\end{remark}
\begin{proposition}\label{hor}
We have
\begin{equation*}
\image(\res_\Sc\circ \Eis^k _q)\cong \image(\rho_{m,n}^0)\otimes\Hc(T/Z)^{\xi-1-q}.
\end{equation*}
\begin{proof}
We prove the proposition with $\C$-coefficients and all of the following constructions will respect the chosen $\Q$-structures. By \Cref{boundary_residueII}
\begin{equation*}
\res_\Sc\circ \Eis^k(\varphi)=\rho_{m,n}^0(\varphi)\cdot\frac{\Norm(\omega)^m}{(m!)^\xi} \vol_P\otimes \vol_V^{*n+1}\in 
H^\bullet(\partial\Sc,\Sym^k\Hc^\prime\otimes \mu^{\otimes n+1})\otimes \mathfrak H^\bullet _\C.
\end{equation*}
Set
\begin{equation*}
\vol_{T_2}=\frac{(-1)^{\frac{\xi(\xi-1)}{2}}(-2)^{\xi}\sqrt{d_F}}{\xi\kappa_F}\sum_{i=1} ^\xi (-1)^{i-1}\frac{dt_{2}}{t_{2}}_{\left\langle \xi\right\rangle\setminus i},\ \vol_{U}=\frac{dx_1\wedge...\wedge dx_n}{\sqrt{d_F}}
\end{equation*}
and get $\vol_P=\vol_{T_2}\wedge \vol_{U}$.
Dividing $\res_\Sc\circ \Eis^k(\varphi)$ by $\frac{\Norm(\omega)^m}{(m!)^\xi}\vol_U\otimes\vol_V^{*n+1}$ we get
\begin{equation*}
\rho_{m,n}^0(\varphi)\cdot \vol_{T_2}\in
\Ind_{B(\A_f)\times \pi_0(B(\R))}^{G(\A_f)\times \pi_0(G(\R))}\bigoplus_{\phi:\ type(\phi)=\gamma_{m,n+1},\tilde{\phi}_{f|Z(\R)=1}}\C \cdot\phi \cdot \Hc^\bullet (T/Z)\otimes \mathfrak H^\bullet_\C
\end{equation*}
So what remains to be shown is that $\mathfrak H^{\bullet,*}_\C\rightarrow \Hc^\bullet (T/Z)_\C$, $v \mapsto \vol_{T_2}(v)$ is surjective. Clearly, $\Hc^\bullet (T/Z)\otimes \mathfrak H^{\bullet}_\C=H^\bullet(T(\R)^0/T(\Q)\cap (K_f\cdot T(\R)^0,\C)$ and we may interpret $\vol_{T_2}$ as rational volume form in $t_2$-direction. Therefore 
\begin{equation*}
\vol_{T_2}=c^\prime\cdot \frac{dt_{2}}{t_2}_{\left\langle \xi-1\right\rangle} \in H^\bullet(T(\R)^0/T(\Q)\cap (K_f\cdot T(\R)^0,\C)
\end{equation*}
for a $c^\prime\in \C^{\times}$. We take $\frac{dr}{r}_i$, $i=1,...,\xi-1$, as basis for $\mathfrak H_\C ^1$ and $\frac{dy}{y}_i$, $i=1,...,\xi-1$, as basis for $\Hc^1 (T/Z)_\C$. We have $\frac{dt_{2,i}}{t_{2,i}}=\frac{1}{2}\left(\frac{dy}{y}_i+\frac{dr}{r}_i\right)$, as $t_{2,i}=\sqrt{y_ir_i}$, and therefore
\begin{equation*}
 \vol_{T_2}=c^\prime\cdot \frac{1}{2}\left(\frac{dy}{y}_1+\frac{dr}{r}_1\right)\wedge...\wedge \frac{1}{2}\left(\frac{dy}{y}_{\xi-1}+\frac{dr}{r}_{\xi-1}\right)=
\end{equation*}
\begin{equation*}
 \frac{c^\prime}{2^{\xi-1}}\cdot\sum_{I\subset \left\langle \xi-1\right\rangle}\sgn(I)\frac{dy}{y}_I\wedge\frac{dr}{r}_{\left\langle \xi-1\right\rangle\setminus I}.
\end{equation*}
If we take the dual basis $V(r)_i\in \mathfrak H_\C ^{1\ *}$ of $\frac{dr}{r}_i$, $i=1,...,\xi-1$, we may see
\begin{equation*}
 \vol_{T_2}(V(r)_J)= \frac{c^\prime}{2^{\xi-1}}\cdot \sgn(\left\langle \xi-1\right\rangle\setminus J)\frac{dy}{y}_{\left\langle \xi-1\right\rangle\setminus J},\ J\subset \left\langle \xi-1\right\rangle.
\end{equation*}
In particular, $\mathfrak H^{\bullet\ *}_\C\rightarrow \Hc^\bullet (T/Z)_\C$, $v \mapsto \vol_{T_2}(v)$, is surjective. 
\end{proof}
\end{proposition}
\begin{remark}\label{im_res_Eis}
The last proposition shows that for any $\varphi\in \Sc(V(\A_f),\mu^{\otimes n}\otimes \C)^0$
\begin{equation*}
\rho_{m,n}^0(\varphi)\frac{\Norm(\omega)^m}{(m!)^\xi}\otimes \vol_V^{*n+1}\otimes\frac{d\tau\wedge d\overline \tau}{\overline\tau-\tau}_{I^c}\wedge d\overline\tau_{I},\ \emptyset\neq I\subset\left\langle \xi\right\rangle,
\end{equation*}
is in the image of $\sum_q \res_\Sc\circ \Eis^k_{q}$ with $\C$-coefficients.
\end{remark}

\section{Comparison with Harder's Eisenstein classes}

The reason why we talk about polylogarithmic Eisenstein classes is that these classes are actually Eisenstein cohomology classes in the sense of Harder. This is what we prove now.\\ 
\begin{theorem}\label{Harder}
Define $\tilde{H}^\bullet(\Sc,\Sym^k\Hc^\prime\otimes \mu^{\otimes n+1}):=$
\begin{equation*}
\ker\left(\res_{\Sc}:H^\bullet(\Sc,\Sym^k\Hc^\prime\otimes \mu^{\otimes n+1})\rightarrow H^\bullet(\partial\Sc,\Sym^k\Hc^\prime\otimes \mu^{\otimes n+1})\right).
\end{equation*}
One has a $G(\A_f)\times \pi_0(G(\R))$-equivariant operator 
\begin{equation*}
\Eis_{Harder}:\image(\res_\Sc)\rightarrow H^\bullet(\Sc,\Sym^k\Hc^\prime\otimes \mu^{\otimes n+1}),
\end{equation*}
which is a section for $\res_{\Sc}$. We denote the image of the operator $\Eis_{Harder}$ by
\begin{equation*}
H^\bullet_{\Eis}(\Sc,\Sym^k\Hc^\prime\otimes \mu^{\otimes n+1})
\end{equation*}
and get a decomposition of $G(\A_f)\times \pi_0(G(\R))$-modules
\begin{equation*}
H^\bullet(\Sc,\Sym^k\Hc^\prime\otimes \mu^{\otimes n+1})=\tilde{H}^\bullet(\Sc,\Sym^k\Hc^\prime\otimes \mu^{\otimes n+1})\oplus H^\bullet_{\Eis}(\Sc,\Sym^k\Hc^\prime\otimes \mu^{\otimes n+1})
\end{equation*}
\begin{proof}
\cite{Ha1} Theorem 2.
\end{proof}
\end{theorem}
The aim of this section is to prove the following
\begin{theorem}\label{pol_Eis}
\begin{equation*}
\image(\Eis^k_q)\subset H^\bullet_{\Eis}(\Sc,\Sym^k\Hc^\prime\otimes \mu^{\otimes n+1})
\end{equation*}
\end{theorem}
We generated $\mathfrak H^\bullet _\C\subset H^\bullet(\Mc_K,\C)$ by the invariant differential forms $\frac{dr}{r}_i$, $i\in \left\langle \xi-1\right\rangle$. Let $V(r)_i$ be the invariant vector fields dual to $\frac{dr}{r}_i$ and set
\begin{equation*}
V(r)_I:=V(r)_{i_1}\wedge...\wedge V(r)_{i_q},\ I=\left\{i_1,...,i_{q}\right\}\subset \left\langle \xi-1\right\rangle.
\end{equation*} 
\begin{lemma}
To prove \Cref{pol_Eis} it suffices to show
\begin{equation*}
\Eis^k_q(\varphi)(V(r)_I)\in H^\bullet_{\Eis}(\Sc,\Sym^k\Hc^\prime\otimes \mu^{\otimes n+1}\otimes \C),
\end{equation*}
for all $I\subset\left\langle \xi-1\right\rangle$ and
\begin{equation*}
\varphi(v,g)=f(v)\cdot\eta(\det(g))\left(\left\|\det(g)\right\|_f\sgn(\Norm(\det(g)))\right)^{n},
\end{equation*}
for $f\in\Sc(V(\A_f),\C)^0$ and $\eta:\prod_{\nu|\infty}\R_{>0}\backslash\I_{F}/F^{\times}\rightarrow \C$ a continuous character.
\begin{proof}
To prove \Cref{pol_Eis} we may extend coefficients to $\C$. The $\varphi$ described above generate $\Sc(V(\A_f),\mu^{\otimes n}\otimes\C)^0$. As $(V(r)_I)_{I\subset\left\langle \xi-1\right\rangle}$ is a basis of $\mathfrak H_\C ^{\bullet *}$, the lemma follows.
\end{proof}
\end{lemma}
From now on we suppose that $\varphi$ is of the form above. For $\emptyset \neq I\subset \left\langle \xi\right\rangle$ we define $\tilde{\Theta}_{I}$ as the form
\begin{equation*}
\frac{(-2i)^{-\xi(m+2)}\Norm(e^{-(m+1)i\theta})e^{-i\theta}_{I}e^{i\theta}_{I^c}}{\Norm(t_{2}) ^{m+2}\left(\left\| \det g_f\right\|_f \sgn(\Norm(\det g_\infty))^{-1}\right)^{-n}}\Norm(\omega)^m\frac{d\tau\wedge d\overline\tau}{\overline\tau-\tau}_{I^c}\wedge  d \overline\tau_{I}\otimes \vol_V^{*n+1}
\end{equation*}
and set $\widetilde{\Eis}^k _{I}(\varphi,\chi):=$
\begin{equation*}
\lim{s\rightarrow 0}\sum_{\gamma\in G(\Q)/B(\Q))}\int_{t\in Z(\A)}
\phi(\hat{f})^m _{I}(tg\gamma e^1,g)\eta(\det(g)) \chi(\det(g)t)\left\|\det(g)t\right\|^{m+2+2s} d^{\times}t\cdot \tilde{\Theta}_{I}.
\end{equation*}
\begin{lemma}
To prove \Cref{pol_Eis} it suffices to show that the differential form $\sum_{\chi \in\widehat{Cl_F ^{K}}(m)}\widetilde{\Eis}^k _{I}(\varphi,\chi)$ defines a class in $H^\bullet_{\Eis}(\Sc,\Sym^k\Hc^\prime\otimes \mu^{\otimes n+1}\otimes \C)$ for all $\emptyset\neq I\subset\left\langle \xi\right\rangle$.
\begin{proof}
$\Eis^k_q(\varphi)(V(r)_I)$ are $\C$-linear combinations of $\sum_{\chi \in\widehat{Cl_F ^{K}}(m)}\widetilde{\Eis}^k _{I}(\varphi,\chi)$.
\end{proof}
\end{lemma}
To prove this last step we have to rewrite our classes in $(\mathfrak g,K)$-cohomology, as Harder constructed his operator $\Eis_{Harder}$ in this setting.
\begin{remark}\label{g_K}
The local system $\Sym^k \Hc^\prime\otimes \mu ^{\otimes n+1}$ corresponds to the $G(\Q)$-representation $\Sym^k _\Q V(\Q)\otimes \vol_V^{*n+1}$ on which $G(\Q)$ acts by 
\begin{equation*}
\rho(g)v_1 ...v_k \vol_V^{*n+1}:= \Norm(\det(g))^{n+1}gv_1... gv_k \vol_V ^{*n+1},
\end{equation*}
for $v_1,...,v_n\in V(\Q)$, $g\in G(\Q)$ and $gv_i$ the standard action by matrix multiplication. Set 
\begin{equation*}
\mathfrak g:=Lie(G(\R))=T_1 G(\R),\ \mathfrak k:=Lie(K_\infty)=T_1K_{\infty}.
\end{equation*}
Of course, we have $\mathfrak g=\bigoplus_{\nu|\infty}Lie(G_0(F_\nu))$ and $\mathfrak k=\bigoplus_{\nu|\infty}Lie(\R_{>0}\cdot SO(2))$. Moreover, $\mathfrak g$ and $\mathfrak k$ are $F_\R$-modules. We get the isomorphism $\Psi$
\begin{equation*}
\Gamma(\Sc_K,\Omega_{\Sc_K} \otimes_\Q \Sym^n\Hc^\prime\otimes\mu^{\otimes n+1})=
\end{equation*}
\begin{equation*}
H^0(G(\Q),\Gamma(K\backslash G(\A),\Omega_{K\backslash G(\A)} \otimes_\Q \Sym^nV(\Q)\otimes \vol_V^{*n+1}))\rightarrow 
\end{equation*}
\begin{equation*}
Hom_{K_\infty}(\bigwedge ^\bullet\mathfrak g/\mathfrak k,\Cc^\infty(K_f\backslash G(\A)/G(\Q))\otimes \Sym^n V(\Q)\otimes \vol_V^{*n+1})
\end{equation*}
as described in \cite{B-W} chapter VII, proposition 2.5. Given a differential form $\eta$ on the left hand side and $X\in \bigwedge ^\bullet\mathfrak g/\mathfrak k$ this means explicitly
\begin{equation*}
\Psi(\eta)(X)(g):=\rho(g_\infty)\eta(p(g))(dp_{g}dr_{g,1} X)\text{ where}
\end{equation*}
\begin{equation*}
r_{g}:K_f\backslash G(\A)\rightarrow K_f\backslash G(\A),\ x\mapsto xg,\ \text{and }
p:K_f \backslash G(\A)\rightarrow  K\backslash G(\A)
\end{equation*}
is the canonical map. Unlike \cite{B-W} we consider right translation instead of left translation due to the fact that $G(\Q)$ acts from the right. We can make everything explicit by giving a basis for $\mathfrak g/\mathfrak k\otimes \C$. We take 
\begin{equation*}
P_{\pm}:=\frac{1}{2}\begin{pmatrix}1& \pm i\\ \pm i & -1\end{pmatrix}=\frac{1}{2}\begin{pmatrix}1 & 0\\ 0& -1 \end{pmatrix} \pm\frac{1}{2}\begin{pmatrix}0 & 1\\ 1& 0\end{pmatrix}\otimes  i
\end{equation*}
as a basis for $Lie(GL_2(\R))/Lie(\R_{>0}SO(2))_\C$, compare \cite{Ha4} p. 130. 
We have for the adjoint representation
\begin{equation*}
Ad(k)P_{\pm}=kP_{\pm} k^{-1}=(a\mp ic)^2P_{\pm}, \ k=\begin{pmatrix}a& -c\\c&a\end{pmatrix}\in SO(2).
\end{equation*}
We have $\H_{\pm}:=\C\setminus \R$ and the isomorphism
\begin{equation*}
\R^\times SO(2)\backslash GL_2(\R)\rightarrow \H_{\pm},\ g\mapsto z=ig=\frac{b+id}{a+ic}=\frac{w_2}{w_1},
\end{equation*}
\begin{equation*}
\text{ for }
g=\begin{pmatrix}a&b\\c&d\end{pmatrix},\ w_1:=a+ic,\ w_2:=b+id,
\end{equation*}
and we get
\begin{equation*}
dr_{g,1}P_+=(w_1\frac{\partial}{\partial\overline w_1}+w_2\frac{\partial}{\partial \overline w_2})_{|(w_1,w_2)=(1,i)g},\ 
dr_{g,1}P_-=(\overline w_1\frac{\partial}{\partial  w_1}+\overline w_2\frac{\partial}{\partial w_2})_{|(w_1,w_2)=(1,i)g},
\end{equation*}
\begin{equation*}
d\overline z(p(g))(dr_{g,1}P_-)=0,\  d\overline z(p(g))(dr_{g,1}P_+)=2i\frac{\det(g)}{(a-ic)^2}=2i\frac{t_2}{t_1}e^{i2\theta}\text{ and}
\end{equation*}
\begin{equation*}
dz(p(g))(dr_{g,1}P_+)=0,\  dz(p(g))(dr_{g,1}P_-)=-2i\frac{\det(g)}{(a+ic)^2}=-2i\frac{t_2}{t_1}e^{-i2\theta}
\end{equation*}
where we used the Iwasawa decomposition of $g\in GL_2(\R)$, see \Cref{Iwasawa}.
Finally, consider $\omega_0:=(e^1\overline z-e^2)\in \Cc^\infty(\H_{\pm},\C^2)$ and rewrite it using the $GL_2(\R)$-action on $\C^2$ as
\begin{equation*}
\omega_0=g^{-1}(\frac{-i\det(g)}{a-ic}(e^1-ie^2))=g^{-1}(-it_2e^{i\theta}(e^1-ie^2)).
\end{equation*}
\end{remark}
Remember that we had
\begin{equation*}
\varphi(v,g)=f(v)\cdot\eta(\det(g))\left(\left\|\det(g)\right\|_f\sgn(\Norm(\det(g)))\right)^{n},
\end{equation*}
for $f\in\Sc(V(\A_f))^0$ and $\eta:\prod_{\nu|\infty}\R_{>0}\backslash\I_{F}/F^{\times}\rightarrow \C^{\times}$ a continuous character. For $\chi\in \widehat{Cl_F ^K}(m)$
we consider the algebraic Hecke character
\begin{equation*} 
\phi:T(\A)/T(\Q)\rightarrow \C^{\times},\ (t_1,t_1)\mapsto \eta(t_1t_2)\chi(t_2)\left\|t_2\right\|^{m+n+2}\left\|t_1\right\|^n.
\end{equation*}
$\phi$ is of type $\gamma_{m,n+1}$ and $\tilde\phi_{f|Z(\R)}=1$, therefore $\phi$ is a character contributing to the cohomology of the boundary.
Write
\begin{equation*}
g=kb,\ k\in SL_2(\hat\Oc)\cdot K^1,\ b=\begin{pmatrix}t_1& t_1 x\\ 0 & t_2\end{pmatrix}\in B(\A).
\end{equation*}
Given $\eta$, $\chi$ and $f$ as above we have well-defined functions $F_{I,\infty} ^m(\chi,s,g):=$
\begin{equation*}
|\Norm(\frac{t_2}{t_1})|^s\eta_\infty(\det(g)) \Norm(-it_2e^{i\theta})^m \Norm(2i\frac{t_2}{t_1})(e^{i2\theta})_{I}\Norm(\det(g))^{n+1}\sgn(\Norm(\det(g)))^{n}
\end{equation*}
on $G(\R)$, $F_{I,f} ^m(\chi,s,g_f):=$
\begin{equation*}
\frac{\Gamma(m+2)^\xi\left\|\frac{t_2}{t_1}\right\|_f ^s\eta_f(\det(g))\int_{t \in Z(\A_f)}\hat{f}(t g_fe_1)\chi_f(\det(g_f)t)\left\|\det(g_f)t\right\|_f ^{2+m}d^{\times}t}{(-2\pi i)^{\xi(m+2)}\left(\left\|\det(g_f)\right\|_f \right)^{-n}}
\end{equation*}
on $G(\A_f)$ and $F_{I} ^m(\chi,s,\ ):=F_{I,f} ^m(\chi,s,\ )\otimes F_{I,\infty} ^m(\chi,s,\ )$ on $G(\A)$.
In other words,
\begin{equation*}
F_{I} ^m(\chi,s,g):=F_{I,f} ^m(\chi,s,g_f)\cdot F_{I,\infty} ^m(\chi,s,g_\infty).
\end{equation*}
Now we need the Harish-Chandra modules $V^* _{\phi\left\|\frac{t_2}{t_1}\right\|^s}$ as defined in \cite{Ha1} 3.3 and p. 79.
\begin{lemma}
We define
\begin{equation*}  
\omega_{I} ^m(\varphi,\chi,s)\in Hom_{K_\infty}(\bigwedge^\bullet\mathfrak g /\mathfrak k,V^* _{\phi\left\|\frac{t_2}{t_1}\right\|^s}\otimes \Sym^k_\C V(\C)\otimes \vol_V^{*n+1})
\end{equation*}
by 
\begin{equation*}
\left(P_{-}\wedge P_{+}\right)_{I^c}\wedge P_{+ I}\mapsto F_{I} ^m(\chi,s,\ )\Norm(e_1-e_2 i)^m \vol_V^{*n+1}
\end{equation*}
and all elements not collinear with $\left(P_{-}\wedge P_{+}\right)_{I^c}\wedge P_{+ I}$ mapping to zero.
\begin{proof}
We have to show that $\omega_{I} ^m(\varphi,\chi,s)$ actually defines a form as claimed. For $b=\begin{pmatrix}t_1& t_1 x\\ 0 & t_2\end{pmatrix}\in B(\A)$ and $g\in G(\A)$ we have
\begin{equation*}
F_{I} ^m(\chi,s,gb)=\left\|\frac{t_2}{t_1}\right\|^s \phi(t_1,t_2)F_{I} ^m(\chi,s,g).
\end{equation*}
$K_f$ and $K_\infty$-finiteness of $F_{I} ^m(\chi,s, )$ are clear and therefore $F_{I} ^m(\chi,s,\ )\in V^* _{\phi\left\|\frac{t_2}{t_1}\right\|^s}$. So what remains to be shown is that 
\begin{equation*}
\left(P_{-}\wedge P_{+}\right)_{I^c}\wedge P_{+ I}\mapsto F_{I} ^m(\chi,s,\ )\Norm(e_1-e_2 i)^m \vol_V^{*n+1}
\end{equation*}
is $K_\infty$-equivariant.
Consider 
\begin{equation*}
k=t\begin{pmatrix}x& -y\\ y& x\end{pmatrix}\in K_\infty\text{ with }t\in Z(\R)^0\text{ and }\begin{pmatrix}x& -y\\ y& x\end{pmatrix}_\nu=\begin{pmatrix}\cos\theta_\nu & -\sin \theta_\nu \\ \sin\theta_\nu & \cos\theta_\nu\end{pmatrix}.
\end{equation*}
We have
\begin{equation*}
Ad(k^{-1})\left(P_{-}\wedge P_{+}\right)_{I^c}\wedge P_{+ I}=Ad\begin{pmatrix}x& -y\\ y& x\end{pmatrix}^{-1}\left(P_{-}\wedge P_{+}\right)_{I^c}\wedge P_{+ I}=
\end{equation*}
\begin{equation*}
e^{2i\theta}_{I}\left(P_{-}\wedge P_{+}\right)_{I^c}\wedge P_{+ I}
\end{equation*}
and on the other hand
\begin{equation*}  
k^{-1}\omega_{I} ^m(\varphi,\chi,s)(\left(P_{-}\wedge P_{+}\right)_{I^c}\wedge P_{+ I})(kg)=k^{-1}F_{I} ^m(\chi,s,kg )\Norm(e_1-e_2 i)^m \vol_V^{* n+1}=
\end{equation*}
\begin{equation*}
\Norm(t^{-1})^{2(n+1)}\Norm(t^{-1}e^{-i\theta})^m F_{I} ^m(\chi,s,kg )\Norm(e_1-e_2 i)^m\vol_V^{* n+1}=
\end{equation*}
\begin{equation*}
\Norm(t^{-1})^{2(n+1)}\Norm(t^{-1}e^{-i\theta})^mF_{I,f} ^m(\chi,s,g_f )F_{I,\infty} ^m(s,kg_\infty ) \Norm(e_1-e_2 i)^m\vol_V^{* n+1}=
\end{equation*}
\begin{equation*}
e^{2i\theta} _{I}F_{I} ^m(\chi,s,g )\Norm(e_1-e_2 i)^m\vol_V^{* n+1}
\end{equation*}
by the definition of $F_{I,\infty} ^m(\chi,s,kg )$ and it follows that $\omega_{I} ^m(\varphi,\chi,s)$ is (right)-$K_\infty$-equivariant. 
\end{proof}
\end{lemma}
Harder defines
\begin{equation*}  
\Eis^*(\omega_{I} ^m(\varphi,\chi,s))\in Hom_{K_\infty}(\bigwedge^\bullet\mathfrak g /\mathfrak k,\Cc^\infty(K_f\backslash G(\A)/G(\Q))\otimes \Sym^k_\C V(\C)\otimes \vol_V^{*n+1})
\end{equation*} by
\begin{equation*}
\left(P_{-}\wedge P_{+}\right)_{I^c}\wedge P_{+ I}\mapsto\sum_{\gamma\in G(\Q)/B(\Q)} F_{I} ^m(\chi,s,g\gamma)\Norm(e_1-e_2 i)^m \vol_V^{*n+1}
\end{equation*}
and all elements not collinear with $\left(P_{-}\wedge P_{+}\right)_{I^c}\wedge P_{+ I}$ mapping to zero, see \cite{Ha1} (4.2.2). 
\begin{proposition}\label{comparison}
\begin{equation*}  
\sum_{\chi \in\widehat{Cl_F ^{K}}(m)}\lim{s\rightarrow 0}\Eis^*(\omega_{I} ^m(\varphi,\chi,s))=\Eis_{Harder}(\rho_{m,n}^0(\varphi)\Norm(\omega)^m\frac{d\tau\wedge d\overline \tau}{\overline\tau-\tau}_{I^c}\wedge d\overline\tau_{I}\otimes \vol_V^{*n+1}),
\end{equation*}
\begin{equation*}  
\lim{s\rightarrow 0}\Eis^*(\omega_{I} ^m(\varphi,\chi,s))=\Psi(\widetilde{\Eis}^k_{I}(\varphi,\chi))
\end{equation*}
and \Cref{pol_Eis} is proved.
\begin{proof}
We have
\begin{equation*}
\rho_{m,n}^0(\varphi)\Norm(\omega)^m\frac{d\tau\wedge d\overline \tau}{\overline\tau-\tau}_{I^c}\wedge d\overline\tau_{I}\otimes \vol_V^{*n+1}=
\end{equation*}
\begin{equation*}
\sum_{\chi \in\widehat{Cl_F ^{K}}(m)}F_{I,f} ^m(\chi,0,g)\eta_\infty(\det(g))\sgn(\Norm(\det(g)))^n \Norm(\omega)^m\frac{d\tau\wedge d\overline\tau}{\overline\tau-\tau}_{I^c}\wedge  d \overline\tau_I\otimes \vol_V^{*n+1}.
\end{equation*}
Following \Cref{g_K} and \cite{Ha1} p.79 and p.80 we see that 
\begin{equation*}  
\sum_{\chi \in\widehat{Cl_F ^{K}}(m)}\lim{s\rightarrow 0}\Eis^*(\omega_{I} ^m(\varphi,\chi,s))=\Eis_{Harder}(\rho_{m,n}^0(\varphi)\Norm(\omega)^m\frac{d\tau\wedge d\overline\tau}{\overline\tau-\tau}_{I^c}\wedge d\overline\tau_{I}\otimes \vol_V^{*n+1})
\end{equation*}
by the very definition of Harder. The right hand side of the equation is defined, as we know by \Cref{im_res_Eis} that
\begin{equation*}  
\rho_{m,n}^0(\varphi)\Norm(\omega)^m\frac{d\tau\wedge d\overline \tau}{\overline\tau-\tau}_{I^c}\wedge d\overline\tau_{I}\otimes \vol_V^{*n+1}\in \image(\res_\Sc).
\end{equation*}
So what remains to be shown is $\lim{s\rightarrow 0}\Eis^*(\omega_{I} ^m(\varphi,\chi,s))=\Psi(\widetilde{\Eis}^k_{I}(\varphi,\chi))$. We start with the calculation of $\Psi(\widetilde{\Eis}^k_{I}(\varphi,\chi))$. Using \Cref{g_K} we see that it is determined by
\begin{equation*}
\lim{s\rightarrow 0}\sum_{\gamma\in G(\Q)/B(\Q))}\int_{t\in Z(\A)}
\phi(\hat{f})^m _{I}(tg\gamma e^1)\eta(\det(g)) \chi(\det(g)t)\left\|\det(g)t\right\|^{m+2+2s} d^{\times}t=
\end{equation*}
\begin{equation*}
\lim{s\rightarrow 0}\sum_{\gamma\in G(\Q)/B(\Q))}\int_{t\in Z(\A)}
\phi(\hat{f})^m _{I}(tg\gamma e^1)\eta(\det(g\gamma)) \chi(\det(g\gamma)t)\left\|\det(g\gamma)t\right\|^{m+2+2s} d^{\times}t=
\end{equation*}
\begin{equation*}
\lim{s\rightarrow 0}\sum_{\gamma\in G(\Q)/B(\Q))}\int_{t\in Z(\A_f)}
\hat{f}(tg_f\gamma e^1)\eta_f(\det(g_f\gamma)) \chi_f(\det(g_f\gamma)t)\left\|\det(g_f\gamma)t\right\|_f^{m+2+2s} d^{\times}t\cdot
\end{equation*}
\begin{equation*}
\int_{t\in Z(\R)}
\phi(\hat{f})^m_{I,\infty}(tg_\infty\gamma e^1)\eta_\infty (\det(g_\infty\gamma)) \chi_\infty(\det(g_\infty\gamma)t)|\Norm(\det(g_\infty\gamma)t)|^{m+2+2s} d^{\times}t
\end{equation*}
Let us set 
\begin{equation*}
F_\infty(s,g):=\int_{t\in Z(\R)}\phi(\hat{f})^m_{I,\infty}(tg_\infty e^1)\eta_\infty (\det(g_\infty)) \chi_\infty(\det(g_\infty)t)|\Norm(\det(g_\infty)t)|^{m+2+2s} d^{\times}t,
\end{equation*}
\begin{equation*}
F_f(s,g):=\int_{t\in Z(\A_f)}\hat{f}(tg_f e^1)\eta_f(\det(g_f)) \chi_f(\det(g_f)t)\left\|\det(g_f)t\right\|_f^{m+2+2s} d^{\times}t
\end{equation*}
for the moment. Consider 
\begin{equation*}
G(g_f):=(g_f^{-1}K_fg_f\cdot G(\R) ^0)\cap G(\Q),\ B(g_f)=:G(g_f)\cap B(\A).
\end{equation*}
$G(g_f)$ acts on $G(\Q)/B(\Q)$ by left translation and we decompose the latter into $G(g_f)$ orbits. Now we may write our sum above as
\begin{equation*}
\lim{s\rightarrow 0}\sum_{\alpha\in G(g_f)\backslash G(\Q)/B(\Q)}\sum_{\gamma\in G(g_f\alpha)/B(g_f\alpha)}F_f(s,g\alpha\gamma)F_\infty(s,g\alpha\gamma)
\end{equation*}
$G(g_f)\backslash G(\Q)/B(\Q)$ is finite, as $G(\Z)\backslash G(\Q)/B(\Q)\cong Cl_F$ is the class group of $F$ and the groups $G(g_f)$ and $G(\Z)$ are commensurable. Moreover, $G(g_f)$ does not affect $F_f(s,\ )$: By construction of $\Eis^k$ we have chosen the level $K_f$ for $f$ and $\eta$ such that $f(kv)=f(v)$ and $\eta(\det(k))=1$ for all $k\in K_f$ and $v\in V(\A_f)$, see \Cref{Eis^k}. Therefore we may write our sum as
\begin{equation*}
\sum_{\alpha\in G(g_f)\backslash G(\Q)/B(\Q)}F_f(0,g\alpha)
\lim{s\rightarrow 0} \sum_{\gamma\in G(g_f \alpha)/B(g_f \alpha)}F_\infty (s,g\alpha\gamma)
\end{equation*}
Let us take care of $F_\infty$. We write 
\begin{equation*}
g_\infty=kb=k t_1\begin{pmatrix}1& x \\ 0& y\end{pmatrix}\text{ and }\alpha\gamma =\begin{pmatrix}
a& b\\ c& d
\end{pmatrix}.
\end{equation*}
We get
\begin{equation*}
F_\infty(s,g\alpha\gamma)=\Norm(e^{i\theta})^{m+1}e^{i\theta}_{ I}e^{-i\theta}_{I^c}\eta_\infty (\det(g_\infty \alpha\gamma)) \Norm(t_2)^{m+2} |\Norm(t_2)|^{2s}\cdot
\end{equation*}
\begin{equation*}
\frac{\Norm(\det(\alpha\gamma))^{m+2}|\Norm(\det(\alpha\gamma))|^{2s}\Norm(a+c\tau)^{m+1}(a+c\tau)_{I}(a+c\overline\tau)_{I^c}\Gamma(m+2+s)^\xi}{\pi^{\xi(m+2)+{\xi s}}|\Norm(a+c\tau)|^{2(m+2+s)}}.
\end{equation*}
Using $|\Norm(\det(\gamma))|=1$ we calculate
\begin{equation*}
\lim{s\rightarrow 0} \sum_{\gamma\in G(g_f \alpha)/B(g_f \alpha)}F_\infty(s,g\alpha\gamma)=
\frac{\Gamma(m+2)^\xi}{\pi^{\xi(m+2)}}\Norm(e^{i\theta})^{m+1}e^{i\theta}_{I}e^{-i\theta}_{I^c}\Norm(t_2)^{m+2}\cdot
\end{equation*}
\begin{equation*}
\lim{s\rightarrow 0} \sum_{\gamma\in G(g_f \alpha)/B(g_f \alpha)}\frac{\eta_\infty (\det(g \alpha\gamma))\Norm(\det(\alpha\gamma ))^{m+2} \Norm(a+c\tau)^{m+1}(a+c\tau)_{I}(a+c\overline\tau)_{I^c}}{|\Norm(a+c\tau)|^{2(m+2+s)}}=
\end{equation*}
\begin{equation*}
\lim{s\rightarrow 0} \sum_{\gamma\in G(g_f \alpha)/B(g_f \alpha)}\frac{\eta_\infty (\det(\alpha\gamma))|\Norm(\det(\alpha\gamma))|^{s}\Norm(\det(\alpha\gamma ))^{m+2} \Norm(a+c\tau)^{m+1}(a+c\tau)_{I}(a+c\overline\tau)_{I^c}}{|\Norm(a+c\tau)|^{2(m+2+s)}}\cdot
\end{equation*}
\begin{equation*}
\frac{\Gamma(m+2)^\xi}{\pi^{\xi(m+2)}}\eta_\infty(\det(g))\Norm(e^{i\theta})^{m+1}e^{i\theta}_{I}e^{-i\theta}_{I^c}\Norm(t_2)^{m+2}|\Norm(\frac{t_2}{t_1})|^s
\end{equation*}
Looking carefully at the formulas in \Cref{g_K} and recalling 
\begin{equation*}
\left\|\det(g_f\alpha)\right\|^{-n}_f|\Norm(\det(\alpha\gamma))|^{-n}=\left\|\det(g_f)\right\|_f^{-n}
\end{equation*}
we may conclude that $\Psi(\widetilde{\Eis}_{I}^k(\varphi,\chi))$ is the form that maps
$\left(P_{-}\wedge P_{+}\right)_{I^c}\wedge P_{+ I}$ to
\begin{equation*}
\lim{s\rightarrow 0}\sum_{\alpha\in G(g_f)\backslash G(\Q)/B(\Q)}F_{I,f}^m(\chi,0,g_f \alpha)
\end{equation*}
\begin{equation*}
\sum_{\gamma\in G(g_f \alpha)/B(g_f \alpha)}\frac{\eta_\infty (\det(\alpha\gamma))|\Norm(\det(\alpha\gamma))|^s \Norm(\det(\alpha\gamma ))^{m+2} \Norm(a+c\tau)^{m+1}(a+c\tau)_{I}(a+c\overline\tau)_{I^c}}{|\Norm(\det(\alpha\gamma))|^{-n}|\Norm(a+c\tau)|^{2(m+2+s)}}\cdot
\end{equation*}
\begin{equation*}
\sgn(\Norm(\det(g)))^{n}\eta_\infty(\det(g))|\Norm(\frac{t_2}{t_1})|^s \Norm(-it_2e^{i\theta}(e_1-ie_2))^m \Norm(2i\frac{t_2}{t_1})(e^{i2\theta})_{I}\Norm(\det(g))^{n+1} \vol_V^{*n+1}
\end{equation*}
and all the elements not collinear with $\left(P_{-}\wedge P_{+}\right)_{I^c}\wedge P_{+ I}$ to zero.
The last sum is nothing else but
\begin{equation*}
\lim{s\rightarrow 0}\sum_{\gamma \in G(\Q)/B(\Q)}F_{I,f}^m(\chi,s,g_f \gamma)F_{I,\infty}^m(\chi,s,g)\Norm(e_1-e_2 i)^m \vol_V^{*n+1}\cdot
\end{equation*}
\begin{equation*}
\frac{\eta_\infty (\det(\alpha\gamma))\Norm(\det(\alpha\gamma ))^{m+2}|\Norm(\det(\alpha\gamma))|^s \Norm(a+c\tau)^{m+1}(a+c\tau)_{I}(a+c\overline\tau)_{I^c}}{|\Norm(\det(\alpha\gamma))|^{-n}|\Norm(a+c\tau)|^{2(m+2+s)}}
\end{equation*}
Finally, 
\begin{equation*}
t_2(g\alpha\gamma)=\frac{\det(\alpha\gamma)}{|a+c \tau|}t_2(g)\text{ and } e^{i\theta(g\alpha\gamma)}=\frac{a+c \tau}{|a+c \tau|}e^{i\theta(g)}
\end{equation*}
and therefore $F_{I,\infty}^m(\chi,s,g \alpha\gamma)$ equals
\begin{equation*}
F_{I,\infty}^m(\chi,s,g)\cdot
\frac{\eta_\infty (\det(\alpha\gamma))\Norm(\det(\alpha\gamma ))^{m+2}|\Norm(\det(\alpha\gamma))|^s \Norm(a+c\tau)^{m+1}(a+c\tau)_{I}(a+c\overline\tau)_{I^c}}{|\Norm(\det(\gamma))|^{-n}|\Norm(a+c\tau)|^{2(m+2+s)}}
\end{equation*}
Consequently $\Psi(\widetilde{\Eis}_{I}^k(\varphi,\chi))$ is the form that maps
\begin{equation*}
\left(P_{-}\wedge P_{+}\right)_{I^c}\wedge P_{+ I}\mapsto
\lim{s\rightarrow 0}\sum_{\gamma \in G(\Q)/B(\Q)}F_{I}^m(\chi,s,g \gamma)\Norm(e_1-e_2 i)^m \vol_V^{*n+1}.
\end{equation*}
This form is exactly $\lim{s\rightarrow 0}\Eis^*(\omega_{I}^m(\varphi,\chi,s))$.
\end{proof}
\end{proposition}

\section{The image of the polylogarithmic Eisenstein operator}

The polylogarithmic Eisenstein classes are built up by Eisenstein series associated to Schwartz-Bruhat functions on $V(\A)$, while Harder's Eisenstein series are associated to functions on $G(\A)$, which are induced from algebraic Hecke characters on $T(\A)/T(\Q)$. Eisenstein series associated to induced functions on $G(\A)$ can be represented as finite sums of Eisenstein series coming from Schwartz-Bruhat functions on $V(\A)$, see \cite{J-Z} Lemma. So we have a good reason to believe that the image of the polylogarithmic Eisenstein operator is quite big. In this section we determine $\image(\Eis^k_q)\subset H^\bullet_{\Eis}(\Sc,\Sym^k\Hc^\prime\otimes \mu^{\otimes n+1})$ completely. As
 \begin{equation*}
\res_\Sc:H^\bullet_{\Eis}(\Sc,\Sym^k\Hc^\prime\otimes \mu^{\otimes n+1})\rightarrow H^\bullet(\partial \Sc,\Sym^k\Hc^\prime\otimes  \mu^{\otimes n+1})
\end{equation*}
is an isomorphism of $G(\A_f)\times \pi_0(G(\R))$-modules onto its image, it suffices to examine the image of $\res_\Sc\circ \Eis^k_q$.
This last operator has already been calculated in \Cref{boundary_residueII} and we know by \Cref{hor} that it suffices to understand the image of the horospherical map $\rho_{m,n}^0$. As we may multiply our functions by
\begin{equation*}
\left\|\det(g)\right\|_f ^{-n}\sgn(\Norm(\det(g)))^{-n},\ n\in \Z,
\end{equation*}
it even suffices to understand the image of $\rho_{m,0}^0$. But first we want to understand a more general map.
\begin{definition}
Define the \textit{horospherical map} for $m\geq 0$
\begin{equation*}
\rho_{m}:\Sc(V(\A_f),\mu^{\otimes 0}\otimes \C)\rightarrow \Ind_{B(\A_f)\times \pi_0(B(\R))}^{G(\A_f)\times \pi_0(G(\R))}\bigoplus_{\phi:\ type(\phi)=\gamma_{m,1},\tilde{\phi}_{f|Z(\R)=1}}\C \cdot\tilde\phi_f
\end{equation*}
by 
\begin{equation*}
\rho_{m}(f)(g):= 
\sum_{\chi\in \widehat{Cl_F ^K}(m)}\frac{\Gamma(m+2)^{\xi}}{(-2\pi i)^{\xi(m+2)}}\int_{t\in Z(\A_f)}\hat{f} (tg_fe^1,g) \chi(\det(g_f)t)\left\|\det(g_f)t\right\|^{m+2} d^{\times}t ,
\end{equation*}
if $f$ is $K_f$-invariant.
\end{definition}
\subsection{Surjectivity of the horospherical map}

 \begin{proposition}\label{hor_surjective}
The horospherical map $\rho_m$ is surjective.

\begin{proof}
The idea how to prove this proposition is already in \cite{J-Z} Lemma. Nevertheless, we give the proof for the sake of completeness. Set for $N\in \N$
\begin{equation*}
U_N:=\ker(\hat\Oc^{\times}\rightarrow (\Oc/N\Oc) ^{\times}),\ Cl_F ^{(N)}:=U_N \prod_{\nu|\infty}\R_{>0}\backslash \I_F/F^{\times}.
\end{equation*}
Let $\phi$ be an algebraic Hecke character as above. We may write
\begin{equation*}
\phi(t_1,t_2)=\eta(t_1t_2)\chi^\prime(t_2)\left\|t_2\right\|^{m+2}, (t_1,t_2)\in T(\A),
\end{equation*}
for characters $\eta,\chi^\prime \in \widehat{Cl_F ^{(N)}}$ for some $N\in \N$.
Then
\begin{equation*}
\tilde\phi_f(t_1,t_2)=\eta(t_1t_2)\chi^\prime(t_2)\left\|t_2\right\|_f^{m+2}\sgn(\Norm(t_2))^{m+2}, (t_1,t_2)\in T(\A).
\end{equation*}
$\tilde\phi_{f|Z(\R)}=1$ means $\chi^\prime(t)=\sgn(\Norm(t))^{m+2},\ t\in F_\R ^{\times}$, in other words $\chi\in Cl_F ^{(N)}(m):=Cl_F ^{K_N}(m)$. Now take 
\begin{equation*}
\psi\in \Ind_{B(\A_f)\times \pi_0(B(\R))}^{G(\A_f)\times \pi_0(G(\R))}\C \cdot\tilde\phi_f.
\end{equation*}
By choosing $N$ big enough we may also assume that $\psi(kx)=\psi(x)$ holds for all $x\in G(\A_f)\times \pi_0(G(\R))$ and $k\in K_N$.
Define $\overline \psi$ by 
\begin{equation*}
\overline \psi(x)=\eta(\det(x))\psi(\det(x)^{-1}\cdot x),\ x\in G(\A_f)\times \pi_0(G(\R)).
\end{equation*}
As $\det(K_N)=U_N$, we still have $\overline\psi(kx)=\overline\psi(x)$ for all $x\in G(\A_f)\times \pi_0(G(\R))$ and $k\in K_N$, but now we have for $b=\begin{pmatrix}t_1&t_1 u\\ 0& t_2\end{pmatrix}\in B(\A_f)\times \pi_0(B(\R))$
\begin{equation*}
\overline\psi(xb)=\overline\psi(x)\chi^\prime(t_1^{-1})\left\|t_1 ^{-1}\right\|_f^{m+2}\sgn(\Norm(t_1^{-1}))^{m+2}=
\overline\psi(x)\chi^\prime_f(t_1)^{-1}\left\|t_1 \right\|_f^{-(m+2)}.
\end{equation*}
This means that $\overline \psi$ is a function
\begin{equation*}
K_N\backslash G(\A_f)\times \pi_0(G(\R))/P(\A_f)\times \pi_0(P_\infty)\rightarrow \C.
\end{equation*}
We restrict $\overline \psi $ to $G(\hat\Z)$ and obtain
\begin{equation*}
\overline\psi:K_N\backslash G(\hat \Z)/P(\hat\Z)=G(\Z/N\Z)/P(\Z/NZ)\rightarrow \C.
\end{equation*}
We consider the embedding $i:G(\A_f)/P(\A_f)\hookrightarrow V(\A_f),\ x\mapsto xe^1$, and make $\overline \psi $ to a function on $V(\A_f)$ by extending it by zero outside $G( \A_f)/P(\A_f)$. When we restrict $\overline\psi$ to $V(\hat\Z)$, it even induces a function on $V(\Z/N\Z)$ and may also be interpreted as a Schwartz-Bruhat function $s\overline\psi\in\Sc(V(\A_f),\C)$. Explicitly $s\overline\psi$ can be described as follows: Take a complete set of representatives $x\in G(\hat\Z)$ of $G( \Z/N\Z)/P(\Z/N\Z)$. Then
\begin{equation*}
s\overline\psi:=\sum_{x\,\text{mod}N\in G( \Z/N\Z)/P(\Z/N\Z)}\overline \psi(x)\chi_{xe^1+NV(\hat\Z)}=\otimes_{\nu\nmid \infty}s\overline\psi_\nu
\end{equation*}
where the local functions $s\overline\psi_\nu:F_\nu\to\C$ are defined by $s\overline\psi_\nu:=\chi_{\Oc_\nu ^2}$, if $\nu\nmid N$, and 
\begin{equation*}
 s\overline\psi_\nu:=\sum_{x_\nu\in G_0(\Oc_\nu/N\Oc_\nu)/P_0(\Oc_\nu/N\Oc_\nu)}
\overline\psi(x_\nu)\chi_{x_\nu e^1+N\Oc_\nu^2}=\overline\psi\cdot \chi_{G_0(\Oc_\nu)/P_0(\Oc_\nu)},\text{ if }\nu| N.
\end{equation*}
Take a full set of representatives $u\in \I_{F}$ of $Cl_F ^{(N)}$. We define $\varphi\in \Sc(V(\A_f),\mu^{\otimes 0}\otimes \C)$ by its Fourier transform with respect to $v$
\begin{equation*}
\hat\varphi(v,g):=\sum_u \chi^\prime_f(u)\left\|u\right\|_f ^{m+2}\eta(\det(g))s\overline\psi(uv).
\end{equation*}
Let us calculate $\rho_m(\varphi)$. We have by definition 
\begin{equation*}
\rho_m(\varphi)(g)=
\sum_{\chi\in \widehat{Cl_F ^{(N)}}(m)}\frac{\Gamma(m+2)^{\xi}}{(-2\pi i)^{\xi(m+2)}}
\int_{t\in Z(\A_f)}\hat{\varphi} (tg_fe^1,g) \chi(\det(g_f)t)\left\|\det(g_f)t\right\|^{m+2} d^{\times}t =
\end{equation*}
\begin{equation*}
\sum_u\sum_{\chi\in \widehat{Cl_F ^{(N)}}(m)}\chi^\prime_f(u)\left\|u\right\|_f ^{m+2}\eta(\det(g))\frac{\Gamma(m+2)^{\xi}}{(-2\pi i)^{\xi(m+2)}}\cdot
\end{equation*}
\begin{equation*}
\int_{t\in Z(\A_f)}s\overline\psi(utg_fe^1) \chi(\det(g_f)t)\left\|\det(g_f)t\right\|^{m+2} d^{\times}t. 
\end{equation*}
So we need to understand the last integral. Analogous to $s\overline \psi$ we define the Schwartz-Bruhat function $s\chi^{\prime-1}\in \Sc(\A_{F,f},\C)$. Extend $\chi^{\prime-1}$ by zero to a function on $\A_{F,f}$ and choose again a full set of representatives $x\in \hat\Oc^{\times}$ of $(\Oc/N\Oc)^{\times}=U_N\backslash\hat\Oc^{\times}$. Then
\begin{equation*}
s\chi^{\prime-1}:=\sum_{x\in (\Oc/N\Oc)^{\times}} \chi^{\prime-1}(x)\chi_{x+N\hat\Oc}=\otimes_{\nu\mid\infty}s\chi_\nu^{\prime-1}
\end{equation*}
and the local functions on $F_\nu$ are defined by
\begin{equation*}
s\chi_\nu^{\prime-1}:=\chi_{\Oc_\nu}, \text{ if }\nu\nmid N,\text{ and } s\chi_\nu^{\prime-1}:=\sum_{x_\nu\in (\Oc_\nu/N\Oc_\nu)^{\times}}\chi^{\prime -1}(x_\nu)\chi_{x_\nu+N\Oc_\nu}=\chi^{\prime-1}\cdot\chi_{\Oc_\nu ^{\times}}, \text{ if }\nu|N.
\end{equation*}
Write $g_f=xb$ for $x\in G(\hat\Z)$ and $b\in B(\A_f)$ and see 
\begin{equation*}
s\overline \psi(tug_fe_1)=s\overline\psi(xtut_1e^1)=\prod_{\nu\nmid N}\chi_{\Oc_\nu^2}((xtut_1e^1)_\nu)\prod_{\nu|N} \overline\psi((xtut_1)_\nu)\chi_{G_0(\Oc_\nu)/P_0(\Oc_\nu)}((xtut_1e^1)_\nu).
\end{equation*}
As $x\in G(\hat\Z)$, we may write this as 
\begin{equation*}
\prod_{\nu\nmid N}\chi_{\Oc_\nu}((tut_1)_\nu)\prod_{\nu|N} \overline\psi((xtut_1)_\nu)\chi_{\Oc_\nu ^{\times}}((tut_1)_\nu).
\end{equation*}
Now 
\begin{equation*}
\overline\psi((xtut_1)_\nu)=\left\|(tut_1)_\nu\right\|_f^{-(m+2)}\chi^{\prime-1}((tut_1)_\nu)\overline\psi(x_\nu)=\chi^{\prime-1}((tut_1)_\nu)\overline\psi(x_\nu) \text{ for } (tut_1)_\nu\in \Oc_\nu ^{\times}
\end{equation*}
and we get using $\overline\psi_{|K_N}=1$ 
\begin{equation*}
\prod_{\nu\nmid N}\chi_{\Oc_\nu}((tut_1)_\nu)\prod_{\nu|N} \overline\psi(x_\nu)\chi^{\prime-1}((tut_1)_\nu)\chi_{\Oc_\nu ^{\times}}((tut_1)_\nu)=
\end{equation*}
\begin{equation*}
\overline\psi(x)\prod_{\nu\nmid N}\chi_{\Oc_\nu}((tut_1)_\nu)\prod_{\nu|N} \chi^{\prime-1}((tut_1)_\nu)\chi_{\Oc_\nu ^{\times}}((tut_1)_\nu)=\overline\psi(x)s\chi^{\prime-1}(tut_1).
\end{equation*}
We may write our integral as
\begin{equation*}
\int_{t\in Z(\A_f)}s\overline\psi(utg_fe^1) \chi(\det(g_f)t)\left\|\det(g_f)t\right\|^{m+2} d^{\times}t= 
\end{equation*}
\begin{equation*}
\overline\psi (x)\int_{t\in Z(\A_f)}s\chi^{\prime-1} (utt_1) \chi(\det(g_f)t)\left\|\det(g_f)t\right\|^{m+2} d^{\times}t= 
\end{equation*}
\begin{equation*}
\overline\psi (x)\chi_f (u t_1)^{-1}\left\|ut_1\right\|_f^{-(m+2)}\chi_f(\det(g_f))\left\|\det(g_f)\right\|_f^{m+2}\int_{t\in Z(\A_f)}s\chi^{\prime-1} (t) \chi(t)\left\|t\right\|^{m+2} d^{\times}t
\end{equation*}
Remembering $\overline\psi(g)=\overline\psi (x)\chi^\prime _f(t_1)^{-1}\left\|t_1\right\|_f ^{-(m+2)}$ and 
\begin{equation*}
\overline\psi(g)=\eta(\det(g))\psi(g\det(g)^{-1})=\psi(g)\eta(\det(g))^{-1}\chi_f^\prime(\det(g)^{-1})\left\|\det(g)\right\|_f ^{-(m+2)}
\end{equation*}
we see now that $\frac{(-2\pi i)^{\xi(m+2)}}{\Gamma(m+2)^{\xi}}\rho_m(\varphi)(g)$ equals
\begin{equation*}
\psi (g)\sum_u\sum_{\chi\in \widehat{Cl_F ^{(N)}}(m)}\chi^\prime_f\cdot\chi_f^{-1}(u t_1 \det(g)^{-1})\cdot\int_{t\in Z(\A_f)}s\chi^{\prime-1} (t) \chi(t)\left\|t\right\|^{m+2} d^{\times}t=
\end{equation*}
\begin{equation*}
\psi (g)\sum_u\sum_{\chi\in \widehat{Cl_F ^{(N)}}(m)}\chi^\prime\cdot\chi^{-1}(u t_1 \det(g)^{-1})\cdot\int_{t\in Z(\A_f)}s\chi^{\prime-1} (t) \chi(t)\left\|t\right\|^{m+2} d^{\times}t.
\end{equation*}
Of course, $\sum_u\chi^\prime\cdot\chi^{-1}(u t_1 \det(g)^{-1})=0$, whenever $\chi\neq \chi^\prime$, and therefore
\begin{equation*}
\rho_m(\varphi)=\frac{|Cl_F^{(N)}|\Gamma(m+2)^{\xi}}{(-2\pi i)^{\xi(m+2)}}\int_{t\in Z(\A_f)}s\chi^{\prime-1} (t) \chi^\prime(t)\left\|t\right\|^{m+2} d^{\times}t\cdot\psi.
\end{equation*}
To complete our proof we need to show that
\begin{equation*}
\int_{t\in Z(\A_f)}s\chi^{\prime-1} (t) \chi^\prime(t)\left\|t\right\|^{m+2} d^{\times}t\neq 0.
\end{equation*}
But following \cite{T} p.342 f.f. and using $s\chi^{\prime-1}=\otimes_{\nu\nmid \infty}s\chi_\nu^{\prime-1}$ we calculate
\begin{equation*}
\int_{t\in Z(\A_f)}s\chi^{\prime-1} (t) \chi^\prime(t)\left\|t\right\|^{m+2} d^{\times}t=\frac{1}{\sqrt{d_F}}\prod_{\nu\nmid N\infty}\frac{1}{1-\chi^\prime(\pi_\nu)\Norm(\mathfrak p_\nu)^{-(m+2)}}
\end{equation*}
an convergent Euler product and therefore non-zero.
\end{proof}
\end{proposition}
\begin{remark}\label{preimage}
We keep the notations from above.
Let us set
\begin{equation*}
\Lambda_N(\chi,m+2)=\frac{|Cl_f ^{(N)}|\Gamma(m+2)^{\xi}}{(-2\pi i)^{\xi(m+2)}}\int_{t\in Z(\A_f)}s\chi^{\prime-1} (t) \chi^\prime(t)\left\|t\right\|^{m+2} d^{\times}t
\end{equation*}
For \textit{any} set of representatives $\left\{u\right\}$ the function $\varphi_{\left\{u\right\}}\in \Sc(V(\A_f),\mu^{\otimes 0}\otimes\C)$ defined by its Fourier transform
\begin{equation*}
\hat\varphi_{\left\{u\right\}}(v,g):=\Lambda_N(\chi,m+2)^{-1}\sum_u \chi^\prime_f(u)\left\|u\right\|_f ^{m+2}\eta(\det(g))s\overline\psi(uv)
\end{equation*}
is a preimage of $\psi$ under $\rho_m$.
\end{remark}
Let us set $T^1:=\ker(\det:T\rightarrow Res_{\Oc/\Z}\mathbb G_m)$.
\begin{lemma}
\begin{equation*}
\Ind_{B(\A_f)\times \pi_0(B(\R))}^{G(\A_f)\times \pi_0(G(\R))}\bigoplus_{\phi:\ type(\phi)=\gamma_{m,1},\tilde{\phi}_{f|Z(\R)=1}}\C \cdot\tilde\phi_f
\end{equation*}
is in the image of $\rho_{m,0} ^0$, if $\phi_{|T^1(\A)}\neq \left\|\ \cdot\ \right\|^2$.
\begin{proof}
Write $\phi(t_1,t_2)=\eta(t_1t_2)\chi(t_2)\left\|t_2\right\|^{m+2}$. The assumption $\phi_{|T^1(\A)}\neq \left\|\ \cdot\ \right\|^2$ tells us that $m\geq1$ or $\chi\neq1$ and if $m=0$ we know $\chi_{|F_\R^{\times}}=1$. So there is an $u_0\in \I_{F}$ such that $\chi_f(u_0)\left\|u_0\right\|_f^m\neq 1$. We calculate
\begin{equation*}
\chi_f(u_0)\left\|u_0\right\|_f^m\int_{v\in V(\A_f)}\hat\varphi_{\left\{u\right\}}(v,g)dv=
\end{equation*}
\begin{equation*}
\Lambda_N(\chi,m+2)^{-1}\sum_u \chi_f(u_0u)\left\|uu_0\right\|_f ^{m}\left\|u\right\|_f ^{2}\eta(\det(g))\int_{v\in V(\A_f)}s\overline\psi(uv)dv
\end{equation*}
\begin{equation*}
\Lambda_N(\chi,m+2)^{-1}\sum_u \chi_f(u_0u)\left\|u_0u\right\|_f ^{m}\left\|u\right\|_f ^{2}\eta(\det(g))\int_{v\in V(\A_f)}s\overline\psi(u_0u(u_0^{-1}v))dv
\end{equation*}
\begin{equation*}
\Lambda_N(\chi,m+2)^{-1}\sum_u \chi_f(u_0u)\left\|u_0u\right\|_f ^{m}\left\|uu_0\right\|_f ^{2}\eta(\det(g))\int_{v\in V(\A_f)}s\overline\psi(u_0uv)dv=
\end{equation*}
\begin{equation*}
\int_{v\in V(\A_f)}\hat\varphi_{\left\{u_0u\right\}}(v,g)dv.
\end{equation*}
In other words,
\begin{equation*}
\chi_f(u_0)^{-1}\left\|u_0\right\|_f^{-m}\hat\varphi_{\left\{u_0u\right\}}-\hat\varphi_{\left\{u\right\}}\in \Sc(V(\A_f),\mu^{\otimes 0}\otimes\C)^0
\end{equation*}
and by Fourier inversion
\begin{equation*}
\chi_f(u_0)^{-1}\left\|u_0\right\|_f^{-m}\varphi_{\left\{u_0u\right\}}-\varphi_{\left\{u\right\}}\in \Sc(V(\A_f),\mu^{\otimes 0}\otimes\C)^0.
\end{equation*}
Now we see
\begin{equation*}
\rho_{m,0}^0\left(\frac{\chi_f(u_0)^{-1}\left\|u_0\right\|_f^{-m}\varphi_{\left\{u_0u\right\}}-\varphi_{\left\{u\right\}}}{\chi_f(u_0)^{-1}\left\|u_0\right\|_f^{-m}-1}\right)=
\rho_m\left(\frac{\chi_f(u_0)^{-1}\left\|u_0\right\|_f^{-m}\varphi_{\left\{u_0u\right\}}-\varphi_{\left\{u\right\}}}{\chi_f(u_0)^{-1}\left\|u_0\right\|_f^{-m}-1}\right)=
\end{equation*}
\begin{equation*}
\frac{\chi_f(u_0)^{-1}\left\|u_0\right\|_f^{-m}\rho_m(\varphi_{\left\{u_0u\right\}})-\rho_m(\varphi_{\left\{u\right\}})}{\chi_f(u_0)^{-1}\left\|u_0\right\|_f^{-m}-1}=
\frac{\chi_f(u_0)^{-1}\left\|u_0\right\|_f^{-m}\psi-\psi}{\chi_f(u_0)^{-1}\left\|u_0\right\|_f^{-m}-1}=\psi
\end{equation*}
\end{proof}
\end{lemma}
\subsection{The spherical functions}\label{spherical}
So what is left to be considered are
\begin{equation*}
\Ind_{B(\A_f)\times \pi_0(B(\R))}^{G(\A_f)\times \pi_0(G(\R))}\C \cdot\tilde\phi_f,\ \phi_{|T_1(\A)}=\left\|\ \cdot\ \right\|^2,
\end{equation*}
in particular $m=0$.
We choose the unique right-invariant Haar measure $ds$ on $SL_2(\A_{F,f})$ with $ds(SL_2(\hat\Oc))=1$. The group $SL_2(\A_{F,f})$ is unimodular. 
\begin{definition}
Let $\phi:T(\A)/T(\Q)\to \C^{\times}$ be an algebraic Hecke character of type $\gamma_{0,n}$ with $\phi_{|T^1(\A)}=\left\|\ \cdot \ \right\|^2$. Then we have the \textit{spherical function} 
\begin{equation*}
S(\phi)\in \Ind_{B(\A_f)\times \pi_0(B(\R))}^{G(\A_f)\times \pi_0(G(\R))}\C \cdot\tilde\phi_f
\end{equation*}
defined by
\begin{equation*}
S(\phi)(g):=\tilde\phi_f(b):=\tilde\phi_f(t_1,t_2),\ \text{if }
g=xb,\ x\in SL_2(\hat\Oc),\ b=\begin{pmatrix}t_1& t_1u\\ 0& t_2\end{pmatrix}\in B(\A)
\end{equation*}
\end{definition}
\begin{remark}
If we have a non-trivial $K_N$-invariant $\psi\in \Ind_{B(\A_f)\times \pi_0(B(\R))}^{G(\A_f)\times \pi_0(G(\R))}\C \cdot\tilde\phi_f$, then $\tilde\phi_f$ already has to be $K_N\cap T(\A_f)$-invariant and, in particular, the function $S(\phi)$ is $K_N$-invariant. 
\end{remark}
\begin{lemma}
Let $\phi:T(\A)\rightarrow \C^{\times}$ be an algebraic Hecke character of type $\gamma_{0,n}$ with $\phi_{|T^1(\A)}=\left\|\ \cdot \ \right\|^2$. Then we have a linear projection operator 
\begin{equation*}
\Psi_\phi:\Ind_{B(\A_f)\times \pi_0(B(\R))}^{G(\A_f)\times \pi_0(G(\R))}\C \cdot\tilde\phi_f\rightarrow \C\cdot S(\phi) 
\end{equation*}
defined by
\begin{equation*}
\Psi_\phi(\psi)(g):=\int_{s\in SL_2(\hat\Oc)}\psi(sg)ds
\end{equation*}
\begin{proof}
Write $g=xb$, $x\in SL_2(\hat\Oc)$, as above. Then
\begin{equation*}
\Psi_\phi(\psi)(g):=\int_{s\in SL_2(\hat\Oc)}\psi(sg)ds=\int_{s\in SL_2(\hat\Oc)}\psi(sxb)ds=
\end{equation*}
\begin{equation*}
\tilde\phi_f(b)\int_{s\in SL_2(\hat\Oc)}\psi(sx)ds=S(\phi)(g)\int_{s\in SL_2(\hat\Oc)}\psi(s)ds=S(\phi)(g)\Psi_\phi(\psi)(1).
\end{equation*}
\end{proof}
\end{lemma}
\begin{remark}
$\Psi_\phi(S(\phi))=S(\phi)$, so $\Psi_\phi$ is a split epimorphism, and $\psi\in \ker(\Psi_\phi)$ if and only if
\begin{equation*}
\int_{s\in SL_2(\hat\Oc)}\psi(s)ds=\Psi_\phi(\psi)(1)=0.
\end{equation*}
\end{remark}
\begin{lemma}
Let $\phi:T(\A)\rightarrow \C^{\times}$ be an algebraic Hecke character of type $\gamma_{0,1}$ with $\tilde\phi_{f|Z(\R)} =1$ and  $\phi_{|T^1(\A)}=\left\|\ \cdot \ \right\|^2$.
Then $\ker(\Psi_\phi)$ is in the image of $\rho_{0,0}^0$.
\begin{proof}
Let $\psi\in \ker(\Psi_\phi)$ be given. We suppose $\psi$ is left $K_N$-invariant and $\phi$ is $U_N$ invariant. 
For any set of representatives $\left\{u\right\}\subset \I_F$ of $Cl_F ^{(N)}$ the function $\varphi_{\left\{u\right\}}\in \Sc(V(\A_f),\mu^{\otimes 0}\otimes\C)$ defined by its Fourier transform
\begin{equation*}
\hat\varphi_{\left\{u\right\}}(v,g):=\Lambda_N(\chi,m+2)^{-1}\sum_u \chi^\prime_f(u)\left\|u\right\|_f ^{m+2}\eta(\det(g))s\overline\psi(uv)
\end{equation*}
is a preimage of $\psi$ under $\rho_m$. To prove the lemma we show that $\varphi_{\left\{u\right\}}\in \Sc(V(\A_f),\mu^{\otimes 0}\otimes\C)^0$. By Fourier inversion we have to show that $\hat\varphi_{\left\{u\right\}}\in \Sc(V(\A_f),\mu^{\otimes 0}\otimes\C)^0$ and therefore it suffices to prove
\begin{equation*}
\int_{v\in V(\A_f)}s\overline\psi(uv)dv=0.
\end{equation*}
Recall how we defined $s\overline\psi$ in the course of the proof of \cref{hor_surjective} and calculate
\begin{equation*}
\int_{v\in V(\A_f)}s\overline\psi(uv)dv=\left\|u\right\|_f ^{-2}\int_{v\in V(\A_f)}s\overline\psi(v)dv=\left\|u\right\|_f ^{-2}dv(NV(\hat\Z))\sum_{v\in V(\Z/N\Z)}s\overline\psi(v)=
\end{equation*}
\begin{equation*}
\left\|u\right\|_f ^{-2}dv(NV(\hat\Z))\sum_{x\in G(\Z/N\Z)/P(\Z/N\Z)}\overline\psi(x)=
\end{equation*}
\begin{equation*}
\left\|u\right\|_f ^{-2}dv(NV(\hat\Z))\sum_{s\in SL_2(\Oc/N\Oc)/U_0(\Oc/N\Oc)}\overline\psi(s)=\frac{\left\|u\right\|_f ^{-2}dv(NV(\hat\Z))}{|U_0(\Oc/N\Oc)|}\sum_{s\in SL_2(\Oc/N\Oc)}\overline\psi(s).
\end{equation*}
If we set $K_N\cap SL_2(\hat\Oc):=K_N ^{(1)}$, we may write this as
\begin{equation*}
\frac{\left\|u\right\|_f ^{-2}dv(NV(\hat\Z))}{|U_0(\Oc/N\Oc)|}\sum_{s\in K_N^{(1)}\backslash SL_2(\hat\Oc)}\overline\psi(s)=
\end{equation*}
\begin{equation*}
\frac{\left\|u\right\|_f ^{-2}dv(NV(\hat\Z))}{|U_0(\Oc/N\Oc)|ds(K_N ^{(1)})}\int_{s\in SL_2(\hat\Oc)}\overline\psi(s)ds=\frac{\left\|u\right\|_f ^{-2}dv(NV(\hat\Z))}{|U_0(\Oc/N\Oc)|ds(K_N ^1)}\int_{s\in SL_2(\hat\Oc)}\psi(s)ds=0,
\end{equation*}
as $\psi_{|SL_2(\hat\Oc)}=\overline\psi_{|SL_2(\hat\Oc)}$.
\end{proof}
\end{lemma}
Now we may define the operator $\Psi_{m,n}$
\begin{equation*}
\Ind_{B(\A_f)\times \pi_0(B(\R))}^{G(\A_f)\times \pi_0(G(\R))}\bigoplus_{\phi:\ type(\phi)=\gamma_{m,n},\tilde{\phi}_{f|Z(\R)=1}}\C \cdot\tilde\phi_f\rightarrow  
\bigoplus_{\phi:\ type(\phi)=\gamma_{m,n},\\ \tilde{\phi}_{f|Z(\R)=1},\phi_{|T^1(\A)}=\left\|\ \cdot\ \right\|^2}\C\cdot S(\phi), 
\end{equation*}
where $\Psi_{m,n}$ is zero on those summands with $\phi_{|T^1(\A)}\neq\left\|\ \cdot\ \right\|^2$ and $\Psi_{m,n}$ equals $\Psi_\phi$ defined above on those summands with $\phi_{|T^1(\A)}=\left\|\ \cdot\ \right\|^2$.
\begin{remark}
$\Psi_{m,n}$ is defined over $\overline\Q$. To see this consider a function
\begin{equation*}
\psi\in \Ind_{B(\A_f)\times \pi_0(B(\R))}^{G(\A_f)\times \pi_0(G(\R))}\overline\Q \cdot\tilde\phi_f\ \text{and}\ \Psi_\phi(\psi)=\int_{s\in SL_2(\hat\Oc)}\psi(s)ds\cdot S(\phi).
\end{equation*}
There is some compact open subgroup $K_f\subset SL_2(\hat\Oc)$ such that 
\begin{equation*}
\int_{s\in SL_2(\hat\Oc)}\psi(s)ds=\sum_{s\in K_f\backslash SL_2(\hat\Oc)}\psi(s)ds(K_f).
\end{equation*}
As the index of $K_f\subset SL_2(\hat\Oc)$ is finite and $ds(SL_2(\hat\Oc))=1$, we have $ds(K_f)\in \Q^{\times}$, and therefore 
\begin{equation*}
\int_{s\in SL_2(\hat\Oc)}\psi(s)ds\text{ and } S(\phi)(g)\in \overline \Q^{\times}\text{ for }g\in G(\A).
\end{equation*}
So the function $\Psi_\phi(\psi)$ has its values in $\overline\Q$.
Moreover, $\Psi_{m,n}$ is Galois-equivariant. This follows from
\begin{equation*}
 \int_{s\in SL_2(\hat\Oc)}\sigma \cdot\psi(s)ds=\sum_{s\in K_f\backslash SL_2(\hat\Oc)}\sigma\cdot\psi(s)ds(K_f)=\sigma\left(\int_{s\in SL_2(\hat\Oc)}\psi(s)ds\right),\ \sigma\in \text{Gal}(\overline\Q/\Q),
\end{equation*}
and we conclude that the vector space $\ker(\Psi_{m,n})$ has a natural $\Q$ structure which we denote by  $\ker(\Psi_{m,n|\Q})$.
\end{remark}

\begin{corollary}\label{ker(Psi)}
$\ker(\Psi_{m,n+1|\Q})$ is in the image of $\rho_{m,n}^0$.
\end{corollary}

\subsection{The long exact sequence for the cohomology of the boundary}

The next step is to show that $\ker(\Psi_{m,n+1|\Q})$ is exactly the image of $\rho_{m,n}^0$. To do so we use a long exact cohomology sequence relating the cohomology with compact supports, the usual cohomology of $\Sc_K$ and the cohomology of the boundary.

We denote again by $\overline{\Sc_K}$ the Borel-Serre compactification of $\Sc_K$. Consider
$j:\Sc_K\to \overline{\Sc_K}$ the open inclusion and $i:\overline{\Sc_K}\setminus\Sc_K\to \overline{\Sc_K}$ the closed inclusion of the boundary. For any abelian sheaf $\Fc$ on $\overline{\Sc_K}$ we have the exact sequence 
\begin{equation*}
 0\to j_!j^{-1}\Fc\to\Fc\to i_* i^{-1}\Fc\to 0.
\end{equation*}
If we apply $H^p_c(\overline{\Sc_K},\ )=H^p(\overline{\Sc_K},\ )$ to this sequence, we obtain a long exact cohomology sequence
\begin{equation*}
\cdots\to H^p_c(\overline{\Sc_K},j_!\Q)\to H^p(\overline{\Sc_K},\Q)\to H^p(\overline{\Sc_K}\setminus\Sc_K,\Q)\to H^{p+1}_c(\overline{\Sc_K},j_!\Q)\to\cdots
\end{equation*}
$H^p_c(\overline{\Sc_K},j_!\ )=H^p_c(\Sc_K,\ )$, the space $\overline{\Sc_K}\setminus\Sc_K$ is homotopy equivalent to $\partial \Sc_K$ and $\overline{\Sc_K}$ is homotopy equivalent to $\Sc_K$ by \cite{Ha1} 2.1, so we get the long exact cohomology sequence 
\begin{equation*}
\cdots\to H^p_c(\Sc_K,\Q)\to H^p(\Sc_K,\Q)\to H^p(\partial\Sc_K,\Q)\to H^{p+1}_c(\Sc_K,\Q)\to\cdots
\end{equation*}
\begin{lemma}
We have an exact sequence
\begin{equation*}
H^{2\xi-1}(\Sc_{K_N},\Q)\stackrel{\res_{\Sc_{K_N}}}{\longrightarrow} H^{2\xi-1}(\partial\Sc_{K_N},\Q)\rightarrow H^{2\xi}_c(\Sc_{K_N},\Q)\rightarrow 0
\end{equation*}
and $\image(\res_{\Sc_{K_N}})$ has codimension $|Cl_F ^{K_N}|=|Cl_F ^{(N)}|$ in $ H^{2\xi-1}(\partial\Sc_{K_N},\Q)$.
\begin{proof}
We start with the exact sequence
\begin{equation*}
H^{2\xi-1}(\Sc_{K_N},\Q)\stackrel{\res_{\Sc_{K_N}}}{\longrightarrow} H^{2\xi-1}(\partial\Sc_{K_N},\Q)\rightarrow H^{2\xi}_c(\Sc_{K_N},\Q)\rightarrow H^{2\xi}(\Sc_{K_N},\Q)
\end{equation*}
We use Poincaré duality to see
\begin{equation*}
H^{2\xi}_c(\Sc_{K_N},\C)=H^0(\Sc_{K_N},\C),\ H^{2\xi}(\Sc_{K_N},\C)=H^0 _c(\Sc_{K_N},\C)=0.
\end{equation*}
The latter is zero, because $\Sc_{K_N}$ is not compact. The first group is non-zero, we even conclude
\begin{equation*}
dim_\Q H^{2\xi}_c(\Sc_{K_N},\Q)=dim_\C H^0(\Sc_{K_N},\C)=|Cl_F ^{K_N}|=|Cl_F ^{(N)}|.
\end{equation*}
So we have the exact sequence
\begin{equation*}
H^{2\xi-1}(\Sc_{K_N},\Q)\stackrel{\res_{\Sc_{K_N}}}{\longrightarrow} H^{2\xi-1}(\partial\Sc_{K_N},\Q)\rightarrow H^{2\xi}_c(\Sc_{K_N},\Q)\rightarrow 0
\end{equation*}
and the codimension of $\image(\res_{\Sc_{K_N}})$ in $H^{2\xi-1}(\partial\Sc_{K_N},\Q)$ equals 
\begin{equation*}
dim_\Q H^{2\xi}_c(\Sc_{K_N},\Q)=|Cl_F ^{(N)}|.
\end{equation*}
\end{proof}
\end{lemma}
This lemma tells us that $\image(\Eis_0 ^0)$ cannot give more than a subspace of codimension $|Cl_F ^{(N)}|$ in the cohomology of the boundary in cohomological degree $2\xi-1$. In particular, $\rho_{0,-1}^0$ cannot be surjective.

\subsection{Determination of the image}

Now we can completely determine the image of our polylogarithmic Eisenstein operator.
\begin{lemma}\label{im_hor}
$\image(\rho_{m,n}^0)=\ker(\Psi_{m,n+1|\Q})$
\begin{proof}
It suffices to prove the case $m=0$ and $n=-1$ with complex coefficients.
Consider the epimorphism (the map is split)
\begin{equation*}
\Psi^{K_N}_{0,0}=\Psi_{0,0}:\left(\Ind_{B(\A_f)\times \pi_0(B(\R))}^{G(\A_f)\times \pi_0(G(\R))}\bigoplus_{\phi:\ type(\phi)=\gamma_{0,0},\tilde{\phi}_{f|Z(\R)=1}}\C \cdot\tilde\phi_f\right)^{K_N}\rightarrow 
\end{equation*}
\begin{equation*}
\left(\bigoplus_{\phi:\ type(\phi)=\gamma_{0,0},\tilde{\phi}_{f|Z(\R)=1},\phi_{|T^1(\A)}=\left\|\ \cdot\ \right\|^2}\C\cdot S(\phi)\right)^{K_N}
\end{equation*}
Any $\phi $ on the right-hand side is of the form
\begin{equation*}
\phi(t_1,t_2)=\eta(t_1t_2)\left\|\frac{t_2}{t_1}\right\|
\end{equation*}
and $\eta$ has to be trivial on $U_N$.
Therefore $\eta\in \widehat{Cl_F ^{K_N}}$ and the number of different $\phi$ on the right-hand side is exactly $|Cl_F^{K_N}|$. So the codimension of $\ker(\Psi^{K_N})$ inside 
\begin{equation*}
\left(\Ind_{B(\A_f)\times \pi_0(B(\R))}^{G(\A_f)\times \pi_0(G(\R))}\bigoplus_{\phi:\ type(\phi)=\gamma_{0,0},\tilde{\phi}_{f|Z(\R)=1}}\C \cdot\tilde\phi_f\right)^{K_N}
\end{equation*}
is $|Cl_F^{K_N}|$. By \Cref{ker(Psi)} we know that $\ker(\Psi^{K_N}_{0,0})\subset \image(\rho_{0,-1}^0)^{K_N}$, but 
\begin{equation*}
\image(\rho_{0,-1}^0)^{K_N}\subset \image(\res_{\Sc_{K_N}}) 
\end{equation*}
and $\image(\res_{\Sc_{K_N}})$ has already codimension $|Cl_F^{K_N}|$ inside
\begin{equation*}
H^{2\xi-1}(\partial\Sc_{K_N},\C)\cong\left(\Ind_{B(\A_f)\times \pi_0(B(\R))}^{G(\A_f)\times \pi_0(G(\R))}\bigoplus_{\phi:\ type(\phi)=\gamma_{0,0},\tilde{\phi}_{f|Z(\R)=1}}\C \cdot\tilde\phi_f\right)^{K_N}.
\end{equation*}
Therefore
\begin{equation*}
\ker(\Psi^{K_N}_{0,0})=\image(\rho_{0,-1}^0)^{K_N}= \image(\res_{\Sc_{K_N}})
\end{equation*}
and consequently $\image(\rho_{0,-1}^0)= \ker(\Psi_{0,0})$.
\end{proof}
\end{lemma}
\begin{theorem}\label{im_Eis}
Let 
\begin{equation*}
\Eis^k_q:\Sc(V(\A_f),\mu^{\otimes n})\otimes \mathfrak H^{q*}\rightarrow H^{2\xi-1-q}(\Sc,\Sym^k\Hc^\prime\otimes \mu^{\otimes n+1})
\end{equation*}
be the polylogarithmic Eisenstein operator defined in \Cref{Eis^k _q}. The image of this operator is isomorphic to $\ker(\Psi_{m,n+1|\Q})\otimes\Hc(T/Z)^{\xi-1-q}$.
\begin{proof}
By \Cref{pol_Eis} and \Cref{Harder} we know that $\res_\Sc: \image{\Eis^k_q}\rightarrow \image(\res_\Sc\circ \Eis^k _q)$
is an isomorphism. By \Cref{hor} we have the isomorphism
\begin{equation*}
 \image(\res_\Sc\circ \Eis^k _q)\cong \image(\rho_{m,n}^0)\otimes\Hc(T/Z)^{\xi-1-q}.
\end{equation*}
\Cref{im_hor} tells us $\ker(\Psi_{m,n+1|\Q})=\image(\rho_{m,n}^0)$, from which the theorem follows.
\end{proof}
\begin{corollary}
If $k>0$, then $\image(\Eis^k_q)=H^{2\xi-1-q}_{\Eis}(\Sc,\Sym^k\Hc^\prime\otimes \mu^{\otimes n+1})$. Moreover, $\image(\Eis^0_0)=H^{2\xi-1}_{\Eis}(\Sc,\mu^{\otimes n+1})$ and $\image(\Eis^0_q)^{K_N}\subset H^{2\xi-1-q}_{\Eis}(\Sc_{K_N},\mu^{\otimes n+1})$ is a subspace of codimension $|Cl_F ^{K_N}|$, if $q>0$.
\begin{proof}
This is just \Cref{im_Eis} together with \cite{Ha1} Theorem 2. For the calculation of the codimension see also the proof of \Cref{im_hor}.
\end{proof} 
\end{corollary}
\end{theorem}
\cleardoublepage
\addcontentsline{toc}{chapter}{Bibliography}
\bibliographystyle{apalike2}
\bibliography{diss_graf_arxiv}

\end{document}